\documentclass[12pt]{article}

\usepackage{amsmath}
\usepackage{amssymb}
\usepackage{graphicx}
\usepackage{epsfig}

\usepackage{longtable}

\usepackage{array}

\usepackage[margin=2cm]{geometry}
\newtheorem{theorem}{Theorem}
\newtheorem{lemma}[theorem]{Lemma}

\def \Z {{\mathbb{Z}}}

\def \a {\alpha}
\def \b {\beta}

\def \deg {{\rm deg}}
\def \proof {\noindent{\bf Proof}\quad}
\def \lan#1{\langle #1\rangle_n}

\def \floor#1{\lfloor#1\rfloor}

\def \cms {\overline{M_s}}

\def \Nbd {{\rm Nbd}}

\def \md#1{{\,({\rm mod}\ #1)}}
\def \cay {{\rm Cay}}

\def\endmark{\hskip 2em$\square$\par}
\def \qed {\hfill\endmark}

\title{\bf Decompositions of complete graphs into cycles of arbitrary lengths}

\author{
Darryn Bryant
\thanks{
The University of Queensland,
Department of Mathematics,
Qld 4072, Australia.
\texttt{db@maths.uq.edu.au},
}
\and
Daniel Horsley
\thanks{
School of Mathematical Sciences,
Monash University,
Vic 3800, Australia.
\texttt{danhorsley@gmail.com},
}
\and
William Pettersson
\thanks{
The University of Queensland,
Department of Mathematics,
Qld 4072, Australia.
\texttt{william.pettersson@gmail.com}
}
}

\date{ }

\begin{document}
\maketitle\thispagestyle{empty}

\begin{abstract}
 We show that the complete graph on $n$ vertices can be decomposed into $t$ cycles of specified
lengths $m_1,\ldots,m_t$ if and only if $n$ is odd, $3\leq m_i\leq n$ for $i=1,\ldots,t$, and
$m_1+\cdots+m_t=\binom n2$. We also show that the complete graph on $n$ vertices can be decomposed
into a perfect matching and $t$ cycles of specified lengths $m_1,\ldots,m_t$ if and only if $n$ is
even, $3\leq m_i\leq n$ for $i=1,\ldots,t$, and $m_1+\ldots+m_t=\binom n2-\frac n2$.
\end{abstract}

\section{Introduction}

A condensed version of this paper was published in \cite{BryHorPet}. Here, we include
supplementary data and some proofs which were omitted from \cite{BryHorPet} (the proofs which were
omitted generally run along similar lines to other included proofs). For example, proofs for all
cases of Lemmas \ref{Sofsize3-8}, \ref{3s4sandHamsforKn-S} and \ref{5sandHamsKn-S} were not given
in \cite{BryHorPet}, and neither were the decompositions in Sections \ref{appendix} and
\ref{appendix2}.

A decomposition of a graph $K$ is a set of subgraphs of $K$ whose edge sets partition the edge set
of $K$. In 1981, Alspach \cite{Als} asked whether it is possible to decompose the complete graph on
$n$ vertices, denoted $K_n$, into $t$ cycles of specified lengths $m_1,\ldots,m_t$ whenever the
obvious necessary conditions are satisfied; namely that $n$ is odd, $3\leq m_i\leq n$, and
$m_1+\cdots+m_t=\binom n2$. He also asked whether it is possible to decompose $K_n$ into a perfect
matching and $t$ cycles of specified lengths $m_1,\ldots,m_t$ whenever $n$ is even, $3\leq m_i\leq
n$, and $m_1+\cdots+m_t=\binom n2-\frac n2$. Again, these conditions are obviously necessary.

Here we solve Alspach's problem by proving the following theorem.

\begin{theorem}\label{mainthm}
There is a decomposition $\{G_1,\ldots,G_t\}$ of $K_n$ in which $G_i$ is an $m_i$-cycle for
$i=1,\ldots,t$ if and only if $n$ is odd, $3\leq m_i\leq n$ for $i=1,\ldots,t$, and
$m_1+\cdots+m_t=\frac{n(n-1)}2$. There is a decomposition $\{G_1,\ldots,G_t,I\}$ of $K_n$ in which
$G_i$ is an $m_i$-cycle for $i=1,\ldots,t$ and $I$ is a perfect matching if and only if $n$ is
even, $3\leq m_i\leq n$ for $i=1,\ldots,t$, and $m_1+\cdots+m_t=\frac{n(n-2)}2$.
\end{theorem}

Let $K$ be a graph and let $M=(m_1,\ldots,m_t)$ be a list of integers with $m_i\geq 3$ for
$i=1,\ldots,t$. If each vertex of $K$ has even degree, then an $(M)$-decomposition of $K$ is a
decomposition $\{G_1,\ldots,G_t\}$ such that $G_i$ is an $m_i$-cycle for $i=1,\ldots,t$. If each
vertex of $K$ has odd degree, then an $(M)$-decomposition of $K$ is a decomposition
$\{G_1,\ldots,G_t,I\}$ such that $G_i$ is an $m_i$-cycle for $i=1,\ldots,t$ and $I$ is a perfect
matching in $K$.

We say that a list $(m_1,\ldots,m_t)$ of integers is {\em $n$-admissible} if $3\leq
m_1,\ldots,m_t\leq n$ and $m_1+\cdots+m_t=n\lfloor\frac{n-1}{2}\rfloor$. Note that
$n\lfloor\frac{n-1}{2}\rfloor=\binom n2$ if $n$ is odd, and $n\lfloor\frac{n-1}{2}\rfloor=\binom
n2-\frac n2$ if $n$ is even. Thus, we can rephrase Alspach's question as follows. {\it Prove that
for each $n$-admissible list $M$, there exists an $(M)$-decomposition of $K_n$.}

A decomposition of $K_n$ into $3$-cycles is equivalent to a {\em Steiner triple system of order
$n$}, and a decomposition of $K_n$ into $n$-cycles is a {\em Hamilton decomposition}. Thus, the
work of Kirkman \cite{Kir} and Walecki (see \cite{Als3,Luc}) from the 1800s addresses Alspach's
problem in the cases where $M$ is of the form $(3,3,\ldots,3)$ or $(n,n,\ldots,n)$. The next
results on Alspach's problem appeared in the 1960s \cite{Kot,Ros1,Ros2}, and a multitude of results
have appeared since then. Many of these focused on the case of decompositions into cycles of
uniform length \cite{AlsMar,AlsVar,Bel,BerHuaSot,Hag,HofLinRod,Jac,RosHua}, and a complete solution
in this case was eventually obtained \cite{AlsGav,Saj}.

There have also been many papers on the case where the lengths of the cycles in the decomposition
may vary. In recent work \cite{BryHorPack,BryHorLong,BryHorAsym}, the first two authors have made
progress by developing methods introduced in \cite{BryHorPSTS} and \cite{BryHorMae}. In
\cite{BryHorLong}, Alspach's problem is settled in the case where all the cycle lengths are greater
than about $\frac n2$, and in \cite{BryHorAsym} the problem is completely settled for sufficiently
large odd $n$. Earlier results for the case of cycles of varying lengths can be found in
\cite{AdaBryKho1,AdaBryKho2,Bal1,BryKhoFu,BryMae,Hag,HeiHorRos,Jor}. See \cite{Bry2} for a survey
on Alspach's problem, and see \cite{BryRod} for a survey on cycle decompositions generally.

The analogous problems on decompositions of complete graphs into matchings, stars or paths have all
been completely solved, see \cite{Bar}, \cite{LinShy} and \cite{Bry3} respectively. It is also
worth mentioning that the easier problems in which each $G_i$ is required only to be a closed trail
of length $m_i$ or each $G_i$ is required only to be a $2$-regular graph of order $m_i$ have been
solved in \cite{Bal2}, \cite{BryHorMae} and \cite{BryMae2}. Decompositions of complete multigraphs
into cycles are considered in \cite{BryHorMaeSmi}.

Balister \cite{Bal1} has verified by computer that Theorem \ref{mainthm} holds for $n\leq 14$, and
we include this result as a lemma for later reference.

\begin{lemma}[\cite{Bal1}]\label{SmallCases}
Theorem \ref{mainthm} holds for $n\leq 14$.
\end{lemma}

Our proof of Theorem \ref{mainthm} relies heavily on the reduction of Alspach's problem
obtained in \cite{BryHorAsym}, see Theorem \ref{ReductionTheorem} below.
Throughout the paper, we use
the notation $\nu_i(M)$ to denote the number of occurrences
of $i$ in a given list $M$.

\vspace{0.5cm}

\noindent{\bf Definition}\quad A list $M$ is an \emph{$n$-ancestor} list if it is $n$-admissible
and satisfies
\begin{itemize}
\item [(1)] $\nu_6(M)+\nu_7(M)+\cdots+\nu_{n-1}(M)\in\{0,1\}$;
\item [(2)] if $\nu_5(M)\geq 3$, then $2\nu_4(M)\leq n-6$;
\item [(3)] if $\nu_5(M)\geq 2$, then $3\nu_3(M)\leq n-10$;
\item [(4)] if $\nu_4(M)\geq 1$ and $\nu_5(M)\geq 1$, then $3\nu_3(M)\leq n-9$;
\item [(5)] if $\nu_4(M)\geq 1$, then $\nu_i(M)=0$ for $i\in\{n-2,n-1\}$; and
\item [(6)] if $\nu_5(M)\geq 1$, then $\nu_i(M)=0$ for $i\in\{n-4,n-3,n-2,n-1\}$.
\end{itemize}

\vspace{0.5cm}

Thus, an $n$-ancestor list is of the form
$$(3,3,\ldots,3,4,4,\ldots,4,5,5,\ldots,5,k,n,n,\ldots,n)$$
where $k$ is either absent or in the range $6 \leq k \leq n-1$,
and there are additional constraints involving the number of occurrences of
cycle lengths in the list.
The following theorem was proved in \cite{BryHorAsym}.

\begin{theorem}\label{ReductionTheorem}{\rm (\cite{BryHorAsym}, Theorem 4.1)}
For each positive integer $n$, if there exists an $(M)$-decomposition of $K_n$ for each
$n$-ancestor list $M$, then there exists an $(M)$-decomposition of $K_n$ for each $n$-admissible
list $M$.
\end{theorem}

Our goal is to construct an $(M)$-decomposition of $K_n$ for each $n$-ancestor list $M$.
We split this problem into two cases:
the case where $\nu_n(M)\geq 2$ and
the case where $\nu_n(M)\leq 1$.
In particular, we prove the following two results.

\begin{lemma}\label{MainLemmaforatleast2hams}
If $M$ is an $n$-ancestor list with $\nu_n(M)\geq 2$, then there is an $(M)$-decomposition of
$K_n$.
\end{lemma}

\proof
See Section \ref{atleast2hamssection}.
\qed

\begin{lemma}\label{MainLemmaforatmost1ham}
If Theorem \ref{mainthm} holds for $K_{n-3}$, $K_{n-2}$ and $K_{n-1}$,
then there is an $(M)$-decomposition of $K_n$ for each $n$-ancestor list
$M$ satisfying $\nu_n(M)\leq 1$.
\end{lemma}

\proof
The case $\nu_n(M)=0$ is proved in Section \ref{zerohamsection}
(see Lemma \ref{Mainlemmaforzerohams})
and the
case $\nu_n(M)=1$ is proved in Section \ref{onehamsection}
(see Lemma \ref{Mainlemmaforoneham}).
\qed

\vspace{0.5cm}

Lemmas
\ref{MainLemmaforatleast2hams}
and
\ref{MainLemmaforatmost1ham}
allow us to
prove our main result using induction on $n$.

\vspace{0.5cm}

\noindent{\bf Proof of Theorem \ref{mainthm}}\quad The proof is by induction on $n$. By Lemma
\ref{SmallCases}, Theorem \ref{mainthm} holds for $n \leq 14$. So let $n\geq 15$ and assume Theorem
\ref{mainthm} holds for complete graphs having fewer than $n$ vertices. By Theorem
\ref{ReductionTheorem}, it suffices to prove the existence of an $(M)$-decomposition of $K_n$ for
each $n$-ancestor list $M$. Lemma \ref{MainLemmaforatleast2hams} covers each $n$-ancestor list $M$
with $\nu_n(M)\geq 2$, and using the inductive hypothesis, Lemma \ref{MainLemmaforatmost1ham}
covers those with $\nu_n(M)\leq 1$. \qed

\section{Notation}

We shall sometimes use superscripts to specify the number of occurrences of a particular integer in
a list. That is, we define $(m_1^{\a_1},\ldots,m_t^{\a_t})$ to be the list comprised of $\a_i$
occurrences of $m_i$ for $i=1,\ldots,t$. Let $M=(m_1^{\a_1},\ldots,m_t^{\a_t})$ and let
$M'=(m_1^{\b_1},\ldots,m_t^{\b_t})$, where $m_1,\ldots,m_t$ are distinct. Then $(M,M')$ is the list
$(m_1^{\a_1+\b_1},\ldots,m_t^{\a_t+\b_t})$ and, if $0\leq \b_i\leq \a_i$ for $i=1,\ldots,t$, $M-M'$
is the list $(m_1^{\a_1-\b_1},\ldots, m_t^{\a_t-\b_t})$.

Let $\Gamma$ be a finite group and let $S$ be a subset of $\Gamma$ such that
the identity of $\Gamma$ is not in $S$ and such that the inverse of any element of
$S$ is also in $S$. The {\em Cayley graph} on $\Gamma$ with {\em connection set} $S$,
denoted $\cay(\Gamma,S)$, has the elements of $\Gamma$ as its vertices and there is an edge
between vertices $g$ and $h$ if and only if $g=hs$ for some $s\in S$.

A Cayley graph on a cyclic group is called a {\em circulant graph}. For any graph with vertex set
$\Z_n$, we define the {\it length} of an edge $xy$ to be $x-y$ or $y-x$, whichever is in
$\{1,\ldots,\floor{\frac n2}\}$. It is convenient to be able to describe the connection set of a
circulant graph on $\Z_n$ by listing only one of $s$ and $n-s$. Thus, we use the following
notation. For any subset $S$ of $\Z_n\setminus\{0\}$ such that $s\in S$ and $n-s\in S$ implies
$n=2s$, we define $\lan{S}$ to be the Cayley graph $\cay(\Z_n,S\cup -S)$.

Let $m\in\{3,4,5\}$ and let $D=\{a_1,\ldots,a_m\}$ where $a_1,\ldots,a_m$ are positive integers. If
there is a partition $\{D_1,D_2\}$ of $D$ such that $\sum D_1-\sum D_2=0$, then $D$ is called a
\emph{difference $m$-tuple}. If there is a partition $\{D_1,D_2\}$ of $D$ such that $\sum D_1-\sum
D_2=0\md n$, then $D$ is called a \emph{modulo $n$ difference $m$-tuple}. Clearly, any difference
$m$-tuple is also a modulo $n$ difference $m$-tuple for all $n$. We may use the terms difference
triple, quadruple and quintuple respectively rather than $3$-tuple, $4$-tuple and $5$-tuple. For
$m\in\{3,4,5\}$, it is clear that if $D$ is a difference $m$-tuple, then there is an
$(m^n)$-decomposition of $\lan{D}$ for all $n\geq 2\max(D)+1$, and that if $D$ is a modulo $n$
difference $m$-tuple, then there is an $(m^n)$-decomposition of $\lan{D}$.

We denote the complete graph with vertex set $V$ by $K_V$ and the complete bipartite graph with
parts $U$ and $V$ by $K_{U,V}$. If $G$ and $H$ are graphs then $G-H$ is the graph with vertex set
$V(G) \cup V(H)$ and edge set $E(G) \setminus E(H)$. If $G$ and $H$ are graphs whose vertex sets
are disjoint then $G \vee H$ is the graph with vertex set $V(G) \cup V(H)$ and edge set $E(G) \cup
E(H) \cup \{xy:x\in V(G),y\in V(H)\}$. A cycle with $m$ edges is called an \emph{$m$-cycle} and is
denoted $(x_1,\ldots,x_m)$, where $x_1,\ldots,x_m$ are the vertices of the cycle and
$x_1x_2,\ldots,x_{m-1}x_m,x_mx_1$ are the edges. A path with $m$ edges is called an \emph{$m$-path}
and is denoted $[x_0,\ldots,x_m]$, where $x_0,\ldots,x_m$ are the vertices of the path and
$x_0x_1,\ldots,x_{m-1}x_m$ are the edges. A graph is said to be \emph{even} if every vertex of the
graph has even degree and is said to be \emph{odd} if every vertex of the graph has odd degree.

A \emph{packing} of a graph $K$ is a decomposition of some subgraph $G$ of $K$,
and the graph $K-G$ is called the \emph{leave} of the packing.
An \emph{$(M)$-packing} of $K_n$ is an $(M)$-decomposition of
some subgraph $G$ of $K_n$ such that $G$ is an even graph if $n$ is odd and $G$ is an odd graph if
$n$ is even
(recall that an $(M)$-decomposition of an odd graph contains a perfect matching).
Thus, the leave of an $(M)$-packing of $K_n$ is an even graph and, like an
$(M)$-decomposition of $K_n$, an $(M)$-packing of $K_n$ contains a perfect matching if and only if
$n$ is even. A decomposition of a graph into Hamilton cycles is called a \emph{Hamilton
decomposition}.

\section{The case of at least two Hamilton cycles}\label{atleast2hamssection}

The purpose of this section is to prove Lemma \ref{MainLemmaforatleast2hams} which states that
there is an $(M)$-decomposition of $K_n$ for each $n$-ancestor list $M$ with $\nu_n(M)\geq 2$. We
first give a general outline of this proof. Theorem \ref{mainthm} has been proved in the case where
$M = (3^a,n^b)$ for some $a,b \geq 0$ \cite{BryMae}, so we will restrict our attention to ancestor
lists which are not of this form. The basic construction involves decomposing $K_n$ into $\lan{S}$
and $K_n-\lan{S}$ where, for some $x\leq 8$, the connection set $S$ is either $\{1,\ldots,x\}$ or
$\{1,\ldots,x-1\} \cup \{x+1\}$ so that $\sum S$ is even. We partition any given $n$-ancestor list
$M$ into two lists $M_s$ and $\overline{M_s}=M-M_s$, and construct an $(M_s)$-decomposition of
$\lan{S}$ and an $(\cms)$-decomposition of $K_n-\lan{S}$. This yields the desired
$(M)$-decomposition of $K_n$. Taking $S=\{1,\ldots,x-1\} \cup \{x+1\}$, rather than
$S=\{1,\ldots,x\}$, is necessary when $1+\cdots+x$ is odd as many desired cycle decompositions of
$\lan{\{1,\ldots,x\}}$ do not exist when $1+\cdots+x$ is odd, see \cite{BryMar}.

If $M=(3^{\alpha_3n+\beta_3},4^{\alpha_4n+\beta_4},5^{\alpha_5n+\beta_5},k^\gamma,n^\delta)$ where
$\alpha_i \geq 0$ and $0\leq \b_i\leq n-1$ for $i \in \{3,4,5\}$, $6 \leq k \leq n-1$, $\gamma \in
\{0,1\}$, and $\delta \geq 2$, then we usually choose
$M_s=(3^{\beta_3},4^{\beta_4},5^{\beta_5},k^\gamma)$. However, if this would result in $\sum M_s$
being less than $4n$, then we sometimes adjust this definition slightly. We always choose $M_s$
such that $\sum M_s$ is at most $8n$,
which explains why we have $|S|\leq 8$.

Our $(M_s)$-decompositions of $\lan{S}$ will be constructed using adaptations of techniques used in
\cite{BryMar} and \cite{BrySch}. We construct our $(\cms)$-decompositions of $K_n-\lan{S}$ using a
combination of difference methods and results on Hamilton decompositions of circulant graphs. In
general, we split the problem into the case $\nu_5(M)\leq 2$ and the case $\nu_5(M)\geq 3$. In the
former case it will follow from our choice of $M_s$ that $\cms=(3^{tn},4^{qn},n^h)$ for some $t,q,h
\geq 0$ and in the latter case it will follow from our choice of $M_s$ that $\cms=(5^{rn},n^h)$ for
some $r,h \geq 0$.

The precise definition of $M_s$ is given in Lemma \ref{MsLemma}, which details the properties that
we require of our partition of $M$ into $M_s$ and $\cms$, and establishes its existence. The
definition includes several minor technicalities in order to deal with complications and exceptions
that arise in the above-described approach. Throughout the remainder of this section, for a given
$n$-ancestor list $M$ such that $\nu_n(M)\geq 2$ and $M \neq (3^a,n^b)$ for any $a,b \geq 0$, we
shall use the notation $M_s$ and $\cms$ to denote the lists constructed in the proof of Lemma
\ref{MsLemma}. If $\nu_n(M)\leq 1$ or $M = (3^a,n^b)$ for some $a,b \geq 0$, then $M_s$ and $\cms$
are not defined.

\begin{lemma}\label{MsLemma}
If $M$ is any $n$-ancestor list such that $\nu_n(M)\geq 2$ and $M \neq (3^a,n^b)$ for any $a,b \geq
0$, then there exists a partition of $M$ into sublists $M_s$ and $\cms$ such that
\begin{itemize}
    \item[(1)]
$\sum M_s\in\{2n,3n,\ldots,8n\}$ and $\sum M_s\neq 8n$ when $\nu_5(M)\leq 2$;
    \item[(2)]
if $\sum M_s=2n$, then $\nu_n(M_s)=1$ and $\cms=(n^h)$ for some $h \geq 1$;
    \item[(3)]
if $\sum M_s=3n$, then $\nu_n(M_s)\in\{0,1\}$ and $\cms=(n^h)$ for some $h \geq 1$;
    \item[(4)]
if $\sum M_s\in\{4n,5n,\ldots,8n\}$ and $\nu_5(M)\geq 3$, then $\nu_n(M_s)=0$ and
$\cms=(5^{rn},n^h)$ for some $r\geq 0$, $h\geq 2$;
    \item[(5)]
if $\sum M_s\in\{4n,5n,\ldots,7n\}$ and $\nu_5(M)\leq 2$, then $\nu_n(M_s)=0$ and
$\cms=(3^{tn},4^{qn},n^h)$ for some $t,q\geq 0$, $h\geq 2$; and
    \item [(6)]
$M_s \neq (3^{\frac{5n}3})$.
\end{itemize}
\end{lemma}

\proof Let $M$ be an $n$-ancestor list. The conditions of the lemma imply $n \geq 7$.
We will first define a list $M_e$ which in many cases will serve as $M_s$, but will sometimes
need to be adjusted slightly.

If
$$M=(3^{\alpha_3n+\beta_3},4^{\alpha_4n+\beta_4},5^{\alpha_5n+\beta_5},k^\gamma,n^\delta)$$
where $\alpha_i \geq 0$ and $0\leq \b_i\leq n-1$ for $i \in \{3,4,5\}$, $6 \leq k \leq n-1$,
$\gamma \in \{0,1\}$, and $\delta \geq 2$, then
$$M_e=(3^{\beta_3},4^{\beta_4},5^{\beta_5},k^\gamma).$$
It is clear from the definition of $n$-ancestor list that if we take $M_s=M_e$, then (4) and (5)
are satisfied.

We now show that $\sum M_e \in \{0,n,2n,\ldots,8n\}$, and that $\sum M_e \neq 8n$ when
$\nu_5(M)\leq 2$. Noting that $\sum M_e \leq 3\beta_3+4\beta_4+5\beta_5+(n-1)$ and separately
considering the cases $\nu_5(M) \geq 3$, $\nu_5(M) \in \{1,2\}$ and $\nu_5(M)=0$, it is routine to
use the definition of $(M)$-ancestor lists to show that $\sum M_e < 9n$, and that $\sum M_e < 8n$
when $\nu_5(M)\leq 2$. Thus, because it follows from $\sum M=n\floor{\frac{n-1}2}$ and the
definition of $M_e$ that $n$ divides $\sum M_e$, we have that $\sum M_e \in \{0,n,2n,\ldots,8n\}$,
and that $\sum M_e \neq 8n$ when $\nu_5(M)\leq 2$.

If $\sum M_e\in\{4n,5n,6n,7n,8n\}$, then we let $M_s=M_e$. If $\sum M_e\in\{0,n,2n,3n\}$, then we
define $M_s$ by
$$
M_s=
\begin{cases}
(M_e,4^n) &\mbox{if $\alpha_4>0$;}\\
(M_e,5^n) &\mbox{if $\alpha_4=0$ and $\alpha_5>0$;}\\
(M_e,3^n) &\mbox{if $\alpha_4=\alpha_5=0$ and $\alpha_3>0$;}\\
(M_e,n)   &\mbox{if $\alpha_3=\alpha_4=\alpha_5=0$ and $\sum M_e\in\{n,2n\}$;} \\
M_e       &\mbox{otherwise.} \\
\end{cases}
$$
Using the definition of $M_s$ and the fact that $M$ is an $n$-ancestor list
with $M \neq (3^a,n^b)$ for any $a,b \geq 0$, it is routine to check
that $M_s$ satisfies (1)--(6).
\qed

\vspace{0.5cm}

Before proving Lemma \ref{MainLemmaforatleast2hams}, we need a number of preliminary lemmas. The
first three give us the necessary decompositions of $\lan{S}$ where $S=\{1,\ldots,x\}$ or
$S=\{1,\ldots,x-1\} \cup \{x+1\}$ for some $x\leq 8$. Lemma \ref{S=12} was proven independently in
\cite{Bry} and \cite{Rod}, and is a special case of Theorem 5 in \cite{BryMar}. Lemmas
\ref{Sofsize3-8} and \ref{S=123includingham} will be proved in Section \ref{decomposeS}.

\begin{lemma}\label{Sofsize3-8}
If
$$S\in\{\{1,2,3\},\{1,2,3,4\},\{1,2,3,4,6\},\{1,2,3,4,5,7\},\{1,2,3,4,5,6,7\},\{1,2,3,4,5,6,7,8\}\},$$
$n\geq 2\max(S)+1$, and $M=(m_1,\ldots,m_t,k)$ is any list satisfying $m_i\in\{3,4,5\}$ for
$i=1,\ldots,t$, $3\leq k\leq n$, and $\sum M=|S|n$, then there is an $(M)$-decomposition of
$\lan{S}$, except possibly when
\begin{itemize}
\item $S=\{1,2,3,4,6\}$, $n\equiv 3\md 6$ and $M=(3^{\frac{5n}3})$; or
\item $S=\{1,2,3,4,6\}$, $n\equiv 4\md 6$ and $M=(3^{\frac{5n-5}3},5)$.
\end{itemize}
\end{lemma}

\proof See Section \ref{decomposeS}. \qed

\vspace{0.5cm}

\begin{lemma}\label{S=12}{\rm (\cite{Bry,Rod})}
If $n\geq 5$ and $M=(m_1,\ldots,m_t,n)$ is any list satisfying $m_i\in\{3,\ldots,n\}$ for
$i=1,\ldots,t$, and $\sum M=2n$, then there is an $(M)$-decomposition of $\lan{\{1,2\}}$.
\end{lemma}

\begin{lemma}\label{S=123includingham}
If $n\geq 7$ and $M=(m_1,\ldots,m_t,k,n)$ is any list satisfying $m_i\in\{3,4,5\}$ for
$i=1,\ldots,t$, $3\leq k\leq n$, and $\sum M=3n$, then there is an $(M)$-decomposition of
$\lan{\{1,2,3\}}$.
\end{lemma}

\proof See Section \ref{decomposeS}. \qed

\vspace{0.5cm}

We now present the lemmas which give us the necessary decompositions of $K_n-\lan{S}$. Lemma
\ref{hamdecompxtonover2} was proved in \cite{BryMae} where it was used to prove Theorem
\ref{mainthm} in the case where $M = (3^a,n^b)$ for some $a,b \geq 0$. Lemmas
\ref{3s4sandHamsforKn-S} and \ref{5sandHamsKn-S} give our main results on decompositions of
$K_n-\lan{S}$. Lemma \ref{3s4sandHamsforKn-S} is for the case $\nu_5(M)\leq 2$ and Lemma
\ref{5sandHamsKn-S} is for the case $\nu_5(M)\geq 3$.

\begin{lemma}{\rm (\cite{BryMae}, Lemma 3.1)}\label{hamdecompxtonover2}
If $1\leq h\leq\lfloor\frac{n-1}2\rfloor$, then there is an $(n^h)$-decomposition of
$K_n-\lan{\{1,\ldots,\floor{\frac{n-1}2}-h\}}$.
\end{lemma}

\begin{lemma}\label{3s4sandHamsforKn-S}
If $S\in\{\{1,2,3,4\},\{1,2,3,4,6\},\{1,2,3,4,5,7\},\{1,2,3,4,5,6,7\}\}$ and $n\geq 2\max(S) +1$,
$t\geq 0$, $q\geq 0$ and $h\geq 2$ are integers satisfying $3t+4q+h=\lfloor\frac{n-1}2\rfloor-|S|$,
then there is a $(3^{tn},4^{qn},n^h)$-decomposition of $K_n-\lan{S}$, except possibly when $h=2$,
$S=\{1,2,3,4,5,6,7\}$ and
\begin{itemize}
\item $n\in\{25,26\}$ and $t=1$; or
\item $n=31$ and $t=2$.
\end{itemize}
\end{lemma}

\proof See Section \ref{decomposeKminusS}. \qed

\begin{lemma}\label{5sandHamsKn-S}
If $S\in\{\{1,2,3,4\},\{1,2,3,4,6\},\{1,2,3,4,5,7\},\{1,2,3,4,5,6,7\},\{1,2,3,4,5,6,7,8\}\}$ and
$n\geq 2\max(S)+1$, $r\geq 0$ and $h\geq 2$ are integers satisfying
$5r+h=\lfloor\frac{n-1}2\rfloor-|S|$, then there is a $(5^{rn},n^h)$-decomposition of
$K_n-\lan{S}$.
\end{lemma}

\proof See Section \ref{decomposeKminusS}. \qed

\vspace{0.5cm}

We also need Lemmas \ref{specialcase1} and \ref{specialcase2} below to deal with cases arising
from the possible exceptions in Lemmas \ref{Sofsize3-8} and \ref{3s4sandHamsforKn-S}
respectively. To prove Lemma \ref{specialcase1} we use the following special case of Lemma 2.8 in
\cite{BryHorAsym}.

\begin{lemma}\label{chainlemma}
If there exists an $(M,4^2)$-decomposition of $K_n$ in which there are two $4$-cycles intersecting
in exactly one vertex, then there exists an $(M,3,5)$-decomposition of $K_n$.
\end{lemma}

\begin{lemma}\label{specialcase1}
If $M$ is an $n$-ancestor list such that $\nu_n(M)\geq 2$, $M_s=(3^{\frac{5n-5}3},5)$ and $n\equiv
4\md 6$ then there is an $(M)$-decomposition of $K_n$.
\end{lemma}

\proof
We will construct an $(\cms,3^{\frac{5n-8}3},4^2)$-decomposition of $K_n$ in which two $4$-cycles
intersect in exactly one vertex. The required $(M)$-decomposition of $K_n$ can then be obtained by
applying Lemma \ref{chainlemma}.

By Lemma \ref{3s4sandHamsforKn-S} there is an $(\cms)$-decomposition of $K_n-\lan{\{1,2,3,4,6\}}$,
so it suffices to construct a $(3^{\frac{5n-8}3},4^2)$-decomposition of $\lan{\{1,2,3,4,6\}}$ in
which the two $4$-cycles intersect in exactly one vertex for all $n\equiv 4\md 6$ with $n\geq 16$
(note that the conditions of the lemma imply $n\geq 16$). The union of the following two sets of
cycles gives such a decomposition.
$$\{\ (0,4,2,6),(2,3,5,8),(1,5,7),(3,4,7),(3,6,9),(4,5,6)\ \}$$
$$
\begin{array}{l}
  \{\ (x+6i,y+6i,z+y6)\ : i\in\{0,\ldots,\textstyle\frac{n-10}6\},\
  (x,y,z)\in\{(4,8,10),(5,9,11),(6,8,12),(6,7,10), \\ \quad(7,11,13),(7,8,9),
(9,12,15),(9,10,13),(10,11,12),(8,11,14)\ \}
\end{array}
$$
\qed

\begin{lemma}\label{specialcase2}
If $M$ is an $n$-ancestor list such that $\nu_5(M)\leq 2$, $\nu_n(M)=2$, $\sum M_s=7n$, and
\begin{itemize}
\item $n=25$ and $\nu_3(\cms)=25$;
\item $n=26$ and $\nu_3(\cms)=26$; or
\item $n=31$ and $\nu_3(\cms)=62$;
\end{itemize}
then there is an $(M)$-decomposition of $K_n$.
\end{lemma}

\proof We begin by showing that it is possible to partition $M_s$ into two lists $M_s^1$ and
$M_s^2$ such that $\sum M_s^1=3n$ and $\sum M_s^2=4n$. If $\nu_3(M_s)\geq n$ or $\nu_4(M_s)\geq n$,
then clearly such a partition exists. Otherwise, $\nu_n(M_s)=0$ by Property (5) of Lemma
\ref{MsLemma}, and so by the definition of $n$-ancestor list and the hypotheses of this lemma, we
have that
$$7n=\textstyle{\sum M_s}\leq 3\nu_3(M_s)+4\nu_4(M_s)+10+(n-1).$$
It is routine to check, using $3\nu_3(M_s)\leq 3n-3$ and $4\nu_4(M_s)\leq 4n-4$, that
$\nu_4(M_s)\geq \frac{3n-6}{4}$ and $\nu_3(M_s)\geq \frac{2n-5}{3}$. Thus for $n=25$, $n=26$ and
$n=31$, we can choose $M_s^1=(3,4^{18})$, $M_s^1=(3^2,4^{18})$, and $M_s^1=(3^3,4^{21})$
respectively. This yields the desired partition of $M_s$.

For $n=25$ we note that $\lan{\{1,2,3\}}\cong\lan{\{2,4,6\}}$ (with $x\mapsto 2x$ being an
isomorphism) and $\lan{\{1,2,3,4\}}\cong\lan{\{1,7,8,9\}}$ (with $x\mapsto 8x$ being an
isomorphism). Since $\{3,10,12\}$ is a modulo $25$ difference triple and
$\langle\{5,11\}\rangle_{25}$ has a Hamilton decomposition (by a result of Bermond et al
\cite{BerFavMah}, see Lemma \ref{BFMTheorem}), this gives us a decomposition of $K_{25}$ into a
copy of $\langle\{1,2,3\}\rangle_{25}$, a copy of $\langle\{1,2,3,4\}\rangle_{25}$, twenty-five
$3$-cycles and two Hamilton cycles. By Lemma \ref{Sofsize3-8}, there is an $(M_s^1)$-decomposition
of $\langle\{1,2,3\}\rangle_{25}$ and an $(M_s^2)$-decomposition of
$\langle\{1,2,3,4\}\rangle_{25}$, and this gives us the required $(M)$-decomposition of $K_{25}$.

For $n=26$ we note that $\lan{\{1,2,3,4\}}\cong\lan{\{5,6,10,11\}}$ (with $x\mapsto 5x$ being an
isomorphism). Since $\{4,8,12\}$ is a difference triple and $\langle\{7,9\}\rangle_{26}$ has a
Hamilton decomposition (by a result of Bermond et al \cite{BerFavMah}, see Lemma \ref{BFMTheorem}),
this gives us a decomposition of $K_{26}$ into a copy of $\langle\{1,2,3\}\rangle_{26}$, a copy of
$\langle\{1,2,3,4\}\rangle_{26}$, twenty-six $3$-cycles and two Hamilton cycles. By Lemma
\ref{Sofsize3-8}, there is an $(M_s^1)$-decomposition of $\langle\{1,2,3\}\rangle_{26}$ and an
$(M_s^2)$-decomposition of $\langle\{1,2,3,4\}\rangle_{26}$, and this gives us the required
$(M)$-decomposition of $K_{26}$.

For $n=31$ we note that $\lan{\{1,2,3,4\}}\cong\lan{\{4,8,12,15\}}$ (with $x\mapsto 4x$ being an
isomorphism). Since $\{5,6,11\}$ is a difference triple, $\{7,10,14\}$ is a modulo $31$ difference
triple, and $\langle\{9,13\}\rangle_{31}$ has a Hamilton decomposition (by a result of Bermond et
al \cite{BerFavMah}, see Lemma \ref{BFMTheorem}), this gives us a decomposition of $K_{31}$ into a
copy of $\langle\{1,2,3\}\rangle_{31}$, a copy of $\langle\{1,2,3,4\}\rangle_{31}$, sixty-two
$3$-cycles and two Hamilton cycles. By Lemma \ref{Sofsize3-8}, there is an $(M_s^1)$-decomposition
of $\langle\{1,2,3\}\rangle_{31}$ and an $(M_s^2)$-decomposition of
$\langle\{1,2,3,4\}\rangle_{31}$, which yields required $(M)$-decomposition of $K_{31}$. \qed

\vspace{0.5cm}

We can now prove Lemma \ref{MainLemmaforatleast2hams} which states that if $M$ is an $n$-ancestor
list with $\nu_n(M)\geq 2$, then there is an $(M)$-decomposition of $K_n$.

\vspace{0.5cm}

\noindent{\bf Proof of Lemma \ref{MainLemmaforatleast2hams}}\quad If $M = (3^a,n^b)$ for some $a,b
\geq 0$, then we can use the main result from \cite{BryMae} to obtain an $(M)$-decomposition of
$K_n$, so we can assume that $M \neq (3^a,n^b)$ for any $a,b \geq 0$. By Lemma \ref{SmallCases} we
can assume that $n \geq 15$. Partition $M$ into $M_s$ and $\cms$.
The proof splits into cases according to the value of $\sum M_s$,
which by Lemma \ref{MsLemma} is in $\{2n,3n,\ldots,8n\}$.

\noindent{\bf Case 1}\quad Suppose that $\sum M_s=2n$. In this case, from Property (2) of Lemma
\ref{MsLemma} we have $\nu_n(M_s)=1$ and $\cms=(n^h)$ for some $h \geq 1$. The required
decomposition of $K_n$ can be obtained by combining an $(M_s)$-decomposition of $\lan{\{1,2\}}$
(which exists by Lemma \ref{S=12}) with a Hamilton decomposition of $K_n-\lan{\{1,2\}}$ (which
exists by Lemma \ref{hamdecompxtonover2}).

\noindent{\bf Case 2}\quad Suppose that $\sum M_s=3n$. In this case, from Property (3) of Lemma
\ref{MsLemma} we have $\nu_n(M_s)\in\{0,1\}$ and $\cms=(n^h)$ for some $h \geq 1$. The required
decomposition of $K_n$ can be obtained by combining an $(M_s)$-decomposition of $\lan{\{1,2,3\}}$
(which exists by Lemma \ref{Sofsize3-8} or \ref{S=123includingham}) with a Hamilton decomposition
of $K_n-\lan{\{1,2,3\}}$ (which exists by Lemma \ref{hamdecompxtonover2}).

\noindent{\bf Case 3}\quad Suppose that $\sum M_s \in \{4n,5n,6n,7n,8n\}$ and $\nu_5(M)\geq 3$. In
this case, from Property (4) of Lemma \ref{MsLemma} we have $\nu_n(M_s)=0$ and $\cms=(5^{rn},n^h)$
for some $r\geq 0$, $h\geq 2$, and we also have $3\nu_3(M)\leq n-10$ from the definition of
$n$-ancestor list. We let
$$S\in\{\{1,2,3,4\},\{1,2,3,4,6\},\{1,2,3,4,5,7\},\{1,2,3,4,5,6,7\},\{1,2,3,4,5,6,7,8\}\}$$
such that $|S|=\frac{1}{n}\sum M_s$ and obtain the required decomposition of $K_n$ by combining an
$(M_s)$-decomposition of $\lan{S}$ (which exists by Lemma \ref{Sofsize3-8}), with an
$(\cms)$-decomposition of $K_n-\lan{S}$ (which exists by Lemma \ref{5sandHamsKn-S}). Note that the
condition $3\nu_3(M)\leq n-10$ implies that the required $(M_s)$-decomposition of $\lan{S}$ is not
among the listed possible exceptions in Lemma \ref{Sofsize3-8}. Note also that the condition $n
\geq 2\max(S)+1$ required in Lemmas \ref{Sofsize3-8} and \ref{5sandHamsKn-S} is easily seen to be
satisfied because $n \geq 15$ and $\sum M_s \leq n\lfloor\frac{n-1}{2}\rfloor$.

\noindent{\bf Case 4}\quad Suppose that $\sum M_s \in \{4n,5n,6n,7n,8n\}$ and $\nu_5(M)\leq 2$. In
this case we have $\nu_n(M_s)=0$ and $\cms=(3^{tn},4^{qn},n^h)$ for some $t,q\geq 0$, $h\geq 2$
(see Property (5) in Lemma \ref{MsLemma}), and $\sum M_s\neq 8n$ (see Property (1) in Lemma
\ref{MsLemma}). We let
$$S\in\{\{1,2,3,4\},\{1,2,3,4,6\},\{1,2,3,4,5,7\},\{1,2,3,4,5,6,7\}\}$$
such that $|S|=\frac{1}{n}{\sum M_s}$. If Lemma \ref{Sofsize3-8} gives us an $(M_s)$-decomposition
of $\lan{S}$ and Lemma \ref{3s4sandHamsforKn-S} gives us an $(\cms)$-decomposition of
$K_n-\lan{S}$, then we have the required decomposition of $K_n$.
The condition $n\geq 2\max(S)+1$ required in Lemmas \ref{Sofsize3-8} and \ref{3s4sandHamsforKn-S}
is satisfied because $n\geq 15$.
This leaves only the cases
arising from the possible exceptions in Lemma \ref{Sofsize3-8} and Lemma \ref{3s4sandHamsforKn-S},
and these are covered by Lemmas \ref{specialcase1} and \ref{specialcase2} respectively. \qed

\section{The case of no Hamilton cycles}\label{zerohamsection}

In this section we prove that Lemma \ref{MainLemmaforatmost1ham} holds in the case $\nu_n(M)=0$. In
this case, for $n \geq 15$, one of $\nu_3(M)$, $\nu_4(M)$ and $\nu_5(M)$ must be sizable, and the
proof splits into three cases accordingly. Each of these three cases splits into subcases according
to whether $n$ is even or odd. In each case we construct the required decomposition of $K_n$ from a
suitable decomposition of $K_{n-1}$ or $K_{n-2}$.

\subsection{Many $3$-cycles and no Hamilton cycles}

In Lemma \ref{nOddNoHamMany3s} we construct the required decompositions of complete graphs of odd
order and in Lemma \ref{nEvenNoHamMany3s} we construct the required decompositions of complete
graphs of even order.

\begin{lemma}\label{nOddNoHamMany3s}
If $n$ is odd, Theorem \ref{mainthm} holds for $K_{n-1}$, and $(M,3^{\frac{n-1}{2}})$ is an
$n$-ancestor list with $\nu_n(M) = 0$, then there is an $(M,3^{\frac{n-1}{2}})$-decomposition of
$K_n$.
\end{lemma}

\proof By Lemma \ref{SmallCases} we can assume that $n \geq 15$. Let $U$ be a vertex set with
$|U|=n-1$, let $\infty$ be a vertex not in $U$, and let $V=U\cup\{\infty\}$. Since
$(M,3^{\frac{n-1}{2}})$ is an $n$-ancestor list with $\nu_n(M) = 0$, it follows that $M$ is
$(n-1)$-admissible. Thus, by assumption there is an $(M)$-decomposition $\mathcal D$ of $K_U$. Let
$I$ be the perfect matching in $\mathcal{D}$. Then
$$\mathcal{D} \cup \mathcal{D}_1$$
is an $(M,3^{\frac{n-1}{2}})$-decomposition of $K_{V}$, where $\mathcal{D}_1$ is a
$(3^{\frac{n-1}2})$-decomposition of $K_{\{\infty\}} \vee I$. \qed

\begin{lemma}\label{nEvenNoHamMany3s}
If $n$ is even, Theorem \ref{mainthm} holds for $K_{n-1}$, and $(M,3^{\frac{n-2}{2}})$ is an
$n$-ancestor list with $\nu_n(M) = 0$, then there is an $(M,3^{\frac{n-2}{2}})$-decomposition of
$K_n$.
\end{lemma}

\proof By Lemma \ref{SmallCases} we can assume that $n \geq 16$. Let $U$ be a vertex set with
$|U|=n-1$, let $\infty$ be a vertex not in $U$, and let $V=U\cup\{\infty\}$. Since
$(M,3^{\frac{n-2}{2}})$ is an $n$-ancestor list with $\nu_n(M) = 0$, it follows that $(M,n-2)$ is
$(n-1)$-admissible and so by assumption there is an $(M,n-2)$-decomposition $\mathcal{D}$ of $K_U$.
Let $C$ be an $(n-2)$-cycle in $\mathcal{D}$, let $\{I,I_1\}$ be a decomposition of $C$ into two
matchings, and let $x$ be the vertex in $U \setminus V(C)$. Then
$$(\mathcal{D} \setminus \{C\}) \cup \{I + \infty x\} \cup \mathcal{D}_1$$
is an $(M,3^{\frac{n-2}{2}})$-decomposition of $K_{V}$, where $\mathcal{D}_1$ is a
$(3^{\frac{n-2}2})$-decomposition of $K_{\{\infty\}} \vee I_1$. \qed

\subsection{Many $4$-cycles and no Hamilton cycles}

\begin{lemma}\label{nOddNoHamMany4s}
If $n$ is odd, Theorem \ref{mainthm} holds for $K_{n-2}$, and $(M,4^{\frac{n+1}{2}})$ is an
$n$-ancestor list with $\nu_n(M) = 0$, then there is an $(M,4^{\frac{n+1}{2}})$-decomposition of
$K_n$.
\end{lemma}

\proof By Lemma \ref{SmallCases} we can assume that $n \geq 15$. Let $U$ be a vertex set with
$|U|=n-2$, let $\infty_1$ and $\infty_2$ be distinct vertices not in $U$, and let
$V=U\cup\{\infty_1,\infty_2\}$. Since $(M,4^{\frac{n+1}{2}})$ is an $n$-ancestor list with
$\nu_n(M)=0$, it follows from (5) in the definition of ancestor lists that any cycle length in $M$
is at most $n-3$. Thus, it is easily seen that $(M,5)$ is $(n-2)$-admissible and by assumption
there is an $(M,5)$-decomposition $\mathcal{D}$ of $K_{U}$.

Let $C$ be a
$5$-cycle in $\mathcal{D}$ and let $x$, $y$ and $z$ be vertices of $C$ such that $x$ and $y$ are
adjacent in $C$ and $z$ is not adjacent to either $x$ or $y$ in $C$. Then
$$(\mathcal{D} \setminus \{C\}) \cup \mathcal{D}_1 \cup \mathcal{D}_2$$
is an $(M,4^{\frac{n+1}{2}})$-decomposition of $K_V$, where
\begin{itemize}
    \item
$\mathcal{D}_1$ is a $(4^{\frac{n-5}2})$-decomposition of
$K_{\{\infty_1,\infty_2\},U \setminus \{x,y,z\}}$; and
    \item
$\mathcal{D}_2$ is a $(4^3)$-decomposition of $K_{\{\infty_1,\infty_2\},\{x,y,z\}} \cup
[\infty_1,\infty_2] \cup C$.
\end{itemize}
These decompositions are straightforward to construct.
\qed

\begin{lemma}\label{nEvenNoHamMany4s}
If $n$ is even, Theorem \ref{mainthm} holds for $K_{n-2}$, and $(M,4^{\frac{n-2}{2}})$ is an
$n$-ancestor list with $\nu_n(M) = 0$, then there is an $(M,4^{\frac{n-2}{2}})$-decomposition of
$K_n$.
\end{lemma}

\proof By Lemma \ref{SmallCases} we can assume that $n \geq 16$. Let $U$ be a vertex set with
$|U|=n-2$, let $\infty_1$ and $\infty_2$ be distinct vertices not in $U$, and let
$V=U\cup\{\infty_1,\infty_2\}$. Since $(M,4^{\frac{n-2}{2}})$ is an $n$-ancestor list with
$\nu_n(M)=0$, it follows from (5) in the definition of ancestor lists that any cycle length in $M$
is at most $n-3$. Thus, it is easily seen that $M$ is $(n-2)$-admissible and by assumption there is
an $(M)$-decomposition $\mathcal{D}$ of $K_{U}$. Let $I$ be the perfect matching in $\mathcal{D}$.
Then
$$(\mathcal{D} \setminus \{I\}) \cup \{I + \infty_1\infty_2\} \cup \mathcal{D}_1$$
is an $(M,4^{\frac{n-2}{2}})$-decomposition of $K_V$, where $\mathcal{D}_1$ is a
$(4^{\frac{n-2}2})$-decomposition of $K_{\{\infty_1,\infty_2\},U}$. \qed

\subsection{Many $5$-cycles and no Hamilton cycles}

We will make use of the following lemma in this subsection and in Subsection
\ref{many5sOneHamSection}.

\begin{lemma}\label{5CycleTrick}
If $G$ is a $3$-regular graph which contains a perfect matching and $\infty$ is a vertex not in
$V(G)$, then there is a decomposition of $K_{\{\infty\}} \vee G$ into $\frac{1}{2}|V(G)|$
$5$-cycles.
\end{lemma}

\proof Let $I$ be a perfect matching in $G$. Then $G-I$ is a $2$-regular graph on the vertex set
$V(G)$ and hence it can be given a coherent orientation $O$. Let
$$\mathcal{D}=\{(\infty,a,b,c,d): bc \in E(I) \mbox{ and } (b,a),(c,d) \in E(O)\}$$
be a set of (undirected) $5$-cycles. Because $O$ contains exactly one arc directed from each vertex
of $V(G)$, $|\mathcal{D}|=|E(I)|=\frac{1}{2}|V(G)|$ and each edge of $G$ appears in exactly one
cycle in $\mathcal{D}$. Further, because $O$ contains exactly one arc directed to each vertex of
$V(G)$, each edge of $K_{\{\infty\},V}$ appears in exactly one cycle in $\mathcal{D}$. Thus
$\mathcal{D}$ is a decomposition of $K_{\{\infty\}} \vee G$ into $\frac{1}{2}|V(G)|$ $5$-cycles.
\qed

\begin{lemma}\label{nOddNoHamMany5s}
If $n$ is odd, Theorem \ref{mainthm} holds for $K_{n-1}$, and $(M,5^{\frac{n-1}{2}})$ is an
$n$-ancestor list with $\nu_n(M) = 0$, then there is an $(M,5^{\frac{n-1}{2}})$-decomposition of
$K_n$.
\end{lemma}

\proof By Lemma \ref{SmallCases} we can assume that $n \geq 15$. Let $U$ be a vertex set with
$|U|=n-1$, let $\infty$ be a vertex not in $U$, and let $V=U\cup\{\infty\}$. Since the list
$(M,n-1)$ is easily seen to be $(n-1)$-admissible, by assumption there is an
$(M,n-1)$-decomposition $\mathcal{D}$ of $K_U$. Let $C$ be an $(n-1)$-cycle in $\mathcal{D}$ and
let $I$ be the perfect matching in $\mathcal{D}$. Then
$$(\mathcal{D} \setminus \{C,I\}) \cup \mathcal{D}_1$$
is an $(M,5^{\frac{n-1}{2}})$-decomposition of $K_V$, where $\mathcal{D}_1$ is a
$(5^{\frac{n-1}2})$-decomposition of $K_{\{\infty\}} \vee (C \cup I)$ (this exists by Lemma
\ref{5CycleTrick}). \qed

\begin{lemma}\label{nEvenNoHamMany5s}
If $n$ is even, Theorem \ref{mainthm} holds for $K_{n-2}$, and $(M,5^{n-2})$ is an $n$-ancestor
list with $\nu_n(M) = 0$, then there is an $(M,5^{n-2})$-decomposition of $K_n$.
\end{lemma}

\proof By Lemma \ref{SmallCases} we can assume that $n \geq 16$. Let $U$ be a vertex set with
$|U|=n-2$, let $\infty_1$ and $\infty_2$ be distinct vertices not in $U$, and let
$V=U\cup\{\infty_1,\infty_2\}$. Since $(M,5^{n-2})$ is an $n$-ancestor list with $\nu_n(M)=0$, it
follows from (6) in the definition of ancestor lists that any cycle length in $M$ is at most $n-5$.
Thus, it is easily seen that the list $(M,(n-2)^3)$ is $(n-2)$-admissible and by assumption there
is an $(M,(n-2)^3)$-decomposition $\mathcal{D}$ of $K_U$. Let $C_1$, $C_2$ and $C_3$ be distinct
$(n-2)$-cycles in $\mathcal{D}$ and let $I$ be the perfect matching in $\mathcal{D}$. Let
$\{I_1,I_2\}$ be a decomposition of $C_3$ into two perfect matchings. Then
$$(\mathcal{D} \setminus \{C_1,C_2,C_3,I\}) \cup \{I + \infty_1\infty_2\} \cup \mathcal{D}_1 \cup \mathcal{D}_2$$
is an $(M,5^{n-2})$-decomposition of $K_V$, where for $i=1,2$, $\mathcal{D}_i$ is a
$(5^{\frac{n-2}2})$-decomposition of $K_{\{\infty_i\}} \vee (C_i \cup I_i)$ (these exist by Lemma
\ref{5CycleTrick}). \qed

\subsection{Proof of Lemma \ref{MainLemmaforatmost1ham} in the case of no Hamilton cycles}

\begin{lemma}\label{Mainlemmaforzerohams}
If Theorem \ref{mainthm} holds for $K_{n-1}$ and $K_{n-2}$,
then there is an
$(M)$-decomposition of $K_n$ for each $n$-ancestor list $M$ satisfying $\nu_n(M)=0$.
\end{lemma}

\proof By Lemma \ref{SmallCases} we can assume that $n \geq 15$. If there is a cycle length in $M$
which is at least $6$ and at most $n-1$, then let $k$ be this cycle length. Otherwise let $k=0$. We
deal separately with the case $n$ is odd and the case $n$ is even.

\noindent{\bf Case 1}\quad Suppose that $n$ is odd. Since $n\geq 15$ and
$3\nu_3(M)+4\nu_4(M)+5\nu_5(M)+k =\frac{n(n-1)}{2}$, it can be seen that either $\nu_3(M) \geq
\frac{n-1}{2}$, $\nu_4(M) \geq \frac{n+1}{2}$ or $\nu_5(M) \geq \frac{n-1}{2}$. If $\nu_3(M) \geq
\frac{n-1}{2}$, then the result follows by Lemma \ref{nOddNoHamMany3s}. If $\nu_4(M) \geq
\frac{n+1}{2}$, then the result follows by Lemma \ref{nOddNoHamMany4s}. If $\nu_5(M) \geq
\frac{n-1}{2}$, then the result follows by Lemma \ref{nOddNoHamMany5s}.

\noindent{\bf Case 2}\quad Suppose that $n$ is even. Since $n\geq 16$,
$3\nu_3(M)+4\nu_4(M)+5\nu_5(M)+k =\frac{n(n-2)}{2}$ and $k \leq n-1$, it can be seen that either
$\nu_3(M) \geq \frac{n-2}{2}$, $\nu_4(M) \geq \frac{n-2}{2}$ or $\nu_5(M) \geq n-2$. (To see this
consider the cases $\nu_5(M) \geq 3$ and $\nu_5(M) \leq 2$ separately and use the definition of
$n$-ancestor list.) If $\nu_3(M) \geq \frac{n-2}{2}$, then the result follows by Lemma
\ref{nEvenNoHamMany3s}. If $\nu_4(M) \geq \frac{n-2}{2}$, then the result follows by Lemma
\ref{nEvenNoHamMany4s}. If $\nu_5(M) \geq n-2$, then the result follows by Lemma
\ref{nEvenNoHamMany5s}. \qed

\section{The case of exactly one Hamilton cycle}\label{onehamsection}

In this section we prove that Lemma \ref{MainLemmaforatmost1ham} holds in the case $\nu_n(M)=1$.
Again in this case, for $n \geq 15$, one of $\nu_3(M)$, $\nu_4(M)$ and $\nu_5(M)$ must be sizable,
and the proof splits into cases accordingly. The case in which $\nu_3(M)$ is sizable further splits
according to whether $\nu_4(M) \geq 1$, $\nu_5(M) \geq 1$, or $\nu_4(M) = \nu_5(M) = 0$. We first
require some preliminary definitions and results.

\subsection{Preliminaries}

Let $\mathcal{P}$ be an $(M)$-packing of $K_n$, let $\mathcal{P}'$ be an
$(M')$-packing of $K_n$ and let $S$ be a subset of $V(K_n$). We say that
$\mathcal{P}$ and $\mathcal{P}'$ are {\em equivalent on $S$} if we can write
$\{G \in \mathcal{P} : V(G) \cap S \neq \emptyset\} = \{G_1,\ldots,G_t\}$ and
$\{G \in \mathcal{P}' : V(G) \cap S \neq \emptyset\}  = \{G'_1,\ldots,G'_t\}$
such that
\begin{itemize}
  \item for $i \in \{1,\ldots,t\}$, $G_i$ is isomorphic to $G'_i$;
  \item for each $x \in S$ and for $i \in \{1,\ldots,t\}$, $x \in V(G_i)$ if
    and only if $x \in V(G'_i)$; and
  \item for all distinct $x,y \in S$ and for $i \in \{1,\ldots,t\}$, $xy \in
    E(G_i)$ if and only if $xy \in E(G'_i)$.
\end{itemize}


The following lemma is from \cite{BryHorLong}. It encapsulates a key \emph{edge swapping} technique which
was used in many of the proofs in \cite{BryHorAsym}, and which we shall make use of in this
section.

\begin{lemma}{\rm (\cite{BryHorLong}), Lemma 2.1)}\label{PathSwitch}
Let $n$ be a positive integer, let $M$ be a list of integers, let $\mathcal{P}$ be an $(M)$-packing
of $K_n$ with a leave, $L$ say, let $\alpha$ and $\beta$ be vertices of $L$, let $\pi$ be the
transposition $(\alpha\beta)$, and let $Z=Z(\mathcal{P},\alpha,\beta)=(\Nbd_L(\alpha)\cup
\Nbd_L(\beta)) \setminus ((\Nbd_L(\alpha)\cap \Nbd_L(\beta)) \cup \{\alpha,\beta\})$. Then there
exists a partition of the set $Z$ into pairs such that for each pair $\{u,v\}$ of the partition,
there exists an $(M)$-packing of $K_n$, $\mathcal{P}'$ say, with a leave, $L'$ say, which differs
from $L$ only in that $\alpha u$, $\alpha v$, $\beta u$ and $\beta v$ are edges in $L'$ if and only
if they are not edges in $L$. Furthermore, if $\mathcal{P} = \{C_1,\ldots,C_t\}$ ($n$ odd) or
$\mathcal{P} = \{I,C_1,\ldots,C_t\}$ ($n$ even) where $C_1,\ldots,C_t$ are cycles and $I$ is a
perfect matching, then $\mathcal{P}' = \{C'_1,\ldots,C'_t\}$ ($n$ odd) or $\mathcal{P}' =
\{I',C'_1,\ldots,C'_t\}$ ($n$ even) where for $i=1,\ldots,t$, $C'_i$ is a cycle of the same length
as $C_i$ and $I'$ is a perfect matching such that
\begin{itemize}
    \item
either $I'=I$ or $I'=\pi(I)$;
    \item
for $i = 1,\ldots,t$ if neither $\alpha$ nor $\beta$ is in $V(C_i)$ then $C'_i=C_i$;
    \item
for $i = 1,\ldots,t$ if exactly one of $\alpha$ and $\beta$ is in $V(C_i)$ then either
$C'_i=C_i$ or $C'_i=\pi(C_i)$; and
    \item
for $i = 1,\ldots,t$ if both $\alpha$ and $\beta$ are in $V(C_i)$ then $C'_i \in
\{C_i,\pi(C_i),\pi(P_i)\cup P^{\dag}_i,P_i \cup \pi(P^{\dag}_i)\}$ where $P_i$ and $P^{\dag}_i$
are the two paths in $C_i$ which have endpoints $\alpha$ and $\beta$.
\end{itemize}
\end{lemma}
We say that $\mathcal{P}'$ is the $(M)$-packing obtained from $\mathcal{P}$ by performing the
$(\alpha,\beta)$-switch with origin $u$ and terminus $v$ (we could equally call $v$ the origin and
$u$ the terminus). For our purposes
here,
it is important to note that $\mathcal{P}'$ is equivalent to $\mathcal{P}$ on $V(L) \setminus
\{\alpha,\beta\}$.

We will also make use of three lemmas from \cite{BryHorAsym}. The original version of Lemma
\ref{EquitableOneStep} (Lemma 2.15 in \cite{BryHorAsym}) does not include the claim that $\mathcal{P}'$ is equivalent to
$\mathcal{P}$ on $V(L) \setminus \{a,b\}$, this follows directly from the fact that the proof uses
only $(a,b)$-switches.

\begin{lemma}\label{EquitableOneStep}
Let $n$ be a positive integer and let $M$ be a list of integers. Suppose that there exists an
$(M)$-packing $\mathcal{P}$ of $K_n$ with a leave $L$ which contains two vertices $a$ and $b$ such
that $\deg_L(a)+2 \leq \deg_L(b)$. Then there exists an $(M)$-packing $\mathcal{P}'$ of $K_n$,
which is equivalent to $\mathcal{P}$ on $V(L) \setminus \{a,b\}$, and which has a leave $L'$ such
that $\deg_{L'}(a)=\deg_L(a)+2$, $\deg_{L'}(b)=\deg_L(b)-2$ and $\deg_{L'}(x)=\deg_L(x)$ for all $x
\in V(L) \setminus \{a,b\}$. Furthermore,
\begin{itemize}
    \item [(i)]
if $a$ and $b$ are adjacent in $L$, then $L'$ has the same number of non-trivial components as
$L$;
    \item [(ii)]
if $\deg_L(a)=0$ and $b$ is not a cut-vertex of $L$, then $L'$ has the same number of
non-trivial components as $L$; and
    \item [(iii)]
if $\deg_L(a)=0$, then either $L'$ has the same number of non-trivial components as $L$, or
$L'$ has one more non-trivial component than $L$.
\end{itemize}
\end{lemma}

Similarly, the original versions of Lemmas \ref{OneDeg4} and \ref{TwoDeg4}
(Lemmas 2.14 and 2.11 respectively in \cite{BryHorAsym}) did not include the
claims that the final decompositions are equivalent to the initial packings on $V\setminus U$.
However, these claims can be seen to hold as the proofs of the lemmas given in \cite{BryHorAsym}
require switching only on vertices of positive degree in the leave, with one exception which we
discuss shortly. The lemmas below each contain the additional hypothesis that $\deg_{L}(x)=0$ for
all $x \in V \setminus U$, and this ensures all the switches are on vertices of $U$ and hence that
the final decomposition is equivalent to the initial packing on $V\setminus U$.

The exception mentioned above occurs in the proof of the original version of Lemma \ref{TwoDeg4}
where a switch on a vertex of degree $0$ in the leave is required when $3\in\{m_1,m_2\}$. We can
ensure this switch is on a vertex in $U$ because we have the additional hypothesis that
$\deg_{L}(x)=0$ for some $x \in U$ when $3\in\{m_1,m_2\}$. This additional hypothesis also allows
us to omit the hypothesis, included in the original version of Lemma \ref{TwoDeg4}, that the size
of the leave be at most $n+1$, because in the proof this was used only to ensure the existence of a
vertex of degree $0$ in the leave when $3\in\{m_1,m_2\}$. Thus the modified versions stated below
hold by the proofs presented in \cite{BryHorAsym}.

\begin{lemma}\label{OneDeg4}
Let $V$ be a vertex set and let $U$ be a subset of $V$. Let $M$ be a list of integers and let $k$,
$m_1$ and $m_2$ be positive integers such that $m_1,m_2 \geq \max(\{3,k+1\})$. Suppose that there
exists an $(M)$-packing $\mathcal{P}$ of $K_V$ with a leave $L$ of size $m_1+m_2$ such that
$\Delta(L)=4$, exactly one vertex of $L$ has degree $4$, $L$ has exactly $k$ non-trivial
components, $L$ does not have a decomposition into two odd cycles if $m_1$ and $m_2$ are both even,
and $\deg_{L}(x)=0$ for all $x \in V \setminus U$. Then there exists an $(M,m_1,m_2)$-decomposition
of $K_V$ which is equivalent to $\mathcal{P}$ on $V \setminus U$.
\end{lemma}

\begin{lemma}\label{TwoDeg4}
Let $V$ be a vertex set and let $U$ be a subset of $V$. Let $M$ be a list of integers and let $m_1$
and $m_2$ be integers such that $m_1,m_2 \geq 3$. Suppose that there exists $(M)$-packing
$\mathcal{P}$ of $K_V$ with a leave $L$ of size $m_1+m_2$ such that $\Delta(L)=4$, exactly two
vertices of $L$ have degree $4$, $L$ has exactly one non-trivial component, $\deg_{L}(x)=0$ for all
$x \in V \setminus U$, and $\deg_{L}(x)=0$ for some $x \in U$ if $3\in\{m_1,m_2\}$. Then there
exists an $(M,m_1,m_2)$-decomposition of $K_V$ which is equivalent to $\mathcal{P}$ on $V \setminus
U$.
\end{lemma}

We also require Lemma \ref{SmallishCases}, which deals with some small order cases.

\begin{lemma}\label{SmallishCases}
Let $n$ be an integer such that $n \in \{15,16,17,18,19,20,22,24,26\}$ and let $M$ be an
$n$-ancestor list such that $\nu_n(M) = 1$, $\nu_5(M) \geq 3$, and $\nu_4(M) \geq 2$ if $n=24$.
Then there is an $(M)$-decomposition of $K_n$.
\end{lemma}

\proof If there is a cycle length in $M$ which is at least $6$ and at most $n-1$, then let $k$ be
this cycle length. Otherwise let $k=0$. Note that $3\nu_3(M)+4\nu_4(M)+5\nu_5(M) + k + n =
n\lfloor\frac{n-1}{2}\rfloor$ and that, because $M$ is an $n$-ancestor list with $\nu_5(M) \geq 3$,
it follows that $3\nu_3(M) \leq n-10$, $2\nu_4(M) \leq n-6$, and $k \leq n-5$.

Using this, it is routine to check that if $n=15$ then $M$ must be one of 12 possible lists and if
$n=16$ then $M$ must be one of 26 possible lists. In each of these cases we have constructed an
$(M)$-decomposition of $K_n$ by computer search.


If $n \in \{17,18,19,20,22,24,26\}$, then we partition $\{1,\ldots,\lfloor\frac{n}{2}\rfloor\}$
into sets $S_1$, $S_2$ and $S_3$ according to the following table.
$$
\begin{array}{|c|c|c|c|}
\hline
n&S_1&S_2&S_3 \\
\hline
17&\{1,2,3,4,5,6,7\}&\emptyset&\{8\}\\
18&\{1,2,3,4,5,6,7\}&\emptyset&\{8,9\}\\
19&\{1,2,3,4,5,6,7,8\}&\emptyset&\{9\}\\
20&\{1,2,3,4,5,6,7,8\}&\emptyset&\{9,10\}\\
22&\{1,2,3,4\}&\{6,7,8,9,10\}&\{5,11\}\\
24&\{1,2,3\}&\{4,6,7,8,11\}&\{5,9,10,12\}\\
26&\{1,2,3,4,5,7\}&\{6,8,9,10,11\}&\{12,13\}\\
\hline
\end{array}
$$

Using $3\nu_3(M) \leq n-10$, $2\nu_4(M) \leq n-6$ and $k \leq n-5$, it is routine to check that
$\nu_5(M) \geq n$ when $n \in \{22,26\}$, and that $\nu_5(M) \geq n+8$ when $n=24$. By Lemma
\ref{Sofsize3-8}, there is an $(M')$-decomposition of $\langle S_1 \rangle_n$, where $M=(M',n)$
when $n \in \{17,18,19,20\}$, $M=(M',5^n,n)$ when $n \in \{22,26\}$, and $M=(M',4^2,5^{n+8},n)$
when $n=24$. For $n\in\{22,24,26\}$, it is easy to see that $S_2$ is a modulo $n$ difference
$5$-tuple, and so there is a $(5^n)$-decomposition of $\langle S_2 \rangle_n$. For
$n\in\{17,18,19,20,22,26\}$, there is an $(n)$-decomposition of the graph $\langle S_3 \rangle_n$,
as it is either an $n$-cycle or a connected $3$-regular Cayley graph on a cyclic group, and the
latter are well known to contain a Hamilton cycle, see \cite{CheQui}. For $n=24$,
$\lan{\{5,9,10\}}\cong\lan{\{1,2,3\}}$ (with $x \mapsto 5x$ being an isomorphism). Thus, by Lemma
\ref{S=123includingham} there is a $(4^2,5^8,24)$-decomposition of $\langle \{5,9,10,12\}
\rangle_{24}$ (as $\langle \{12\} \rangle_{24}$ is a perfect matching). Combining these
decompositions of $\lan{S_1}$, $\lan{S_2}$ and $\lan{S_3}$ gives us the required
$(M)$-decomposition of $K_n$. \qed

\subsection{Many $3$-cycles, one Hamilton cycle, and at least one $4$- or $5$-cycle}

In Lemmas \ref{nOddOneHamMany3sAtLeastOne4} and \ref{nOddOneHamMany3sAtLeastOne5} we construct the
required decompositions of complete graphs of odd order in the cases where the decomposition
contains at least one $4$-cycle or at least one $5$-cycle, respectively. In Lemmas
\ref{nEvenOneHamMany3sAtLeastOne4} and \ref{nEvenOneHamMany3sAtLeastOne5} we construct the required
decompositions of complete graphs of even order in the cases where the decomposition contains at
least one $4$-cycle or at least one $5$-cycle, respectively. These results are proved by
constructing the required decomposition of $K_n$ from a suitable decomposition of $K_{n-1}$,
$K_{n-2}$ or $K_{n-3}$.

\begin{lemma}\label{nOddOneHamMany3sAtLeastOne4}
If $n$ is odd, Theorem \ref{mainthm} holds for $K_{n-1}$, and $(M,3^{\frac{n-5}{2}},4,n)$ is an
$n$-ancestor list with $\nu_n(M) = 0$, then there is an $(M,3^{\frac{n-5}{2}},4,n)$-decomposition
of $K_n$.
\end{lemma}

\proof By Lemma \ref{SmallCases} we can assume that $n \geq 15$. Let $U$ be a vertex set with
$|U|=n-1$, let $\infty$ be a vertex not in $U$, and let $V=U\cup\{\infty\}$. Since
$(M,3^{\frac{n-5}{2}},4,n)$ is an $n$-ancestor list with $\nu_n(M) = 0$, it follows that $(M,n-2)$
is $(n-1)$-admissible and by assumption there is an $(M,n-2)$-decomposition $\mathcal{D}$ of $K_U$.

Let $H$ be an $(n-2)$-cycle in $\mathcal{D}$, let $I$ be the perfect matching in $\mathcal{D}$, and
let $[w,x,y,z]$ be a path in $I \cup H$ such that $w \notin V(H)$, $wx,yz \in E(I)$ and $xy \in
E(H)$. Then
$$(\mathcal{D} \setminus \{I,H\}) \cup \{H',(\infty,x,y,z)\} \cup \mathcal{D}_1$$
is an $(M,3^{\frac{n-5}{2}},4,n)$-decomposition of $K_{V}$, where
\begin{itemize}
    \item
$H'=(H-[x,y]) \cup [x,w,\infty,y]$; and
    \item
$\mathcal{D}_1$ is a $(3^{\frac{n-5}2})$-decomposition of $K_{\{\infty\},U \setminus
\{w,x,y,z\}} \cup (I - \{wx,yz\})$. \qed
\end{itemize}

\begin{lemma}\label{nOddOneHamMany3sAtLeastOne5}
If $n$ is odd, Theorem \ref{mainthm} holds for $K_{n-1}$, and $(M,3^{\frac{n-5}{2}},5,n)$ is an
$n$-ancestor list with $\nu_n(M) = 0$, then there is an $(M,3^{\frac{n-5}{2}},5,n)$-decomposition
of $K_n$.
\end{lemma}

\proof By Lemma \ref{SmallCases} we can assume that $n \geq 15$. Let $U$ be a vertex set with
$|U|=n-1$, let $\infty$ be a vertex not in $U$, and let $V=U\cup\{\infty\}$. Since
$(M,3^{\frac{n-5}{2}},5,n)$ is an $n$-ancestor list with $\nu_n(M) = 0$, it follows that $(M,n-1)$
is $(n-1)$-admissible and so by assumption there is an $(M,n-1)$-decomposition $\mathcal{D}$ of
$K_U$.

Let $H$ be an $(n-1)$-cycle in $\mathcal{D}$, let $I$ be the perfect matching in $\mathcal{D}$, and
let $[w,x,y,z]$ be a path in $I \cup H$ such that $wx,yz \in E(I)$ and $xy \in E(H)$. Then
$$(\mathcal{D} \setminus \{I,H\}) \cup \{H',(\infty,w,x,y,z)\} \cup \mathcal{D}_1$$
is an $(M,3^{\frac{n-5}{2}},5,n)$-decomposition of $K_{V}$, where
\begin{itemize}
    \item
$H'=(H-[x,y]) \cup [x,\infty,y]$; and
    \item
$\mathcal{D}_1$ is a $(3^{\frac{n-5}2})$-decomposition of $K_{\{\infty\},U \setminus
\{w,x,y,z\}} \cup (I-\{wx,yz\})$. \qed
\end{itemize}

\begin{lemma}\label{nEvenOneHamMany3sAtLeastOne4}
If $n$ is even, Theorem \ref{mainthm} holds for $K_{n-3}$, and $(M,3^{\frac{3n-14}{2}},4,n)$ is an
$n$-ancestor list with $\nu_n(M) = 0$, then there is an $(M,3^{\frac{3n-14}{2}},4,n)$-decomposition
of $K_n$.
\end{lemma}

\proof By Lemma \ref{SmallCases} we can assume that $n \geq 16$. Let $U$ be a vertex set with
$|U|=n-3$, let $\infty_1$, $\infty_2$ and $\infty_3$ be distinct vertices not in $U$, and let
$V=U\cup\{\infty_1,\infty_2,\infty_3\}$. Since $(M,3^{\frac{3n-14}{2}},4,n)$ is an $n$-ancestor
list with $\nu_n(M) = 0$, it follows from (5) in the definition of ancestor lists that any cycle
length in $M$ is at most $n-3$. Thus, it is easily seen that $(M,(n-4)^2,n-3)$ is
$(n-3)$-admissible and so by assumption there is an $(M,(n-4)^2,n-3)$-decomposition $\mathcal{D}$
of $K_U$.

Let $C_1$ and $C_2$ be distinct $(n-4)$-cycles in $\mathcal{D}$,
let $H$ be an $(n-3)$-cycle in $\mathcal{D}$,\
let $\{I,I_1\}$ be a decomposition of $C_1$ into two matchings,
let $\{I_2,I_3\}$ be a decomposition of $C_2$ into two matchings,
let $w$ be the vertex in $U\setminus V(C_1)$,
and let $[x,y,z]$ be a path in $H \cup I_3$ such that $x \notin V(C_2)$, $xy \in
E(H)$ and $yz \in E(I_3)$ (possibly $w \in \{x,y,z\}$).
Then
$$(\mathcal{D} \setminus \{H,C_1,C_2\}) \cup \{I + \{\infty_1w,\infty_2\infty_3\},H',(\infty_3,x,y,z)\} \cup \mathcal{D}_1 \cup \mathcal{D}_2 \cup \mathcal{D}_3$$
is an $(M,3^{\frac{3n-14}{2}},4,n)$-decomposition of $K_{V}$, where
\begin{itemize}
    \item
$H'=(H-[x,y]) \cup [x,\infty_2,\infty_1,\infty_3,y]$;
    \item
for $i=1,2$, $\mathcal{D}_i$ is a $(3^{\frac{n-4}2})$-decomposition of $K_{\{\infty_i\}} \vee
I_i$; and
    \item
$\mathcal{D}_3$ is a $(3^{\frac{n-6}2})$-decomposition of $K_{\{\infty_3\},U \setminus
\{x,y,z\}} \cup (I_3-yz)$. \qed
\end{itemize}

\begin{lemma}\label{nEvenOneHamMany3sAtLeastOne5}
If $n$ is even, Theorem \ref{mainthm} holds for $K_{n-1}$, and $(M,3^{\frac{n-6}{2}},5,n)$ is an
$n$-ancestor list with $\nu_n(M) = 0$, then there is an $(M,3^{\frac{n-6}{2}},5,n)$-decomposition
of $K_n$.
\end{lemma}

\proof By Lemma \ref{SmallCases} we can assume that $n \geq 16$. Let $U$ be a vertex set with
$|U|=n-1$, let $\infty$ be a vertex not in $U$, and let $V=U\cup\{\infty\}$. Since
$(M,3^{\frac{n-6}{2}},5,n)$ is an $n$-ancestor list with $\nu_n(M) = 0$, it follows that
$(M,n-2,n-1)$ is $(n-1)$-admissible and so by assumption there is an $(M,n-2,n-1)$-decomposition
$\mathcal{D}$ of $K_U$.

Let $H$ be an $(n-1)$-cycle in $\mathcal{D}$, let $C$ be an $(n-2)$-cycle in $\mathcal{D}$, let
$\{I,I_1\}$ be a decomposition of $C$ into two matchings, let $[w,x,y,z]$ be a path in $I_1 \cup H$
such that $wx,yz \in E(I_1)$ and $xy \in E(H)$, and let $v$ be the vertex in $U \setminus V(C)$.
Then
$$(\mathcal{D} \setminus \{C,H\}) \cup \{I + v\infty,
H', (\infty,w,x,y,z)\} \cup \mathcal{D}_1$$ is an $(M,3^{\frac{n-6}{2}},5,n)$-decomposition
of $K_{V}$, where
\begin{itemize}
    \item
$H'=(H-[x,y]) \cup [x,\infty,y]$; and
    \item
$\mathcal{D}_1$ is a $(3^{\frac{n-6}2})$-decomposition of $K_{\{\infty\},U \setminus
\{v,w,x,y,z\}} \cup (I-\{wx,yz\})$. \qed
\end{itemize}

\subsection{Many $3$-cycles, one Hamilton cycle and no $4$- or $5$-cycles}

We show that the required decompositions exist in Lemma \ref{OneHamNo4sOr5s}. We first require
three preliminary lemmas. These results are proved using the edge swapping techniques mentioned
previously.

\begin{lemma}\label{33kChain}
Let $n$ and $k$ be positive integers, and let $M$ be a list of integers. If there exists an
$(M)$-packing of $K_n$ whose leave has a decomposition into two $3$-cycles, $T_1$ and $T_2$, and a
$k$-cycle $C$ such that $|V(T_1) \cap V(T_2)| = 1$, $|V(T_2) \cap V(C)| = 1$ and $V(T_1) \cap V(C)
= \emptyset$, then there exists an $(M,3,k+3)$-decomposition of $K_n$.
\end{lemma}

\proof Let $\mathcal{P}$ be an $(M)$-packing of $K_n$ which satisfies the conditions of the lemma
and let $L$ be its leave. Let $[w,x,y,z]$ be a path in $L$ such that $w \in V(T_1) \setminus
V(T_2)$, $V(T_1) \cap V(T_2) = \{x\}$, $V(T_2) \cap V(C) = \{y\}$, and $z \in V(C) \setminus
V(T_2)$. Let $\mathcal{P}'$ be the $(M)$-packing of $K_n$ obtained from $\mathcal{P}$ by performing
the $(w,z)$-switch $S$ with origin $x$. If the terminus of $S$ is $y$, then the leave of
$\mathcal{P}'$ has a decomposition into the $3$-cycle $T_2$ and a $(k+3)$-cycle. Otherwise the
terminus of $S$ is not $y$ and the leave of $\mathcal{P}'$ has a decomposition into the $3$-cycle
$(x,y,z)$ and a $(k+3)$-cycle. In either case we complete the proof by adding these cycles to
$\mathcal{P}'$. \qed

\begin{lemma}\label{Two3sAndAk}
Let $n$ and $k$ be positive integers such that $k \leq n-4$, and let $M$ be a list of integers.
Suppose that there exists an $(M)$-packing of $K_n$ whose leave has a decomposition into two
$3$-cycles, $T_1$ and $T_2$, and a $k$-cycle $C$ such that $|V(T_1) \cap V(C)| \leq 1$, $|V(T_2)
\cap V(C)| \leq 1$, and $|V(T_1) \cap V(T_2)| = 1$ if $k = n-4$. Then there exists an
$(M,3,k+3)$-decomposition of $K_n$.
\end{lemma}

\proof Let $\mathcal{P}$ be an $(M)$-packing of $K_n$ which satisfies the conditions of the lemma
and let $L$ be its leave.

\noindent{\bf Case 1}\quad Suppose that $\Delta(L)=2$. Then $T_1$, $T_2$ and $C$ are pairwise
vertex disjoint. Let $x \in V(T_1)$ and $y \in V(C)$ and let $z$ be a neighbour in $T_1$ of $x$.
Let $\mathcal{P}'$ be the $(M)$-packing of $K_n$ obtained from $\mathcal{P}$ by performing the
$(x,y)$-switch $S$ with origin $z$, and let $L'$ be the leave of $\mathcal{P}'$. Then either the
non-trivial components of $L'$ are a $3$-cycle and a $(k+3)$-cycle or $\Delta(L')=4$, exactly one
vertex of $L'$ has degree $4$, and $L'$ has exactly two nontrivial components. In the former case
we can add these cycles to $\mathcal{P}'$ to complete the proof. In the latter case we can apply
Lemma \ref{OneDeg4} to complete the proof.

\noindent{\bf Case 2}\quad Suppose that $\Delta(L)=4$. If exactly one vertex of $L$ has degree $4$,
then $L$ has exactly two nontrivial components and we can complete the proof by applying Lemma
\ref{OneDeg4}. Thus we can assume that $L$ has at least two vertices of degree $4$. We can further
assume that $\mathcal{P}$ does not satisfy the conditions of Lemma \ref{33kChain}, for otherwise we
can complete the proof by applying Lemma \ref{33kChain}. Noting that $|V(T_1) \cap V(T_2)| \in
\{0,1\}$, it follows that $|V(T_1) \cap V(C)| = 1$ and $|V(T_2) \cap V(C)| = 1$. Because $k \leq
n-4$ and $|V(T_1) \cap V(T_2)| = 1$ if $k = n-4$, there is an isolated vertex $z$ in $L$. Let $w$
be the vertex in $V(T_1) \cap V(C)$, let $x$ and $y$ be the neighbours in $T_1$ of $w$. Let
$\mathcal{P}'$ be the $(M)$-packing of $K_n$ obtained from $\mathcal{P}$ by performing the
$(w,z)$-switch $S$ with origin $x$, and let $L'$ be the leave of $\mathcal{P}'$. If the terminus of
$S$ is not $y$, then $L'$ has a decomposition into a $(k+3)$-cycle and a $3$-cycle, and we complete
the proof by adding these cycles to $\mathcal{P}'$. Otherwise the terminus of $S$ is $y$ and either
$\Delta(L')=4$, exactly one vertex of $L'$ has degree $4$, and $L'$ has exactly two nontrivial
components (this occurs when $|V(T_1) \cap V(T_2)| = 0$) or $\mathcal{P}'$ satisfies the conditions
of Lemma \ref{33kChain} (this occurs when $|V(T_1) \cap V(T_2)| = 1$). Thus we can complete the
proof by applying Lemma \ref{OneDeg4} or Lemma \ref{33kChain}.

\noindent{\bf Case 3}\quad  Suppose that $\Delta(L) \geq 6$. In this case, exactly one vertex of
$L$ has degree $6$ and every other vertex of $L$ has degree at most $2$. Let $w$ be the vertex of
degree $6$ in $L$, let $x$ be a neighbour in $T_2$ of $w$, and let $y$ and $z$ be the neighbours in
$T_1$ of $w$. Let $\mathcal{P}'$ be the $(M)$-packing of $K_n$ obtained from $\mathcal{P}$ by
performing the $(w,x)$-switch $S$ with origin $y$, and let $L'$ be the leave of $\mathcal{P}'$. If
the terminus of $S$ is $z$, then $\mathcal{P}'$ satisfies the conditions of Lemma \ref{33kChain}
and we complete the proof by applying Lemma \ref{33kChain}. Otherwise the terminus of $S$ is not
$z$ and $L'$ has a decomposition into the $3$-cycle $T_2$ and a $(k+3)$-cycle, and we complete the
proof by adding these cycles to $\mathcal{P}'$. \qed

\begin{lemma}\label{Add3Tok}
Let $n$, $k$ and $t$ be positive integers such that $3 \leq k \leq n-4$. If there exists a
$(3^t,k,n)$-decomposition of $K_n$, then there exists a $(3^{t-1},k+3,n)$-decomposition of $K_n$.
\end{lemma}

\proof By Lemma \ref{SmallCases} we can assume that $n \geq 15$. Let $V$ be a vertex set with
$|V|=n$, let $\mathcal{D}$ be a $(3^t,k,n)$-decomposition of $K_V$, and let $C$ be a $k$-cycle in
$\mathcal{D}$. Let $U = V \setminus V(C)$. The $n$-cycle in $\mathcal D$ contains at most $|U|-1$
edges of $K_U$ (as the subgraph of the $n$-cycle induced by $U$ is a forest). Also, if $n$ is even,
then the perfect matching in $\mathcal{D}$ contains at most $\lfloor\frac{1}{2}|U|\rfloor$ edges of
$K_U$. The proof now splits into two cases depending on whether $k=n-4$.

\noindent{\bf Case 1}\quad Suppose that $k \leq n-5$. Then $|U| \geq 5$ and, by the comments in the
preceding paragraph, the $3$-cycles in $\mathcal{D}$ contain at least four edges of $K_U$. Thus
there are distinct $3$-cycles $T_1,T_2\in\mathcal{D}$ such that each contains at least one edge of
$K_U$. We can remove $C$, $T_1$ and $T_2$ from $\mathcal{D}$ and apply Lemma \ref{Two3sAndAk} to
the resulting packing to complete the proof.

\noindent{\bf Case 2}\quad Suppose that $k = n-4$. Then $|U| = 4$, the $n$-cycle in $\mathcal D$
contains at most three edges of $K_U$ and the perfect matching in $\mathcal{D}$ contains at most
two edges of $K_U$. This leaves at least one edge of $K_U$ which occurs in a $3$-cycle
$T_1\in\mathcal{D}$. Let $U=\{u_1,u_2,u_3,u_4\}$ and let $H$ be the $n$-cycle in $\mathcal{D}$.
This case now splits into two subcases depending on whether $V(T_1) \cap V(C) = \emptyset$.

\noindent{\bf Case 2a}\quad Suppose that $V(T_1) \cap V(C) = \emptyset$. Then we can assume without
loss of generality that $T_1=(u_1,u_2,u_3)$. If any of the three edges $u_1u_4,u_2u_4,u_3u_4$ is in
a $3$-cycle $T_2\in\mathcal{D}$, then we can remove $C$, $T_1$ and $T_2$ from $\mathcal{D}$ and
apply Lemma \ref{Two3sAndAk} to the resulting packing to complete the proof. Thus, we assume there
is no such $3$-cycle in $\mathcal D$. Without loss of generality, it follows that $n$ is even, that
$u_1u_4$ is an edge of the perfect matching in $\mathcal{D}$, and that $u_2u_4,u_3u_4\in E(H)$. Let
$z$ be a vertex in $C$ which is not adjacent in $H$ to a vertex in $U$ (such a vertex exists as
$n\geq 15$ implies $|V(C)|\geq 11$ and there are only at most four vertices of $C$ which are
adjacent in $H$ to vertices in $U$).

Now let $\mathcal{P}'$ be the $(3^{t-1},n)$-packing of $K_V$ obtained from $\mathcal{D} \setminus
\{C,T_1\}$ by performing the $(u_1,z)$-switch $S$ with origin $u_2$. If the terminus of $S$ is not
$u_3$, then the only non-trivial component in the leave of $\mathcal{P}'$ is a $(k+3)$-cycle and we
can complete the proof by adding this cycle to $\mathcal{P}'$. Otherwise, the terminus of $S$ is
$u_3$ and the only non-trivial component in the leave of $\mathcal{P}'$ is $(u_2,u_3,z) \cup C$.
Furthermore, since $z$ is not adjacent in $H$ to a vertex in $U$, the final dot point in Lemma
\ref{PathSwitch} guarantees that neither $u_1u_2$ nor $u_1u_3$ is an edge of the $n$-cycle in
$\mathcal{P}'$. Since $u_1u_2$ and $u_1u_3$ cannot both be edges of the perfect matching in
$\mathcal{P}'$, this means that one of them must be in a $3$-cycle $T'_2\in\mathcal{P}'$. Thus, we
can remove $T'_2$ from $\mathcal{P}'$ and apply Lemma \ref{Two3sAndAk} to the resulting packing to
complete the proof.

\noindent{\bf Case 2b}\quad Suppose that $|V(T_1) \cap V(C)| = 1$. Let $T_1=(x,y,z)$ with $x\in
V(C)$ and $y,z\in U$, and let $w\in U\setminus\{y,z\}$. Let $\mathcal{P}'$ be the
$(3^{t-1},n)$-packing of $K_V$ obtained from $\mathcal{D} \setminus \{C,T_1\}$ by performing the
$(w,x)$-switch $S$ with origin $y$, and let $L'$ be the leave of $\mathcal{P}'$. If the terminus of
$S$ is not $z$, then the only non-trivial component in the leave of $\mathcal{P}'$ is a
$(k+3)$-cycle and we can complete the proof by adding this cycle to $\mathcal{P}'$. Otherwise the
terminus of $S$ is $z$, and the only non-trivial components in the leave of $\mathcal{P}'$ are $C$
and $(w,y,z)$. By adding these cycles to $\mathcal{P}'$ we obtain a $(3^t,k,n)$-decomposition of
$K_V$ which contains a $3$-cycle and an $(n-4)$-cycle which are vertex disjoint, and we can proceed
as we did in Case 2a. \qed

\vspace{0.5cm}

We are now ready to prove the main result of this subsection.

\begin{lemma}\label{OneHamNo4sOr5s}
If $n$, $k$ and $t$ are positive integers such that $3 \leq k \leq n-1$, Theorem \ref{mainthm}
holds for $K_{n-3}$, $K_{n-2}$ and $K_{n-1}$, and $(3^t,k,n)$ is an $n$-ancestor list, then there
is a $(3^t,k,n)$-decomposition of $K_n$.
\end{lemma}

\proof By Lemma \ref{SmallCases} we can assume that $n \geq 15$. Let $r\in\{3,4,5\}$ such that $r
\equiv k \md 3$. It suffices to find a $(3^{t+\frac{k-r}{3}},r,n)$-decomposition of $K_n$, since we
can then obtain a $(3^t,k,n)$-decomposition of $K_n$ by repeatedly applying Lemma \ref{Add3Tok}
($\frac{k-r}3$ times). If $r=3$ then the existence of a $(3^{t+\frac{k-r}{3}},r,n)$-decomposition
of $K_n$ follows from the main result of \cite{BryMae}, so we may assume $r \in \{4,5\}$. Thus, the
existence of the required $(3^{t+\frac{k-r}{3}},r,n)$-decomposition of $K_n$ is given by one of
Lemmas \ref{nOddOneHamMany3sAtLeastOne4}, \ref{nOddOneHamMany3sAtLeastOne5},
\ref{nEvenOneHamMany3sAtLeastOne4} and \ref{nEvenOneHamMany3sAtLeastOne5}, provided that
$t+\frac{k-r}{3}\geq\frac{3n-14}2$ (the number of $3$-cycles in the decompositions given by Lemma
\ref{nEvenOneHamMany3sAtLeastOne4} is at least $\frac{3n-14}2$ and the number is smaller for the
other three lemmas for $n\geq 15$). However, it follows from $3t+k+n=n\lfloor\frac{n-1}{2}\rfloor$,
$k \leq n-1$ and $n \geq 15$ that $t \geq \frac{3n-14}{2}$. \qed

\subsection{Many $4$-cycles and one Hamilton cycle}

In Lemma \ref{nOddOneHamMany4s} we construct the required decompositions of complete graphs of odd
order and in Lemma \ref{nEvenOneHamMany4s} we construct the required decompositions of complete
graphs of even order. In each case we construct the required decomposition of $K_n$ from a suitable
decomposition of $K_{n-2}$.

\begin{lemma}\label{nOddOneHamMany4s}
If $n$ is odd, Theorem \ref{mainthm} holds for $K_{n-2}$, and $(M,4^{\frac{n-3}{2}},n)$ is an
$n$-ancestor list with $\nu_n(M) = 0$, then there is an $(M ,4^{\frac{n-3}{2}},n)$-decomposition of
$K_n$.
\end{lemma}

\proof By Lemma \ref{SmallCases} we can assume that $n \geq 15$. Let $U$ be a vertex set with
$|U|=n-2$, let $\infty_1$ and $\infty_2$ be distinct vertices not in $U$, and let
$V=U\cup\{\infty_1,\infty_2\}$. Since $(M,4^{\frac{n-3}{2}},n)$ is an $n$-ancestor list with
$\nu_n(M) = 0$, it follows from (5) in the definition of ancestor lists that any cycle length in
$M$ is at most $n-3$. Thus, it is easily seen that $(M,n-3)$ is $(n-2)$-admissible and so by
assumption there is an $(M,n-3)$-decomposition $\mathcal{D}$ of $K_{U}$.

Let $H$ be an
$(n-3)$-cycle in $\mathcal{D}$, let $z$ be the vertex in $U \setminus V(H)$,
and let $x$ and $y$ be adjacent vertices in $H$. Then
$$(\mathcal{D} \setminus \{H\}) \cup \{H',(\infty_1,y,x,\infty_2)\} \cup \mathcal{D}_1$$
is an $(M,4^{\frac{n-3}{2}},n)$-decomposition of $K_{V}$, where
\begin{itemize}
    \item
$H'=(H-[x,y]) \cup [x,\infty_1,z,\infty_2,y]$; and
    \item
$\mathcal{D}_1$ is a $(4^{\frac{n-5}2})$-decomposition of $K_{\{\infty_1,\infty_2\},U \setminus
\{x,y,z\}}$. \qed
\end{itemize}

\begin{lemma}\label{nEvenOneHamMany4s}
If $n$ is even, Theorem \ref{mainthm} holds for $K_{n-2}$, and $(M,4^{\frac{n-2}{2}},n)$ is an
$n$-ancestor list with $\nu_n(M) = 0$, then there is an $(M,4^{\frac{n-2}{2}},n)$-decomposition of
$K_n$.
\end{lemma}

\proof By Lemma \ref{SmallCases} we can assume that $n \geq 16$. Let $U$ be a vertex set with
$|U|=n-2$, let $\infty_1$ and $\infty_2$ be distinct vertices not in $U$, and let
$V=U\cup\{\infty_1,\infty_2\}$. Since $(M,4^{\frac{n-2}{2}},n)$ is an $n$-ancestor list with
$\nu_n(M) = 0$, it follows from (5) in the definition of ancestor lists that any cycle length in
$M$ is at most $n-3$. Thus, it is easily seen that $(M,3,n-3)$ is $(n-2)$-admissible and so by
assumption there is an $(M,3,n-3)$-decomposition $\mathcal{D}$ of $K_{U}$.

Let $H$ be an $(n-3)$-cycle in $\mathcal{D}$, let $C$ be a $3$-cycle in $\mathcal{D}$, and let $I$
be the perfect matching in $\mathcal{D}$. Let $z$ be the vertex in $U \setminus V(H)$, let $w$ and
$x$ be distinct vertices in $V(C) \cap V(H)$, let $u$ be the vertex in $V(C)\setminus\{w,x\}$
(possibly $u=z$), and let $y$ be a vertex adjacent to $x$ in $H$. Then
$$(\mathcal{D} \setminus \{I,C,H\}) \cup \{I + \infty_1\infty_2,H',(\infty_1,y,x,w),(\infty_2,x,u,w)\} \cup \mathcal{D}_1$$
is an $(M,4^{\frac{n-2}{2}},n)$-decomposition of $K_{V}$, where
\begin{itemize}
    \item
$H'=(H-[x,y]) \cup [x,\infty_1,z,\infty_2,y]$; and
    \item
$\mathcal{D}_1$ is a $(4^{\frac{n-6}2})$-decomposition of $K_{\{\infty_1,\infty_2\},U \setminus
\{w,x,y,z\}}$. \qed
\end{itemize}

\subsection{Many $5$-cycles and one Hamilton cycle}\label{many5sOneHamSection}

In Lemma \ref{nOddOneHamMany5s} we construct the required decompositions of complete graphs of odd
order and in Lemma \ref{nEvenOneHamMany5s} we construct the required decompositions of complete
graphs of even order. We first require four preliminary lemmas.

\begin{lemma}\label{EvenGraphDecomp}
Every even graph has a decomposition into cycles such that any two cycles in the decomposition
share at most two vertices.
\end{lemma}

\proof It is well known that every even graph has a decomposition into cycles. Let $G$ be an even
graph. Amongst all decompositions of $G$ into cycles, let $\mathcal{D}$ be one with a maximum
number of cycles. We claim that any pair of cycles in $\mathcal{D}$ shares at most two vertices.
Suppose otherwise. That is, there are distinct cycles $A$ and $B$ in $\mathcal{D}$ and distinct
vertices $x$, $y$ and $z$ of $G$ such that $\{x,y,z\} \subseteq V(A) \cap V(B)$. Let $A = A_1 \cup
A_2$ and $B = B_1 \cup B_2$, where $A_1$, $A_2$, $B_1$ and $B_2$ are paths from $x$ to $y$ such
that $z \in V(A_2)$ and $z \in V(B_2)$. Then it is easy to see that $A_1 \cup B_1$ and $A_2 \cup
B_2$ are both nonempty even graphs. For $i=1,2$, let $\mathcal{D}_i$ be a decomposition of $A_i
\cup B_i$ into cycles, and note that $|\mathcal{D}_2| \geq 2$ because $\deg_{A_2 \cup B_2}(z)=4$.
Then $(\mathcal{D}\setminus\{A,B\}) \cup \mathcal{D}_1 \cup \mathcal{D}_2$ is a decomposition of
$G$ into cycles which contains more cycles than $\mathcal{D}$, contradicting our definition of
$\mathcal{D}$. \qed

\begin{lemma}\label{SubgraphNotLasso}
Let $V$ be a vertex set and let $U$ be a subset of $V$ such that $|U| \geq 10$. Let $M$ be a list
of integers, let $m \in \{3,4,5\}$, and let $\mathcal{P}$ be an $(M)$-packing of $K_V$ with a leave
$L$ such that $\deg_{L}(x)=0$ for all $x \in V \setminus U$. If there exists a spanning even
subgraph $G$ of $L$ such that at least one of the following holds,
\begin{itemize}
    \item[(i)]
$\Delta(G)=4$, exactly one vertex of $G$ has degree $4$, $G$ has at most two nontrivial
components, $|E(G)| \geq m+3$, and $G$ does not have a decomposition into two odd cycles if
$m=4$;
    \item[(ii)]
$\Delta(G)=4$, exactly two vertices of $G$ have degree $4$, $G$ has exactly one nontrivial
component, $|E(G)| \geq m+3$, and $\deg_{G}(x)=0$ for some $x \in U$ if $m=3$;
    \item[(iii)]
$m=4$, $G$ has exactly one nontrivial component, $G$ has a decomposition into three cycles each
pair of which intersect in exactly one vertex, and $\deg_{G}(x)=0$ for some $x \in U$; or
    \item[(iv)]
$m=5$, $\Delta(G) \geq 4$, and $G$ has a decomposition into three cycles such that
any two intersect in at most two vertices, and such that any two which intersect
have lengths adding to $6$ or $7$;
\end{itemize}
then there exists $(M,m)$-packing $\mathcal{P}'$ of $K_V$ which is equivalent to $\mathcal{P}$ on
$V \setminus U$.
\end{lemma}

\proof Suppose that there is a spanning even subgraph $G$ of $L$ which satisfies one of (i), (ii),
(iii) or (iv). Because $L$ and $G$ are both even graphs, it follows that $L-G$ is an even graph and
hence has a decomposition $\mathcal{A}=\{A_1,\ldots,A_t\}$ into cycles. Let $a_i=|V(A_i)|$ for $i =
1,\ldots,t$ and let $M^{\dag}=(a_1,\ldots,a_t)$. So $\mathcal{P} \cup \mathcal{A}$ is an
$(M,M^{\dag})$-packing of $K_{V}$ which is equivalent to $\mathcal{P}$ on $V \setminus U$. The
leave of $\mathcal{P} \cup \mathcal{A}$ is $G$.  Let $e=|E(G)|$.

If we can produce an $(M,M^{\dag},m,e-m)$-decomposition $\mathcal{D}$ of $K_V$ which is equivalent
to $\mathcal{P}$ on $V \setminus U$, then there will be cycles in $\mathcal{D}$ with lengths
$a_1,\ldots,a_t,e-m$ whose vertex sets are subsets of $U$, and we can complete the proof by
removing these cycles from $\mathcal{D}$. So it suffices to find such a decomposition. The proof
now splits into cases.

\noindent{\bf Case 1}\quad Suppose that $G$ satisfies (i). Then we can apply Lemma \ref{OneDeg4} to
obtain the required decomposition.

\noindent{\bf Case 2}\quad Suppose that $G$ satisfies (ii). Then we can apply Lemma \ref{TwoDeg4}
to obtain the required decomposition. The only non-trivial thing to check is that there is an $x\in
U$ with $\deg_G(x)=0$ when $m\in\{4,5\}$ and $e-m=3$. In this case we have $e \in \{7,8\}$ and,
because $G$ is even and $|U| \geq 10$, there is indeed an $x\in U$ with $\deg_G(x)=0$.

\noindent{\bf Case 3}\quad Suppose that $G$ satisfies (iii). Either exactly one vertex of $G$ has
degree $6$ and every other vertex of $G$ has degree at most $2$, or exactly three vertices of $G$
have degree $4$ and every other vertex of $G$ has degree at most $2$. In the former case we apply
Lemma \ref{EquitableOneStep}, choosing $b$ to be the vertex of degree $6$ in $G$ and $a$ to be a
neighbour in $G$ of $b$. In the latter case we apply Lemma \ref{EquitableOneStep}, choosing $b$ to
be a vertex of degree $4$ in $G$ and $a$ to be a vertex in $U$ which has degree $0$ in $G$. In
either case we obtain an $(M,M^{\dag})$-packing $\mathcal{P}'$ of $K_{V}$, which is equivalent to
$\mathcal{P}$ on $V \setminus U$, with a leave $G'$ such that $\Delta(G')=4$, exactly two vertices
of $G'$ have degree $4$, $G'$ has exactly one nontrivial component, and $|E(G')| \geq 9$. Thus we
can apply Lemma \ref{TwoDeg4} to obtain the required decomposition.

\noindent{\bf Case 4}\quad Suppose that $G$ satisfies (iv). Since $\Delta(G)\geq 4$ there is at
least one pair of intersecting cycles in any cycle decomposition of $G$. Thus, there exists a
decomposition $\{B_1,B_2,B_3\}$ of $G$ into three cycles such that $|V(B_1) \cap V(B_2)| \in
\{1,2\}$ and $|E(B_1)|+|E(B_2)| \in \{6,7\}$.

\noindent{\bf Case 4a}\quad Suppose that $B_3$ is a component of $G$. If $|V(B_1) \cap V(B_2)| =
1$, then we can apply Lemma \ref{OneDeg4} to obtain the required decomposition, so we may assume
that $|V(B_1) \cap V(B_2)| = 2$. Let $x \in V(B_1) \cap V(B_2)$, let $y \in V(B_3)$, and let $z$ be
a neighbour in $G$ of $x$. Let $\mathcal{P}'$ be the $(M,M^{\dag})$-packing of $K_n$ obtained from
$\mathcal{P}$ by performing the $(x,y)$-switch $S$ with origin $z$ and let $G'$ be its leave. Then
$\mathcal{P}'$ is equivalent to $\mathcal{P}$ on $V \setminus U$, $G'$ has exactly one nontrivial
component, $\Delta(G')=4$, exactly two vertices of $G'$ have degree $4$, and $|E(G')| \geq 9$. Thus
we can apply Lemma \ref{TwoDeg4} to obtain the required decomposition.

\noindent{\bf Case 4b}\quad Suppose that $B_3$ is not a component of $G$. Then $B_3$ intersects
with $B_1$ or $B_2$ and so $|E(B_3)| \in \{3,4\}$. Thus, $e \in \{9,10,11\}$ as we have
$|E(B_1)|+|E(B_2)| \in \{6,7\}$. Let $\mathcal{P}'$ be the $(M,M^{\dag})$-packing of $K_n$ obtained
from $\mathcal{P}$ by repeatedly applying Lemma \ref{EquitableOneStep}, each time choosing $b$ to
be a vertex maximum degree in the leave and $a$ to be a vertex in $U$ of degree $0$ in the leave,
until the leave has maximum degree $4$ and has exactly one vertex of degree $4$ (a suitable choice
for $a$ will exist each time since $e \leq 11$ and $|U| \geq 10$). Let $G'$ be the leave of
$\mathcal{P}'$. Then $\mathcal{P}'$ is equivalent to $\mathcal{P}$ on $V \setminus U$,
$\Delta(G')=4$, exactly one vertex of $G'$ has degree $4$, $|E(G')| \in \{9,10,11\}$, and $G'$ has
at most two components (because $|E(G')| \leq 11$). Thus we can apply Lemma \ref{OneDeg4} to obtain
the required decomposition. \qed

\begin{lemma}\label{NotLasso}
Let $V$ be a vertex set and let $U$ be a subset of $V$ such that $|U| \geq 10$. Let $m \in
\{3,4,5\}$, let $M$ be a list of integers, and let $\mathcal{P}$ be an $(M)$-packing of $K_V$ with
a leave $L$ such that $|E(L)| \geq |U|+m$ and $\deg_{L}(x)=0$ for all $x \in V \setminus U$. Then
there exists an $(M,m)$-packing of $K_V$ which is equivalent to $\mathcal{P}$ on $V \setminus U$.
\end{lemma}

\proof Since $L$ is an even graph, Lemma \ref{EvenGraphDecomp} guarantees that there is a
decomposition $\mathcal{D}$ of $L$ such that any pair of cycles in $\mathcal{D}$ intersect in at
most two vertices. Let $e=|E(L)|$. Since $e \geq |U|+m \geq 13$ it follows that $\mathcal{D}$
contains at least three cycles. Also, since $e>|U|$, there is at least one pair of intersecting
cycles in $\mathcal D$. We now consider separately the cases $m=3$, $m=4$ and $m=5$.

\noindent{\bf Case 1}\quad Suppose that $m=3$. We can assume that there are no $3$-cycles in
$\mathcal{D}$ (otherwise we can simply add one to $\mathcal{P}$ to complete the proof). Let $C_1$,
$C_2$ and $C_3$ be distinct cycles in $\mathcal{D}$ such that $C_1$ and $C_2$ intersect. If
$|V(C_1) \cap V(C_2)| = 1$, then we can apply Lemma \ref{SubgraphNotLasso} (i) (with
$E(G)=E(C_1\cup C_2$)) to complete the proof, so we may assume that $|V(C_1) \cap V(C_2)| = 2$. If
$|E(C_1)|+|E(C_2)| \leq |U|+1$, then there is at least one vertex of $U$ that is not in $V(C_1)\cup
V(C_2)$ and we can apply Lemma \ref{SubgraphNotLasso} (ii) (with $E(G)=E(C_1\cup C_2$)) to complete
the proof. Thus, we may assume $|E(C_1)|+|E(C_2)|=|U|+2$, and it follows from this that $V(C_1)
\cup V(C_2)=U$. This means that $V(C_3)\subseteq V(C_1)\cup V(C_2)$. Thus, since we have
$V(C_3)\geq 4$, $V(C_1)\cap V(C_3)\leq 2$ and $V(C_2)\cap V(C_3)\leq 2$, it follows that
$V(C_3)=4$, $V(C_1)\cap V(C_3)=2$ and $V(C_2)\cap V(C_3)=2$. We can assume without loss of
generality that $|E(C_1)| \leq |E(C_2)|$ and hence that $|E(C_2)| \geq 5$ (since $|U| \geq 10$).
This means that there is at least one vertex of $U$ that is in neither $C_1$ nor $C_3$, and so we
can apply Lemma \ref{SubgraphNotLasso} (ii) (with $E(G)=E(C_1\cup C_3$)) to complete the proof.

\noindent{\bf Case 2}\quad Suppose that $m=4$. If two cycles in $\mathcal{D}$ intersect in exactly
two vertices, then we can apply Lemma \ref{SubgraphNotLasso} (ii) (with the edges of $G$ being the
edges of two such cycles) to complete the proof. So we may assume that any two cycles in
$\mathcal{D}$ intersect in at most one vertex. Let $\{C_1,C_2\}$ be a pair of intersecting cycles
in $\mathcal D$ such that $|E(C_1\cup C_2)|\leq |E(C_i\cup C_j)|$ for any pair $\{C_i,C_j\}$ of
intersecting cycles in $\mathcal D$. If there is a cycle in $\mathcal{D}$ which is vertex disjoint
from $C_1 \cup C_2$, then we can apply Lemma \ref{SubgraphNotLasso} (i) (with the edges of $G$
being the edges of $C_1$, $C_2$ and this cycle) to complete the proof. If there is a cycle in
$\mathcal{D}$ which intersects with exactly one of $C_1$ and $C_2$, then we can apply Lemma
\ref{SubgraphNotLasso} (ii) (with the edges of $G$ being the edges of $C_1$, $C_2$ and this cycle)
to complete the proof. So we may assume that every cycle in $\mathcal{D} \setminus \{C_1,C_2\}$
intersects (in exactly one vertex) with $C_1$ and with $C_2$. Let $C_3$ be a shortest cycle in
$\mathcal{D} \setminus \{C_1,C_2\}$ and note that $|V(C_i)| \leq |V(C_3)|$ for $i=1,2$ by our
definition of $C_1$ and $C_2$. If $V(C_1 \cup C_2 \cup C_3) \neq U$, then we can apply Lemma
\ref{SubgraphNotLasso} (iii) (with $E(G)=E(C_1 \cup C_2 \cup C_3)$) to complete the proof.
Otherwise $V(C_1 \cup C_2 \cup C_3) = U$ which means that $|V(C_1)|+|V(C_2)|+|V(C_3)|\in
\{|U|+2,|U|+3\}$. However, we have $e \geq |U|+4$ and so there is a cycle $C_4 \in \mathcal{D}
\setminus \{C_1,C_2,C_3\}$. Thus $C_4$ is a $3$-cycle (as $V(C_1 \cup C_2 \cup C_3) = U$ and $C_4$
intersects each of $C_1$, $C_2$ and $C_3$ in exactly one vertex). It then follows from the
minimality of $C_3$ and from $|V(C_i)| \leq |V(C_3)|$ for $i=1,2$ that $C_1$, $C_2$ and $C_3$ are
also $3$-cycles. Since $|U| \geq 10$, this is a contradiction and the result is proved.

\noindent{\bf Case 3}\quad Suppose that $m=5$. Let $C_1$, $C_2$ and $C_3$ be three cycles in
$\mathcal{D}$ such that $C_1$ and $C_2$ intersect. If there are a pair of cycles in
$\{C_1,C_2,C_3\}$ which intersect and whose lengths add to at least $8$, then the union of this
pair of cycles has at least $m+3$ edges and we can apply Lemma \ref{SubgraphNotLasso} (i) or Lemma
\ref{SubgraphNotLasso} (ii) (with the edges of $G$ being the edges of this pair of cycles) to
complete the proof. Otherwise we can apply Lemma \ref{SubgraphNotLasso} (iv) to complete the proof.
\qed

\begin{lemma}\label{StartDecomp}
Let $u$ and $k$ be integers such that $u$ is even, $u \geq 16$ and $6 \leq k \leq u-1$, let $U$ be
a vertex set such that $|U| = u$, and let $x$ and $y$ be distinct vertices in $U$. Then there
exists a packing of $K_U$ with a perfect matching, a $u$-cycle, a $(u-1)$-path from $x$ to $y$, a
$k$-cycle, three $(u-2)$-cycles each having vertex set $U \setminus \{x,y\}$, and a $2$-path from
$x$ to $y$.
\end{lemma}

\proof Let $U = \mathbb{Z}_{u-3} \cup \{\infty,x,y\}$.  For $i=0,\ldots,5$, let
$$H_i=(\infty,i,i+1,i+(u-4),i+2,i+(u-5),\ldots,i+\tfrac{u-6}{2},i+\tfrac{u}{2},i+\tfrac{u-4}{2},i+\tfrac{u-2}{2}),$$
and let
$$I = [\tfrac{u-6}{2},\tfrac{u-2}{2}] \cup [\tfrac{u-8}{2},\tfrac{u}{2}] \cup \cdots \cup [0,u-4] \cup [\infty,\tfrac{u-4}{2}],$$
so that $\{I,H_0,\ldots,H_5\}$ is a packing of $\mathbb{Z}_{u-3} \cup \{\infty\}$ with one perfect
matching and six $(u-2)$-cycles (recall that $u \geq 16$). Then
$$\{I + xy,(H_0 - [\infty,0,1]) \cup [\infty,x,0,y,1], (H_1 - [1,2]) \cup [1,x] \cup [2,y], P \cup [a,x,b], H_3, H_4, H_5, [x,c,y]\}$$
is the required packing, where $P$ is a $(k-2)$-path in $H_2$ with endpoints $a$ and $b$ such that
$a,b\in\mathbb{Z}_{u-3}\setminus\{0,1\}$ ($P$ exists as there are $u-2$ distinct paths of length
$k-2$ in $H_2$, and at most six having $\infty$, $0$ or $1$ as an endpoint), and $c$ is any vertex
in $\mathbb{Z}_{u-3} \setminus\{0,1,2,a,b\}$. \qed

\begin{lemma}\label{nOddOneHamMany5s}
If $n$ is odd and $(M,5^{\frac{3n-11}{2}},n)$ is an $n$-ancestor list with $\nu_n(M) = 0$, then
there is an $(M,5^{\frac{3n-11}{2}},n)$-decomposition of $K_n$.
\end{lemma}

\proof By Lemma \ref{SmallCases} (for $n\leq 13$) and Lemma \ref{SmallishCases} (for
$n\in\{15,17\}$) we may assume that $n \geq 19$ (Lemma \ref{SmallishCases} can indeed be applied as
$\nu_5(M,5^{\frac{3n-11}{2}},n)\geq 3$ when $n\in\{15,17\}$). Let $U$ be a vertex set with
$|U|=n-3$, let $x$ and $y$ be distinct vertices in $U$, let $\infty^{\dag}$, $\infty_1$ and
$\infty_2$ be distinct vertices not in $U$, and let $V=U\cup\{\infty^{\dag},\infty_1,\infty_2\}$.

Since $(M,5^{\frac{3n-11}{2}},n)$ is an $n$-ancestor list with $\nu_n(M) = 0$, it follows from (6)
in the definition of ancestor lists that any cycle length in $M$ is at most $n-5$. If there is a
cycle length in $M$ which is at least $6$, then let $k$ be this cycle length. Otherwise let $k=0$
(so $k\in\{0\} \cup \{6,\ldots,n-5\}$). By Lemma \ref{StartDecomp}, there exists a packing
$\mathcal{P}$ of $K_U$ with a perfect matching $I$, an $(n-4)$-path $P_1$ from $x$ to $y$, a
$k$-cycle (if $k \neq 0$), three $(n-5)$-cycles $C_1$, $C_2$ and $C_3$ each having vertex set $U
\setminus \{x,y\}$, and a $2$-path $P_2$ from $x$ to $y$. Let $\{I_1,I_2\}$ be a decomposition of
$C_3$ into two matchings.

Let
$$\mathcal{P'}=(\mathcal{P} \setminus \{I,P_1,P_2,C_1,C_2,C_3\}) \cup \{P_1 \cup [x,\infty_1,\infty^{\dag},\infty_2,y], P_2 \cup [x,\infty_2,\infty_1,y]\} \cup \mathcal{D} \cup \mathcal{D}_1 \cup \mathcal{D}_2,$$
where
\begin{itemize}
    \item
$\mathcal{D}$ is a $(3^{\frac{n-3}2})$-decomposition of $K_{\{\infty^{\dag}\}} \vee I$;
    \item
for $i=1,2$, $\mathcal{D}_i$ is a $(5^{\frac{n-5}2})$-decomposition of
$K_{\{\infty_i\}} \vee (C_i \cup I_i)$ (this exists by Lemma \ref{5CycleTrick}).
\end{itemize}
Then $\mathcal{P}'$ is a $(3^{\frac{n-3}2},5^{n-4},k,n)$-packing of $K_{V}$ (a
$(3^{\frac{n-3}2},5^{n-4},n)$-packing of $K_{V}$ if $k=0$) such that
\begin{itemize}
    \item[(i)]
$\frac{n-3}{2}$ $3$-cycles in $\mathcal{P}'$ contain the vertex $\infty^{\dag}$;
    \item[(ii)]
$\infty^{\dag}\infty_1$ and $\infty^{\dag}\infty_2$ are edges of the $n$-cycle in
$\mathcal{P}'$; and
    \item[(iii)]
$\infty^{\dag}$, $\infty_1$ and $\infty_2$ all have degree $0$ in the leave of $\mathcal{P}'$.
\end{itemize}

Since $(M,5^{\frac{3n-11}{2}},n)$ is an $n$-ancestor list with $\nu_n(M)=0$, it can be seen that by
beginning with $\mathcal{P}'$ and repeatedly applying Lemma \ref{NotLasso} we can obtain an
$(M,3^{\frac{n-3}{2}},5^{n-4},n)$-packing of $K_{V}$ which is equivalent to $\mathcal{P}'$ on $V
\setminus U$. Note that the leave of this packing has $n-3$ edges. Thus, by then repeatedly
applying Lemma \ref{EquitableOneStep} we can obtain an $(M,3^{\frac{n-3}{2}},5^{n-4},n)$-packing
$\mathcal{P}''$ of $K_{V}$ which is equivalent to $\mathcal{P}'$ on $V \setminus U$ and whose leave
$L''$ has the property that $\deg_{L''}(x)=0$ for each $x \in \{\infty^{\dag},\infty_1,\infty_2\}$
and $\deg_{L''}(x)=2$ for each $x \in U$. Because $\mathcal{P}''$ is equivalent to $\mathcal{P}'$
on $V \setminus U$, it follows from (i) and (ii) that there is a set $\mathcal{T}$ of
$\frac{n-3}{2}$ $3$-cycles in $\mathcal{P}''$ each of which contains the vertex $\infty^{\dag}$ and
two vertices in $U$. Let $T$ be the union of the $3$-cycles in $\mathcal{T}$. Then
$$(\mathcal{P}'' \setminus \mathcal{T}) \cup \mathcal{D}''$$
is an $(M,5^{\frac{3n-11}{2}},n)$-decomposition of $K_{V}$ where $\mathcal{D}''$ is a
$(5^{\frac{n-3}2})$-decomposition of $T \cup L''$ (this exists by Lemma \ref{5CycleTrick}, noting
that $E(T \cup L'') = E(K_{\{\infty^{\dag}\}} \vee G)$ for some $3$-regular graph $G$ on vertex set
$U$ which contains a perfect matching). \qed

\begin{lemma}\label{nEvenOneHamMany5s}
If $n$ is even and $(M,5^{2n-9},n)$ is an $n$-ancestor list with $\nu_n(M) = 0$, then there is an
$(M,5^{2n-9},n)$-decomposition of $K_n$.
\end{lemma}

\proof By Lemma \ref{SmallCases} (for $n\leq 14$) and Lemma \ref{SmallishCases} (for
$n\in\{16,18\}$) we may assume that $n \geq 20$ (Lemma \ref{SmallishCases} can indeed be applied as
$\nu_5(M,5^{2n-9},n)\geq 3$ when $n\in\{16,18\}$). Let $U$ be a vertex set with $|U|=n-4$, let $x$
and $y$ be distinct vertices in $U$, let $\infty^{\dag}_1$, $\infty^{\dag}_2$, $\infty_1$ and
$\infty_2$ be distinct vertices not in $U$, and let
$V=U\cup\{\infty^{\dag}_1,\infty^{\dag}_2,\infty_1,\infty_2\}$.

Since $(M,5^{2n-9},n)$ is an $n$-ancestor list with $\nu_n(M) = 0$, it follows from (6) in the
definition of ancestor lists that any cycle length in $M$ is at most $n-5$. If there is a cycle
length in $M$ which is at least $6$ then let $k$ be this cycle length. Otherwise let $k=0$ (so
$k\in\{0\}\cup\{6,\ldots,n-5\}$). By Lemma \ref{StartDecomp}, there exists a packing $\mathcal{P}$
of $K_U$ with a perfect matching $I$, an $(n-4)$-cycle $B$, an $(n-5)$-path $P_1$ from $x$ to $y$,
a $k$-cycle (if $k \neq 0$), and three $(n-6)$-cycles $C_1$, $C_2$ and $C_3$ each having vertex set
$U \setminus \{x,y\}$, and a $2$-path $P_2$ from $x$ to $y$. Let $\{I^{\dag}_1,I^{\dag}_2\}$ be a
decomposition of $B$ into two matchings and $\{I_1,I_2\}$ be a decomposition of $C_3$ into two
matchings.

Let
$$
\begin{array}{ll}
\mathcal{P'}=&(\mathcal{P} \setminus \{B,P_1,P_2,C_1,C_2,C_3\}) \cup \\
&\{I + \{\infty_1\infty_2^\dag,\infty_2\infty_1^\dag\},
P_1 \cup [x,\infty_1,\infty^{\dag}_1,\infty^{\dag}_2,\infty_2,y],
P_2 \cup [x,\infty_2,\infty_1,y]\} \cup \\ &\mathcal{D}^{\dag}_1 \cup \mathcal{D}^{\dag}_2 \cup \mathcal{D}_1 \cup \mathcal{D}_2,
\end{array}
$$
where
\begin{itemize}
    \item
for $i=1,2$, $\mathcal{D}^{\dag}_i$ is a $(3^{\frac{n-4}{2}})$-decomposition of
$K_{\{\infty^{\dag}_i\}} \vee I^{\dag}_i$;
    \item
for $i=1,2$, $\mathcal{D}_i$ is a $(5^{\frac{n-6}2})$-decomposition of $K_{\{\infty_i\}} \vee (C_i \cup I_i)$ (this exists by Lemma \ref{5CycleTrick}).
\end{itemize}
Then $\mathcal{P}'$ is a $(3^{n-4},5^{n-5},k,n)$-packing of $K_V$ (a $(3^{n-4},5^{n-5},n)$-packing
of $K_V$ if $k=0$) such that
\begin{itemize}
    \item[(i)]
for $i=1,2$, $\frac{n-4}{2}$ $3$-cycles in $\mathcal{P}'$ contain the vertex $\infty^{\dag}_i$;
    \item[(ii)]
$\infty^{\dag}_1\infty_1$, $\infty^{\dag}_1\infty^{\dag}_2$ and $\infty^{\dag}_2\infty_2$ are
edges of the $n$-cycle in $\mathcal{P}'$;
    \item[(iii)]
$\infty^{\dag}_1\infty_2$ and $\infty^{\dag}_2\infty_1$ are edges of the perfect matching in
$\mathcal{P}'$; and
    \item[(iv)]
$\infty^{\dag}_1$, $\infty^{\dag}_2$, $\infty_1$ and $\infty_2$ all have degree $0$ in the
leave of $\mathcal{P}'$.
\end{itemize}

Since $(M,5^{2n-9},n)$ is an $n$-ancestor list with $\nu_n(M)=0$, by beginning with $\mathcal{P}'$
and repeatedly applying Lemma \ref{NotLasso} we can obtain an $(M,3^{n-4},5^{n-5},n)$-packing of
$K_{V}$, which is equivalent to $\mathcal{P}'$ on $V \setminus U$. Note that the leave of this
packing has $2n-8$ edges. Thus, by then repeatedly applying Lemma \ref{EquitableOneStep} we can
obtain an $(M,3^{n-4},5^{n-5},n)$-packing $\mathcal{P}''$ of $K_{V}$ which is equivalent to
$\mathcal{P}'$ on $V \setminus U$ and whose leave $L''$ has the property that $\deg_{L''}(x)=0$ for
$x \in \{\infty^{\dag}_1,\infty^{\dag}_2,\infty_1,\infty_2\}$ and $\deg_{L''}(x)=4$ for all $x \in
U$. By Petersen's Theorem \cite{Pet}, $L''$ has a decomposition $\{H_1,H_2\}$ into two $2$-regular
graphs, each with vertex set $U$. Because $\mathcal{P}''$ is equivalent to $\mathcal{P}'$ on $V
\setminus U$, it follows from (i), (ii) and (iii) that, for $i=1,2$ there is a set $\mathcal{T}_i$
of $\frac{n-4}{2}$ $3$-cycles in $\mathcal{P}''$ each of which contains the vertex
$\infty^{\dag}_i$ and two vertices in $U$. For $i=1,2$, let $T_i$ be the union of the $3$-cycles in
$\mathcal{T}_i$. Then
$$(\mathcal{P}'' \setminus (\mathcal{T}_1 \cup \mathcal{T}_2)) \cup \mathcal{D}''_1 \cup \mathcal{D}''_2$$
is an $(M,5^{2n-9},n)$-decomposition of $K_{V}$ where, for $i=1,2$, $\mathcal{D}''_i$ is a
$(5^{\frac{n-4}2})$-decomposition of $T_i \cup H_i$ (these decompositions exist by Lemma
\ref{5CycleTrick}, noting that for $i=1,2$, $E(T_i \cup H_i) = E(K_{\{\infty^{\dag}_i\}} \vee G)$
for some $3$-regular graph $G$ with vertex set $U$ that contains a perfect matching). \qed

\subsection{Proof of Lemma \ref{MainLemmaforatmost1ham} in the case of one Hamilton cycle}

\begin{lemma}\label{Mainlemmaforoneham}
If Theorem \ref{mainthm} holds for $K_{n-1}$, $K_{n-2}$ and $K_{n-3}$,
then there is an
$(M)$-decomposition of $K_n$ for each $n$-ancestor list $M$ satisfying $\nu_n(M)=1$.
\end{lemma}

\proof By Lemma \ref{SmallCases} we can assume that $n \geq 15$. If there is a cycle length in $M$
which is at least $6$ and at most $n-1$ then let $k$ be this cycle length. Otherwise let $k=0$.
We deal separately with the case $n$ is odd and the case $n$ is even.

\noindent{\bf Case 1}\quad Suppose that $n$ is odd. Since $n\geq 15$ and
$3\nu_3(M)+4\nu_4(M)+5\nu_5(M)+k+n =\frac{n(n-1)}{2}$, it can be seen that either
\begin{itemize}
    \item[(i)]
$n \in \{15,17,19\}$ and $\nu_5(M) \geq 3$;
    \item[(ii)]
$\nu_3(M) \geq \frac{n-5}{2}$;
    \item[(iii)]
$\nu_4(M) \geq \frac{n-3}{2}$; or
    \item[(iv)]
$\nu_5(M) \geq \frac{3n-11}{2}$.
\end{itemize}
(To see this consider the cases $\nu_5(M) \geq 3$ and $\nu_5(M) \leq 2$ separately and use the
definition of $n$-ancestor list.) If (i) holds, then the result follows by Lemma
\ref{SmallishCases}. If (ii) holds, then the result follows by one of Lemmas
\ref{nOddOneHamMany3sAtLeastOne4}, \ref{nOddOneHamMany3sAtLeastOne5} or \ref{OneHamNo4sOr5s}. If
(iii) holds, then the result follows by Lemma \ref{nOddOneHamMany4s}. If (iv) holds, then the
result follows by Lemma \ref{nOddOneHamMany5s}.

\noindent{\bf Case 2}\quad Suppose that $n$ is even. Since $n\geq 16$ and
$3\nu_3(M)+4\nu_4(M)+5\nu_5(M)+k+n =\frac{n(n-2)}{2}$, it can be seen that either
\begin{itemize}
    \item[(i)]
$n \in \{16,18,20,22,24,26\}$, $\nu_5(M) \geq 3$, and $\nu_4(M)\geq 2$ if $n=24$;
    \item[(ii)]
$\nu_3(M) \geq \frac{3n-14}{2}$;
    \item[(iii)]
$\nu_4(M) \geq \frac{n-2}{2}$; or
    \item[(iv)]
$\nu_5(M) \geq 2n-9$.
\end{itemize}
(To see this consider the cases $\nu_5(M) \geq 3$ and $\nu_5(M) \leq 2$ separately and use the
definition of $n$-ancestor list.) If (i) holds, then the result follows by Lemma
\ref{SmallishCases}. If (ii) holds, then the result follows by one of Lemmas
\ref{nEvenOneHamMany3sAtLeastOne4},
\ref{nEvenOneHamMany3sAtLeastOne5} or
\ref{OneHamNo4sOr5s}
(note that
$\frac{3n-14}{2} \geq \frac{n-6}{2}$ for $n \geq 16$). If (iii) holds, then the result follows by
Lemma \ref{nEvenOneHamMany4s}. If (iv) holds, then the result follows by Lemma
\ref{nEvenOneHamMany5s}. \qed

\section{Decompositions of $\lan{S}$}\label{decomposeS}\label{sec:konly}

Our general approach to constructing decompositions of $\lan{S}$ follows the approach used in
\cite{BryMar} and \cite{BrySch}. For each connection set $S$ in which we are interested, we define
a graph $J_n$ for each positive integer $n$ such that there is a natural bijection between $E(J_n)$
and $E(\lan{S})$, and such that $\lan{S}$ can be obtained from $J_n$ by identifying a small number
(approximately $|S|$) of pairs of vertices. Thus, decompositions of $J_n$ yield decompositions of
$\lan{S}$.

The key property of the graph $J_n$ is that it can be decomposed into a copy of $J_{n-y}$ and a
copy of $J_y$ for any positive integer $y$ such that $1 \leq y < n$, and this facilitates the
construction of desired decompositions of $J_n$ for arbitrarily large $n$ from decompositions of
$J_i$ for various small values of $i$. For example, in the case $S=\{1,2,3\}$ we define $J_n$ by
$V(J_n)=\{0,\ldots,n+2\}$ and $E(J_n)=\{\{i,i+1\},\{i+1,i+3\},\{i,i+3\}\}:i=0,\ldots n-1\}$. It is
straightforward to construct a $(3)$-decomposition of $J_1$, a $(4,5)$-decomposition of $J_3$, a
$(4^3)$-decomposition of $J_4$ and a $(5^3)$-decomposition of $J_5$. Moreover, since $J_n$
decomposes into $J_{n-y}$ and $J_y$, it is easy to see that these decompositions can be combined to
produce an $(M)$-decomposition of $J_n$ for any list $M=(m_1,\ldots,m_t)$ satisfying $\sum M=3n$
and $m_i\in\{3,4,5\}$ for $i=1,\ldots,t$. For all $n\geq 7$, an $(M)$-decomposition of
$\lan{\{1,2,3\}}$ can be obtained from an $(M)$-decomposition of $J_n$ by identifying vertex $i$
with vertex $i+n$ for $i=0,1,2$.

In what follows, this general approach is modified to allow for the construction of decompositions
which, in addition to cycles of lengths $3$, $4$ and $5$, contain one arbitrarily long cycle or, in
the case $S=\{1,2,3\}$, one arbitrarily long cycle and one Hamilton cycle. The constructions used
to prove Lemma \ref{Sofsize3-8} proceed in a similar fashion for each
connection set $S$.

\subsection{Proof of Lemma \ref{Sofsize3-8}}
\label{sec:Sofsize3-8}
In this section we prove Lemma \ref{Sofsize3-8}, which we restate here for convenience.

\noindent{\bf Lemma \ref{Sofsize3-8}}\quad
{\em
If
$$S\in\{\{1,2,3\},\{1,2,3,4\},\{1,2,3,4,6\},\{1,2,3,4,5,7\},\{1,2,3,4,5,6,7\},\{1,2,3,4,5,6,7,8\}\},$$
$n\geq 2\max(S)+1$, and $M=(m_1,\ldots,m_t,k)$ is any list satisfying $m_i\in\{3,4,5\}$ for
$i=1,\ldots,t$, $3\leq k\leq n$, and $\sum M=|S|n$, then there is an $(M)$-decomposition of
$\lan{S}$, except possibly when
\begin{itemize}
\item $S=\{1,2,3,4,6\}$, $n\equiv 3\md 6$ and $M=(3^{\frac{5n}3})$; or
\item $S=\{1,2,3,4,6\}$, $n\equiv 4\md 6$ and $M=(3^{\frac{5n-5}3},5)$.
\end{itemize}
}

Let
$$\mathcal S=\{\{1,2,3\},\{1,2,3,4\},\{1,2,3,4,6\},
\{1,2,3,4,5,7\},\{1,2,3,4,5,6,7\},\{1,2,3,4,5,6,7,8\}\}.$$
We shall show that the required decompositions exist for each $S \in \mathcal
S$ separately.

%
%

Our proof is essentially inductive and requires a large number of specific base
decompositions, and these are given in the appendix. Some of the
constructions could possibly have been completed using a smaller number of base
decompositions, but since these were found using a computer search, we decided
to keep the inductive steps themselves as simple as possible, at the cost of
requiring a larger number of base decompositions.

\subsubsection{$S=\{1,2,3\}$}
In this section we show the existence of required decompositions for the case
$S=\{1,2,3\}$ in Lemma \ref{Sofsize3-8}.
We first define $J_n$ by
\[
  E(J_n)=\{\{i,i+1\},\{i+1,i+3\},\{i,i+3\}: i=0,\ldots,n-1\}\}
\]
and $V(J_n)=\{0,\ldots,n+2\}$.  We note the following basic properties of
$J_n$.  For a list of integers $M$, an $(M)$-decomposition of $J_n$ will be
denoted by $J_n\rightarrow M$.
\begin{itemize}
    \item
For $n\geq 7$, if for each $i\in\{0,1,2\}$ we identify vertex $i$ of $J_n$ with vertex
$i+n$ of $J_n$ then the resulting graph is $\lan{\{1,2,3\}}$. This means that for $n\geq
7$, we can obtain an $(M)$-decomposition of $\lan{\{1,2,3\}}$ from a decomposition
$J_n\rightarrow M$, provided that for each $i\in\{0,1,2\}$, no cycle in the decomposition
of $J_n$ contains both vertex $i$ and vertex $i+n$.
    \item
For any integers $y$ and $n$ such that $1 \leq y < n$, the graph $J_n$ is the
union of $J_{n-y}$ and the graph obtained from $J_y$ by applying the
permutation $x\mapsto x+(n-y)$.  Thus, if there is a decomposition
$J_{n-y}\rightarrow M$ and a decomposition $J_y\rightarrow M'$, then there is a
decomposition $J_n\rightarrow M,M'$. We will call this construction, and the
similar constructions that follow, \emph{concatenations}.
\end{itemize}

\begin{lemma}\label{123:Jgraph345}
If $n$ is a positive integer and $M=(m_1,\ldots,m_t)$ is a list such that $\sum M=3n$,
$m_i\in\{3,4,5\}$ for $i=1,\ldots,t$, then there is a decomposition $J_n\rightarrow M$.
\end{lemma}

\proof Since $J_1$ is a $3$-cycle, the result holds trivially for $n=1$, so let
$n\geq 2$ and suppose by induction that the result holds for each integer $n'$
in the range $1 \leq n' <n$. The following decompositions
are given in Table \ref{tab:123_decomps} in the appendix.
$$
\begin{array}{lllllll}
J_1 \rightarrow 3 &\quad& J_3 \rightarrow 4,5 &\quad& J_4 \rightarrow 4^3 &\quad& J_5 \rightarrow 5^3 \end{array}
$$

It is routine to check that if $M$ satisfies the hypotheses of the lemma, then
$M$ can be written as $M=(X,Y)$ where $J_y \rightarrow Y$ is one of the
decompositions above and $X$ is some (possibly empty) list. If $X$ is empty,
then we are finished immediately. If $X$ is nonempty then we can obtain a
decomposition $J_n\rightarrow M$ by concatenating a decomposition
$J_{n-y}\rightarrow X$ (which exists by our inductive hypothesis) with a
decomposition $J_y\rightarrow Y$.  \qed

\begin{lemma}\label{123:Jgraphk=6-10}
For $6 \leq k \leq 10$, if $n \geq k$ is an integer and $M=(m_1,\ldots,m_t)$ is
a list such that $\sum M + k = 3n$ and $m_i\in\{3,4,5\}$ for $i=1,\ldots,t$,
then there is a decomposition $J_n\rightarrow M,k$ such that the $k$-cycle is
incident upon vertices $\{0,1,\ldots,k-1\}$.
\end{lemma}
\proof
First we note that Table
\ref{tab:123_6-10} in the appendix lists a number of decompositions containing a $k$-cycle
$C_k$ for some $6 \leq k \leq 10$ such that $V(C_k) \subseteq \{0,\ldots,k-1\}$.
For each $k$, it is easy to use the value of $k\md{3}$ to check that for $n\geq k$
and any $M$ that satisfies the hypotheses of the lemma we can write $M$ as
$(X,Y)$ where $J_x \rightarrow X,k$ is one of the decompositions in Table
\ref{tab:123_6-10},
and $Y$ is some (possibly empty) list with the property $\sum Y = 3y$ for some
integer $y$. If $Y$ is empty we are done, else Lemma \ref{123:Jgraph345} gives
us the existence of a decomposition $J_y \rightarrow Y$ and the required
decomposition can be obtained by concatenation of $J_x \rightarrow X,k$ with
$J_y \rightarrow Y$.

The $k$-cycle in the resulting decomposition will be incident on vertices
$\{0,1,\ldots,k-1\}$ as this property held in the decomposition $J_x
\rightarrow X,k$ and thus the resulting decomposition satisfies the condition
in the lemma.
\qed
\vspace{0.3cm}
Let $M=(m_1,\ldots,m_t)$ be a list of integers with $m_i\geq 3$ for $i=1,\ldots,t$. A decomposition
$\{G_1,\ldots,G_t,C\}$ of $J_n$ such that
\begin{itemize}
    \item
$G_i$ is an $m_i$-cycle for $i=1,\ldots,t$; and
    \item
$C$ is a $k$-cycle such that $V(C)=\{n-k+3,\ldots,n+2\}$ and
$\{n,n+2\} \in E(C)$;
\end{itemize}
will be denoted $J_n\rightarrow M,k^*$.

In Lemma \ref{123:Jgraphk=n} we will form new decompositions of graphs $J_n$ by
concatenating decompositions of $J_{n-y}$ with decompositions of graphs $J^+_y$
which we will now define. For $y \in \{3,\ldots,8\}$, the graph obtained from
$J_y$ by adding the edge $\{0,2\}$ will be denoted $J^+_y$. Let
$M=(m_1,\ldots,m_t)$ be a list of integers with $m_i\geq 3$ for $i=1,\ldots,t$.
A decomposition $\{G_1,\ldots,G_t,A\}$ of $J^+_y$ such that
\begin{itemize}
\item $G_i$ is an $m_i$-cycle for $i=1,\ldots,t$;
\item $A$ is a path from $0$ to $2$ such that $1\notin V(A)$; and
\item $|E(A)|=l+1$;
\end{itemize}
will be denoted $J^+_y \rightarrow M,l^+$. Moreover, if $l=y$ and $\{n,n+2\}
\in E(A)$ , then the decomposition will be denoted $J^+_y\rightarrow M,y^{+*}$.

For $y \in \{3,\ldots,8\}$ and $n > y$, the graph $J_n$ is the union of the
graph obtained from $J_{n-y}$ by deleting the edge $\{n-y,n-y+2\}$ and the
graph obtained from $J^+_y$ applying the permutation $x\mapsto x+(n-y)$. It
follows that if there is a decomposition $J_{n-y}\rightarrow M,k^*$ and a
decomposition $J^+_y\rightarrow M',l^+$, then there is a decomposition $J_n
\rightarrow M,M',k+l$. The edge $\{n,n+2\}$  of the $k$-cycle in the
decomposition of $J_{n-y}$ is replaced by the path in the decomposition of
$J^+_y$ to form the $(k+l)$-cycle in the new decomposition. Similarly, if there
is a decomposition $J_{n-y}\rightarrow M,k^*$ and a decomposition
$J^+_y\rightarrow M',y^{+*}$, then there is a decomposition $J_n\rightarrow
M,M',(k+y)^*$.

\begin{lemma}\label{123:Jgraphk=11-16}
For $11 \leq k \leq 16$, if $n \geq k$ is an integer and $M=(m_1,\ldots,m_t)$ is
a list such that $\sum M + k = 3n$ and $m_i\in\{3,4,5\}$ for $i=1,\ldots,t$,
then there is a decomposition $J_n\rightarrow M,k^*$.
\end{lemma}
\proof
First we note the existence of decompositions of the form $J_x \rightarrow
X,k^*$ listed in Table \ref{tab:123_11-16} in the appendix.
For each $k$, it is routine to use the value of $k\md{3}$ to check that for
$n\geq k$ and any $M$ that satisfies the hypotheses of the lemma we can write
$M$ as $(X,Y)$ where $J_x \rightarrow X,k^*$ is one of the decompositions in
Table \ref{tab:123_11-16}, and $Y$ is some (possibly empty) list with the
property $\sum Y = 3y$ for some integer $y$. If $Y$ is empty we are done, else
Lemma
\ref{123:Jgraph345} gives us the existence of a decomposition $J_y
\rightarrow Y$ and the required decomposition can be obtained by concatenation
of $J_y \rightarrow Y$ with $J_x \rightarrow X,k^*$.

\qed

\begin{lemma}\label{123:Jgraphk=n}
If $n \geq 11$ is an integer and $M=(m_1,\ldots,m_t)$ is a list such that $\sum M=2n$ and
$m_i\in\{3,4,5\}$ for $i=1,\ldots,t$, then there is a decomposition
$J_n\rightarrow M,n^*$.
\end{lemma}

\proof Lemma \ref{123:Jgraphk=11-16} shows that the result holds for $11\leq n
\leq 16$. So let $n \geq 17$ and suppose by induction that the result holds for
each integer $n'$ in the range $11\leq n'<n$. The following decompositions can
be seen in Table \ref{tab:123_decomps} in the appendix.
$$
\begin{array}{lllll}
J^+_5 \rightarrow 5^2,5^{+*} &\quad& J^+_6 \rightarrow 4^3,6^{+*} &\quad& J^+_6 \rightarrow 3^4,6^{+*}
\end{array}
$$

It is routine to check, using $\sum M = 2n \geq 34$, that if $M$ satisfies the
hypotheses of the lemma, then $M$ can be written as
$M=(X,Y)$ where $J_y^+ \rightarrow Y,y^{+*}$ is one of the decompositions above
and $X$ is some nonempty list. We can obtain a decomposition $J_n\rightarrow
M,n^*$ by concatenating a decomposition $J_{n-y}\rightarrow X,(n-y)^*$ (which
exists by our inductive hypothesis, since $n-y \geq n-6 \geq 11$) with a
decomposition $J^+_y \rightarrow Y,y^{+*}$.
\qed

\begin{lemma}\label{123:Jgraphk<n}
If $n$ and $k$ are integers such that $6 \leq k \leq n$ and
$M=(m_1,\ldots,m_t)$ is a list such that $\sum M=3n-k$ and $m_i\in\{3,4,5\}$
for $i=1,\ldots,t$, then there is a decomposition $J_n\rightarrow M,k$.
Furthermore, for $n\geq 7$ all cycles in this decomposition have the property
that for $i\in\{0,1,2\}$ no cycle is incident upon both vertex $i$ and vertex
$n+i$.
\end{lemma}

\proof We first note that if $n\geq 7$ it is clear that any $3$-, $4$- or
$5$-cycle in such a decomposition cannot be incident on two vertices $i$ and
$i+n$ for any $i\in\{0,1,2\}$. As such, Lemma \ref{123:Jgraphk=6-10} shows that
the result holds for all $n$ with $6 \leq k \leq 10$, so in the following we
deal only with $k \geq 11$.

Lemma \ref{123:Jgraphk=11-16} shows that the result holds
for all $n$ with $11\leq k \leq 16$ with the additional property that the
$k$-cycle is not incident upon any vertex in $\{0,1,2\}$, and Lemma
\ref{123:Jgraphk=n} shows that the result holds for all $n=k$ with the same
property on the $k$-cycle. We can therefore assume that $17 \leq k \leq n-1$,
so let $n \geq 18$ and suppose by induction that the result holds for each
positive integer $n'$ in the range $6\leq n'< n$ with the additional property
that the $k$-cycle is not incident upon any vertex in $\{0,1,2\}$.

The following decompositions exist by Lemma \ref{123:Jgraph345}.
$$
\begin{array}{lllllll}
J_1 \rightarrow 3 &\quad& J_3 \rightarrow 4,5 &\quad& J_4 \rightarrow 4^3 &\quad& J_5 \rightarrow 5^3
\end{array}
$$

\noindent{\bf Case 1}\quad Suppose that $k \leq n-5$. Then it is routine to
check, using $\sum M = 3n-k \geq 2n+5 \geq 41$, that $M=(X,Y)$ where $J_y
\rightarrow Y$ is one of the decompositions above and $X$ is some nonempty
list. We can obtain a decomposition $J_n\rightarrow M,k$ by concatenating a
decomposition $J_{n-y}\rightarrow X,k$ (which exists by our inductive
hypothesis, since $k \leq n-6 \leq n-y$) with a decomposition $J_y\rightarrow
Y$. Since $n \geq 17$ it is clear that any $3$-, $4$- or $5$-cycle in this
decomposition having a vertex in $\{0,1,2\}$ has no vertex in
$\{n,n+1,n+2\}$, and by our inductive hypothesis the same holds for the
$k$-cycle.

\noindent{\bf Case 2}\quad Suppose that $n-4 \leq k \leq n-1$. In a similar
manner to Case 1, we can obtain the required decomposition $J_n\rightarrow M,k$
if $M=(X,3)$ for some list $X$, if $k \in \{n-4,n-3\}$ and $M=(X,4,5)$ for some
list $X$, and if $k = n-4$ and $M=(X,4^3)$ for some list $X$. So we may assume
that none of these hold. Additionally, we can construct the
decomposition $J_{18} \rightarrow 5^8,14$ as the concatenation of $J_{13}
\rightarrow 5^5,14^*$ with $J_5 \rightarrow 5^3$ (both given in Table
\ref{tab:123_decomps} in the appendix) as it cannot be constructed by the
method shown below, so in the following also note that we do not consider this
decomposition.

Given this, using $\sum M = 3n-k \geq 2n+1 \geq 37$, it is routine to check
that the required decomposition $J_n\rightarrow M,k$ can be obtained using one
of the concatenations given in the table below (note that, since $\nu_3(M)=0$,
in each case we can deduce the given value of $\nu_5(M) \md 2$ from $\sum M =
3n-k$). The decompositions in the third column exist by Lemma
\ref{123:Jgraphk=n} (since $k\geq 17$), and the decompositions listed in the
last column are shown in Table \ref{tab:123_decomps} in the appendix.

\begin{center}
\begin{tabular}{|l|l|l|l|}
  \hline
  $k$ & $\nu_5(M) \md 2$ & first decomposition & second decomposition \\
  \hline \hline
  $n-4$ & 0 & $J_{n-8} \rightarrow (M-(5^4)),(n-8)^*$ & $J^+_8 \rightarrow 5^4,4^+$  \\
  \hline
  $n-3$ & 1 & $J_{n-6} \rightarrow (M-(5^3)),(n-6)^*$ & $J^+_6 \rightarrow 5^3,3^+$  \\
  \hline
  $n-2$ & 0 & $J_{n-3} \rightarrow (M-(4^2)),(n-3)^*$ & $J^+_3 \rightarrow 4^2,1^+$   \\
  && $J_{n-4} \rightarrow (M-(5^2)), (n-4)^*$ & $J^+_4 \rightarrow 5^2,2^+$ \\
  \hline
  $n-1$ & 1 & $J_{n-4} \rightarrow (M-(4,5)),(n-4)^*$ & $J^+_4 \rightarrow 4,5,3^+$   \\
  && $J_{n-7} \rightarrow (M-(5^3)),(n-7)^*$ & $J^+_7 \rightarrow 5^3,6^+$ \\
  \hline
\end{tabular}
\end{center}
Since $n \geq 18$ it is clear that any $3$-, $4$- or $5$-cycle in this
decomposition having a vertex in $\{0,1,2\}$ has no vertex in
$\{n,\ldots,n+2\}$, and by the definition of the decompositions given in the
third column the $k$-cycle has no vertex in $\{0,1,2\}$, so these
decompositions do have the required properties.

\qed

\vspace{0.5cm}

\begin{lemma}\label{Sofsize3-8-123}
If $S=\{1,2,3\}$, $n\geq 7$, and $M=(m_1,\ldots,m_t,k)$ is any list satisfying
$m_i\in\{3,4,5\}$ for $i=1,\ldots,t$, $3\leq k\leq n$, and $\sum M=3n$, then
there is an $(M)$-decomposition of $\lan{S}$.
\end{lemma}

\proof
As noted
above, for $n\geq 7$ we can obtain an $(M)$-decomposition of $\lan{\{1,2,3\}}$
from an $(M)$-decomposition of $J_n$, provided that for each $i\in\{0,1,2\}$,
no cycle contains both vertex $i$ and vertex $i+n$. Thus, for $S=\{1,2,3\}$,
the required result follows by Lemma \ref{123:Jgraph345} for $k\in\{3,4,5\}$
and by Lemma \ref{123:Jgraphk<n} for $6\leq k\leq n$. \qed

\subsubsection{$S=\{1,2,3,4\}$}
In this section we show the existence of required decompositions for the case $S=\{1,2,3,4\}$ in Lemma \ref{Sofsize3-8}.
We first define $J_n$ by
\[
  E(J_n)=\{\{i+2,i+3\},\{i+2,i+4\},\{i,i+3\},\{i,i+4\}: i=0,\ldots,n-1\}\}
\]
and $V(J_n)=\{0,\ldots,n+3\}$.  We note the following basic properties of
$J_n$.  For a list of integers $M$, an $(M)$-decomposition of $J_n$ will be
denoted by $J_n\rightarrow M$.
\begin{itemize}
    \item
For $n\geq 9$, if for each $i\in\{0,1,2,3\}$ we identify vertex $i$ of $J_n$ with
vertex $i+n$ of $J_n$ then the resulting graph is $\lan{\{1,2,3,4\}}$. This means
that for $n\geq 9$, we can obtain an $(M)$-decomposition of $\lan{\{1,2,3,4\}}$
from a decomposition $J_n\rightarrow M$, provided that for each
$i\in\{0,1,2,3\}$, no cycle in the decomposition of $J_n$ contains both vertex
$i$ and vertex $i+n$.
    \item
For any integers $y$ and $n$ such that $1 \leq y < n$, the graph $J_n$ is the
union of $J_{n-y}$ and the graph obtained from $J_y$ by applying the
permutation $x\mapsto x+(n-y)$.  Thus, if there is a decomposition
$J_{n-y}\rightarrow M$ and a decomposition $J_y\rightarrow M'$, then there is a
decomposition $J_n\rightarrow M,M'$. We will call this construction, and the
similar constructions that follow, \emph{concatenations}.
\end{itemize}

\begin{lemma}\label{1234:Jgraph345}
If $n$ is a positive integer and $M=(m_1,\ldots,m_t)$ is a list such that $\sum
M=4n$, $m_i\in\{3,4,5\}$ for $i=1,\ldots,t$, then there is a decomposition
$J_n\rightarrow M$.
\end{lemma}

\proof Since $J_1$ is a $4$-cycle, the result holds trivially for $n=1$, so let
$n\geq 2$ and suppose by induction that the result holds for each integer $n'$
in the range $1 \leq n' <n$. The following decompositions can be seen in Table
\ref{tab:1234_decomps} in the appendix.
$$
\begin{array}{lllllll}
J_1 \rightarrow 4 &\quad& J_2 \rightarrow 3,5 &\quad& J_3 \rightarrow 3^4 &\quad& J_5 \rightarrow 5^4 \end{array}
$$

It is routine to check that if $M$ satisfies the hypotheses of the lemma, then
$M$ can be written as $M=(X,Y)$ where $J_y \rightarrow Y$ is one of the
decompositions above and $X$ is some (possibly empty) list. If $X$ is empty,
then we are finished immediately. If $X$ is nonempty then we can obtain a
decomposition $J_n\rightarrow M$ by concatenating a decomposition
$J_{n-y}\rightarrow X$ (which exists by our inductive hypothesis) with a
decomposition $J_y\rightarrow Y$.  \qed

\begin{lemma}\label{1234:Jgraphk=6-8}
For $6 \leq k \leq 8$, if $n \geq k$ is an integer and $M=(m_1,\ldots,m_t)$ is
a list such that $\sum M + k = 4n$ and $m_i\in\{3,4,5\}$ for $i=1,\ldots,t$,
then there is a decomposition $J_n\rightarrow M,k$ such that the $k$-cycle is
incident upon vertices $\{0,1,\ldots,k-1\}$.
\end{lemma}
\proof
First we note that Table \ref{tab:1234_6-8} in the appendix lists a number of
decompositions containing a $k$-cycle $C_k$ for some $6 \leq k \leq 8$ such
that $V(C_k) \subseteq \{0,\ldots,k-1\}$.

For each $k$, it is easy to use the value of $k\md{4}$ to check that for $n\geq k$
and any $M$ that satisfies the hypotheses of the lemma we can write $M$ as
$(X,Y)$ where $J_x \rightarrow X,k$ is one of the decompositions in Table
\ref{tab:1234_6-8} and
$Y$ is some (possibly empty) list with the property $\sum Y = 4y$ for some
integer $y$. If $Y$ is empty we are done, else Lemma \ref{1234:Jgraph345} gives us the
existence of a decomposition $J_y \rightarrow Y$ and the required decomposition
can be obtained by concatenation of $J_x \rightarrow X,k$ with $J_y \rightarrow
Y$.

The $k$-cycle in the resulting decomposition will be incident on vertices
$\{0,1,\ldots,k\}$ as this property held in the decomposition $J_x \rightarrow
X,k$ and thus the resulting decomposition satisfies the condition in the lemma.
\qed

\vspace{0.3cm}

Let $M=(m_1,\ldots,m_t)$ be a list of integers with $m_i\geq 3$ for $i=1,\ldots,t$. A decomposition
$\{G_1,\ldots,G_t,C\}$ of $J_n$ such that
\begin{itemize}
    \item
$G_i$ is an $m_i$-cycle for $i=1,\ldots,t$; and
    \item
$C$ is a $k$-cycle such that $V(C)=\{n-k+4,\ldots,n+3\}$ and
$\{\{n,n+2\},\{n+1,n+3\}\} \subseteq E(C)$;
\end{itemize}
will be denoted $J_n\rightarrow M,k^*$.

In Lemma \ref{1234:Jgraphk=n} we will form new decompositions of graphs $J_n$ by
concatenating decompositions of $J_{n-y}$ with decompositions of graphs $J^+_y$
which we will now define. For $y \in \{1,\ldots,7\}$, the graph obtained from
$J_y$ by adding the edges $\{0,2\}$ and $\{1,3\}$ will be denoted $J^+_y$. Let
$M=(m_1,\ldots,m_t)$ be a list of integers with $m_i\geq 3$ for $i=1,\ldots,t$.
A decomposition $\{G_1,\ldots,G_t,A_1,A_2\}$ of $J^+_y$ such that
\begin{itemize}
\item $G_i$ is an $m_i$-cycle for $i=1,\ldots,t$;
\item $A_1$ and $A_2$ are vertex-disjoint paths, one from $0$ to $2$ and one from $1$ to $3$; and
\item $|E(A_1)|+|E(A_2)|=l+2$;
\end{itemize}
will be denoted $J^+_y \rightarrow M,l^+$. Moreover, if $l=y$ and
$\{\{n,n+2\},\{n+1,n+3\}\} \subseteq E(A)$ , then the decomposition will be
denoted $J^+_y\rightarrow M,y^{+*}$.

For $y \in \{1,\ldots,7\}$ and $n > y$, the graph $J_n$ is the union of the
graph obtained from $J_{n-y}$ by deleting the edges $\{n-y,n-y+2\}$ and
$\{n-y+1,n-y+3\}$, and the graph obtained from $J^+_y$ applying the permutation
$x\mapsto x+(n-y)$. It follows that if there is a decomposition
$J_{n-y}\rightarrow M,k^*$ and a decomposition $J^+_y\rightarrow M',l^+$, then
there is a decomposition $J_n \rightarrow M,M',k+l$. The edges $\{n,n+2\}$ and
$\{n+1,n+3\}$ of the $k$-cycle in the decomposition of $J_{n-y}$ are replaced
by the two paths in the decomposition of $J^+_y$ to form the $(k+l)$-cycle in
the new decomposition. Similarly, if there is a decomposition
$J_{n-y}\rightarrow M,k^*$ and a decomposition $J^+_y\rightarrow M',y^{+*}$,
then there is a decomposition $J_n\rightarrow M,M',(k+y)^*$.

\begin{lemma}\label{1234:Jgraphk=9-13}
For $9 \leq k \leq 13$, if $n \geq k$ is an integer and $M=(m_1,\ldots,m_t)$ is
a list such that $\sum M + k = 4n$ and $m_i\in\{3,4,5\}$ for $i=1,\ldots,t$,
then there is a decomposition $J_n\rightarrow M,k^*$.
\end{lemma}
\proof
First we note the existence of decompositions of the form $J_x \rightarrow
X,k^*$ listed in Table \ref{tab:1234_9-13} in the appendix.
For each $k$, it is routine to use the value of $k\md{4}$ to check that for $n\geq k$
and any $M$ that satisfies the hypotheses of the lemma we can write $M$ as
$(X,Y)$ where $J_x \rightarrow X,k^*$ is one of the decompositions in Table
\ref{tab:1234_9-13}, and $Y$ is some (possibly empty) list with the property
$\sum Y = 4y$ for some integer $y$. If $Y$ is empty we are done, else Lemma
\ref{1234:Jgraph345} gives us the existence of a decomposition $J_y \rightarrow
Y$ and the required decomposition can be obtained by concatenation of $J_y
\rightarrow Y$ with $J_x \rightarrow X,k^*$.

\qed

\begin{lemma}\label{1234:Jgraphk=n}
If $n \geq 9$ is an integer and $M=(m_1,\ldots,m_t)$ is a list such that $\sum
M=3n$ and $m_i\in\{3,4,5\}$ for $i=1,\ldots,t$, then there is a decomposition
$J_n\rightarrow M,n^*$.
\end{lemma}

\proof Lemma \ref{1234:Jgraphk=9-13} shows that the result holds for $9\leq n
\leq 13$. So let $n \geq 14$ and suppose by induction that the result holds for
each integer $n'$ in the range $9\leq n'<n$. The following decompositions are
given in Table \ref{tab:1234_decomps} in the appendix.
\[
\begin{array}{lllllll}
J_{1}^+ \rightarrow  3,1^{+*} &\quad& J_{3}^+ \rightarrow  4,5,3^{+*} &\quad&
J_{4}^+ \rightarrow  4^{3},4^{+*} &\quad& J_{5}^+ \rightarrow  5^{3},5^{+*}
\end{array}
\]

It is routine to check, using $\sum M = 3n \geq 42$, that $M$ can be written as
$M=(X,Y)$ where $J_y^+ \rightarrow Y,y^{+*}$ is one of the decompositions above
and $X$ is some nonempty list. We can obtain a decomposition $J_n\rightarrow
M,n^*$ by concatenating a decomposition $J_{n-y}\rightarrow X,(n-y)^*$ (which
exists by our inductive hypothesis, since $n-y \geq n-5 \geq 9$) with a
decomposition $J^+_y \rightarrow Y,y^{+*}$.
\qed

\begin{lemma}\label{1234:Jgraphk<n}
If $n$ and $k$ are integers such that $6 \leq k \leq n$ and
$M=(m_1,\ldots,m_t)$ is a list such that $\sum M=4n-k$ and $m_i\in\{3,4,5\}$
for $i=1,\ldots,t$, then there is a decomposition $J_n\rightarrow M,k$.
Furthermore, for $n\geq 9$ all cycles in this decomposition have the property
that for $i\in\{0,1,2,3\}$ no cycle is incident upon both vertex $i$ and vertex
$n+i$.
\end{lemma}

\proof We first note that if $n\geq 9$ it is clear that any $3$-, $4$- or
$5$-cycle in such a decomposition cannot be incident on two vertices $i$ and
$i+n$ for any $i\in\{0,1,2,3\}$. As such, Lemma \ref{1234:Jgraphk=6-8} shows
that the result
holds for all $n$ with $6 \leq k \leq 8$, so in the following we deal only with
$k \geq 9$.

Lemma \ref{1234:Jgraphk=9-13} shows that the result holds
for all $n$ with $9\leq k \leq 13$ with the additional property that the
$k$-cycle is not incident upon any vertex in $\{0,1,2,3\}$, and Lemma
\ref{1234:Jgraphk=n} shows that the result holds for all $n=k$ with the same
property on the $k$-cycle. We can therefore assume that $14 \leq k \leq n-1$,
so let $n \geq 15$ and suppose by induction that the result holds for each
positive integer $n'$ in the range $6\leq n'< n$ with the additional property
that the $k$-cycle is not incident upon any vertex in $\{0,1,2,3\}$.

The following decompositions exist by Lemma \ref{1234:Jgraph345}.
$$
\begin{array}{lllllll}
J_1 \rightarrow 4 &\quad& J_2 \rightarrow 3,5 &\quad& J_3 \rightarrow 3^4 &\quad& J_5 \rightarrow 5^4
\end{array}
$$

\noindent{\bf Case 1}\quad Suppose that $k \leq n-5$. Then it is routine to
check, using $\sum M = 4n-k \geq 3n+5 \geq 50$, that $M=(X,Y)$ where $J_y
\rightarrow Y$ is one of the decompositions above and $X$ is some nonempty
list. We can obtain a decomposition $J_n\rightarrow M,k$ by concatenating a
decomposition $J_{n-y}\rightarrow X,k$ (which exists by our inductive
hypothesis, since $k \leq n-6 \leq n-y$) with a decomposition $J_y\rightarrow
Y$. Since $n \geq 15$ it is clear that any $3$-, $4$- or $5$-cycle in this
decomposition having a vertex in $\{0,1,2,3\}$ has no vertex in
$\{n,n+1,n+2,n+3\}$, and by our inductive hypothesis the same holds for the
$k$-cycle.

\noindent{\bf Case 2}\quad Suppose that $n-4 \leq k \leq n-1$. In a similar
manner to Case 1, we can obtain the required decomposition $J_n\rightarrow M,k$
if $M=(X,4)$ for some list $X$, if $k \in \{n-4,n-3,n-2\}$ and $M=(X,3,5)$ for some
list $X$, and if $k \in\{n-4,n-3\}$ and $M=(X,3^4)$ for some list $X$. So we
may assume that none of these hold.

Given this, using $\sum M = 4n-k \geq 3n+1 \geq 46$, it is routine to check
that the required decomposition $J_n\rightarrow M,k$ can be obtained using one
of the concatenations given in the table below (note that, since $\nu_4(M)=0$,
in each case we can deduce the given value of $\nu_5(M) \md 3$ from $\sum M =
3n-k$). The decompositions in the third column exist by Lemma
\ref{1234:Jgraphk=n} (since $k\geq 14$), and the decompositions listed in the
last column are shown in Table \ref{tab:1234_decomps} in the appendix.

\begin{center}
\begin{tabular}{|l|l|l|l|}
  \hline
  $k$ & $\nu_5(M) \md 3$ & first decomposition & second decomposition \\
  \hline \hline
  $n-4$ & 2 & $J_{n-7} \rightarrow (M-(5^5)),(n-7)^*$ & $J^+_7 \rightarrow 5^5,3^+$  \\
  \hline
  $n-3$ & 0 & $J_{n-4} \rightarrow (M-(5^3)),(n-4)^*$ & $J^+_4 \rightarrow 5^3,1^+$  \\
  \hline
  $n-2$ & 1 & $J_{n-6} \rightarrow (M-(5^4)),(n-6)^*$ & $J^+_6 \rightarrow 5^4,4^+$   \\
  \hline
  $n-1$ & 2 & $J_{n-3} \rightarrow (M-(5^2)),(n-3)^*$ & $J^+_2 \rightarrow 5^2,2^+$   \\
  \hline
\end{tabular}
\end{center}
For $n \geq 9$ it is clear that any $3$-, $4$- or $5$-cycle in this
decomposition having a vertex in $\{0,1,2,3\}$ has no vertex in
$\{n,\ldots,n+3\}$, and by the definition of the decompositions given in the
third column the $k$-cycle has no vertex in $\{0,1,2,3\}$, so these
decompositions do have the required properties.

\qed

\vspace{0.5cm}

\begin{lemma}\label{Sofsize3-8-1234}
If $S=\{1,2,3,4\}$, $n\geq 9$, and $M=(m_1,\ldots,m_t,k)$ is any list
satisfying $m_i\in\{3,4,5\}$ for $i=1,\ldots,t$, $3\leq k\leq n$, and $\sum
M=4n$, then there is an $(M)$-decomposition of $\lan{S}$.
\end{lemma}

\proof
As noted
earlier, for $n\geq 9$ we can obtain an $(M)$-decomposition of $\lan{\{1,2,3,4\}}$
from an $(M)$-decomposition of $J_n$, provided that for each $i\in\{0,1,2,3\}$,
no cycle contains both vertex $i$ and vertex $i+n$. Thus, for $S=\{1,2,3,4\}$,
the required result follows by Lemma \ref{1234:Jgraph345} for $k\in\{3,4,5\}$
and by Lemma \ref{1234:Jgraphk<n} for $6\leq k\leq n$. \qed

\subsubsection{$S=\{1,2,3,4,6\}$}\label{sec:12346}
In this section we show the existence of required decompositions for the case
$S=\{1,2,3,4,6\}$ in Lemma \ref{Sofsize3-8}.
We first define $J_n$ by
\[
E(J_n)=\{\{i+2,i+3\},\{i+2,i+4\},\{i+3,i+6\},\{i,i+4\},\{i,i+6\}:i=0,\ldots,n-1\}
\]
and $V(J_n)=\{0,\ldots,n+5\}$.  We
note the following basic properties of $J_n$.

For a list of integers $M$, an $(M)$-decomposition of $J_n$ will be denoted by $J_n\rightarrow M$.
We note the following basic properties of $J_n$.
\begin{itemize}
    \item
For $n\geq 13$, if for each $i\in\{0,1,2,3,4,5\}$ we identify vertex $i$ of $J_n$ with vertex
$i+n$ of $J_n$ then the resulting graph is $\lan{\{1,2,3,4,6\}}$. This means that for $n\geq
13$, we can obtain an $(M)$-decomposition of $\lan{\{1,2,3,4,6\}}$ from a decomposition
$J_n\rightarrow M$, provided that for each $i\in\{0,1,2,3,4,5\}$, no cycle in the decomposition
of $J_n$ contains both vertex $i$ and vertex $i+n$.
    \item
For any integers $y$ and $n$ such that $1 \leq y < n$, the graph $J_n$ is the union of
$J_{n-y}$ and the graph obtained from $J_y$ by applying the permutation $x\mapsto x+(n-y)$.
Thus, if there is a decomposition $J_{n-y}\rightarrow M$ and a decomposition $J_y\rightarrow
M'$, then there is a decomposition $J_n\rightarrow M,M'$. We will call this construction, and
the similar constructions that follow, \emph{concatenations}.
\end{itemize}

\begin{lemma}\label{12346:Jgraph345}
If $n$ is a positive integer and $M=(m_1,\ldots,m_t)$ is a list such that $\sum
M=5n$, $m_i\in\{3,4,5\}$ for $i=1,\ldots,t$, and $M \notin \mathcal{E}$ where
$$\mathcal{E}=\{(3^2,4)\} \cup \{(3^{5i}):\mbox{$i\geq 1$ is odd}\} \cup
\{(3^{5i},5):\mbox{$i\geq 1$ is odd}\},$$ then there is a decomposition
$J_n\rightarrow M$.
\end{lemma}

\proof We have verified by computer search and concatenation that the result
holds for $n\leq 10$. So assume $n\geq 11$ and suppose by induction that the
result holds for each integer $n'$ in the range $1 \leq n' < n$.
The following decompositions are given in
Table \ref{tab:12346_decomps} in the appendix.
$$
\begin{array}{llllllllllll}
J_1 \rightarrow 5 &\quad& J_4 \rightarrow 4^5 &\quad& J_6 \rightarrow 3^{10} &\quad& J_5 \rightarrow 3^7,4 &\quad& J_4 \rightarrow 3^4,4^2 &\quad& J_3 \rightarrow 3,4^3
\end{array}
$$

It is routine to check that if $M$ satisfies the hypotheses of the lemma, then
$M$ can be written as $M=(X,Y)$ where $J_y \rightarrow Y$ is one of the
decompositions above and $X$ is some (possibly empty) list. If $X$ is
empty, then we are finished immediately. If $X$ is nonempty and $X \notin
\mathcal{E}$, then we can obtain a decomposition $J_n\rightarrow M$ by
concatenating a decomposition $J_{n-y}\rightarrow X$ (which exists by our
inductive hypothesis) with a decomposition $J_y\rightarrow Y$. Thus, we can
assume $X \in \mathcal{E}$.
But since $n\geq 11$ and $\sum Y \leq 30$, we have $\sum X\geq 25$ which
implies $X\in \{(3^{5i}):\mbox{$i\geq 3$ is odd}\} \cup
\{(3^{5i},5):\mbox{$i\geq 3$ is odd}\}$.

It follows that $M=(3^{10},X')$ for some
nonempty list $X'\notin \mathcal{E}$ (because $M \notin \mathcal{E}$) and we
can obtain a decomposition $J_n\rightarrow M$ by concatenating a decomposition
$J_{n-6}\rightarrow X'$ (which exists by our inductive hypothesis) with a
decomposition $J_6\rightarrow 3^{10}$.
\qed
\vspace{0.5cm}

\begin{lemma}\label{12346:Jgraphk=6-10}
For $6 \leq k \leq 10$, if $n \geq k$ is an integer and $M=(m_1,\ldots,m_t)$ is
a list such that $\sum M + k = 5n$ and $m_i\in\{3,4,5\}$ for $i=1,\ldots,t$,
then there is a decomposition $J_n\rightarrow M,k$.
Furthermore, for $n\geq 13$ any cycle in this decomposition having a vertex in
$\{0,\ldots,5\}$ has no vertex in $\{n,\ldots,n+5\}$.
\end{lemma}
\proof
First we note that Table \ref{tab:12346_6-10} in the appendix lists a number of
decompositions containing a $k$-cycle $C_k$ for some $6 \leq k \leq 10$ such
that $V(C_k) \subseteq \{0,\ldots,12\}$.
For each $k$, it is routine to check that for $n\geq k$
and any $M$ that satisfies the hypotheses of the lemma we can write $M$ as
$(X,Y)$ where $J_x \rightarrow X,k$ is one of the decompositions in Table
\ref{tab:12346_6-10},
and $Y$ is some (possibly empty) list with the property $\sum Y = 5y$ for some
integer $y$. If $Y$ is empty we are done, else Lemma \ref{12346:Jgraph345}
gives us the existence of a decomposition $J_y \rightarrow Y$ and the required
decomposition can be obtained by concatenation of $J_y \rightarrow Y$ with $J_x
\rightarrow X,k$.  All decompositions in the list have the property that the
$k$-cycle is only incident upon some subset of the vertices $\{0,\ldots,12\}$ ,
and the resulting decomposition after concatenation will still have this
property.  Additionally it is simple to check that every decomposition used has
the property that no $3$-, $4$- or $5$-cycle contains two vertices $v_a$ and
$v_b$ such that $|v_a - v_b| \geq 8$, so this gives the required
decompositions.
\qed

\vspace{0.3cm}

Let $M=(m_1,\ldots,m_t)$ be a list of integers with $m_i\geq 3$ for $i=1,\ldots,t$. A decomposition
$\{G_1,\ldots,G_t,C\}$ of $J_n$ such that
\begin{itemize}
    \item
$G_i$ is an $m_i$-cycle for $i=1,\ldots,t$; and
    \item
$C$ is a $k$-cycle such that $V(C)=\{n-k+6,\ldots,n+5\}$ and
$\{\{n,n+3\},\{n+1,n+4\},\{n+2,n+5\}\} \subseteq E(C)$;
\end{itemize}
will be denoted $J_n\rightarrow M,k^*$.

In Lemma \ref{12346:Jgraphk=n} we will form new decompositions of graphs $J_n$
by concatenating decompositions of $J_{n-y}$ with decompositions of graphs
$J_y^+$ which we will now define. For $y \in \{4,\ldots,11\}$, the graph obtained
from $J_y$ by adding the three edges $\{0,3\}$, $\{1,4\}$ and $\{2,5\}$ will be
denoted $J^+_y$. Let $M=(m_1,\ldots,m_t)$ be a list of integers with $m_i\geq
3$ for $i=1,\ldots,t$. A decomposition $\{G_1,\ldots,G_t,A_1,A_2,A_3\}$ of
$J^+_y$ such that
\begin{itemize}
\item $G_i$ is an $m_i$-cycle for $i=1,\ldots,t$;
\item $A_1$, $A_2$ and $A_3$ are vertex-disjoint paths, one from $0$ to $3$, one from $1$ to $4$, and one from $2$ to
$5$; and
\item $|E(A_1)|+|E(A_2)|+|E(A_3)|=l+3$;
\end{itemize}
will be denoted $J^+_y \rightarrow M,l^+$. Moreover, if $l=y$ and
$\{\{n,n+3\},\{n+1,n+4\},\{n+2,n+5\}\} \subseteq E(A_1) \cup E(A_2) \cup
E(A_3)$, then the decomposition will be denoted $J^+_y\rightarrow M,y^{+*}$.

For $y \in \{4,\ldots,11\}$ and $n > y$, the graph $J_n$ is the union of the
graph obtained from $J_{n-y}$ by deleting the edges in
$\{\{n-y,n-y+3\},\{n-y+1,n-y+4\},\{n-y+2,n-y+5\}\}$ and the graph obtained from
$J^+_y$ applying the permutation $x\mapsto x+(n-y)$. It follows that if there
is a decomposition $J_{n-y}\rightarrow M,k^*$ and a decomposition
$J^+_y\rightarrow M',l^+$, then there is a decomposition $J_n \rightarrow
M,M',k+l$. The three edges $\{n,n+3\}$, $\{n+1,n+4\}$ and $\{n+2,n+5\}$ of the
$k$-cycle in the decomposition of $J_{n-y}$ are replaced by the three paths in
the decomposition of $J^+_y$ to form the $(k+l)$-cycle in the new
decomposition. Similarly, if there is a decomposition $J_{n-y}\rightarrow
M,k^*$ and a decomposition $J^+_y\rightarrow M',y^{+*}$, then there is a
decomposition $J_n\rightarrow M,M',(k+y)^*$.

\begin{lemma}\label{12346:Jgraphk=11-16}
For $11 \leq k \leq 16$, if $n \geq k$ is an integer and $M=(m_1,\ldots,m_t)$ is
a list such that $\sum M + k = 3n$ and $m_i\in\{3,4,5\}$ for $i=1,\ldots,t$,
then there is a decomposition $J_n\rightarrow M,k^*$.
\end{lemma}
\proof
First we note the existence of the decompositions given in Table
\ref{tab:12346_11-16} in the appendix. For $11 \leq k \leq 14$ it is routine to
check that for $n \geq k$ and any $M$ that satisfies the hypotheses of the
lemma we can write $M$ as $(X,Y)$ where $J_x \rightarrow X,k^*$ is one of the
decompositions in Table \ref{tab:12346_11-16}, and $Y$ is some (possibly empty)
list with the properties that $\sum Y = 5y$ for some integer $y$ and $Y$ is not
one of the exceptions to Lemma \ref{12346:Jgraph345}. If $Y$ is
empty, we are done, else Lemma \ref{12346:Jgraph345} gives us the existence of
a decomposition $J_y \rightarrow Y$ and the required decomposition can be
obtained by concatenation of $J_y \rightarrow Y$ with $J_x \rightarrow X,k^*$.

For $15 \leq k \leq 16$ we note the existence of the following
decompositions, given in Table \ref{tab:12346_decomps} in the appendix.
\[
\begin{array}{lllllll}
J_{4}^+ \rightarrow  4^{4},+4^{*} &\quad& J_{4}^+ \rightarrow
3^{2},5^{2},+4^{*} &\quad&
J_{4}^+ \rightarrow  3,4^{2},5,+4^{*} &\quad& J_{4}^+ \rightarrow
3^{4},4,+4^{*} \\
J_5^+ \rightarrow 5^4,+5^{*} &\quad& J_5^+ \rightarrow 3^5,+5^{*}
\end{array}
\]

For $n\geq k$, if $M$ can be written as $(X,Y)$ such that $J_y^+ \rightarrow
Y,+y$ is a decomposition in the above list and $k-y \geq 11$, then the
required decomposition can be obtained by concatenation of a decomposition of
$J_x \rightarrow X,(k-y)^*$ (which exists by the argument above) with the
decomposition $J_y^+ \rightarrow Y,+y$.

For $k\geq 15$ we therefore assume that $Y \notin M$ for any $Y \in\{(4^4),
(3^2,5^2), (3,4^2,5), (3^4,4)\}$ and for $k=16$ we add the additional
assumptions that $5^4 \notin M$ and $3^5,5 \notin M$.

Given $n\geq k$ and a list $M$ that satisfies these assumptions (where
applicable) and the conditions of the lemma, it is routine to check that we can
write $M$ as $(X,Y)$ where $J_x \rightarrow X,k^*$ is one of the decompositions
listed in Table \ref{tab:12346_11-16} in the appendix, and $Y$ is some
(possibly empty) list with the property $\sum Y = 5y$ for some integer $y$. If
$Y$ is empty we are done, else Lemma \ref{12346:Jgraph345} gives us the
existence of a decomposition $J_y \rightarrow Y$ and the required decomposition
can be obtained by concatenation of $J_y \rightarrow Y$ with $J_x \rightarrow
X,k^*$.  \qed

\begin{lemma}\label{12346:Jgraphk=n}
If $n \geq 11$ is an integer and $M=(m_1,\ldots,m_t)$ is a list such that $\sum
M=4n$ and $m_i\in\{3,4,5\}$ for $i=1,\ldots,t$, then there is a decomposition
$J_n\rightarrow M,n^*$.
\end{lemma}

\proof
Lemma \ref{12346:Jgraphk=11-16} shows that the result holds for $11\leq n \leq
16$. So let $n \geq 17$ and suppose by induction that the result holds for each
integer $n'$ in the range $11\leq n'<n$. The
following decompositions are given in Table \ref{tab:12346_decomps} in the
appendix.
$$
\begin{array}{lllll}
J^+_4 \rightarrow 4^4,4^{+*} &\quad& J^+_5 \rightarrow 5^4,5^{+*} &\quad& J^+_6 \rightarrow 3^8,6^{+*}
\end{array}
$$

It is routine to check, using $\sum M = 4n \geq 68$, that $M$ can be written as
$M=(X,Y)$ where $J_y^+ \rightarrow Y,y^{+*}$ is one of the decompositions above
and $X$ is some nonempty list. We can obtain a decomposition $J_n\rightarrow
M,n^*$ by concatenating a decomposition $J_{n-y}\rightarrow X,(n-y)^*$ (which
exists by our inductive hypothesis, since $n-y \geq n-6 \geq 11$) with a
decomposition $J^+_y \rightarrow Y,y^{+*}$. \qed

\begin{lemma}\label{12346:Jgraphk<n}
If $n$ and $k$ are integers such that $6 \leq k \leq n$ and
$M=(m_1,\ldots,m_t)$ is a list such that $\sum M=5n-k$ and $m_i\in\{3,4,5\}$
for $i=1,\ldots,t$, then there is a decomposition $J_n\rightarrow M,k$.
Furthermore, for $n\geq 13$ any cycle in this decomposition having a vertex in
$\{0,\ldots,5\}$ has no vertex in $\{n,\ldots,n+5\}$.
\end{lemma}

\proof Lemmas \ref{12346:Jgraphk=6-10} and \ref{12346:Jgraphk=11-16} show that
the result holds for $6 \leq k \leq 16$, and Lemma \ref{12346:Jgraphk=n} shows
that it holds for $k=n$, so we can also assume that $17 \leq k \leq n-1$.
Let $n \geq 18$ and suppose by induction that the result holds for each positive
integer $n'$ in the range $6\leq n'< n$.
The following decompositions exist by Lemma \ref{12346:Jgraph345}.
$$
\begin{array}{lllllll}
J_1 \rightarrow 5^1 &\quad& J_3 \rightarrow 3^14^3 &\quad& J_4 \rightarrow 4^5 &\quad& J_6 \rightarrow 3^{10}
\end{array}
$$

\noindent{\bf Case 1}\quad Suppose that $k \leq n-6$. Then it is routine to
check, using $\sum M = 5n-k \geq 4n+6 \geq 78$, that $M=(X,Y)$ where $J_y
\rightarrow Y$ is one of the decompositions above and $X$ is some nonempty
list. We can obtain a decomposition $J_n\rightarrow M,k$ by concatenating a
decomposition $J_{n-y}\rightarrow X,k$ (which exists by our inductive
hypothesis, since $k \leq n-6 \leq n-y$) with a decomposition $J_y\rightarrow
Y$. As this concatenation does not change the $k$-cycle, this decomposition has
the desired properties.

\noindent{\bf Case 2}\quad Suppose that $n-5 \leq k \leq n-1$. In a similar
manner to Case 1, we can obtain the required decomposition $J_n\rightarrow M,k$
if $M=(X,5)$ for some list $X$, if $k \in \{n-5,n-4,n-3\}$ and $M=(X,3,4^3)$
for some list $X$, and if $k \in \{n-5,n-4\}$ and $M=(X,4^5)$ for some list
$Y$. So we may assume that none of these hold. Additionally, the following two
decompositions are noted here.
$$
\begin{array}{cc}
J_{18} \rightarrow 3^{23},4,17 &
J_{21} \rightarrow 3^{28},4,17
\end{array}
$$
The first of these decompositions is given in Table \ref{tab:12346_decomps} in
the appendix. The decomposition $J_{21} \rightarrow 3^{28},4,17$ can be
obtained by concatenation of $J_{12} \rightarrow 3^{16},12^*$ (which exists by
Lemma \ref{12346:Jgraphk=n}) with $J_9^+ \rightarrow 3^{12},4,5^+$ (also given
in Table \ref{tab:12346_decomps}). As such, we do not consider these
decompositions in what follows.

Given the exceptions above, using $\sum M = 5n-k \geq 4n+1 \geq 73$, it is routine
to check that the required decomposition $J_n\rightarrow M,k$ can be obtained
using one of the concatenations given in the table below (note that, since
$\nu_5(M)=0$, in each case we can deduce the given value of $\nu_3(M) \md 4$
from $\sum M = 5n-k$). The decompositions in the third column exist by Lemma
\ref{12346:Jgraphk=n} (since $k\geq 17$), and the decompositions
listed in the last column are given in Table \ref{tab:12346_decomps}.

\begin{center}
\begin{tabular}{|l|l|l|l|}
  \hline
  $k$ & $\nu_3(M) \md 4$ & first decomposition & second decomposition \\
  \hline \hline
  $n-5$ & 3 & $J_{n-10} \rightarrow (M-(3^{15})),(n-10)^*$ & $J^+_{10} \rightarrow 3^{15},5^+$  \\
  \hline
  $n-4$ & 0 & $J_{n-11} \rightarrow (M-(3^{16})),(n-11)^*$ & $J^+_{11} \rightarrow 3^{16},7^+$  \\
  \hline
  $n-3$ & 1 & $J_{n-9} \rightarrow (M-(3^{13})),(n-9)^*$ & $J^+_9 \rightarrow 3^{13},6^+$  \\
  \hline
  $n-2$ & 2 & $J_{n-7} \rightarrow (M-(3^{10})),(n-7)^*$ & $J^+_7 \rightarrow 3^{10},5^+$   \\
  && $J_{n-5} \rightarrow (M-(3^2,4^4)), (n-5)^*$ & $J^+_5 \rightarrow 3^2,4^4,3^+$ \\
  \hline
  $n-1$ & 3 & $J_{n-8} \rightarrow (M-(3^{11})),(n-8)^*$ & $J^+_8 \rightarrow 3^{11},7^+$   \\
  && $J_{n-4} \rightarrow (M-(3^3,4^2)),(n-4)^*$ & $J^+_4 \rightarrow 3^3,4^2,3^+$ \\
  \hline
\end{tabular}
\end{center}
By the definition of the decompositions given in the third column the vertex
set of the $k$-cycle is some subset of $\{6,\ldots,n+5\}$. \qed

\vspace{0.5cm}

\begin{lemma}\label{Sofsize3-8-12346}
If $S=\{1,2,3,4,6\}$, $n\geq 13$, and $M=(m_1,\ldots,m_t,k)$ is any list
satisfying $m_i\in\{3,4,5\}$ for $i=1,\ldots,t$, $3\leq k\leq n$, and $\sum
M=5n$, then there is an $(M)$-decomposition of $\lan{S}$, except possibly when
\begin{itemize}
\item $n\equiv 3\md 6$ and $M=(3^{\frac{5n}3})$; or
\item $n\equiv 4\md 6$ and $M=(3^{\frac{5n-5}3},5)$.
\end{itemize}
\end{lemma}

\proof  As noted above, for $n\geq
13$ we can obtain an $(M)$-decomposition of $\lan{\{1,2,3,4,6\}}$ from an
$(M)$-decomposition of $J_n$, provided that for each $i\in\{0,\ldots,5\}$, no
cycle contains both vertex $i$ and vertex $i+n$. Thus, for $S=\{1,2,3,4,6\}$,
Lemma \ref{Sofsize3-8} follows by Lemma \ref{12346:Jgraph345} for
$k\in\{3,4,5\}$ and by Lemma \ref{12346:Jgraphk<n} for $6\leq k\leq n$. \qed

\subsubsection{$S=\{1,2,3,4,5,7\}$}
In this section we show the existence of required decompositions for the case
$S=\{1,2,3,4,5,7\}$ in Lemma \ref{Sofsize3-8}.
We first define $J_n$
by
\[
  E(J_n)=\{\{i+6,i+7\},\{i+5,i+7\},\{i+3,i+6\},\{i+3,i+7\},\{i,i+5\},\{i,i+7\}:
i=0,\ldots,n-1\}\}
\]
and $V(J_n)=\{0,\ldots,n+6\}$.  We note the following basic
properties of $J_n$.  For a list of integers $M$, an $(M)$-decomposition of
$J_n$ will be denoted by
$J_n\rightarrow M$.
\begin{itemize}
    \item
For $n\geq 15$, if for each $i\in\{0,1,2,3,4,5,6\}$ we identify vertex $i$ of $J_n$ with
vertex $i+n$ of $J_n$ then the resulting graph is $\lan{\{1,2,3,4,5,7\}}$. This means
that for $n\geq 15$, we can obtain an $(M)$-decomposition of $\lan{\{1,2,3,4,5,7\}}$
from a decomposition $J_n\rightarrow M$, provided that for each
$i\in\{0,1,\ldots,6\}$, no cycle in the decomposition of $J_n$ contains both vertex
$i$ and vertex $i+n$.
    \item
For any integers $y$ and $n$ such that $1 \leq y < n$, the graph $J_n$ is the
union of $J_{n-y}$ and the graph obtained from $J_y$ by applying the
permutation $x\mapsto x+(n-y)$.  Thus, if there is a decomposition
$J_{n-y}\rightarrow M$ and a decomposition $J_y\rightarrow M'$, then there is a
decomposition $J_n\rightarrow M,M'$. We will call this construction, and the
similar constructions that follow, \emph{concatenations}.
\end{itemize}

\begin{lemma}\label{123457:Jgraph345}
If $n$ is a positive integer and $M=(m_1,\ldots,m_t)$ is a list such that $\sum
M=6n$, $m_i\in\{3,4,5\}$ for $i=1,\ldots,t$, and $M\neq (4^3)$ then there
is a decomposition $J_n\rightarrow M$.
\end{lemma}

\proof We first note the existence of the following decompositions, given in
Table \ref{tab:123457_decomps} in the appendix.
$$
\begin{array}{lllllll}
J_{1} \rightarrow  3^{2} &\quad& J_{2} \rightarrow  3,4,5 &\quad&
J_{3} \rightarrow  4^{2},5^{2} &\quad& J_{3} \rightarrow  3,5^{3} \\
J_{3} \rightarrow  3^{2},4^{3} &\quad& J_{4} \rightarrow  4,5^{4} &\quad&
J_{4} \rightarrow  4^{6} &\quad& J_{4} \rightarrow  3,4^{4},5 \\
J_{5} \rightarrow  5^{6} &\quad& J_{5} \rightarrow  4^{5},5^{2} &\quad&
J_{6} \rightarrow  4^{9} &\quad& \end{array}
$$

The only required decomposition of $J_1$ is shown in the table above, so we
may assume $n\geq 2$ and assume by induction that the result holds for any
positive integer $n'$ in the range $1 \leq n' < n$.

It is routine to check that for $n\geq 2$ if $M$ satisfies the hypotheses of
the lemma, then $M$ can be written as $M=(X,Y)$ for some (possibly empty) list
$Y\neq (4^3)$ where $J_x \rightarrow X$ is one of the decompositions listed
above.  If $Y$ is empty, we are done, else we can obtain the required
decomposition by concatenation of $J_y \rightarrow Y$ (which exists by our
inductive hypothesis since $Y\neq (4^3)$) with the decomposition $J_x\rightarrow
X$.  \qed

\begin{lemma}\label{123457:Jgraphk=6}
For $n \geq 4$, if $M=(m_1,\ldots,m_t)$ is
a list such that $\sum M + 6 = 6n$ and $m_i\in\{3,4,5\}$ for $i=1,\ldots,t$,
then there is a decomposition $J_n\rightarrow M,6$ such that the $6$-cycle is
incident upon vertices $\{i,i+1,\ldots,i+5\}$ for some integer $i$.
\end{lemma}
\proof
First we note that Table \ref{tab:123457_6} in the appendix lists a number of
decompositions containing a $6$-cycle $C_6$ such that $V(C_6) \subseteq
\{4,\ldots,9\}$.

It is routine to check that for $n\geq 4$ and any $M$ that satisfies the
hypotheses of the lemma we can write $M$ as $(X,Y)$ where $J_x \rightarrow X,k$
is one of the decompositions in Table \ref{tab:123457_6}, and $Y\neq (4^3)$ is
some (possibly empty) list. If $Y$ is empty we are done, else Lemma
\ref{123457:Jgraph345} gives us the existence of a decomposition $J_y
\rightarrow Y$ and the required decomposition can be obtained by concatenation
of $J_x \rightarrow X,k$ with $J_y \rightarrow Y$.  \qed

\vspace{0.3cm}

Let $M=(m_1,\ldots,m_t)$ be a list of integers with $m_i\geq 3$ for $i=1,\ldots,t$. A decomposition
$\{G_1,\ldots,G_t,C\}$ of $J_n$ such that
\begin{itemize}
    \item
$G_i$ is an $m_i$-cycle for $i=1,\ldots,t$; and
    \item
$C$ is a $k$-cycle such that $V(C)=\{n-k+7,\ldots,n+6\}$ and
$\{n+4,n+6\} \in E(C)$;
\end{itemize}
will be denoted $J_n\rightarrow M,k^*$.

In Lemma \ref{123457:Jgraphk=n} we will form new decompositions of graphs $J_n$ by
concatenating decompositions of $J_{n-y}$ with decompositions of graphs $J^+_y$
which we will now define. For $y \in \{2,\ldots,7\}$, the graph obtained from
$J_y$ by adding the edge $\{4,6\}$ will be denoted $J^+_y$. Let
$M=(m_1,\ldots,m_t)$ be a list of integers with $m_i\geq 3$ for $i=1,\ldots,t$.
A decomposition $\{G_1,\ldots,G_t,A\}$ of $J^+_y$ such that
\begin{itemize}
\item $G_i$ is an $m_i$-cycle for $i=1,\ldots,t$;
\item $A$ is a path from $4$ to $6$ such that $\{0,1,2,3,5\} \cap V(A) =
\emptyset$; and
\item $|E(A)|=l+1$;
\end{itemize}
will be denoted $J^+_y \rightarrow M,l^+$. Moreover, if $l=y$ and
$\{n+4,n+6\} \in E(A)$, then the decomposition will be
denoted $J^+_y\rightarrow M,y^{+*}$.

For $y \in \{2,\ldots,7\}$ and $n > y$, the graph $J_n$ is the union of the
graph obtained from $J_{n-y}$ by deleting the edge $\{n-y+4,n-y+6\}$, and the
graph obtained from $J^+_y$ applying the permutation $x\mapsto x+(n-y)$. It
follows that if there is a decomposition $J_{n-y}\rightarrow M,k^*$ and a
decomposition $J^+_y\rightarrow M',l^+$, then there is a decomposition $J_n
\rightarrow M,M',k+l$. The edge $\{n+4,n+6\}$ of the $k$-cycle in the
decomposition of $J_{n-y}$ are replaced by the two paths in the decomposition
of $J^+_y$ to form the $(k+l)$-cycle in the new decomposition. Similarly, if
there is a decomposition $J_{n-y}\rightarrow M,k^*$ and a decomposition
$J^+_y\rightarrow M',y^{+*}$, then there is a decomposition $J_n\rightarrow
M,M',(k+y)^*$.

\begin{lemma}\label{123457:Jgraphk=7-12}
For $7 \leq k \leq 12$, if $n$ is an integer with $n \geq k$ and $M=(m_1,\ldots,m_t)$ is
a list such that $\sum M + k = 6n$ and $m_i\in\{3,4,5\}$ for $i=1,\ldots,t$,
then there is a decomposition $J_n\rightarrow M,k^*$.
\end{lemma}
\proof

For each $k$, it is routine to use the value of $k\md{6}$ to check that for
$n\geq k$ and any $M$ that satisfies the hypotheses of the lemma we can write
$M$ as $(X,Y)$ where $J_x \rightarrow X,k^*$ is one of the decompositions
in Table \ref{tab:123457_7-12} in the appendix, and $Y \neq (4^3)$ is some
(possibly empty) list. If $Y$ is empty, then we are done, else Lemma
\ref{123457:Jgraph345} gives us the existence of a decomposition $J_y
\rightarrow Y$ and the required decomposition can be obtained by concatenation
of $J_y \rightarrow Y$ with $J_x \rightarrow X,k^*$.  \qed

\begin{lemma}\label{123457:Jgraphk=n}
Given an integer $n \geq 7$, if $M=(m_1,\ldots,m_t)$ is a list such that $\sum
M=5n$ and $m_i\in\{3,4,5\}$ for $i=1,\ldots,t$, then there is a decomposition
$J_n\rightarrow M,n^*$.
\end{lemma}

\proof Lemma \ref{123457:Jgraphk=7-12} shows that the result holds for $7\leq n
\leq 12$. So let $n \geq 13$ and suppose by induction that the result holds for
each integer $n'$ in the range $7\leq n'<n$. The following decompositions are
given in Table \ref{tab:123457_decomps} in the appendix.
\[
\begin{array}{lllllll}
J_{4}^+ \rightarrow  5^4,4^* &\quad& J_{4}^+ \rightarrow  4^{5},4^{*} &\quad&
J_{6}^+ \rightarrow  3^{10},6^{*} &\quad& \end{array}
\]

It is routine to check, using $\sum M = 5n \geq 65$, that $M$ can be written as
$M=(X,Y)$ where $J_y^+ \rightarrow Y,y^{+*}$ is one of the decompositions above
and $X$ is some nonempty list. We can obtain a decomposition $J_n\rightarrow
M,n^*$ by concatenating a decomposition $J_{n-y}\rightarrow X,(n-y)^*$ (which
exists by our inductive hypothesis, since $n-y \geq n-6 \geq 7$) with a
decomposition $J^+_y \rightarrow Y,y^{+*}$.
\qed

\begin{lemma}\label{123457:Jgraphk<n}
If $n$ and $k$ are integers such that $6 \leq k \leq n$ and
$M=(m_1,\ldots,m_t)$ is a list such that $\sum M=6n-k$ and $m_i\in\{3,4,5\}$
for $i=1,\ldots,t$, then there is a decomposition $J_n\rightarrow M,k$.
Furthermore, for $n\geq 15$ all cycles in this decomposition have the property
that for $i\in\{0,1,\ldots,6\}$ no cycle is incident upon both vertex $i$ and vertex
$n+i$.
\end{lemma}

\proof We first note that if $n\geq 15$ it is clear that any $3$-, $4$- or
$5$-cycle in such a decomposition cannot be incident on two vertices $i$ and
$i+n$ for any $i\in\{0,1,\ldots,6\}$. As such, Lemma \ref{123457:Jgraphk=6}
shows that the result holds for all $n$ with $k=6$, so in the following
we deal only with $k \geq 7$.

Lemma \ref{123457:Jgraphk=7-12} shows that the result holds
for all $n\geq k$ with $7\leq k \leq 12$ with the property that the $k$-cycle
is not incident upon any vertex in $\{0,1,\ldots,6\}$, and Lemma
\ref{123457:Jgraphk=n} shows that the result holds for all $n=k$ with the same
property on the $k$-cycle. We can
therefore assume that $13 \leq k \leq n-1$, so let $n \geq 14$ and
suppose by induction that the result holds for each positive integer $n'$ in
the range $6\leq n'< n$ with the property that the $k$-cycle is not incident
upon any vertex in $\{0,1,\ldots,6\}$.

The following decompositions exist by Lemma \ref{123457:Jgraph345}.
$$
\begin{array}{lllllll}
J_1 \rightarrow 3^2 &\quad& J_2 \rightarrow 3,4,5 &\quad& J_3 \rightarrow
4^2,5^2 \\
J_4 \rightarrow 4^6 &\quad& J_4 \rightarrow 4,5^4 &\quad& J_5 \rightarrow 5^6
\end{array}
$$

\noindent{\bf Case 1}\quad Suppose that $k \leq n-5$. Then it is routine to
check, using $\sum M = 6n-k \geq 5n+5 \geq 75$, that $M=(X,Y)$ where $J_y
\rightarrow Y$ is one of the decompositions above and $X$ is some nonempty
list. We can obtain a decomposition $J_n\rightarrow M,k$ by concatenating a
decomposition $J_{n-y}\rightarrow X,k$ (which exists by our inductive
hypothesis, since $k \leq n-5 \leq n-y$) with a decomposition $J_y\rightarrow
Y$. For $n \geq 15$ it is clear that any $3$-, $4$- or $5$-cycle in this
decomposition having a vertex in $\{0,1,\ldots,6\}$ has no vertex in
$\{n,n+1,\ldots,n+6\}$, and by our inductive hypothesis the same holds for the
$k$-cycle.

\noindent{\bf Case 2}\quad Suppose that $n-4 \leq k \leq n-1$. In a similar
manner to Case 1, we can obtain the required decomposition $J_n\rightarrow M,k$
if $M=(X,3^2)$ for some list $X$, if $k \in \{n-4,n-3,n-2\}$ and $M=(X,3,4,5)$
for some list $X$, if $k \in\{n-4,n-3\}$ and $M=(X,4^2,5^2)$ for some list $X$,
and if $k = n-4$ and $M=(X,4^6)$ or $M=(X,4,5^4)$ for some list $X$. So we may
assume that none of these hold.

Given this, using $\sum M = 6n-k \geq 5n+1 \geq 71$, it is routine to check
that the required decomposition $J_n\rightarrow M,k$ can be obtained using one
of the concatenations given in the table below. Note that, since $\nu_3(M) <
2$, there are only two cases to take for $\nu_3(M)$ and in either case we can
deduce the given value of $\nu_4(M) \md 5$ from $\sum M = 6n-k$.  For $k=n-4$
this shows that $\nu_4(M) \geq 1$ so it is routine to see that all required
decompositions for $k=n-4$ can be constructed in the manner described in the
previous paragraph.  The decompositions in the fourth column exist by Lemma
\ref{123457:Jgraphk=n}, and the decompositions listed in the last column are
shown in Table \ref{tab:123457_decomps} in the appendix.

\begin{center}
\begin{tabular}{|l|l|l|l|l|}
  \hline
  $k$ & $\nu_3(M)$ &$\nu_4(M) \md 5$ & first decomposition & second decomposition \\
  \hline \hline
  $n-3$ & 1 & 0 & $J_{n-4} \rightarrow (M-(3,4^5)),(n-4)^*$ & $J^+_4 \rightarrow 3,4^5,1^+$  \\
   & & & $J_{n-4} \rightarrow (M-(3,5^4)),(n-4)^*$ & $J^+_4 \rightarrow 3,5^4,1^+$  \\
  $n-3$ & 0 & 2 & $J_{n-5} \rightarrow (M-(4^7)),(n-5)^*$ & $J^+_5 \rightarrow 4^7,2^+$  \\
  \hline
  $n-2$ & 1 & 1 & $J_{n-5} \rightarrow (M-(3,4^6)),(n-5)^*$ & $J^+_5 \rightarrow 3,4^6,3^+$  \\
  $n-2$ & 0 & 3 & $J_{n-3} \rightarrow (M-(4^3,5)),(n-3)^*$ & $J^+_3 \rightarrow 4^3,5,1^+$  \\
   & & & $J_{n-6} \rightarrow (M-(4^8)),(n-6)^*$ & $J^+_6 \rightarrow 4^8,4^+$  \\
  \hline
  $n-1$ & 1 & 2 & $J_{n-2} \rightarrow (M-(3,4^2)),(n-2)^*$ & $J^+_2 \rightarrow 3,4^2,1^+$   \\
  $n-1$ & 0 & 4 & $J_{n-3} \rightarrow (M-(4^3,5)),(n-3)^*$ & $J^+_3 \rightarrow 4^3,5,1^+$   \\
   & & & $J_{n-7} \rightarrow (M-(4^9)),(n-7)^*$ & $J^+_7 \rightarrow 4^9,6^+$   \\
  \hline
\end{tabular}
\end{center}
For $n \geq 15$ it is clear that any $3$-, $4$- or $5$-cycle in this
decomposition having a vertex in $\{0,1,\ldots,6\}$ has no vertex in
$\{n,\ldots,n+6\}$, and by the definition of the decompositions given in the
fourth column the $k$-cycle has no vertex in $\{0,1,\ldots,6\}$, so these
decompositions do have the required properties.
\qed

\begin{lemma}\label{Sofsize3-8-123457}
If $S=\{1,2,3,4,5,7\}$, $n\geq 15$, and $M=(m_1,\ldots,m_t,k)$ is any list
satisfying $m_i\in\{3,4,5\}$ for $i=1,\ldots,t$, $3\leq k\leq n$, and $\sum
M=6n$, then there is an $(M)$-decomposition of $\lan{S}$.
\end{lemma}

\proof As noted earlier, for $n\geq 15$ we can obtain an $(M)$-decomposition of
$\lan{\{1,2,3,4,5,7\}}$ from an $(M)$-decomposition of $J_n$, provided that for
each $i\in\{0,1,\ldots,6\}$, no cycle contains both vertex $i$ and vertex
$i+n$. Thus, for $S=\{1,2,3,4,5,7\}$, the required result follows by Lemma
\ref{123457:Jgraph345} for $k\in\{3,4,5\}$ and by Lemma \ref{123457:Jgraphk<n}
for $6\leq k\leq n$. \qed

\subsubsection{$S=\{1,2,3,4,5,6,7\}$}
In this section we show the existence of required decompositions for the case
$S = \{1,2,3,4,5,6,7\}$ in Lemma \ref{Sofsize3-8}.
We first define $J_n$
by
\begin{align*}
  E(J_n)=\{&\{i+3,i+4\},\{i+5,i+7\},\{i+3,i+6\},\{i,i+4\},\\ &\{i,i+5\},\{i,i+6\},\{i,i+7\}:
i=0,\ldots,n-1\}\}
\end{align*}
and $V(J_n)=\{0,\ldots,n+6\}$.  We note the following basic
properties of $J_n$.  For a list of integers $M$, an $(M)$-decomposition of
$J_n$ will be denoted by
$J_n\rightarrow M$.
\begin{itemize}
    \item
For $n\geq 15$, if for each $i\in\{0,1,2,3,4,5,6\}$ we identify vertex $i$ of $J_n$ with
vertex $i+n$ of $J_n$ then the resulting graph is $\lan{\{1,2,3,4,5,6,7\}}$. This means
that for $n\geq 15$, we can obtain an $(M)$-decomposition of $\lan{\{1,2,3,4,5,6,7\}}$
from a decomposition $J_n\rightarrow M$, provided that for each
$i\in\{0,1,\ldots,6\}$, no cycle in the decomposition of $J_n$ contains both vertex
$i$ and vertex $i+n$.
    \item
For any integers $y$ and $n$ such that $1 \leq y < n$, the graph $J_n$ is the
union of $J_{n-y}$ and the graph obtained from $J_y$ by applying the
permutation $x\mapsto x+(n-y)$.  Thus, if there is a decomposition
$J_{n-y}\rightarrow M$ and a decomposition $J_y\rightarrow M'$, then there is a
decomposition $J_n\rightarrow M,M'$. We will call this construction, and the
similar constructions that follow, \emph{concatenations}.
\end{itemize}

\begin{lemma}\label{1234567:Jgraph345}
If $n$ is a positive integer and $M=(m_1,\ldots,m_t)$ is a list such that $\sum
M=7n$, $m_i\in\{3,4,5\}$ for $i=1,\ldots,t$, and $M\neq (3^7)$ then there
is a decomposition $J_n\rightarrow M$.
\end{lemma}

\proof The following decompositions are given in
Table \ref{tab:1234567_decomps} in the appendix.

\[
\begin{array}{lllllll}
J_{1} \rightarrow  3,4 &\quad& J_{2} \rightarrow  4,5^{2} &\quad&
J_{2} \rightarrow  3^{3},5 &\quad& J_{3} \rightarrow  4^{4},5 \\
J_{3} \rightarrow  3^{2},5^{3} &\quad& J_{4} \rightarrow  4^{7} &\quad&
J_{4} \rightarrow  3,5^{5} &\quad& J_{5} \rightarrow  5^{7} \\
J_{6} \rightarrow  3^{14} &\quad& J_{9} \rightarrow  3^{21} &\quad&
\end{array}
\]

The only required decomposition of $J_1$ is shown in the table above, so we
may assume $n\geq 2$ and assume by induction that the result holds for any
positive integer $n'$ in the range $1 \leq n' < n$.

It is routine to check that for $n\geq 2$ if $M$ satisfies the hypotheses of
the lemma, then $M$ can be written as $M=(X,Y)$ for some (possibly empty) list
$Y\neq (3^7)$ where $J_x \rightarrow X$ is one of the decompositions listed
above.  If $Y$ is empty, we are done, else we can obtain the required
decomposition by concatenation of $J_y \rightarrow Y$ (which exists by our
inductive hypothesis since $Y\neq (3^7)$) with the decomposition $J_x\rightarrow
X$.  \qed

\begin{lemma}\label{1234567:Jgraphk=6-7}
For $k\in\{6,7\}$ and $n \geq k+1$, if $M=(m_1,\ldots,m_t)$ is
a list such that $\sum M + k = 7n$ and $m_i\in\{3,4,5\}$ for $i=1,\ldots,t$,
then there is a decomposition $J_n\rightarrow M,k$ such that the $k$-cycle is
incident upon vertices $\{i,i+1,\ldots,i+5\}$ for some integer $i$.
\end{lemma}
\proof

First we note that Table \ref{tab:1234567_6-7} in the appendix lists a
number of decompositions required for this lemma. All of these decompositions
contain a $k$-cycle for some $k \in\{6,7\}$ where the $k$-cycle is incident
on some subset of the vertices $\{4,\ldots,k+3\}$.

It is routine to check that for $n\geq 4$ and any $M$ that satisfies the
hypotheses of the lemma we can write $M$ as $(X,Y)$ where $J_x \rightarrow X,k$
is one of the decompositions in Table \ref{tab:1234567_6-7} in the appendix,
and $Y\neq (3^7)$ is some (possibly empty) list. If $Y$ is empty we are done,
else Lemma \ref{1234567:Jgraph345} gives us the existence of a decomposition
$J_y \rightarrow Y$ and the required decomposition can be obtained by
concatenation of $J_x \rightarrow X,k$ with $J_y \rightarrow Y$.  \qed

\vspace{0.3cm}

Let $M=(m_1,\ldots,m_t)$ be a list of integers with $m_i\geq 3$ for $i=1,\ldots,t$. A decomposition
$\{G_1,\ldots,G_t,C\}$ of $J_n$ such that
\begin{itemize}
    \item
$G_i$ is an $m_i$-cycle for $i=1,\ldots,t$; and
    \item
$C$ is a $k$-cycle such that $V(C)=\{n-k+7,\ldots,n+6\}$ and
$\{n+4,n+6\} \in E(C)$;
\end{itemize}
will be denoted $J_n\rightarrow M,k^*$.

In Lemma \ref{1234567:Jgraphk=n} we will form new decompositions of graphs $J_n$ by
concatenating decompositions of $J_{n-y}$ with decompositions of graphs $J^+_y$
which we will now define. For $y \in \{5,\ldots,10\}$, the graph obtained from
$J_y$ by adding the edge $\{4,6\}$ will be denoted $J^+_y$. Let
$M=(m_1,\ldots,m_t)$ be a list of integers with $m_i\geq 3$ for $i=1,\ldots,t$.
A decomposition $\{G_1,\ldots,G_t,A\}$ of $J^+_y$ such that
\begin{itemize}
\item $G_i$ is an $m_i$-cycle for $i=1,\ldots,t$;
\item $A$ is a path from $4$ to $6$ such that $\{0,1,2,3,5\} \cap V(A) =
\emptyset$; and
\item $|E(A)|=l+1$;
\end{itemize}
will be denoted $J^+_y \rightarrow M,l^+$. Moreover, if $l=y$ and
$\{n+4,n+6\} \in E(A)$, then the decomposition will be
denoted $J^+_y\rightarrow M,y^{+*}$.

For $y \in \{5,\ldots,10\}$ and $n > y$, the graph $J_n$ is the union of the
graph obtained from $J_{n-y}$ by deleting the edge $\{n-y+4,n-y+6\}$, and the
graph obtained from $J^+_y$ applying the permutation $x\mapsto x+(n-y)$. It
follows that if there is a decomposition $J_{n-y}\rightarrow M,k^*$ and a
decomposition $J^+_y\rightarrow M',l^+$, then there is a decomposition $J_n
\rightarrow M,M',k+l$. The edge $\{n+4,n+6\}$ of the $k$-cycle in the
decomposition of $J_{n-y}$ are replaced by the two paths in the decomposition
of $J^+_y$ to form the $(k+l)$-cycle in the new decomposition. Similarly, if
there is a decomposition $J_{n-y}\rightarrow M,k^*$ and a decomposition
$J^+_y\rightarrow M',y^{+*}$, then there is a decomposition $J_n\rightarrow
M,M',(k+y)^*$.

\begin{lemma}\label{1234567:Jgraphk=8-17}
For $8 \leq k \leq 17$, if $n$ is an integer with $n \geq k$ and $M=(m_1,\ldots,m_t)$ is
a list such that $\sum M + k = 7n$ and $m_i\in\{3,4,5\}$ for $i=1,\ldots,t$,
then there is a decomposition $J_n\rightarrow M,k^*$.
\end{lemma}
\proof

For each $k$, it is routine to use the value of $k\md{7}$ to check that for
$n\geq k$ and any $M$ that satisfies the hypotheses of the lemma we can write
$M$ as $(X,Y)$ where $J_x \rightarrow X,k^*$ is one of the decompositions
in Table \ref{tab:1234567_8-17} in the appendix, and $Y \neq (3^7)$ is some
(possibly empty) list. If $Y$ is empty, then we are done, else Lemma
\ref{1234567:Jgraph345} gives us the existence of a decomposition $J_y
\rightarrow Y$ and the required decomposition can be obtained by concatenation
of $J_y \rightarrow Y$ with $J_x \rightarrow X,k^*$.
\qed

\begin{lemma}\label{1234567:Jgraphk=n}
Given an integer $n \geq 7$, if $M=(m_1,\ldots,m_t)$ is a list such that $\sum
M=6n$ and $m_i\in\{3,4,5\}$ for $i=1,\ldots,t$, then there is a decomposition
$J_n\rightarrow M,n^*$.
\end{lemma}

\proof Lemma \ref{1234567:Jgraphk=8-17} shows that the result holds for $8\leq n
\leq 17$. So let $n \geq 18$ and suppose by induction that the result holds for
each integer $n'$ in the range $8\leq n'<n$. The following decompositions are
given in Table \ref{tab:1234567_decomps} in the appendix.

\[
\begin{array}{lllllll}
J_{7}^+ \rightarrow  4^{3},5^{6},7^{+*} &\quad& J_{7}^+ \rightarrow
3^{4},5^{6},7^{+*} &\quad&
J_{7}^+ \rightarrow  3^{6},4^{6},7^{+*} &\quad& J_{7}^+ \rightarrow
3^{14},7^{+*} \\
J_{8}^+ \rightarrow  4^{12},8^{+*} &\quad& J_{10}^+ \rightarrow  5^{12},10^{+*}
&\quad&
\end{array}
\]

It is routine to check, using $\sum M = 6n \geq 108$, that $M$ can be written as
$M=(X,Y)$ where $J_y^+ \rightarrow Y,y^{+*}$ is one of the decompositions above
and $X$ is some nonempty list. We can obtain a decomposition $J_n\rightarrow
M,n^*$ by concatenating a decomposition $J_{n-y}\rightarrow X,(n-y)^*$ (which
exists by our inductive hypothesis, since $n-y \geq n-10 \geq 8$) with a
decomposition $J^+_y \rightarrow Y,y^{+*}$.
\qed

\begin{lemma}\label{1234567:Jgraphk<n}
If $n$ and $k$ are integers such that $6 \leq k \leq n$ and
$M=(m_1,\ldots,m_t)$ is a list such that $\sum M=7n-k$ and $m_i\in\{3,4,5\}$
for $i=1,\ldots,t$, then there is a decomposition $J_n\rightarrow M,k$.
Furthermore, for $n\geq 15$ all cycles in this decomposition have the property
that for $i\in\{0,1,\ldots,6\}$ no cycle is incident upon both vertex $i$ and vertex
$n+i$.
\end{lemma}

\proof We first note that if $n\geq 15$ it is clear that any $3$-, $4$- or
$5$-cycle in such a decomposition cannot be incident on two vertices $i$ and
$i+n$ for any $i\in\{0,1,\ldots,6\}$. As such, Lemma \ref{1234567:Jgraphk=6-7}
shows that the result holds for all $n$ with $k\in\{6,7\}$, so in the following
we deal only with $k \geq 8$.

Lemma \ref{1234567:Jgraphk=8-17} shows that the result holds
for all $n\geq k$ with $8\leq k \leq 17$ with the property that the $k$-cycle
is not incident upon any vertex in $\{0,1,\ldots,6\}$, and Lemma
\ref{1234567:Jgraphk=n} shows that the result holds for all $n=k$ with the same
property on the $k$-cycle. We can
therefore assume that $18 \leq k \leq n-1$, so let $n \geq 19$ and
suppose by induction that the result holds for each positive integer $n'$ in
the range $6\leq n'< n$ with the property that the $k$-cycle is not incident
upon any vertex in $\{0,1,\ldots,6\}$.

The following decompositions exist by Lemma \ref{1234567:Jgraph345}.
\[
\begin{array}{lllllll}
J_{1} \rightarrow  3,4 &\quad& J_{2} \rightarrow  4,5^{2} &\quad&
J_{2} \rightarrow  3^{3},5 &\quad& J_{3} \rightarrow  4^{4},5 \\
J_{3} \rightarrow  3^{2},5^{3} &\quad& J_{4} \rightarrow  4^{7} &\quad&
J_{4} \rightarrow  3,5^{5} &\quad& J_{5} \rightarrow  5^{7} \\
J_{6} \rightarrow  3^{14} &\quad&
\end{array}
\]

\noindent{\bf Case 1}\quad Suppose that $k \leq n-6$. Then it is routine to
check, using $\sum M = 7n-k \geq 6n+6 \geq 120$, that $M=(X,Y)$ where $J_y
\rightarrow Y$ is one of the decompositions above and $X$ is some nonempty
list. We can obtain a decomposition $J_n\rightarrow M,k$ by concatenating a
decomposition $J_{n-y}\rightarrow X,k$ (which exists by our inductive
hypothesis, since $k \leq n-6 \leq n-y$) with a decomposition $J_y\rightarrow
Y$. Since $n \geq 18$ it is clear that any $3$-, $4$- or $5$-cycle in this
decomposition having a vertex in $\{0,1,\ldots,6\}$ has no vertex in
$\{n,n+1,\ldots,n+6\}$, and by our inductive hypothesis the same holds for the
$k$-cycle.

\noindent{\bf Case 2}\quad Suppose that $n-5 \leq k \leq n-1$. In a similar
manner to Case 1, we can obtain the required decomposition $J_n\rightarrow M,k$
if $M=(X,Y)$ for some list $X$ where $J_x \rightarrow X$ is one of the
decompositions shown above and $k+x \leq n$. In what follows we take deal with
each case of $k \in\{n-5,n-4,\ldots,n-1\}$ separately, and in each case we
assume that $M$ cannot be written as $(X,Y)$ for any such $X$.

Given this, using $\sum M = 7n-k \geq 6n+1 \geq 115$, it is routine to check
that the required decomposition $J_n\rightarrow M,k$ can be obtained using one
of the concatenations given in the table below. We can use the fact that $\sum
M = 7n -k$ to determine that $\nu_4(M) + 2\nu_5(M) \equiv (n-k)\md{3}$, and this is
also shown in the table (Note that by this, it is routine to see that for
$k=n-5$, $\nu_4(M) + 2\nu_5(M) \equiv 2\md{3}$ and thus all required
decompositions can be constructed in the manner described in the previous
paragraph). The decompositions in the third column exist by Lemma
\ref{1234567:Jgraphk=n}, and the decompositions listed in the last column are
shown in Table \ref{tab:1234567_decomps} in the appendix.

\begin{center}
\begin{tabular}{|l|l|l|l|l|}
  \hline
  $k$ & $\nu_4(M)+2\nu_5(M) \md 3$ & first decomposition & second decomposition \\
  \hline \hline
  $n-4$ & 1 & $J_{n-6} \rightarrow (M-(5^8)),(n-6)^*$ & $J_4^+ \rightarrow 5^8,2^+$  \\
  \hline
  $n-3$ & 0 & $J_{n-5} \rightarrow (M-(3^{11})),(n-5)^*$ & $J^+_5 \rightarrow 3^{11},2^+$  \\
   & & $J_{n-7} \rightarrow (M-(5^9)),(n-7)^*$ & $J^+_7 \rightarrow 5^9,4^+$  \\
  \hline
  $n-2$ & 2 & $J_{n-5} \rightarrow (M-(4^8)),(n-5)^*$ & $J^+_5 \rightarrow 4^8,3^+$  \\
  & & $J_{n-8} \rightarrow (M-(5^{10})),(n-8)^*$ & $J^+_8 \rightarrow 5^{10},6^+$  \\
  \hline
  $n-1$ & 1 & $J_{n-5} \rightarrow (M-(3^2,5^5)),(n-5)^*$ & $J^+_5 \rightarrow 3^2,5^5,4^+$   \\
  & & $J_{n-5} \rightarrow (M-(3^7,5^2)),(n-5)^*$ & $J^+_5 \rightarrow 3^7,5^2,4^+$   \\
  & & $J_{n-6} \rightarrow (M-(4^8,5)),(n-6)^*$ & $J^+_6 \rightarrow 4^8,5,5^+$   \\
  & & $J_{n-8} \rightarrow (M-(4,5^9)),(n-8)^*$ & $J^+_8 \rightarrow 4,5^9,7^+$   \\
  & & $J_{n-9} \rightarrow (M-(5^{11})),(n-9)^*$ & $J^+_9 \rightarrow 5^{11},8^+$   \\
  \hline
\end{tabular}
\end{center}
Since $n \geq 18$ it is clear that any $3$-, $4$- or $5$-cycle in this
decomposition having a vertex in $\{0,1,\ldots,6\}$ has no vertex in
$\{n,\ldots,n+6\}$, and by the definition of the decompositions given in the
fourth column the $k$-cycle has no vertex in $\{0,1,\ldots,6\}$, so these
decompositions do have the required properties.
\qed

\begin{lemma}\label{Sofsize3-8-1234567}
If $S=\{1,2,3\}$, $n\geq 15$, and $M=(m_1,\ldots,m_t,k)$ is any list satisfying
$m_i\in\{3,4,5\}$ for $i=1,\ldots,t$, $3\leq k\leq n$, and $\sum M=7n$, then
there is an $(M)$-decomposition of $\lan{S}$.
\end{lemma}

\proof As noted earlier, for $n\geq 15$ we can obtain an $(M)$-decomposition of
$\lan{\{1,2,3,4,5,6,7\}}$ from an $(M)$-decomposition of $J_n$, provided that
for each $i\in\{0,1,\ldots,6\}$, no cycle contains both vertex $i$ and vertex
$i+n$. Thus, for $S=\{1,2,3,4,5,6,7\}$, the required result follows by Lemma
\ref{1234567:Jgraph345} for $k\in\{3,4,5\}$ and by Lemma
\ref{1234567:Jgraphk<n} for $6\leq k\leq n$. \qed

\subsubsection{$S=\{1,2,3,4,5,6,7,8\}$}
In this section we show the existence of required decompositions for the case
$S=\{1,2,3,4,5,6,7,8\}$ in Lemma \ref{Sofsize3-8}.
We first define $J_n$
by
\begin{align*}E(J_n)=\{&\{i+7,i+8\},\{i+5,i+7\},\{i+5,i+8\},\{i+5,i+9\},
\{i,i+5\}, \\ &\{i+1,i+7\},\{i,i+7\},\{i+1,i+9\} :
i=0,\ldots,n-1\}\}
\end{align*}
and $V(J_n)=\{0,\ldots,n+8\}$.  We note the following basic
properties of $J_n$.  For a list of integers $M$, an $(M)$-decomposition of
$J_n$ will be denoted by
$J_n\rightarrow M$.
\begin{itemize}
    \item
For $n\geq 17$, if for each $i\in\{0,1,\ldots,8\}$ we identify vertex $i$ of $J_n$ with
vertex $i+n$ of $J_n$ then the resulting graph is $\lan{\{1,2,3,4,5,6,7,8\}}$. This means
that for $n\geq 17$, we can obtain an $(M)$-decomposition of $\lan{\{1,2,3,4,5,6,7,8\}}$
from a decomposition $J_n\rightarrow M$, provided that for each
$i\in\{0,1,\ldots,9\}$, no cycle in the decomposition of $J_n$ contains both vertex
$i$ and vertex $i+n$.
    \item
For any integers $y$ and $n$ such that $1 \leq y < n$, the graph $J_n$ is the
union of $J_{n-y}$ and the graph obtained from $J_y$ by applying the
permutation $x\mapsto x+(n-y)$.  Thus, if there is a decomposition
$J_{n-y}\rightarrow M$ and a decomposition $J_y\rightarrow M'$, then there is a
decomposition $J_n\rightarrow M,M'$. We will call this construction, and the
similar constructions that follow, \emph{concatenations}.
\end{itemize}

\begin{lemma}\label{12345678:Jgraph345}
If $n$ is a positive integer and $M=(m_1,\ldots,m_t)$ is a list such that $\sum
M=8n$, $m_i\in\{3,4,5\}$ for $i=1,\ldots,t$, and $M \neq (3^8)$ then there
is a decomposition $J_n\rightarrow M$.
\end{lemma}

\proof We first note the existence of the following decompositions, given in
Table \ref{tab:12345678_decomps} in the appendix.

\[
  \begin{array}{lllllll}
    J_{1} \rightarrow  4^{2} &\quad& J_{1} \rightarrow  3,5 &\quad&
    J_{2} \rightarrow  3^{4},4 &\quad& J_{3} \rightarrow  4,5^{4} \\
    J_{4} \rightarrow  3^{9},5 &\quad& J_{5} \rightarrow  5^{8} &\quad&
    J_{5} \rightarrow  3^{12},4 &\quad& J_{6} \rightarrow  3^{16} \\
    J_{9} \rightarrow  3^{24} &\quad& \end{array}
\]

The only required decompositions of $J_1$ are shown in the table above, so we
may assume $n\geq 2$ and assume by induction that the result holds for any
positive integer $n'$ in the range $1 \leq n' < n$.

It is routine to check that for $n\geq 2$ if $M$ satisfies the hypotheses of
the lemma, then $M$ can be written as $M=(X,Y)$ for some (possibly empty) list
$Y\neq (3^8)$ where $J_x \rightarrow X$ is one of the decompositions listed
above.  If $Y$ is empty, we are done, else we can obtain the required
decomposition by concatenation of $J_y \rightarrow Y$ (which exists by our
inductive hypothesis since $Y\neq (3^8$)) with the decomposition $J_x\rightarrow
X$.  \qed

\begin{lemma}\label{12345678:Jgraphk=6-8}
For $k\in\{6,7,8\}$ and $n \geq 3$, if $M=(m_1,\ldots,m_t)$ is
a list such that $\sum M + k = 8n$ and $m_i\in\{3,4,5\}$ for $i=1,\ldots,t$,
then there is a decomposition $J_n\rightarrow M,k$ such that $V(C_k) \subseteq
\{n,n+1,\ldots,n+9\}$.
\end{lemma}
\proof
First we note that Table \ref{tab:12345678_6-8} in the appendix lists a
number of decompositions required for this lemma. All of these decompositions
contain a $k$-cycle for some $6 \leq k \leq 8$ where the $k$-cycle is incident
on some subset of the vertices $\{n,\ldots,n+9\}$.
 some subset of the
It is routine to check that for $n\geq 3$ and any $M$ that satisfies the
hypotheses of the lemma we can write $M$ as $(X,Y)$ where $J_x \rightarrow X,k$
is one of the decompositions in Table \ref{tab:12345678_6-8}, and $Y\neq (3^8)$
is some (possibly empty) list. If $Y$ is empty we are done, else Lemma
\ref{12345678:Jgraph345} gives us the existence of a decomposition $J_y
\rightarrow Y$ and the required decomposition can be obtained by concatenation
of $J_y \rightarrow Y$ with $J_x \rightarrow X,k$. Since we concatenate with
the $k$-cycle on the right, it is clear that the $k$-cycle is still incident
upon some subset of $\{n,n+1,\ldots,n+9\}$.  \qed

\vspace{0.3cm}

Let $M=(m_1,\ldots,m_t)$ be a list of integers with $m_i\geq 3$ for $i=1,\ldots,t$. A decomposition
$\{G_1,\ldots,G_t,C\}$ of $J_n$ such that
\begin{itemize}
    \item
$G_i$ is an $m_i$-cycle for $i=1,\ldots,t$; and
    \item
$C$ is a $k$-cycle such that $V(C)=\{n-k+9,\ldots,n+8\}$ and
$\{\{n+1,n+5\},\{n+2,n+6\},\{n+3,n+7\},\{n+4,n+8\}\} \subseteq E(C)$;
\end{itemize}
will be denoted $J_n\rightarrow M,k^*$.

In Lemma \ref{12345678:Jgraphk=n} we will form new decompositions of graphs $J_n$ by
concatenating decompositions of $J_{n-y}$ with decompositions of graphs $J^+_y$
which we will now define. For $y \in \{2,\ldots,8\}$, the graph obtained from
$J_y$ by adding the edges $\{1,5\},\{2,6\},\{3,7\},\{4,8\}$ will be denoted $J^+_y$. Let
$M=(m_1,\ldots,m_t)$ be a list of integers with $m_i\geq 3$ for $i=1,\ldots,t$.
A decomposition $\{G_1,\ldots,G_t,A_1,A_2,A_3,A_4\}$ of $J^+_y$ such that
\begin{itemize}
\item $G_i$ is an $m_i$-cycle for $i=1,\ldots,t$;
\item $A_i$ is a path from $i$ to $i+4$ with $0 \notin V(A_i)$ for $i=1,2,3,4$;
\item $V(A_i) \cap V(A_j) = \emptyset$ for $i\neq j$; and
\item $|E(A_1)| + |E(A_2)| + |E(A_3)| + |E(A_4)|=l+4$;
\end{itemize}
will be denoted $J^+_y \rightarrow M,l^+$. Moreover, if $l=y$ and
$
\{\{n+1,n+5\},\{n+2,n+6\},\{n+3,n+7\},\{n+4,n+8\}\} \subseteq E(A_1) \cup
E(A_2) \cup E(A_3) \cup E(A_4)
$,
then the decomposition will be denoted $J^+_y\rightarrow M,y^{+*}$.

For $y \in \{2,\ldots,8\}$ and $n > y$, the graph $J_n$ is the union of the
graph obtained from $J_{n-y}$ by deleting the edges
$\{n-y+1,n-y+5\},\{n-y+2,n-y+6\},\{n-y+3,n-y+7\},\{n-y+4,n-y+8\}$, and the
graph obtained from $J^+_y$ applying the permutation $x\mapsto x+(n-y)$. It
follows that if there is a decomposition $J_{n-y}\rightarrow M,k^*$ and a
decomposition $J^+_y\rightarrow M',l^+$, then there is a decomposition $J_n
\rightarrow M,M',k+l$. The removed edges of the $k$-cycle in the
decomposition of $J_{n-y}$ are replaced by the four paths in the decomposition
of $J^+_y$ to form the $(k+l)$-cycle in the new decomposition. Similarly, if
there is a decomposition $J_{n-y}\rightarrow M,k^*$ and a decomposition
$J^+_y\rightarrow M',y^{+*}$, then there is a decomposition $J_n\rightarrow
M,M',(k+y)^*$.

\begin{lemma}\label{12345678:Jgraphk=9-14}
For $9 \leq k \leq 14$, if $n$ is an integer with $n \geq k$ and $M=(m_1,\ldots,m_t)$ is
a list such that $\sum M + k = 8n$ and $m_i\in\{3,4,5\}$ for $i=1,\ldots,t$,
then there is a decomposition $J_n\rightarrow M,k^*$.
\end{lemma}
\proof

For each $k$, it is routine to use the value of $k\md{8}$ to check that for
$n\geq k$ and any $M$ that satisfies the hypotheses of the lemma we can write
$M$ as $(X,Y)$ where $J_x \rightarrow X,k^*$ is one of the decompositions
given in Table \ref{tab:12345678_9-14} in the appendix, and $Y \neq (3^8)$ is
some (possibly empty) list. If $Y$ is empty, then we are done, else Lemma
\ref{12345678:Jgraph345} gives us the existence of a decomposition $J_y
\rightarrow Y$ and the required decomposition can be obtained by concatenation
of $J_y \rightarrow Y$ with $J_x \rightarrow X,k^*$.

\qed

\begin{lemma}\label{12345678:Jgraphk=n}
Given an integer $n \geq 9$, if $M=(m_1,\ldots,m_t)$ is a list such that $\sum
M=7n$ and $m_i\in\{3,4,5\}$ for $i=1,\ldots,t$, then there is a decomposition
$J_n\rightarrow M,n^*$.
\end{lemma}

\proof Lemma \ref{12345678:Jgraphk=9-14} shows that the result holds for $9\leq n
\leq 14$. So let $n \geq 15$ and suppose by induction that the result holds for
each integer $n'$ in the range $9\leq n'<n$. The following decompositions are
given in Table \ref{tab:12345678_decomps} in the appendix.

\[
  \begin{array}{lllllll}
    J_{4}^+ \rightarrow  4^{7},4^{+*} &\quad& J_{5}^+ \rightarrow  5^{7},5^{+*}
    &\quad&
    J_{6}^+ \rightarrow  3^{14},6^{+*} &\quad&
  \end{array}
\]

It is routine to check, using $\sum M = 7n \geq 105$, that $M$ can be written as
$M=(X,Y)$ where $J_y^+ \rightarrow Y,y^{+*}$ is one of the decompositions above
and $X$ is some nonempty list. We can obtain a decomposition $J_n\rightarrow
M,n^*$ by concatenating a decomposition $J_{n-y}\rightarrow X,(n-y)^*$ (which
exists by our inductive hypothesis, since $n-y \geq n-6 \geq 9$) with a
decomposition $J^+_y \rightarrow Y,y^{+*}$.
\qed

\begin{lemma}\label{12345678:Jgraphk<n}
If $n$ and $k$ are integers such that $6 \leq k \leq n$ and
$M=(m_1,\ldots,m_t)$ is a list such that $\sum M=8n-k$ and $m_i\in\{3,4,5\}$
for $i=1,\ldots,t$, then there is a decomposition $J_n\rightarrow M,k$.
Furthermore, for $n\geq 17$ all cycles in this decomposition have the property
that for $i\in\{0,1,\ldots,9\}$ no cycle is incident upon both vertex $i$ and vertex
$n+i$.
\end{lemma}

\proof We first note that if $n\geq 17$ it is clear that any $3$-, $4$- or
$5$-cycle in such a decomposition cannot be incident on two vertices $i$ and
$i+n$ for any $i\in\{0,1,\ldots,9\}$. As such, Lemma \ref{12345678:Jgraphk=6-8}
shows that the result holds for all $n$ with $k\in\{6,7,8\}$, so in the following
we deal only with $k \geq 9$.

Lemma \ref{12345678:Jgraphk=9-14} shows that the result holds
for all $n\geq k$ with $9\leq k \leq 14$ with the property that the $k$-cycle
is not incident upon any vertex in $\{0,1,\ldots,9\}$, and Lemma
\ref{12345678:Jgraphk=n} shows that the result holds for all $n=k$ with the same
property on the $k$-cycle. We can
therefore assume that $15 \leq k \leq n-1$, so let $n \geq 16$ and
suppose by induction that the result holds for each positive integer $n'$ in
the range $6\leq n'< n$ with the property that the $k$-cycle is not incident
upon any vertex in $\{0,1,\ldots,9\}$.

The following decompositions exist by Lemma \ref{12345678:Jgraph345}.

\[ \begin{array}{lllllll}
    J_{1} \rightarrow  4^{2} &\quad& J_{1} \rightarrow  3,5 &\quad&
    J_{2} \rightarrow  3^{4},4 &\quad& J_{3} \rightarrow  4,5^{4} \\
    J_{4} \rightarrow  3^{9},5 &\quad& J_{5} \rightarrow  5^{8} &\quad&
    J_{5} \rightarrow  3^{12},4 &\quad& J_{6} \rightarrow  3^{16}
  \end{array}
\]

\noindent{\bf Case 1}\quad Suppose that $k \leq n-6$. Then it is routine to
check, using $\sum M = 8n-k \geq 7n+6 \geq 118$, that $M=(X,Y)$ where $J_y
\rightarrow Y$ is one of the decompositions above and $X$ is some nonempty
list. We can obtain a decomposition $J_n\rightarrow M,k$ by concatenating a
decomposition $J_{n-y}\rightarrow X,k$ (which exists by our inductive
hypothesis, since $k \leq n-6 \leq n-y$) with a decomposition $J_y\rightarrow
Y$. Since $n \geq 15$ it is clear that any $3$-, $4$- or $5$-cycle in this
decomposition having a vertex in $\{0,1,\ldots,9\}$ has no vertex in
$\{n,n+1,\ldots,n+9\}$, and by our inductive hypothesis the same holds for the
$k$-cycle.

\noindent{\bf Case 2}\quad Suppose that $n-5 \leq k \leq n-1$. In a similar
manner to Case 1, we can obtain the required decomposition $J_n\rightarrow M,k$
if $M=(X,Y)$ for some list $X$ where $J_x \rightarrow X$ is one of the
decompositions shown above and $k+x \leq n$. We can therefore assume that $M$
cannot be written as $(X,Y)$ for any such list $X$. In particular,
since $J_1 \rightarrow 4^2$ exists we can assume $\nu_4(M) \leq 1$, and since
$J_1 \rightarrow 3,5$ exists we can assume either $\nu_3(M)=0$ or $\nu_5(M)=0$.
As a result, using $\sum M = 8n-k \geq 7n+1 \geq 113$ we have either $\nu_3(M)
\geq 33$ or $\nu_5(M) \geq 20$.

Given this, it is routine to check that the required decomposition
$J_n\rightarrow M,k$ can be obtained using one of the concatenations given in
the table below.

The decompositions in the second column exist by Lemma \ref{12345678:Jgraphk=n}
(since $k\geq 15$), and the decompositions listed in the last column are shown
in Table \ref{tab:12345678_decomps} in the appendix.

\begin{center}
\begin{tabular}{|l|l|l|}
  \hline
  $k$ & first decomposition & second decomposition \\
  \hline \hline
  $n-5$ & $J_{n-7} \rightarrow (M-(3^{18})),(n-7)^*$ & $J_7^+ \rightarrow 3^{14},2^+$  \\
  \hline
  $n-4$ & $J_{n-8} \rightarrow (M-(5^{12})),(n-8)^*$ & $J_8^+ \rightarrow 5^{12},4^+$  \\
        & $J_{n-5} \rightarrow (M-(3^{13})),(n-5)^*$ & $J_5^+ \rightarrow 3^{13},1^+$  \\
  \hline
  $n-3$ & $J_{n-6} \rightarrow (M-(5^{9})),(n-6)^*$ & $J^+_6 \rightarrow 5^{9},3^+$  \\
        & $J_{n-6} \rightarrow (M-(3^{15})),(n-6)^*$ & $J^+_6 \rightarrow 3^{15},3^+$  \\
  \hline
  $n-2$ & $J_{n-4} \rightarrow (M-(5^6)),(n-4)^*$ & $J^+_4 \rightarrow 5^6,2^+$  \\
        & $J_{n-4} \rightarrow (M-(3^{10})),(n-4)^*$ & $J^+_4 \rightarrow 3^{10},2^+$  \\
  \hline
  $n-1$ & $J_{n-2} \rightarrow (M-(5^3)),(n-2)^*$ & $J^+_2 \rightarrow 5^3,1^+$   \\
        & $J_{n-2} \rightarrow (M-(3^5)),(n-2)^*$ & $J^+_2 \rightarrow 3^5,1^+$   \\
  \hline
\end{tabular}
\end{center}
If $n \geq 17$ it is clear that any $3$-, $4$- or $5$-cycle in this
decomposition having a vertex in $\{0,1,\ldots,9\}$ has no vertex in
$\{n,\ldots,n+9\}$, and by the definition of the decompositions given in the
second column the $k$-cycle has no vertex in $\{0,1,\ldots,9\}$, so these
decompositions do have the required properties.
\qed

\begin{lemma}\label{Sofsize3-8-12345678}
If $S=\{1,2,3,4,5,6,7,8\}$, $n\geq 17$, and $M=(m_1,\ldots,m_t,k)$ is any list
satisfying $m_i\in\{3,4,5\}$ for $i=1,\ldots,t$, $3\leq k\leq n$, and $\sum
M=8n$, then there is an $(M)$-decomposition of $\lan{S}$.
\end{lemma}

\proof As noted
earlier, for $n\geq 17$ we can obtain an $(M)$-decomposition of
$\lan{\{1,2,3,4,5,6,7,8\}}$ from an $(M)$-decomposition of $J_n$, provided that for
each $i\in\{0,1,\ldots,8\}$, no cycle contains both vertex $i$ and vertex
$i+n$. Thus, for $S=\{1,2,3,4,5,6,7,8\}$, the required result follows by Lemma
\ref{12345678:Jgraph345} for $k\in\{3,4,5\}$ and by Lemma
\ref{12345678:Jgraphk<n} for $6\leq k\leq n$. \qed

\vspace{0.3cm}

We now prove Lemma \ref{Sofsize3-8}.

\noindent{\bf Lemma \ref{Sofsize3-8}}\quad
{\em
If
$$S\in\{\{1,2,3\},\{1,2,3,4\},\{1,2,3,4,6\},\{1,2,3,4,5,7\},\{1,2,3,4,5,6,7\},\{1,2,3,4,5,6,7,8\}\},$$
$n\geq 2\max(S)+1$, and $M=(m_1,\ldots,m_t,k)$ is any list satisfying $m_i\in\{3,4,5\}$ for
$i=1,\ldots,t$, $3\leq k\leq n$, and $\sum M=|S|n$, then there is an $(M)$-decomposition of
$\lan{S}$, except possibly when
\begin{itemize}
\item $S=\{1,2,3,4,6\}$, $n\equiv 3\md 6$ and $M=(3^{\frac{5n}3})$; or
\item $S=\{1,2,3,4,6\}$, $n\equiv 4\md 6$ and $M=(3^{\frac{5n-5}3},5)$.
\end{itemize}
}

\proof The required decompositions exist
by Lemma \ref{Sofsize3-8-123} for $S=\{1,2,3\}$,
by Lemma \ref{Sofsize3-8-1234} for $S=\{1,2,3,4\}$,
by Lemma \ref{Sofsize3-8-12346} for $S=\{1,2,3,4,6\}$,
by Lemma \ref{Sofsize3-8-123457} for $S=\{1,2,3,4,5,7\}$,
by Lemma \ref{Sofsize3-8-1234567} for $S=\{1,2,3,4,5,6,7\}$, and
by Lemma \ref{Sofsize3-8-12345678} for $S=\{1,2,3,4,5,6,7,8\}$. \qed

\subsection{Proof of Lemma \ref{S=123includingham}}\label{sec:123n}

In this section we prove Lemma \ref{S=123includingham}, which we restate here for convenience.

\vspace{0.5cm}

\noindent{\bf Lemma \ref{S=123includingham}}\quad {\em If $n\geq 7$ and $M=(m_1,\ldots,m_t,k,n)$ is
any list satisfying $m_i\in\{3,4,5\}$ for $i=1,\ldots,t$, $3\leq k\leq n$, and $\sum M=3n$, then
there is an $(M)$-decomposition of $\lan{\{1,2,3\}}$. }

\vspace{0.5cm}

The proof of Lemma \ref{S=123includingham} proceeds along similar lines to the proof of Lemma
\ref{Sofsize3-8}. We make use of the graphs $J^{\{1,2,3\}}_n$ defined in the proof of Lemma
\ref{Sofsize3-8}, which in this subsection we denote by just $J_n$. Recall that $J^{\{1,2,3\}}_n$
is the graph with vertex set $\{0,\ldots,n+2\}$ and edge set
$$\{\{i,i+1\},\{i+1,i+3\},\{i,i+3\}:i =0,\ldots,n-1\}.$$
We first construct decompositions of graphs which are related to the graphs $J_n$, then concatenate
these decompositions to produce decompositions of the graphs $J_n$, and finally identify pairs of
vertices to produce the required decompositions of $\lan{\{1,2,3\}}$.

For $n\geq 1$, the graph obtained from $J_n$ by adding the edges $\{n,n+1\}$ and $\{n+1,n+2\}$ will
be denoted $L_n$. Let $M=(m_1,\ldots,m_t)$ be a list of integers with $m_i\geq 3$ for
$i=1,\ldots,t$. A decomposition $\{G_1,\ldots,G_t,A,B\}$ of $L_n$ such that
\begin{itemize}
\item $G_i$ is an $m_i$-cycle for $i=1,\ldots,t$;
\item $A$ is a path of length $k$ from $n$ to $n+1$;
and
\item $B$ is a path of length $n-1$ from $n+1$ to $n+2$ such that
$0,1,2\notin V(B)$;
\end{itemize}
will be denoted $L_n \rightarrow M,k^+,(n-1)^H$.

In Lemmas \ref{Lgraphsk=n} and \ref{Lgraphsksmaller} we form new decompositions of graphs $L_n$ by
concatenating decompositions of $L_{n-y}$ with decompositions of graphs $P_y$ which we will now
define. For $y \in \{3,4,5,6\}$, the graph obtained from $J_y$ by deleting the edges in
$\{\{0,1\},\{1,2\}\}$ and adding the edges in $\{\{y,y+1\},\{y+1,y+2\}\}$ will be denoted $P_y$.
Let $M=(m_1,\ldots,m_t)$ be a list of integers with $m_i\geq 3$ for $i=1,\ldots,t$. A decomposition
$\{G_1,\ldots,G_t,A_1,A_2,B_1,B_2\}$ of $P_y$ such that
\begin{itemize}
\item
$G_i$ is an $m_i$-cycle for $i=1,2,\ldots,t$;
\item $A_1$ and $A_2$ are vertex-disjoint paths with endpoints $0$, $1$, $y$ and $y+1$,
such that $A_1$ has endpoints $0$ and $y$ or $0$ and $y+1$;
\item $|E(A_1)|+|E(A_2)|=k'$ and $2\notin V(A_1)\cup V(A_2)$;
\item $B_1$ and $B_2$ are vertex-disjoint paths with endpoints $1$, $2$, $y+1$ and $y+2$,
such that $B_1$ has endpoints $1$ and $y+1$ or $1$ and $y+2$; and
\item $|E(B_1)|+|E(B_2)|=y$, and $0\notin V(B_1)\cup V(B_2)$;
\end{itemize}
will be denoted $P_y \rightarrow M,k'^+,y^H$.

For $y \in \{3,4,5,6\}$ and $n >y$, the graph $L_n$ is the union of the graph $L_{n-y}$ and the
graph obtained from $P_y$ by applying the permutation $x\mapsto x+(n-y)$. It follows that if there
is a decomposition $L_{n-y} \rightarrow M,k^+,(n-y-1)^H$ and a decomposition $P_y \rightarrow
M',k'^+,y^H$, then there is a decomposition $L_n \rightarrow M,M',(k+k')^+,(n-1)^H$.

\begin{lemma}\label{Lgraphsk=n}
If $n\geq 2$ is an integer and $M=(m_1,\ldots,m_t)$ is a list such that $\sum M=n+1$,  $M \notin
\{(3^{i}):\mbox{$i$ is even}\}$ and $m_i\in\{3,4,5\}$ for $i=1,\ldots,t$, then there is a
decomposition $L_n\rightarrow M,(n+2)^+,(n-1)^H$.
\end{lemma}

\proof The following decompositions are given in full detail in Table
\ref{tab:123n_decomps} in the appendix, thus verifying the lemma for
$n\in\{2,3,4\}$.
$$\begin{array}{lllll}
L_2 \rightarrow 3,4^+,1^H &\quad&
L_3 \rightarrow 4,5^+,2^H &\quad&
L_4 \rightarrow 5,6^+,3^H \\
P_4 \rightarrow 4,4^+,4^H &\quad&
P_5 \rightarrow 5,5^+,5^H &\quad&
P_6 \rightarrow 3^2,6^+,6^H
\end{array}$$

So let $n\geq 5$ and assume by induction that the result holds for each integer
$n'$ in the range $2\leq n'<n$. It is routine to check that for $n\geq 5$, if
$M$ satisfies the hypotheses of the lemma, then $M$ can be written as $M =
(X,Y)$ where $n-y \geq 2$, $X \notin \{(3^{i}):\mbox{$i$ is even}\}$, and $P_y
\rightarrow Y,y^+,y^H$ is one of the decompositions above. We can obtain the
required decomposition $L_n \rightarrow M,(n+2)^+,(n-1)^H$ by concatenating a
decomposition $L_{n-y}\rightarrow X,(n-y+2)^+,(n-y-1)^H$ (which exists by our
inductive hypothesis) with a decomposition $P_y \rightarrow Y,y^+,y^H$. \qed

\vspace{0.5cm}

\begin{lemma}\label{Lgraphsksmaller}
If $n$ and $k$ are positive integers with $\frac{4n+12}{5} \leq k\leq n+2$, and
$M=(m_1,\ldots,m_t)$ is a list such that $\sum M=2n-k+3$, $M \notin \{(3^{i}):\mbox{$i$ is even}\}$
and $m_i\in\{3,4,5\}$ for $i=1,\ldots,t$, then there is a decomposition $L_n\rightarrow
M,k^+,(n-1)^H$.
\end{lemma}

\proof The proof will be by induction on $j=n-k+2$. For a given $n$ we need to
prove the result for each integer $j$ in the range $0\leq j\leq\frac{n-2}5$.
The case $j=0$ is covered in Lemma \ref{Lgraphsk=n}, so assume $1\leq
j\leq\frac{n-2}5$ and that the result holds for each integer $j'$ in the range
$0\leq j'<j$. Note that, since $\frac{4n+12}{5} \leq k$ and $j\geq 1$, we have
$n \geq 7$. The following decompositions are given in Table
\ref{tab:123n_decomps} in the appendix.
$$\begin{array}{lllll}
P_3 \rightarrow 4,2^+,3^H &\quad&
P_4 \rightarrow 5,3^+,4^H &\quad&
P_5 \rightarrow 3^2,4^+,5^H
\end{array}$$

It is routine to check that for $j\leq \frac{n-2}5$, if $M$ satisfies the hypotheses of the lemma,
then $M$ can be written as $M = (X,Y)$ where $X\not\in\{(3^i): \mbox{$i$ is even}\}$ and $P_y
\rightarrow Y,(y-1)^+,y^H$ is one of the decompositions listed above. A decomposition
$L_{n-y}\rightarrow X,(k-y+1)^+,(n-y-1)^H$ will exist by our inductive hypothesis provided that
$$\tfrac{4(n-y)+12}{5} \leq k-y+1 \leq n-y+2$$
and it is routine to check that this holds using $\frac{4n+12}{5} \leq k$, $j \geq 1$ and $y \in
\{3,4,5\}$. Thus, the required decomposition $L_n\rightarrow M,k^+,(n-1)^H$ can be obtained by
concatenating a decomposition $L_{n-y}\rightarrow X,(k-y+1)^+,(n-y-1)^H$ with a decomposition $P_y
\rightarrow Y,(y-1)^+,y^H$. \qed

\vspace{0.3cm}

Let $(m_1,\ldots,m_t)$ be a list of integers with $m_i\geq 3$ for $i=1,\ldots,t$. A decomposition of $J_n$ into partitions
$\{G_1,\ldots,G_t,H\}$ such that
\begin{itemize}
\item $G_i$ is an $m_i$-cycle for $i=1,\ldots,t$; and
\item $H$ is an $n$-cycle such that $0,1,2 \notin V(H)$ and $\{n,n+2\} \in E(H)$;
\end{itemize}
will be denoted $J_n \rightarrow m_1,\ldots,m_t,n^H$.

In Lemma \ref{JgraphsClosek} we will form decompositions of graphs $J_n$ by concatenating
decompositions of graphs $L_{n-y}$ obtained from Lemma \ref{Lgraphsksmaller} with decompositions of
graphs $Q_y$ which we will now define. For each $y \in \{4,5,6\}$, the graph obtained from $J_y$ by
deleting the edges $\{0,1\}$ and $\{1,2\}$ will be denoted by $Q_y$. Let $M=(m_1,\ldots,m_t)$ be a
list of integers with $m_i \geq 3$ for $i = 1,\ldots,t$. A decomposition $\{G_1,\ldots,G_t,A,B\}$
of $Q_y$ such that
\begin{itemize}
\item $G_i$ is an $m_i$-cycle for $i=1,\ldots,t$;
\item $A$ is a path of length $k'$ from $0$ to $1$ such that
 $\{2,y,y+1,y+2\}\notin V(A)$; and
\item $B$ is a path of length $y+1$ from $1$ to $2$ such that $0\not\in V(B)$ and $\{y,y+2\}
    \in E(B)$;
\end{itemize}
will be denoted $Q_y \rightarrow M,k'^+,(y+1)^H$.

For $y \in \{4,5,6\}$ and $n>y$, the graph $J_n$ is the union of the graph $L_{n-y}$ and the graph
obtained from $Q_y$ by applying the permutation $x\mapsto x+(n-y)$. It follows that if there is a
decomposition $L_{n-y} \rightarrow M,k^+,(n-y-1)^H$ and a decomposition $Q_y \rightarrow
M',k'^+,(y+1)^H$, then there is a decomposition $J_n \rightarrow M,M',k+k',n^H$. Note that, for $y
\in \{4,5,6\}$ and $n-y \geq 3$, no cycle of this decomposition contains both vertex $i$ and vertex
$i+n$ for $i\in\{0,1,2\}$.

\begin{lemma}\label{123:JgraphsClosekSmalln}
If $n$ and $k$ are integers with $6 \leq n \leq 32$, $k \geq 6$ and $n-5 \leq k
\leq n$ and $M=(m_1,\ldots,m_t)$ is a list such that $\sum M=2n-k$,
$m_i\in\{3,4,5\}$ for $i=1,\ldots,t$, then there is a decomposition
$J_n\rightarrow M,k,n^H$ such that for $i\in\{0,1,2\}$ no cycle of the
decomposition contains both vertex $i$ and vertex $i+n$.
\end{lemma}

\proof
First note the existence of the following decompositions, given in full detail
in Table \ref{tab:123n_decomps} in the appendix.
$$
\begin{array}{lllll}
Q_4 \rightarrow 3,2^+,5^H &\quad&
Q_5 \rightarrow 4,3^+,6^H &\quad&
Q_6 \rightarrow 5,4^+,7^H
\end{array}
$$
Let $j=n-k$, then by the conditions of the lemma we have $0 \leq j \leq 5$.
Note that for any $M$ that satisfies the conditions of the lemma, if $M$ can be
written as $M=(X,y-1)$ where $L_{n-y} \rightarrow
X,(k-y+2)^+,(n-y-1)^H$ exists and $Q_y \rightarrow (y-1),(y-2)^+,(y+1)^H$ is
one of the decompositions listed above, then we can construct the required
decomposition by concatenating $L_{n-y} \rightarrow X,(k-y+2)^+,(n-y-1)^H$ with
$Q_y \rightarrow (y-1),(y-2)^+,(y+1)^H$.
We can write $L_{n-y} \rightarrow X,(k-y+2)^+,(n-y-1)^H$ as $L_{n-y}
\rightarrow X,(n-j-y+2)^+,(n-y-1)^H$. For $j=0$, such decompositions exist by
Lemma \ref{Lgraphsk=n} and for $1 \leq j \leq 5$ such decompositions exist for $n \geq
12  + 5j$ by Lemma \ref{Lgraphsksmaller}.
Therefore, in the following we assume $1 \leq j \leq 5$ and $n < 12 + 5i$.

In tables
\ref{tab:123n_concats_1}, \ref{tab:123n_concats_2},
\ref{tab:123n_concats_3}, \ref{tab:123n_concats_4} and
\ref{tab:123n_concats_5}, we list all required decompositions for $j=1,2,3,4$
and $5$ respectively. That is, each table lists all required decompositions of
the form $J_n \rightarrow M,k,n$ where $6 \leq k \leq n<
12 + 5j$ for the given value of $j=n-k$. \qed

\begin{lemma}\label{JgraphsClosek}
If $n$ and $k$ are integers with $n \geq 6$, $k \geq 3$ and $n-5 \leq k \leq n$ and
$M=(m_1,\ldots,m_t)$ is a list such that $\sum M=2n-k$, $m_i\in\{3,4,5\}$ for $i=1,\ldots,t$, then
there is a decomposition $J_n\rightarrow M,k,n^H$ such that for $i\in\{0,1,2\}$ no cycle of the
decomposition contains both vertex $i$ and vertex $i+n$.
\end{lemma}

\proof For $k\in\{3,4,5\}$ the result holds by Lemma \ref{123:Jgraphk<n} (by
letting $n$ in this theorem be $k$ in Lemma \ref{123:Jgraphk=n})
We therefore assume $6 \leq k \leq n$.

For $6 \leq n \leq 32$ the required result holds by Lemma
\ref{123:JgraphsClosekSmalln}.  Given this, we may assume $n \geq 33$. The special
case where $M \in \{(3^i):\mbox{$i$ is odd}\}$ will be dealt with separately in
a moment.

\noindent{\bf Case 1}\quad Suppose that $M \not\in \{(3^i):\mbox{$i$ is
odd}\}$.
The following decompositions are given in full detail in Table
\ref{tab:123n_decomps} in the appendix.
$$
\begin{array}{lllll}
Q_4 \rightarrow 3,2^+,5^H &\quad&
Q_5 \rightarrow 4,3^+,6^H &\quad&
Q_6 \rightarrow 5,4^+,7^H
\end{array}
$$
It is routine to check for $n\geq 33$, if $M$ satisfies the hypotheses of the
lemma (and $M \notin \{(3^i):\mbox{$i$ is odd}\}$), then $M$ can be written as
$M=(X,y-1)$ where $X \notin \{(3^i):\mbox{$i$ is even}\}$ and $Q_y \rightarrow
(y-1),(y-2)^+,(y+1)^H$ is one of the decompositions listed in Lemma
\ref{123:JgraphsClosekSmalln} above. Using $n \geq 33$ and $y\in\{4,5,6\}$,
it can be verified that a decomposition $L_{n-y} \rightarrow
X,(k-y+2)^+,(n-y-1)^H$ exists by Lemma \ref{Lgraphsksmaller}.  Concatenation of
this decomposition with $Q_y \rightarrow (y-1),(y-2)^+,(y+1)^H$ yields the
required decomposition $J_n\rightarrow M,k,n^H$.

\noindent{\bf Case 2}\quad Suppose that $M \in \{(3^i):\mbox{$i$ is odd}\}$. Let
$p=\frac{i-3}2-(n-k)$. We deal separately with the case $n=k$ and the case
$n\in\{k+1,k+2,k+3,k+4,k+5\}$.

\noindent{\bf Case 2a}\quad Suppose that $n=k$. Note that since $n\geq 33$ and $\sum M=3i=2n-k$, we
have $p\geq 4$ when $n=k$. The set of $3$-cycles in the decomposition is the union of the following
two sets.
$$
\begin{array}{c}
\{ (0,1,3),(2,4,5), (n-3,n-2,n-1)\} \\
\{(6j+6,6j+7,6j+8),(6j+9,6j+10,6j+11) : j \in \{0,\ldots,p-1\}\}
\\
\end{array}
$$
The edge set of one $n$-cycle is $E_1\cup E_2\cup E_3$ where
$$
\begin{array}{cl}
E_1=&\{ \{5,3\},\{3,4\},\{4,6\}\},\\
E_2=&\{\{n-4,n-2\},\{n-2,n+1\},\{n+1,n-1\},
\{n-1,n+2\},\{n+2,n\},\{n,n-3\}\}, \\
E_3=&\{\{6j+6, 6j+9\},\{6j+9, 6j+8\},
\{6j+8,6j+11\},\\&
\{6j+5,6j+7\},\{6j+7,6j+10\},\{6j+10,6j+12\}: j \in \{0,\ldots,p-1\}\}.
\end{array}
$$
Note that this $n$-cycle contains the edge $\{n,n+2\}$ and does not contain any of the vertices in
$\{0,1,2\}$. The remaining edges form the edge set of the other $n$-cycle (here $n=k$).

\noindent{\bf Case 2b}\quad Suppose that $n\in\{k+1,k+2,k+3,k+4,k+5\}$. Since
$n\geq 33$, it is routine to verify that for any integers $n$ and $k$ and list
$M\in \{(3^i):\mbox{$i$ is odd}\}$ which satisfy the conditions of the lemma we
have $p\geq 1$, except in the case where $M=(3^{13})$ and $n=k+5$. In this
special case we have $(n,k)=(34,29)$ and we have constructed the decomposition
required in this case explicitly and it is listed in Table
\ref{tab:123n_concats_5} in the appendix.  Thus we can assume $p\geq 1$. Let $l
= 5(n-k)$. The set of $3$-cycles in the decomposition is the union of the
following three sets.
$$
\begin{array}{c}
\{ (0,1,3),(2,4,5), (n-3,n-2,n-1)\}  \\
\{(5j+6,5j+7,5j+8),(5j+9,5j+10,5j+11) : j \in \{0,\ldots,(n-k)-1\} \} \\
\{(6j+l+6,6j+l+7,6j+l+8),(6j+l+9,6j+l+10,6j+l+11) : j \in \{0,\ldots,p-1\}\} \\
\end{array}
$$
The edge set of the $n$-cycle is $E_1 \cup E_2 \cup E_3 \cup E_4$ where
$$
\begin{array}{cl}
E_1=&\{\{5,3\},\{3,4\},\{4,6\}\}, \\
E_2=&\{\{n-4,n-2\},\{n-2,n+1\},\{n+1,n-1\},
     \{n-1,n+2\},\{n+2,n\},\{n,n-3\}\}, \\
E_3=&\{\{5j+6,5j+9\}, \{5j+9, 5j+8\}, \{5j+8, 5j+11\} , \\
&\{5j+5, 5j+7\}, \{5j+7,5j+10\} : j \in \{0,\ldots,(n-k)-1\}\},  \\
E_4=&\{\{6j+l+6, 6j+l+9\},\{6j+l+9, 6j+l+8\}, \\&
\{6j+l+8,6j+l+11\}, \{6j+l+5,6j+l+7\},  \\&
\{6j+l+7,6j+l+10\},\{6j+l+10,6j+l+12\} : j \in \{0,\ldots,p-1\}\}.
\end{array}
$$
Note that this $n$-cycle contains the edge $\{n,n+2\}$ and does not contain any of the vertices in
$\{0,1,2\}$. The remaining edges form the edge set of the cycle of length $k$. \qed

\vspace{0.5cm}

In Lemma \ref{JgraphsAnyk} we will form decompositions of graphs $J_n$ by concatenating
decompositions of $J_{n-y}$ with decompositions of graphs $R_y$ which we will now define. For $y
\in \{5,6\}$, the graph obtained from $J_y$ by adding the edge $\{0,2\}$ will be denoted by $R_y$.
Let $M=(m_1,\ldots,m_t)$ be a list of integers with $m_i \geq 3$ for $i=1,\ldots,t$. A
decomposition $\{G_1,\ldots,G_t,A\}$ of $R_y$ such that
\begin{itemize}
\item $G_i$ is an $m_i$-cycle for $i = 1,\ldots,t$; and
\item $A$ is a path of length $y+1$ from $0$ to $2$ such that $1\not\in V(A)$ and $\{y,y+2\}\in
    E(A)$;
\end{itemize}
will be denoted $R_y \rightarrow M,y^H$.

For $y \in \{5,6\}$ and $n >y$, the graph $J_n$ is the union of the graph obtained from $J_{n-y}$
by removing the edge $\{n-y,n-y+2\}$ and the graph obtained from $R_y$ by applying the permutation
$x\mapsto x+(n-y)$. It follows that if there is a decomposition $J_{n-y}\rightarrow M,k,(n-y)^H$
and a decomposition $R_y \rightarrow M',y^H$, then there is a decomposition $J_n\rightarrow
M,M',k,n^H$. In this construction the edge $\{n-y,n-y+2\}$ in the $(n-y)$-cycle of the
decomposition of $J_{n-y}$ is replaced with the path from the decomposition of $R_y$ to form the
$n$-cycle in the decomposition of $J_n$. Note that, for $y \in \{5,6\}$ and $n-y \geq 3$, no cycle
of the decomposition contains both vertex $i$ and vertex $i+n$ for $i\in\{0,1,2\}$.

\begin{lemma}\label{JgraphsAnyk}
Let $n$ and $k$ be integers with $6 \leq k \leq n$. If $M=(m_1,\ldots,m_t)$ is a list such that
$\sum M=2n-k$ and $m_i\in\{3,4,5\}$ for $i=1,\ldots,t$, then there is a decomposition
$J_n\rightarrow M,k,n^H$ such that for $i\in\{0,1,2\}$ no cycle of the decomposition contains both
vertex $i$ and vertex $i+n$.
\end{lemma}

\proof If $k \geq n-5$, then the result follows by Lemma \ref{JgraphsClosek}, which means the
result holds for $n \leq 11$. We can therefore assume that $k\leq n-6$, $n \geq 12$ and, by
induction, that the result holds for each integer $n'$ in the range $6 \leq n'<n$.

The following decompositions are given in full detail in Table
\ref{tab:123n_decomps} in the appendix.
$$
\begin{array}{lllllll}
R_5 \rightarrow 5^2,5^H &\quad&
R_6 \rightarrow 3,4,5,6^H &\quad&
R_6 \rightarrow 4^3,6^H &\quad&
R_6 \rightarrow 3^4,6^H
\end{array}
$$

It is routine to check that if $M$ satisfies the conditions of the lemma, then $M$ can be written
as $M =(X,Y)$ where $R_y \rightarrow Y,y^H$ is one of the decompositions above. The required
decomposition can be obtained by concatenating a decomposition $J_{n-y} \rightarrow X,k,(n-y)^H$
(which exists by our inductive hypothesis since $k \leq n-6 \leq n-y$) with a decomposition $R_y
\rightarrow Y,y^H$. \qed

\vspace{0.5cm}

\noindent{\bf Proof of Lemma \ref{S=123includingham}}\quad If $k\in\{3,4,5\}$, then the result
follows by Lemma \ref{Sofsize3-8}. So we can assume $k\geq 6$. Since $n\geq 7$, we can obtain an
$(M)$-decomposition of $\lan{\{1,2,3\}}$ from an $(M)$-decomposition of $J_n$ by identifying vertex
$i$ with vertex $i+n$ for each $i\in\{0,1,2\}$, provided that for each $i\in\{0,1,2\}$, no cycle of
our decomposition contains both vertex $i$ and vertex $i+n$. Thus, Lemma \ref{S=123includingham}
follows immediately from Lemma \ref{JgraphsAnyk}. \qed

\section{Decompositions of $K_n-\lan{S}$}\label{decomposeKminusS}

The purpose of this section is to prove Lemmas \ref{3s4sandHamsforKn-S} and \ref{5sandHamsKn-S},
and these proofs are given in Subsections \ref{3s4sHamsProof} and \ref{5sHamsProof} respectively.
In Subsection \ref{hamcycledecompsofcayleygraphs} we present results on Hamilton decompositions of
circulant graphs that we will require.

To prove Lemma \ref{3s4sandHamsforKn-S}, we require a $(3^{tn},4^{qn},n^h)$-decomposition of
$\lan{\overline S}$, where $\overline S=\{1,\ldots,\floor{\frac n2}\}\setminus S$, for almost all
$n$, $t$, $q$ and $h$ satisfying $h\geq 2$, $n\geq 2\max(S)+1$ and
$3t+4q+h=\floor{\frac{n-1}2}-|S|$. To construct this, $\overline S$ will be partitioned into three
subsets $S_1$, $S_2$ and $S_3$ such that there is a $(3^{tn})$-decomposition of $\lan{S_1}$, a
$(4^{qn})$-decomposition of $\lan{S_2}$, and an $(n^h)$-decomposition of $\lan{S_3}$. Our
$(3^{tn})$-decompositions of $\lan{S_1}$ are constructed by partitioning $S_1$ into modulo $n$
difference triples, our $(4^{qn})$-decompositions of $\lan{S_2}$ are constructed by partitioning
$S_2$ into modulo $n$ difference quadruples, and our $(n^h)$-decompositions of $\lan{S_3}$ are
constructed by partitioning $S_3$ into sets of size at most $3$ to yield connected circulant graphs
of degree at most $6$ that are known to have Hamilton decompositions. Lemma \ref{5sandHamsKn-S} is
proved in an analogous manner.

\subsection{Decompositions of circulant graphs into Hamilton cycles}\label{hamcycledecompsofcayleygraphs}

Theorems \ref{BFMTheorem}--\ref{MDeanTheorem} address the open problem of whether every connected
Cayley graph on a finite abelian group has a Hamilton decomposition \cite{Als2}. Note that $\lan S$
is connected if and only if $\gcd(S\cup\{n\})=1$.

\begin{theorem}\label{BFMTheorem}{\rm (\cite{BerFavMah})}
Every connected $4$-regular Cayley graph on a finite abelian group has a decomposition into two
Hamilton cycles.
\end{theorem}

The following theorem is an easy corollary of Theorem \ref{BFMTheorem}.

\begin{theorem}\label{5regCayintoHams}
Every connected $5$-regular Cayley graph on a finite abelian group has a decomposition into two
Hamilton cycles and a perfect matching.
\end{theorem}

\proof Let the graph be $X=\cay(\Gamma,S)$. Since each vertex of $X$ has odd degree, $S$ contains
an element $s$ of order $2$ in $\Gamma$. Let $F$ be the perfect matching of $X$ generated by $s$.
If $\cay(\Gamma,S\setminus\{s\})$ is connected then, as it is also $4$-regular, the result follows
immediately from Theorem \ref{BFMTheorem}. On the other hand, if $\cay(\Gamma,S\setminus\{s\})$ is
not connected, then it consists of two isomorphic connected components, with $x \mapsto sx$ being
an isomorphism. These components are $4$-regular and so by Theorem \ref{BFMTheorem}, each can be
decomposed into two Hamilton cycles. Moreover, since $x \mapsto sx$ is an isomorphism, there exists
a Hamilton decomposition $\{H_1,H'_1\}$ of the first and a Hamilton decomposition $\{H_2,H'_2\}$ of
the second such that there is a pair of vertex-disjoint $4$-cycles $(x_1,y_1,y_2,x_2)$ and
$(x'_1,y'_1,y'_2,x'_2)$ in $X$ with $x_1y_1\in E(H_1)$, $x'_1y'_1\in E(H'_1)$, $x_2y_2\in E(H_2)$,
$x'_2y'_2\in E(H'_2)$, and $x_1x_2,y_1y_2,x'_1x'_2,y'_1y'_2\in E(F)$. It follows that if we let $G$
be the graph with edge set
$$(E(H_1)\cup E(H_2)\cup\{x_1x_2,y_1y_2\})\setminus\{x_1y_1,x_2y_2\}$$
and let $G'$ be the graph with edge set
$$(E(H'_1)\cup E(H'_2)\cup\{x'_1x'_2,y'_1y'_2\})\setminus\{x'_1y'_1,x'_2y'_2\},$$
then $G$ and $G'$ are edge-disjoint Hamilton cycles in $X$. This proves the result.
\qed

\begin{theorem}\label{MDeanTheorem}{\rm (\cite{Dea})}
Every $6$-regular Cayley graph on a group which is generated by an element
of the connection set has a decomposition into three Hamilton cycles.
\end{theorem}

This theorem implies that, for distinct $a,b,c \in \{1,\ldots,\floor{\frac{n-1}{2}}\}$, the graph
$\lan{\{a,b,c\}}$ has a decomposition into three Hamilton cycles if $\gcd(x,n)=1$ for some $x \in
\{a,b,c\}$.

In the next two lemmas we give results similar to that of Lemma \ref{hamdecompxtonover2}, but for
the case where the connection set is of the form $\{x-1\} \cup \{x+1,\ldots,\floor{\frac n2}\}$
rather than $\{x,\ldots,\floor{\frac n2}\}$. Lemma \ref{hamdecompxtonover2withgapnodd} deals with
the case $n$ is odd, and Lemma \ref{hamdecompxtonover2withgapneven} deals with the case $n$ is
even.

\begin{lemma}\label{hamdecompxtonover2withgapnodd}
If $n$ is odd and $1\leq h\leq\frac{n-3}2$, then there is an $(n^h)$-decomposition of
$\langle\{\frac{n-1}2-h\} \cup \{\frac{n-1}2-h+2,\ldots,\frac{n-1}2\}\rangle_n$; except when $h=1$
and $n\equiv 3\md 6$ in which case the graph is not connected.
\end{lemma}

\proof If $h=1$, then the graph is $\lan{\frac{n-3}2}$. If $n \equiv 1,5 \md{6}$, then
$\gcd(\frac{n-3}2,n)=1$ and $\lan{\frac{n-3}2}$ is an $n$-cycle. If $n \equiv 3 \md{6}$, then
$\gcd(\frac{n-3}2,n)=3$ and $\lan{\frac{n-3}2}$ is not connected. Thus the result holds for $h=1$.
In the remainder of the proof we assume $h\geq 2$.

We first decompose $\langle\{\frac{n-1}2-h\} \cup \{\frac{n-1}2-h+2,\ldots,\frac{n-1}2\}\rangle_n$
into circulant graphs by partitioning the connection set, and then decompose the resulting
circulant graphs into $n$-cycles using Theorems \ref{BFMTheorem} and \ref{MDeanTheorem}.

If $h$ is even, then we partition the connection set into pairs by pairing $\frac{n-1}2-h$ with
$\frac{n-1}2$ and partitioning $\{\frac{n-1}2-h+2,\ldots,\frac{n-3}2\}$ into consecutive pairs (if
$h=2$, then our partition is just $\{\{\frac{n-5}2,\frac{n-1}2\}\}$). Each of the resulting
circulant graphs is $4$-regular and connected and thus can be decomposed into two $n$-cycles by
Theorem \ref{BFMTheorem}. If $h$ is odd, then we partition the connection set into the triple $\{
\frac{n-1}2-h,\frac{n-1}2-h+2,\frac{n-1}2 \}$ and consecutive pairs from
$\{\frac{n-1}2-h+3,\ldots,\frac{n-3}2\}$ (if $h=3$, then our partition is just $\{\{
\frac{n-7}2,\frac{n-3}2,\frac{n-1}2 \}\}$).  Since $\gcd(\frac{n-1}2,n)=1$, the graph
$\langle\{\frac{n-1}2-h,\frac{n-1}2-h+2,\frac{n-1}2\}\rangle_n$ can be decomposed into three
$n$-cycles by Theorem \ref{MDeanTheorem}. Any other resulting circulant graphs are $4$-regular and
connected and thus can each be decomposed into two $n$-cycles by Theorem \ref{BFMTheorem}. \qed

\begin{lemma}\label{hamdecompxtonover2withgapneven}
If $n$ is even and $1\leq h\leq\frac{n-4}2$, then there is an $(n^h)$-decomposition of
$\langle\{\frac n2-h-1\} \cup \{\frac{n}{2}-h+1,\ldots,\frac n2\}\rangle_n$; except when $h=1$ and
$n\equiv 0\md 4$ in which case the graph is not connected.
\end{lemma}

\proof If $h=1$, then the graph is $\lan{\{\frac{n-4}{2},\frac n2\}}$. If $n \equiv 2 \md{4}$, then
$\gcd(\frac{n-4}2,n)=1$, $\lan{\frac{n-4}2}$ is an $n$-cycle, and $\lan{\{\frac{n}2\}}$ is a
perfect matching. If $n \equiv 0 \md{4}$, then $\gcd(\frac{n-4}2,\frac n2,n)=2$ and $\lan{\{
\frac{n-4}2,\frac n2\}}$ is not connected. Thus the result holds for $h=1$. In the remainder of the
proof we assume $h\geq 2$.

We first decompose $\langle\{\frac n2-h-1\} \cup \{\frac{n}{2}-h+1,\ldots,\frac n2\}\rangle_n$ into
circulant graphs by partitioning the connection set, and then decompose the resulting circulant
graphs into $n$-cycles using Theorems \ref{BFMTheorem}, \ref{5regCayintoHams} and
\ref{MDeanTheorem}.

If $n\equiv 0\md 4$ and $h$ is even, then we partition the connection set into pairs and the
singleton $\{\frac n2\}$ by pairing $\frac n2-h-1$ with $\frac{n-2}2$ and partitioning $\{\frac
n2-h+1,\ldots,\frac{n-4}2\}$ into pairs of consecutive integers (if $h=2$, then our partition is
just $\{\{\frac n2\},\{\frac{n-6}2,\frac{n-2}2\}\}$). The graph $\lan{\{\frac{n}2\}}$ is a perfect
matching. The other resulting circulant graphs are $4$-regular and connected and thus can each be
decomposed into two $n$-cycles by Theorem \ref{BFMTheorem} (note that $\gcd(\frac{n-2}2,n)=1$).

If $n\equiv 0\md 4$ and $h$ is odd, then we partition the connection set into pairs, the triple
$\{\frac n2-h-1,\frac n2-h+1,\frac{n-2}2\}$ and the singleton $\{\frac n2\}$ by partitioning
$\{\frac n2-h+2,\ldots,\frac{n-4}2\}$ into pairs of consecutive integers (if $h=3$, then our
partition is just $\{\{\frac n2\},\{\frac{n-8}{2},\frac{n-4}{2},\frac{n-2}2\}\}$). The graph
$\lan{\{\frac{n}2\}}$ is a perfect matching and, since $\gcd(\frac{n-2}2,n)=1$, the graph
$\langle\{\frac n2-h-1,\frac n2-h+1,\frac{n-2}2\}\rangle_n$ can be decomposed into three $n$-cycles
using Theorem \ref{MDeanTheorem}.  Any other resulting circulant graphs are $4$-regular and
connected and thus can each be decomposed into two $n$-cycles by Theorem \ref{BFMTheorem}.

If $n\equiv 2\md 4$, then we partition the connection set into pairs, the triple $\{\tfrac
n2-h-1,\tfrac{n-2}2,\tfrac n2\}$ and, when $h$ is odd, the singleton $\{\frac{n-4}2\}$ by
partitioning $\{\frac n2-h+1,\ldots,\frac{n-4}2\}$ into pairs of consecutive integers (when $h$ is
even) or into pairs of consecutive integers and the singleton $\{\frac{n-4}2\}$ (when $h$ is odd).
(Our partition is just $\{\{\tfrac {n-6}{2},\tfrac{n-2}2,\tfrac n2\}\}$ if $h=2$, and just
$\{\{\tfrac {n-8}{2},\tfrac{n-2}2,\tfrac n2\},\{\tfrac {n-4}{2}\}\}$ if $h=3$.) Since
$\gcd(\frac{n-2}2,\frac n2)=1$, the graph $\langle\{\frac n2-h-1,\frac{n-2}2,\frac n2\}\rangle_n$
can be decomposed into two $n$-cycles and a perfect matching using Theorem \ref{5regCayintoHams}.
When $h$ is odd, $\langle\{\frac{n-4}2\}\rangle_n$ is an $n$-cycle (note that
$\gcd(\frac{n-4}2,n)=1$). Any other resulting circulant graphs are $4$-regular and connected and
thus can each be decomposed into two $n$-cycles by Theorem \ref{BFMTheorem}. \qed

\vspace{0.5cm}

\subsection{Proof of Lemma \ref{3s4sandHamsforKn-S}}\label{3s4sHamsProof}

Before we prove Lemma \ref{3s4sandHamsforKn-S}, we require three preliminary lemmas which establish
the existence of various $(4^{qn},n^h)$-decompositions of circulant graphs.

\begin{lemma}\label{4cycledifftrips}
If $S\subseteq\{1,\ldots,\lfloor\frac{n-1}2\rfloor\}$ such that
\begin{itemize}
    \item
$S=\{x+1,\ldots,x+4q\}$ for some $x$;
    \item
$S=\{x\} \cup \{x+2,\ldots,x+4q-1\} \cup \{x+4q+1\}$ for some $x$; or
    \item
$S=\{\frac{n-1}2-4q\} \cup \{\frac{n-1}2-4q+2,\ldots,\frac{n-1}2\}$ where $n$ is odd;
\end{itemize}
then there is a $(4^{qn})$-decomposition of $\lan S$.
\end{lemma}

\proof It is sufficient to partition $S$ into $q$ modulo $n$ difference quadruples. If
$S=\{x+1,\ldots,x+4q\}$, then we partition $S$ into $q$ sets of the form $\{y,y+1,y+2,y+3\}$, each
of which is a difference quadruple. If $S=\{x\} \cup \{x+2,\ldots,x+4q-1\} \cup \{x+4q+1\}$, then
we partition $S$ into $q$ sets of the form $\{y,y+2,y+3,y+5\}$, each of which is a difference
quadruple. If $S=\{\frac{n-1}2-4q\} \cup \{\frac{n-1}2-4q+2,\ldots,\frac{n-1}2\}$ and $n$ is odd,
then we partition $S$ into $q-1$ sets of the form $\{y,y+2,y+3,y+5\}$, each of which is a
difference quadruple, and the set $\{ \frac{n-9}2,\frac{n-5}2,\frac{n-3}2,\frac{n-1}2\}$, which is
a modulo $n$ difference quadruple (note that $\frac{n-5}2+\frac{n-3}2+\frac{n-1}2-\frac{n-9}2=n$).
\qed

\begin{lemma}\label{hamsand4snogap}
If $h$, $q$ and $n$ are non-negative integers with $1\leq 4q+h\leq\floor{\frac{n-1}2}$, then there
is a $(4^{qn},n^h)$-decomposition of $\lan{ \{\floor{\frac{n-1}2}-h-4q+1,\ldots,\floor{\frac n2}\}
}$.
\end{lemma}

\proof If $h=0$ then the result follows immediately by Lemma \ref{4cycledifftrips}, and if $q=0$
then the result follows immediately by Lemma \ref{hamdecompxtonover2}. For $q,h\geq 1$ we partition
the connection set into the set $\{\floor{\frac{n-1}2}-h-4q+1,\ldots,\floor{\frac{n-1}2}-h\}$ and
the set $\{\floor{\frac{n-1}2}-h+1,\ldots,\floor{\frac n2}\}$. Then
$\lan{\{\floor{\frac{n-1}2}-h-4q+1,\ldots,\floor{\frac{n-1}2}-h\}}$ has a $(4^{qn})$-decomposition
by Lemma \ref{4cycledifftrips}, and $\lan{\{\floor{\frac{n-1}2}-h+1,\ldots,\floor{\frac n2}\}}$ has
an $(n^h)$-decomposition by Lemma \ref{hamdecompxtonover2}. \qed

\begin{lemma}\label{hamsand4swithgap}
If $h$, $q$ and $n$ are non-negative integers with $1\leq 4q+h\leq\floor{\frac{n-3}2}$ such that
$n$ is odd when $h=0$ and $n\equiv 1,2,5,6,7,10,11\md{12}$ when $h=1$, then there is a
$(4^{qn},n^h)$-decomposition of $\lan{ \{\floor{\frac{n-1}2}-h-4q\} \cup
\{\floor{\frac{n-1}2}-h-4q+2,\ldots,\floor{\frac n2}\}}$.
\end{lemma}

\proof If $h=0$, then the result follows immediately by Lemma \ref{4cycledifftrips}. If $q=0$, then
the result follows immediately by Lemma \ref{hamdecompxtonover2withgapnodd} ($n$ odd) or Lemma
\ref{hamdecompxtonover2withgapneven} ($n$ even). For $h,q\geq 1$ we partition the connection set
into the set $\{\floor{\frac{n-1}2}-h-4q\} \cup
\{\floor{\frac{n-1}2}-h-4q+2,\ldots,\floor{\frac{n-1}2}-h-1\} \cup \{\floor{\frac{n-1}2}-h+1\}$ and
the set $\{\floor{\frac{n-1}2}-h\} \cup\{\floor{\frac{n-1}2}-h+2,\ldots,\floor{\frac n2}\}$. Then
$\lan{\{\floor{\frac{n-1}2}-h-4q\} \cup
\{\floor{\frac{n-1}2}-h-4q+2,\ldots,\floor{\frac{n-1}2}-h-1\} \cup \{\floor{\frac{n-1}2}-h+1\}}$
has a $(4^{qn})$-decomposition by Lemma \ref{4cycledifftrips}, and $\lan{\{\floor{\frac{n-1}2}-h\}
\cup\{\floor{\frac{n-1}2}-h+2,\ldots,\floor{\frac n2}\}}$ has an $(n^h)$-decomposition by Lemma
\ref{hamdecompxtonover2withgapnodd} ($n$ odd) or Lemma \ref{hamdecompxtonover2withgapneven} ($n$
even). \qed

\vspace{0.3cm}

We now prove Lemma \ref{3s4sandHamsforKn-S}, which we restate here for convenience.

\vspace{0.5cm}

\noindent{\bf Lemma \ref{3s4sandHamsforKn-S} }
{\em
If $S\in\{\{1,2,3,4\},\{1,2,3,4,6\},\{1,2,3,4,5,7\},\{1,2,3,4,5,6,7\}\}$ and $n\geq 2\max(S) + 1$,
$t\geq 0$, $q\geq 0$ and $h\geq 2$ are integers satisfying $3t+4q+h=\lfloor\frac{n-1}2\rfloor-|S|$,
then there is a $(3^{tn},4^{qn},n^h)$-decomposition of $K_n-\lan{S}$, except possibly when $h=2$,
$S=\{1,2,3,4,5,6,7\}$ and
\begin{itemize}
\item $n\in\{25,26\}$ and $t=1$; or
\item $n=31$ and $t=2$.
\end{itemize}

}
\vspace{0.5cm}
\proof
We give the proof of Lemma
\ref{3s4sandHamsforKn-S} for each
$$S\in\{\{1,2,3,4\},\{1,2,3,4,6\},\{1,2,3,4,5,7\},\{1,2,3,4,5,6,7\}\}$$
separately.

\noindent{\bf Case A: ($S = \{1,2,3,4\}$)} \quad
The conditions $h\geq 2$ and $3t+4q+h=\lfloor\frac{n-1}2\rfloor-4$ imply $n\geq 6t+13$. If $t=0$,
then the result follows immediately by Lemma \ref{hamsand4swithgap}. We deal separately with the
three cases $t\in\{1,2,3,4\}$, $t\in\{5,6,7,8\}$, and $t\geq 9$.

\vspace{0.5cm}

\noindent{\bf Case A1:}\quad Suppose that $t\in\{1,2,3,4\}$. The cases $6t+13\leq n\leq 6t+17$
and the cases
$$(n,t)\in\{(30,2),(32,2),(36,3)\}$$
are dealt with first. Since $h\geq 2$, it follows from $3t+4q+h=\floor{\frac{n-1}2}-5$ that in each
of these cases we have $q=0$. Thus, the value of $h$ is uniquely determined by the values of $n$
and $t$. The required decompositions are obtained by partitioning $\{5,\ldots,\floor{\frac
n2}\}$ into $t$ modulo $n$ difference triples and a collection of connection sets for circulant
graphs such that the circulant graphs can be decomposed into Hamilton cycles (or Hamilton cycles
and a perfect matching) using the results in Section \ref{hamcycledecompsofcayleygraphs}. Suitable
partitions are given in the following tables.

\noindent
{\bf t=1:}
$$
\begin{array}{|c|c|c|}
\hline
n&{\rm modulo\ } n&{\rm connection\ sets}\\
&{\rm difference\ triples}&\\
\hline
\hline
19&\{5,6,8\}&\{7,9\}\\
\hline
20&\{5,6,9\}&\{7,8\}\\
\hline
21&\{5,7,9\}&\{6,8\},\{10\}\\
\hline
22&\{5,7,10\}&\{6,8,11\},\{9\}\\
\hline
\end{array}
$$

\noindent
{\bf t=2:}
$$
\begin{array}{|c|c|c|}
\hline
n&{\rm modulo\ } n&{\rm connection\ sets}\\
&{\rm difference\ triples}&\\
\hline
\hline
25&\{5,9,11\},\{7,8,10\}&\{6,12\}\\
\hline

26&\{5,7,12\},\{6,9,11\}&\{8,10,13\}\\
\hline

27&\{5,6,11\},\{8,9,10\}&\{7\},\{12,13\}\\
\hline

28&\{5,8,13\},\{6,10,12\}&\{7,9\},\{11\}\\
\hline

30&\{5,7,12\},\{6,8,14\}&\{9,10\},\{11,13\}\\
\hline

32&\{5,7,12\},\{6,8,14\}&\{9,10\},\{11,13\},\{15\}\\
\hline
\end{array}
$$

\noindent
{\bf t=3:}
$$
\begin{array}{|c|c|c|}
\hline
n&{\rm modulo\ } n&{\rm connection\ sets}\\
&{\rm difference\ triples}&\\
\hline
\hline
31&\{6,8,14\},\{7,11,13\},\{9,10,12\}&\{5,15\}\\
\hline

32&\{5,7,12\},\{6,8,14\},\{9,10,13\}&\{11,15\}\\
\hline

33&\{5,7,12\},\{6,8,14\},\{9,11,13\}&\{10\},\{15,16\}\\
\hline

34&\{5,7,12\},\{6,8,14\},\{10,11,13\}&\{9\},\{15,16\}\\
\hline

36&\{5,7,12\},\{6,8,14\},\{10,11,15\}&\{9,13\},\{16,17\}\\
\hline
\end{array}
$$

\noindent
{\bf t=4:}
$$
\begin{array}{|c|c|c|}
\hline
n&{\rm modulo\ } n&{\rm connection\ sets}\\
&{\rm difference\ triples}&\\
\hline
\hline
37&\{5,10,15\},\{6,7,13\},\{8,9,17\},\{11,12,14\}&\{16,18\}\\
\hline

38&\{5,9,14\},\{6,7,13\},\{8,10,18\},\{11,12,15\}&\{16,17\}\\
\hline

39&\{5,10,15\},\{6,7,13\},\{8,9,17\},\{11,12,16\}&\{14\},\{18,19\}\\
\hline

40&\{5,9,14\},\{6,7,13\},\{8,10,18\},\{11,12,17\}&\{15,16\},\{19\}\\
\hline
\end{array}
$$

\vspace{0.5cm}

We now deal with $n\geq 6t+18$ which implies $4q+h\geq 4$. Define $S_t$ by $S_t=
\{5,\ldots,3t+8\}$ for $t\in\{1,4\}$, and $S_t=\{5,\ldots,3t+7\} \cup \{3t+9\}$
when $t\in\{2,3\}$. The following table gives a partition $\pi_t$ of $S_t$ into difference triples
and a difference quadruple $Q_t$ such that $Q_t$ can be partitioned into two pairs of relatively
prime integers.

$$
\begin{array}{|c|c|}
\hline
t&\pi_t \\
\hline \hline
1&\{\{5,6,11\},\{7,8,9,10\}\}\\
\hline
2&\{\{5,6,11\},\{7,8,15\},\{9,10,12,13\}\}\\
\hline
3&\{\{5,6,11\},\{7,9,16\},\{8,10,18\},\{12,13,14,15\}\}\\
\hline
4&\{\{5,10,15\},\{6,11,17\},\{7,9,16\},\{8,12,20\},\{13,14,18,19\}\}\\
\hline
\end{array}
$$

\vspace{0.5cm}

Thus, $\lan{Q_t}$ can be decomposed into two connected $4$-regular Cayley graphs, which in turn can
be decomposed into Hamilton cycles using Theorem \ref{BFMTheorem}. It follows that there is both a
$(3^{tn},4^n)$-decomposition and a $(3^{tn},n^4)$-decomposition of $\lan{S_t}$. If $q=0$, then we
use the $(3^{tn},n^4)$-decomposition of $\lan{S_t}$ and if $q\geq1$, then we use the
$(3^{tn},4^n)$-decomposition of $\lan{S_t}$. This leaves us needing an $(n^{h-4})$-decomposition of
$K_n-\lan{\{1,2,3,4\}\cup S_t}$ when $q=0$, and a $(4^{(q-1)n},n^h)$-decomposition of
$K_n-\lan{\{1,2,3,4\}\cup S_t}$ when $q\geq1$. Note that $K_n-\lan{\{1,2,3,4\}\cup S_t}$ is
isomorphic to
\begin{itemize}
    \item
$\lan{\{3t+9,\ldots,\floor{\frac n2}\}}$ when $t \in \{1,4\}$; and
    \item
$\lan{\{3t+8\} \cup \{3t+10,\ldots,\floor{\frac n2}\}}$ when $t \in \{2,3\}$.
\end{itemize}

When $t \in \{1,4\}$ the required decomposition exists by Lemma \ref{hamsand4snogap}. When $t
\in \{2,3\}$ and the required number of Hamilton cycles (that is, $h-4$ when $q=0$ and $h$ when
$q\geq1$) is at least $2$, the required decomposition exists by Lemma \ref{hamsand4swithgap}. So we
need to consider only the cases where $q=0$, $h\in\{4,5\}$ and $t\in\{2,3\}$.

Since $3t+4q+h=\floor{\frac{n-1}{2}}-4$, and since we have already dealt with the cases where
$(n,t)\in\{(30,2),(32,2),(36,3)$, this leaves us with only the five cases where
$$(n,t,h) \in \{(29,2,4),(31,2,5),(35,3,4),(37,3,5),(38,3,5)\}.$$ In the cases $(n,t,h) \in \{(29,2,4),(35,3,4)\}$ we
have that $h-4$ (the required number of Hamilton cycles) is $0$ and $n$ is odd, and in the cases
$(n,t,h) \in \{(31,2,5),(37,3,5),(38,3,5)\}$ we have that $h-4$ (the required number of Hamilton cycles) is $1$ and $n
\equiv 1,2,7 \md{12}$. So in all these cases the required decompositions exist by Lemma
\ref{hamsand4swithgap}.

\vspace{0.5cm}

\noindent{\bf Case A2:}\quad Suppose that $t\in\{5,6,7,8\}$. Redefine $S_t$ by
$S_t=\{5,\ldots,3t+6\}$. The following table gives a partition of $S_t$ into
difference triples and a set $R_t$ consisting of a pair of relatively prime
integers. Thus, $\lan{R_t}$ is a connected $4$-regular Cayley graph, and so can
be decomposed into two Hamilton cycles using Theorem \ref{BFMTheorem}. Thus, we
have a $(3^{tn},n^2)$-decomposition of $\lan{S_t}$.
$$
\begin{array}{|c|c|c|}
\hline
t & \mbox{difference triples} & R_t \\
\hline
\hline
5&\{5,12,17\},\{6,13,19\},\{7,14,21\},\{8,10,18\},\{9,11,20\} & \{15,16\} \\
& \{12,13,25\}&  \\
\hline
6&\{5,13,18\},\{6,14,20\},\{7,15,22\},\{8,16,24\},\{9,10,19\},\{11,12,23\} & \{17,21\} \\
\hline
7&\{5,15,20\},\{6,16,22\},\{7,17,24\},\{8,18,26\},\{9,12,21\},\{10,13,23\}, & \{19,27\} \\
& \{11,14,25\} & \\
\hline
8&\{5,17,22\},\{6,18,24\},\{7,19,26\},\{8,13,21\},\{9,16,25\},\{10,20,30\}, & \{27,28\} \\
& \{11,12,23\},\{14,15,29\} &  \\
\hline
\end{array}
$$
Thus, we only require a $(4^{qn},n^{h-2})$-decomposition of $K_n-\lan{\{1,2,3,4\}\cup S_t}$. But
$K_n-\lan{\{1,2,3,4\}\cup S_t}$ is isomorphic to $\lan{\{3t+7,\ldots,\floor{\frac n2}\}}$ and so
this decomposition exists by Lemma \ref{hamsand4snogap}.

\vspace{0.5cm}

\noindent{\bf Case A3:}\quad Suppose that $t\geq 9$. Redefine $S_t$ by
$S_t=\{5,\ldots,3t+4\}$ when $t\equiv 0,1\md 4$, and $S_t=\{5,\ldots,3t+3\}
\cup \{3t+5\}$ when $t\equiv 2,3\md 4$. We now obtain a
$(3^{tn})$-decomposition of $\lan{S_t}$.

For $t\equiv 0,1\md 4$ (respectively $t\equiv
2,3\md 4$), we can obtain a $(3^{tn})$-decomposition of $\lan{S_t}$ by using a
Langford sequence (respectively hooked Langford sequence) of order $t$ and defect $5$, which
exists since $t \geq 9$, to partition $S_t$ into difference triples (see
\cite{Sha,Sim}). So we have a $(3^{tn})$-decomposition of $\lan{S_t}$, and require a
$(4^{qn},n^h)$-decomposition of $K_n-\lan{\{1,2,3,4,6\}\cup S_t}$. Since
$K_n-\lan{\{1,2,3,4\}\cup S_t}$ is isomorphic to
\begin{itemize}
\item $\lan{\{3t+5,\ldots,\floor{\frac n2}\}}$ when $t\equiv 0,1\md 4$; and
\item $\lan{\{3t+4\}\cup\{3t+6,\ldots,\floor{\frac n2}\}}$ when $t\equiv 2,3\md 4$;
\end{itemize}
this decomposition exists by Lemma \ref{hamsand4snogap} or \ref{hamsand4swithgap}. \qed

\vspace{0.3cm}

\noindent{\bf Case B: ($S = \{1,2,3,4,6\}$)} \quad
The conditions $h\geq 2$ and $3t+4q+h=\lfloor\frac{n-1}2\rfloor-5$ imply $n\geq 6t+15$. If $t=0$,
then the result follows immediately by Lemma \ref{hamsand4swithgap}. We deal separately with the
three cases $t\in\{1,2,3,4,5,6\}$, $t\in\{7,8,9,10\}$, and $t\geq 11$.

\vspace{0.5cm}

\noindent{\bf Case B1:}\quad Suppose that $t\in\{1,2,3,4,5,6\}$. The cases $6t+15\leq n\leq 6t+18$
and the cases
$$(n,t)\in\{(38,3),(39,3),(40,3),(44,4),(45,4)\}$$
are dealt with first. Since $h\geq 2$, it follows from $3t+4q+h=\floor{\frac{n-1}2}-5$ that in each
of these cases we have $q=0$. Thus, the value of $h$ is uniquely determined by the values of $n$
and $t$. The required decompositions are obtained by partitioning $\{5\}\cup\{7,\ldots,\floor{\frac
n2}\}$ into $t$ modulo $n$ difference triples and a collection of connection sets for circulant
graphs such that the circulant graphs can be decomposed into Hamilton cycles (or Hamilton cycles
and a perfect matching) using the results in Section \ref{hamcycledecompsofcayleygraphs}. Suitable
partitions are given in the following tables.

\noindent
{\bf t=1:}
$$
\begin{array}{|c|c|c|}
\hline
n&{\rm modulo\ } n&{\rm connection\ sets}\\
&{\rm difference\ triples}&\\
\hline
\hline
21&\{5,7,9\}&\{8,10\}\\
\hline
22&\{5,7,10\}&\{8,9\},\{11\}\\
\hline
23&\{5,8,10\}&\{7,9\},\{11\}\\
\hline
24&\{5,9,10\}&\{7,8\},\{11\},\{12\}\\
\hline
\end{array}
$$

\noindent
{\bf t=2:}
$$
\begin{array}{|c|c|c|}
\hline
n&{\rm modulo\ } n&{\rm connection\ sets}\\
&{\rm difference\ triples}&\\
\hline
\hline
27&\{5,7,12\},\{8,9,10\}&\{11,13\}\\
\hline

28&\{5,7,12\},\{8,9,11\}&\{10,13\},\{14\}\\
\hline

29&\{5,7,12\},\{8,10,11\}&\{9\},\{13,14\}\\
\hline

30&\{5,9,14\},\{8,10,12\}&\{7\},\{11,13\},\{15\}\\
\hline
\end{array}
$$

\noindent
{\bf t=3:}
$$
\begin{array}{|c|c|c|}
\hline
n&{\rm modulo\ } n&{\rm connection\ sets}\\
&{\rm difference\ triples}&\\
\hline
\hline
33&\{5,8,13\},\{7,9,16\},\{10,11,12\}&\{14,15\}\\
\hline

34&\{5,10,15\},\{7,9,16\},\{8,12,14\}&\{11,13\},\{17\}\\
\hline

35&\{5,8,13\},\{7,9,16\},\{10,11,14\}&\{12,15\},\{17\}\\
\hline

36&\{5,8,13\},\{7,9,16\},\{10,12,14\}&\{11,15\},\{17\},\{18\}\\
\hline

38&\{5,12,17\},\{7,9,16\},\{8,10,18\}&\{11,13\},\{14,15\},\{19\}\\
\hline

39&\{5,12,17\},\{7,9,16\},\{8,10,18\}&\{11,13\},\{14,15\},\{19\}\\
\hline

40&\{5,12,17\},\{7,9,16\},\{8,10,18\}&\{11,13\},\{14,15\},\{19\},\{20\}\\
\hline
\end{array}
$$

\noindent
{\bf t=4:}
$$
\begin{array}{|c|c|c|}
\hline
n&{\rm modulo\ } n&{\rm connection\ sets}\\
&{\rm difference\ triples}&\\
\hline
\hline
39&\{5,11,16\},\{7,8,15\},\{9,10,19\},\{12,13,14\}&\{17,18\}\\
\hline

40&\{5,8,13\},\{7,9,16\},\{10,12,18\},\{11,14,15\}&\{17,19\},\{20\}\\
\hline

41&\{5,14,19\},\{7,8,15\},\{9,11,20\},\{12,13,16\}&\{10,17\},\{18\}\\
\hline

42&\{5,10,15\},\{7,11,18\},\{8,9,17\},\{12,14,16\}&\{13\},\{19,20\},\{21\}\\
\hline

44&\{5,11,16\},\{7,13,20\},\{8,10,18\},\{9,12,21\}&\{14,15\},\{17,19\},\{22\}\\
\hline

45&\{5,11,16\},\{7,13,20\},\{8,10,18\},\{9,12,21\}&\{14,15\},\{17,19\},\{22\}\\
\hline
\end{array}
$$

\noindent
{\bf t=5:}
$$
\begin{array}{|c|c|c|}
\hline
n&{\rm modulo\ } n&{\rm connection\ sets}\\
&{\rm difference\ triples}&\\
\hline
\hline
45&\{5,17,22\},\{7,13,20\},\{8,10,18\},\{9,12,21\},\{14,15,16\}&\{11,19\}\\
\hline

46&\{5,17,22\},\{7,13,20\},\{8,10,18\},\{9,12,21\},\{11,16,19\}&\{14,15\},\{23\}\\
\hline

47&\{5,18,23\},\{7,13,20\},\{8,11,19\},\{10,12,22\},\{14,16,17\}&\{9,15\},\{21\}\\
\hline

48&\{5,17,22\},\{7,13,20\},\{8,10,18\},\{9,12,21\},\{11,14,23\}&\{15,16\},\{19\},\{24\}\\
\hline
\end{array}
$$

\noindent
{\bf t=6:}
$$
\begin{array}{|c|c|c|}
\hline
n&{\rm modulo\ } n&{\rm connection\ sets}\\
&{\rm difference\ triples}&\\
\hline
\hline
51&\{5,13,18\},\{7,15,22\},\{8,16,24\}, &\{21,25\}\\
&\{9,10,19\},\{11,12,23\},\{14,17,20\}&\\
\hline

52&\{5,13,18\},\{7,15,22\},\{8,16,24\}, &\{20,25\},\{26\}\\
&\{9,10,19\},\{11,12,23\},\{14,17,21\}&\\
\hline

53&\{5,13,18\},\{7,12,19\},\{8,14,22\}, &\{23\},\{25,26\}\\
&\{9,15,24\},\{10,11,21\},\{16,17,20\}&\\
\hline

54&\{5,17,22\},\{7,8,15\},\{9,10,19\}, &\{21,23\},\{25\},\{27\}\\
&\{11,13,24\},\{12,14,26\},\{16,18,20\}&\\
\hline
\end{array}
$$

\vspace{0.5cm}

We now deal with $n\geq 6t+19$ which implies $4q+h\geq 4$. Define $S_t$ by $S_t=\{5\} \cup
\{7,\ldots,3t+9\}$ for $t\in\{1,2,5,6\}$, and $S_t=\{5\} \cup \{7,\ldots,3t+8\} \cup \{3t+10\}$
when $t\in\{3,4\}$. The following table gives a partition $\pi_t$ of $S_t$ into difference triples
and a difference quadruple $Q_t$ such that $Q_t$ can be partitioned into two pairs of relatively
prime integers.

$$
\begin{array}{|c|c|}
\hline
t&\pi_t \\
\hline \hline
1&\{\{5,7,12\},\{8,9,10,11\}\}\\
\hline
2&\{\{5,9,14\},\{7,8,15\},\{10,11,12,13\}\}\\
\hline
3&\{\{5,9,14\},\{7,10,17\},\{8,11,19\},\{12,13,15,16\}\}\\
\hline
4&\{\{5,9,14\},\{7,13,20\},\{8,11,19\},\{10,12,22\},\{15,16,17,18\}\}\\
\hline
5&\{\{5,14,19\},\{7,13,20\},\{8,10,18\},\{9,15,24\},\{11,12,23\},\{16,17,21,22\}\}\\
\hline
6&\{\{5,15,20\},\{7,16,23\},\{8,14,22\},\{9,12,21\},\{10,17,27\},\{11,13,24\},\\
&\{18,19,25,26\}\}\\
\hline
\end{array}
$$

\vspace{0.5cm}

Thus, $\lan{Q_t}$ can be decomposed into two connected $4$-regular Cayley graphs, which in turn can
be decomposed into Hamilton cycles using Theorem \ref{BFMTheorem}. It follows that there is both a
$(3^{tn},4^n)$-decomposition and a $(3^{tn},n^4)$-decomposition of $\lan{S_t}$. If $q=0$, then we
use the $(3^{tn},n^4)$-decomposition of $\lan{S_t}$ and if $q\geq1$, then we use the
$(3^{tn},4^n)$-decomposition of $\lan{S_t}$. This leaves us needing an $(n^{h-4})$-decomposition of
$K_n-\lan{\{1,2,3,4,6\}\cup S_t}$ when $q=0$, and a $(4^{(q-1)n},n^h)$-decomposition of
$K_n-\lan{\{1,2,3,4,6\}\cup S_t}$ when $q\geq1$. Note that $K_n-\lan{\{1,2,3,4,6\}\cup S_t}$ is
isomorphic to
\begin{itemize}
    \item
$\lan{\{3t+10,\ldots,\floor{\frac n2}\}}$ when $t \in \{1,2,5,6\}$; and
    \item
$\lan{\{3t+9\} \cup \{3t+11,\ldots,\floor{\frac n2}\}}$ when $t \in \{3,4\}$.
\end{itemize}

When $t \in \{1,2,5,6\}$ the required decomposition exists by Lemma \ref{hamsand4snogap}. When $t
\in \{3,4\}$ and the required number of Hamilton cycles (that is, $h-4$ when $q=0$ and $h$ when
$q\geq1$) is at least $2$, the required decomposition exists by Lemma \ref{hamsand4swithgap}. So we
need to consider only the cases where $q=0$, $h\in\{4,5\}$ and $t\in\{3,4\}$.

Since $3t+4q+h=\floor{\frac{n-1}{2}}-5$, and since we have already dealt with the cases where
$(n,t)\in\{(38,3),(39,3),(40,3),(44,4),(45,4)\}$, this leaves us with only the three cases where
$(n,t,h) \in \{(37,3,4),(43,4,4),(46,4,5)\}.$ In the cases $(n,t,h) \in \{(37,3,4),(43,4,4)\}$ we
have that $h-4$ (the required number of Hamilton cycles) is $0$ and $n$ is odd, and in the case
$(n,t,h) = (46,4,5)$ we have that $h-4$ (the required number of Hamilton cycles) is $1$ and $n
\equiv 10 \md{12}$. So in all these cases the required decompositions exist by Lemma
\ref{hamsand4swithgap}.

\vspace{0.5cm}

\noindent{\bf Case B2:}\quad Suppose that $t\in\{7,8,9,10\}$. Redefine $S_t$ by $S_t=\{5\} \cup
\{7,\ldots,3t+7\}$. The following table gives a partition of $S_t$ into difference triples and a
set $R_t$ consisting of a pair of relatively prime integers. Thus, $\lan{R_t}$ is a connected
$4$-regular Cayley graph, and so can be decomposed into two Hamilton cycles using Theorem
\ref{BFMTheorem}. Thus, we have a $(3^{tn},n^2)$-decomposition of $\lan{S_t}$.
$$
\begin{array}{|c|c|c|}
\hline
t & \mbox{difference triples} & R_t \\
\hline
\hline
7&\{5,17,22\},\{7,16,23\},\{8,20,28\},\{9,18,27\},\{10,14,24\},\{11,15,26\}, & \{19,21\} \\
& \{12,13,25\}&  \\
\hline
8&\{5,17,22\},\{7,19,26\},\{8,16,24\},\{9,20,29\},\{10,21,31\},\{11,14,25\}, & \{23,27\} \\
& \{12,18,30\},\{13,15,28\}& \\
\hline
9&\{5,19,24\},\{7,20,27\},\{8,21,29\},\{9,22,31\},\{10,23,33\},\{11,15,26\}, & \{25,34\} \\
& \{12,16,28\},\{13,17,30\},\{14,18,32\} & \\
\hline
10&\{5,21,26\},\{7,22,29\},\{8,23,31\},\{9,24,33\},\{10,25,35\},\{11,17,28\}, & \{32,37\} \\
& \{12,18,30\},\{13,14,27\},\{15,19,34\},\{16,20,36\} &  \\
\hline
\end{array}
$$
Thus, we only require a $(4^{qn},n^{h-2})$-decomposition of $K_n-\lan{\{1,2,3,4,6\}\cup S_t}$. But
$K_n-\lan{\{1,2,3,4,6\}\cup S_t}$ is isomorphic to $\lan{\{3t+8,\ldots,\floor{\frac n2}\}}$ and so
this decomposition exists by Lemma \ref{hamsand4snogap}.

\vspace{0.5cm}

\noindent{\bf Case B3:}\quad Suppose that $t\geq 11$. Redefine $S_t$ by
$S_t=\{5\} \cup \{7,\ldots,3t+5\}$ when $t\equiv 1,2\md 4$, and $S_t=\{5\} \cup
\{7,\ldots,3t+4\} \cup \{3t+6\}$ when $t\equiv 0,3\md 4$. We now obtain a
$(3^{tn})$-decomposition of $\lan{S_t}$. For $11\leq t\leq 52$, we have found
such a decomposition by partitioning $S_t$ into difference triples with the aid
of a computer. These difference triples are shown in Table
\ref{tab:5789_triples} in the appendix.  For $t \geq 53$ and $t\equiv 1,2 \md
4$ (respectively $t\equiv 0,3 \md 4$) we set aside as one difference triple
$\{5,3t,3t+5\}$ (respectively $\{5,3t+1,3t+6\}$) and form the set $S'_t =
\{7,\ldots,3t-1,3t+1,\ldots,3t+4\} = \{7,\ldots,3t+4\} \setminus \{3t\}$
(respectively $S'_t = \{7,\ldots,3t,3t+2,\ldots,3t+4\} = \{7,\ldots,3t+4\}
\setminus \{3t+1\}$).  We can obtain a $(3^{(t-1)n})$-decomposition of
$\lan{S'_t}$ by using an extended Langford sequence of order $t-1$ and defect
$7$ to partition $S'_t$ into difference triples. Since $t \geq 53$, this
sequence exists by Theorem 7.1 in \cite{LinMor} (also see \cite{Sha,Sim}). So
we have a $(3^{tn})$-decomposition of $\lan{S_t}$, and require a
$(4^{qn},n^h)$-decomposition of $K_n-\lan{\{1,2,3,4,6\}\cup S_t}$. Since
$K_n-\lan{\{1,2,3,4,6\}\cup S_t}$ is isomorphic to
\begin{itemize}
\item $\lan{\{3t+6,\ldots,\floor{\frac n2}\}}$ when $t\equiv 1,2\md 4$; and
\item $\lan{\{3t+5\}\cup\{3t+7,\ldots,\floor{\frac n2}\}}$ when $t\equiv 0,3\md 4$;
\end{itemize} this decomposition exists by Lemma
\ref{hamsand4snogap} or \ref{hamsand4swithgap}. \qed


\vspace{0.3cm}

\noindent{\bf Case C: ($S = \{1,2,3,4,5,7\}$)} \quad
The conditions $h\geq 2$ and $3t+4q+h=\lfloor\frac{n-1}2\rfloor-5$ imply $n\geq 6t+17$. If $t=0$,
then the result follows immediately by Lemma \ref{hamsand4swithgap}. We deal separately with the
four cases $t\in\{1,5,6,7,8\}$, $t\in\{2,3,4\}$, $t\in\{9,10,11,12\}$, and $t\geq 11$.

\vspace{0.5cm}

\noindent{\bf Case C1:}\quad Suppose that $t\in\{1,5,6,7,8\}$. The cases $6t+17\leq n\leq 6t+20$
and the cases
$$(n,t)\in\{  (28,1),(52,5),(70,8),(72,8)  \}$$
are dealt with first. Since $h\geq 2$, it follows from $3t+4q+h=\floor{\frac{n-1}2}-6$ that in each
of these cases we have $q=0$. Thus, the value of $h$ is uniquely determined by the values of $n$
and $t$. The required decompositions are obtained by partitioning $\{6\}\cup\{8,\ldots,\floor{\frac
n2}\}$ into $t$ modulo $n$ difference triples and a collection of connection sets for circulant
graphs such that the circulant graphs can be decomposed into Hamilton cycles (or Hamilton cycles
and a perfect matching) using the results in Section \ref{hamcycledecompsofcayleygraphs}. Suitable
partitions are given in the following tables.

\noindent
{\bf t=1:}
$$
\begin{array}{|c|c|c|}
\hline
n&{\rm modulo\ } n&{\rm connection\ sets}\\
&{\rm difference\ triples}&\\
\hline
\hline
23&\{6,8,9\}&\{10,11\}\\
\hline
24&\{6,8,10\}&\{9,11\}\\
\hline
25&\{6,8,11\}&\{9,10\},\{12\}\\
\hline
26&\{6,8,12\}&\{9,10\},\{11\}\\
\hline
27&\{6,8,13\}&\{9,10\},\{11,12\}\\
\hline
28&\{6,10,12\}&\{8,9\},\{11,13\}\\
\hline
\end{array}
$$

\noindent
{\bf t=5:}
$$
\begin{array}{|c|c|c|}
\hline
n&{\rm modulo\ } n&{\rm connection\ sets}\\
&{\rm difference\ triples}&\\
\hline
\hline
47&\{6,13,19\},\{8,15,23\},\{9,11,20\},\{10,12,22\},\{14,16,17\}&\{18,21\}\\
\hline

48&\{6,13,19\},\{8,15,23\},\{9,11,20\},\{10,12,22\},\{14,16,18\}&\{17,21\}\\
\hline

49&\{6,13,19\},\{8,15,23\},\{9,11,20\},\{10,12,22\},\{14,17,18\}&\{16,21\},\{24\}\\
\hline

50&\{6,13,19\},\{8,16,24\},\{9,11,20\},\{10,12,22\},\{15,17,18\}&\{14,21\},\{23\}\\
\hline

52&\{6,13,19\},\{8,16,24\},\{9,11,20\},\{10,12,22\},\{14,17,21\}&\{15,18\},\{23,25\}\\
\hline
\end{array}
$$

\noindent
{\bf t=6:}
$$
\begin{array}{|c|c|c|}
\hline
n&{\rm modulo\ } n&{\rm connection\ sets}\\
&{\rm difference\ triples}&\\
\hline
\hline
53&\{6,15,21\},\{8,17,25\},\{9,13,22\}, &\{20,26\}\\
&\{10,14,24\},\{11,12,23\},\{16,18,19\}&\\
\hline

54&\{6,15,21\},\{8,17,25\},\{9,13,22\}, &\{19,26\}\\
&\{10,14,24\},\{11,12,23\},\{16,18,20\}&\\
\hline

55&\{6,15,21\},\{8,17,25\},\{9,13,22\}, &\{18,26\},\{27\}\\
&\{10,14,24\},\{11,12,23\},\{16,19,20\}&\\
\hline

56&\{6,15,21\},\{8,18,26\},\{9,13,22\}, &\{16,25\},\{27\}\\
&\{10,14,24\},\{11,12,23\},\{17,19,20\}&\\
\hline
\end{array}
$$

\noindent
{\bf t=7:}
$$
\begin{array}{|c|c|c|}
\hline
n&{\rm modulo\ } n&{\rm connection\ sets}\\
&{\rm difference\ triples}&\\
\hline
\hline
59&\{6,16,22\},\{8,15,23\},\{9,19,28\},\{10,17,27\}, &\{25,29\}\\
&\{11,13,24\},\{12,14,26\},\{18,20,21\}&\\
\hline

60&\{6,16,22\},\{8,17,25\},\{9,15,24\},\{10,18,28\}, &\{26,29\}\\
&\{11,12,23\},\{13,14,27\},\{19,20,21\}&\\
\hline

61&\{6,15,21\},\{8,16,24\},\{9,17,26\},\{10,18,28\}, &\{25,29\},\{30\}\\
&\{11,12,23\},\{13,14,27\},\{19,20,22\}&\\
\hline

62&\{6,20,26\},\{8,17,25\},\{9,15,24\},\{10,18,28\}, &\{16,30,31\},\{29\}\\
&\{11,12,23\},\{13,14,27\},\{19,21,22\}&\\
\hline
\end{array}
$$

\noindent
{\bf t=8:}
$$
\begin{array}{|c|c|c|}
\hline
n&{\rm modulo\ } n&{\rm connection\ sets}\\
&{\rm difference\ triples}&\\
\hline
\hline
65&\{6,21,27\},\{8,18,26\},\{9,15,24\},\{10,19,29\}, &\{31,32\}\\
&\{11,17,28\},\{12,13,25\},\{14,16,30\},\{20,22,23\}&\\
\hline

66&\{6,21,27\},\{8,18,26\},\{9,23,32\},\{10,19,29\}, &\{15,31\}\\
&\{11,17,28\},\{12,13,25\},\{14,16,30\},\{20,22,24\}&\\
\hline

67&\{6,27,33\},\{8,18,26\},\{9,23,32\},\{10,19,29\}, &\{15,20\},\{31\}\\
&\{11,17,28\},\{12,13,25\},\{14,16,30\},\{21,22,24\}&\\
\hline

68&\{6,18,24\},\{8,21,29\},\{9,23,32\},\{10,15,25\}, &\{17,31\},\{33\}\\
&\{11,19,30\},\{12,16,28\},\{13,14,27\},\{20,22,26\}&\\
\hline

70&\{6,18,24\},\{9,22,31\},\{10,15,25\},\{11,19,30\}, &\{8,29\},\{32,34,35\}\\
&\{12,16,28\},\{13,14,27\},\{17,20,33\},\{21,23,26\}&\\
\hline

72&\{6,16,22\},\{8,13,21\},\{9,17,26\},\{10,18,28\}, &\{29,33\},\{31,32\},\{35\}\\
&\{11,19,30\},\{12,15,27\},\{14,20,34\},\{23,24,25\}&\\
\hline
\end{array}
$$

\vspace{0.5cm}

We now define $S_t$ by $S_t=\{6\} \cup \{8,\ldots,3t+9\}\cup \{3t+11\}$ for
$t\equiv 0,1 \md 4$, and $S_t=\{6\} \cup \{8,\ldots,3t+9\}$ when $t\equiv
2,3 \md 4$. The following table gives a partition $\pi_t$ of $S_t$ into
difference triples and a difference quadruple $Q_t$ such that $Q_t$ can be
partitioned into two pairs of relatively prime integers.

$$
\begin{array}{|c|c|}
\hline
t&\pi_t \\
\hline \hline
1&\{\{6,8,14\},\{9,10,11,12\}\}\\
\hline
5&\{\{6,14,20\},\{8,15,23\},\{9,17,26\},\{10,12,22\},\{11,13,24\},\{16,18,19,21\}\}\\
\hline
6&\{\{6,13,19\},\{8,16,24\},\{9,17,26\},\{10,18,28\},\\
&\{11,14,25\},\{12,15,27\},\{20,21,22,23\}\}\\
\hline
7&\{\{6,24,30\},\{8,23,31\},\{9,16,25\},\{10,17,27\},\\
&\{11,18,29\},\{12,14,26\},\{13,15,28\},\{19,20,21,22\}\}\\
\hline
8&\{\{6,17,23\},\{8,16,24\},\{9,19,28\},\{10,15,25\},\{11,20,31\},\\
&\{12,21,33\},\{13,22,35\},\{14,18,32\},\{26,27,29,30\}\}\\
\hline
\end{array}
$$

\vspace{0.5cm}

Thus, $\lan{Q_t}$ can be decomposed into two connected $4$-regular Cayley graphs, which in turn can
be decomposed into Hamilton cycles using Theorem \ref{BFMTheorem}. It follows that there is both a
$(3^{tn},4^n)$-decomposition and a $(3^{tn},n^4)$-decomposition of $\lan{S_t}$. If $q=0$, then we
use the $(3^{tn},n^4)$-decomposition of $\lan{S_t}$ and if $q\geq1$, then we use the
$(3^{tn},4^n)$-decomposition of $\lan{S_t}$. This leaves us needing an $(n^{h-4})$-decomposition of
$K_n-\lan{\{1,2,3,4,5,7\}\cup S_t}$ when $q=0$, and a $(4^{(q-1)n},n^h)$-decomposition of
$K_n-\lan{\{1,2,3,4,5,7\}\cup S_t}$ when $q\geq1$. Note that $K_n-\lan{\{1,2,3,4,6\}\cup S_t}$ is
isomorphic to
\begin{itemize}
    \item
$\lan{\{3t+10,\ldots,\floor{\frac n2}\}}$ when $t \in \{1,5,6\}$; and
    \item
$\lan{\{3t+9\} \cup \{3t+11,\ldots,\floor{\frac n2}\}}$ when $t \in \{7,8\}$.
\end{itemize}

When $t \in \{6,7\}$ the required decomposition exists by Lemma
\ref{hamsand4snogap}. When $t \in \{1,5,8\}$ and the required number of
Hamilton cycles (that is, $h-4$ when $q=0$ and $h$ when $q\geq1$) is at least
$2$, the required decomposition exists by Lemma \ref{hamsand4swithgap}. So we
need to consider only the cases where $q=0$, $h\in\{4,5\}$ and $t\in\{1,5,8\}$.

Since $3t+4q+h=\floor{\frac{n-1}{2}}-6$, and since we have already dealt with
the cases where $(n,t)\in\{(27,1),(28,1),(52,5),(70,8),(72,8)\}$, this leaves
us with the cases where $(n,t,h)$ in
$$\{
(27,1,4),(29,1,5),(30,1,5),(51,5,4),(53,5,5),(54,5,5),(69,8,4),(71,8,5)
\}.$$

In the cases $(n,t,h) \in \{(27,1,4),(51,5,4),(69,8,4)\}$ we have that $h-4$
(the required number of Hamilton cycles) is $0$ and $n$ is odd, and in the
other cases we have that $h-4$ (the required number of Hamilton cycles) is $1$
and $n \equiv 5,6,11 \md{12}$. So in all these cases the required
decompositions exist by Lemma \ref{hamsand4swithgap}.

\vspace{0.5cm}
\noindent{\bf Case C2:}\quad Suppose that $t\in\{2,3,4\}$. The cases $6t+17\leq n\leq 6t+24$
are dealt with first. Since $h\geq 2$, it follows from $3t+4q+h=\floor{\frac{n-1}2}-6$ that in each
of these cases we have $q=0$. Thus, the value of $h$ is uniquely determined by the values of $n$
and $t$. The required decompositions are obtained by partitioning $\{6\}\cup\{8,\ldots,\floor{\frac
n2}\}$ into $t$ modulo $n$ difference triples and a collection of connection sets for circulant
graphs such that the circulant graphs can be decomposed into Hamilton cycles (or Hamilton cycles
and a perfect matching) using the results in Section \ref{hamcycledecompsofcayleygraphs}. Suitable
partitions are given in the following tables.

\noindent
{\bf t=2:}
$$
\begin{array}{|c|c|c|}
\hline
n&{\rm modulo\ } n&{\rm connection\ sets}\\
&{\rm difference\ triples}&\\
\hline
\hline
29&\{6,9,14\},\{8,10,11\}&\{12,13\}\\
\hline

30&\{6,8,14\},\{9,10,11\}&\{12,13\}\\
\hline

31&\{6,8,14\},\{9,10,12\}&\{11,13\},\{15\}\\
\hline

32&\{6,8,14\},\{9,10,13\}&\{11,12\},\{15\}\\
\hline

33&\{6,8,14\},\{10,11,12\}&\{9,16\},\{13,15\}\\
\hline

34&\{6,8,14\},\{10,11,13\}&\{9,12\},\{15,16\}\\
\hline

35&\{6,8,14\},\{10,12,13\}&\{9,11\},\{15,16\},\{17\}\\
\hline

36&\{6,8,14\},\{11,12,13\}&\{9,10\},\{15,16\},\{17\}\\
\hline
\end{array}
$$

\noindent
{\bf t=3:}
$$
\begin{array}{|c|c|c|}
\hline
n&{\rm modulo\ } n&{\rm connection\ sets}\\
&{\rm difference\ triples}&\\
\hline
\hline
35&\{6,8,14\},\{9,11,15\},\{10,12,13\}&\{16,17\}\\
\hline

36&\{6,8,14\},\{9,10,17\},\{11,12,13\}&\{15,16\}\\
\hline

37&\{6,8,14\},\{9,11,17\},\{10,12,15\}&\{13,16\},\{18\}\\
\hline

38&\{6,9,15\},\{8,10,18\},\{11,13,14\}&\{12,16,19\},\{17\}\\
\hline

39&\{6,9,15\},\{8,10,18\},\{12,13,14\}&\{11,16\},\{17,19\}\\
\hline

40&\{6,9,15\},\{8,10,18\},\{11,13,16\}&\{12,17\},\{14,19\}\\
\hline

41&\{6,9,15\},\{8,10,18\},\{12,13,16\}&\{11,14\},\{17,19\},\{20\}\\
\hline

42&\{6,9,15\},\{8,10,18\},\{12,14,16\}&\{11,13\},\{17\},\{19,20\}\\
\hline
\end{array}
$$

\noindent
{\bf t=4:}
$$
\begin{array}{|c|c|c|}
\hline
n&{\rm modulo\ } n&{\rm connection\ sets}\\
&{\rm difference\ triples}&\\
\hline
\hline
41&\{6,8,14\},\{9,15,17\},\{10,11,20\},\{12,13,16\}&\{18,19\}\\
\hline

42&\{6,8,14\},\{9,16,17\},\{10,12,20\},\{11,13,18\}&\{15,19\}\\
\hline

43&\{6,11,17\},\{8,12,20\},\{9,10,19\},\{13,14,16\}&\{15,18\},\{21\}\\
\hline

44&\{6,10,16\},\{8,13,21\},\{9,11,20\},\{12,14,18\}&\{15,17\},\{19\}\\
\hline

45&\{6,11,17\},\{8,12,20\},\{9,10,19\},\{13,14,18\}&\{15,16\},\{21,22\}\\
\hline

46&\{6,11,17\},\{8,12,20\},\{9,10,19\},\{13,15,18\}&\{14,16,23\},\{21,22\}\\
\hline

47&\{6,11,17\},\{8,12,20\},\{9,10,19\},\{13,16,18\}&\{14,15\},\{21,22\},\{23\}\\
\hline

48&\{6,11,17\},\{8,12,20\},\{9,10,19\},\{14,16,18\}&\{13,15\},\{21,22\},\{23\}\\
\hline
\end{array}
$$

\vspace{0.5cm}

We now deal with $n\geq 6t+25$ and $(n,t)$ not covered earlier.  We define
$S_t$ by $S_t = \{6\} \cup \{8,\ldots,3t+12\}$.  The following table gives a
partition $\pi_t$ of $S_t$ into difference triples, a set $R_t$ and a
difference quadruple $Q_t$ such that $Q_t$ can be partitioned into two pairs of
relatively prime integers, and $R_t$ is a pair of relatively prime integers.

$$
\begin{array}{|c|c|}
\hline
t&\pi_t \\
\hline \hline
2&\{\{6,12,18\},\{8,9,17\},\{13,14,15,16\},\{10,11\}\}\\
\hline
3&\{\{6,14,20\},\{8,13,21\},\{9,10,19\},\{15,16,17,18\},\{11,12\}\}\\
\hline
4&\{\{6,11,17\},\{8,12,20\},\{9,13,22\},\{10,14,24\},\{15,16,18,19\},\{21,23\}\}\\
\hline
\end{array}
$$

Thus, $\lan{Q_t}$ can be decomposed into two connected $4$-regular Cayley
graphs, which in turn can be decomposed into Hamilton cycles using Theorem
\ref{BFMTheorem}. It follows that there is both a
$(3^{tn},4^n,n^2)$-decomposition and a $(3^{tn},n^6)$-decomposition of
$\lan{S_t}$. If $q=0$, then we use the $(3^{tn},n^6)$-decomposition of
$\lan{S_t}$ and if $q\geq1$, then we use the $(3^{tn},4^n,n^2)$-decomposition
of $\lan{S_t}$. This leaves us needing an $(n^{h-6})$-decomposition of
$K_n-\lan{\{1,2,3,4,5,7\}\cup S_t}$ when $q=0$, and a
$(4^{(q-1)n},n^{h-2})$-decomposition of $K_n-\lan{\{1,2,3,4,5,7\}\cup S_t}$
when $q\geq1$. Note that $K_n-\lan{\{1,2,3,4,6\}\cup S_t}$ is isomorphic to
$\lan{\{3t+13,\ldots,\floor{\frac n2}\}}$ so the required decomposition exists
by Lemma \ref{hamsand4snogap}.
\vspace{0.5cm}

\noindent{\bf Case C3:}\quad Suppose that $t\in\{9,10,11,12\}$. Redefine $S_t$
by $S_t=\{6\} \cup \{8,\ldots,3t+8\}$. The following table gives a partition
of $S_t$ into difference triples and a set $R_t$ consisting of a pair of
relatively prime integers. Thus, $\lan{R_t}$ is a connected $4$-regular Cayley
graph, and so can be decomposed into two Hamilton cycles using Theorem
\ref{BFMTheorem}. Thus, we have a $(3^{tn},n^2)$-decomposition of $\lan{S_t}$.
$$
\begin{array}{|c|c|c|}
\hline
t & \mbox{difference triples} & R_t \\
\hline
\hline
9&\{6,20,26\},\{8,21,29\},\{9,22,31\},\{10,23,33\},\{11,24,35\},&\\
&\{12,16,28\},\{13,17,30\},\{14,18,32\},\{15,19,34\}&\{25,27\}\\
\hline
10&\{6,23,29\},\{8,20,28\},\{9,18,27\},\{10,22,32\},\{11,25,36\},\{12,26,38\},&\\
&\{13,24,37\},\{14,21,35\},\{15,19,34\},\{16,17,33\}&\{30,31\}\\
\hline
11&\{6,25,31\},\{8,24,32\},\{9,21,30\},\{10,19,29\},\{11,26,37\},\{12,27,39\},&\\
&\{13,28,41\},\{14,20,34\},\{15,18,33\},\{16,22,38\},\{17,23,40\}&\{35,36\}\\
\hline
12&\{6,27,33\},\{8,23,31\},\{9,25,34\},\{10,22,32\},\{11,24,35\},&\\
&\{12,28,40\},\{13,29,42\},\{14,30,44\},\{15,26,41\},&\\
&\{16,21,37\},\{17,19,36\},\{18,20,38\}&\{39,43\}\\
\hline
\end{array}
$$
Thus, we only require a $(4^{qn},n^{h-2})$-decomposition of
$K_n-\lan{\{1,2,3,4,5,7\}\cup S_t}$. But $K_n-\lan{\{1,2,3,4,5,7\}\cup S_t}$
is isomorphic to $\lan{\{3t+9,\ldots,\floor{\frac n2}\}}$ and so this
decomposition exists by Lemma \ref{hamsand4snogap}.

\vspace{0.5cm}

\noindent{\bf Case C4:}\quad Suppose that $t\geq 13$. Redefine $S_t$ by
$S_t=\{6\} \cup \{8,\ldots,3t+6\}$ when $t\equiv 2,3\md 4$, and $S_t=\{6\} \cup
\{8,\ldots,3t+5\} \cup \{3t+7\}$ when $t\equiv 0,1\md 4$. We now obtain a
$(3^{tn})$-decomposition of $\lan{S_t}$. For $13\leq t\leq 60$, we have found
such a decomposition by partitioning $S_t$ into difference triples with the aid
of a computer. These difference triples are shown in Table
\ref{tab:689_triples} in the appendix.  For $t \geq 61$ and $t\equiv 2,3 \md 4$
(respectively $t\equiv 0,1 \md 4$) we set aside as one difference triple
$\{6,3t,3t+6\}$ (respectively $\{6,3t+1,3t+7\}$) and form the set $S'_t =
\{8,\ldots,3t-1,3t+1,\ldots,3t+5\} = \{8,\ldots,3t+5\} \setminus \{3t\}$
(respectively $S'_t = \{8,\ldots,3t,3t+2,\ldots,3t+5\} = \{8,\ldots,3t+5\}
\setminus \{3t+1\}$).  We can obtain a $(3^{(t-1)n})$-decomposition of
$\lan{S'_t}$ by using an extended Langford sequence of order $t-1$ and defect
$8$ to partition $S'_t$ into difference triples. Since $t \geq 61$, this
sequence exists by Theorem 7.1 in \cite{LinMor} (also see \cite{Sha,Sim}). So
we have a $(3^{tn})$-decomposition of $\lan{S_t}$, and require a
$(4^{qn},n^h)$-decomposition of $K_n-\lan{\{1,2,3,4,5,7\}\cup S_t}$. Since
$K_n-\lan{\{1,2,3,4,5,7\}\cup S_t}$ is isomorphic to
\begin{itemize}
\item $\lan{\{3t+7,\ldots,\floor{\frac n2}\}}$ when $t\equiv 2,3\md 4$; and
\item $\lan{\{3t+6\}\cup\{3t+8,\ldots,\floor{\frac n2}\}}$ when $t\equiv 0,1\md 4$;
\end{itemize}
this decomposition exists by Lemma \ref{hamsand4snogap} or \ref{hamsand4swithgap}. \qed

\vspace{0.3cm}

\noindent{\bf Case D: ($S = \{1,2,3,4,5,6,7\}$)} \quad
The conditions $h\geq 2$ and $3t+4q+h=\lfloor\frac{n-1}2\rfloor-7$ imply
$n\geq 6t+19$. If $t=0$, then the result follows immediately by Lemma
\ref{hamsand4swithgap}. We deal separately with the five cases $t\in\{1,2\}$,
$t\in\{3,4,5,6\}$, $t\in\{7,8,9,10\}$, $t\in\{11,12,13,14\}$ and $t\geq 15$.

\vspace{0.5cm}

\noindent{\bf Case D1: }\quad Suppose that $t\in\{1,2\}$. The cases $6t+19\leq
n\leq 6t+30$ and the cases $$(n,t)\in\{(38,1),(39,1),(40,1),(44,2),(45,2)\}$$
are dealt with first. Since $h\geq 2$, it follows from
$3t+4q+h=\floor{\frac{n-1}2}-5$ that in each of these cases we have
$q\in\{0,1\}$.
The required decompositions are obtained by partitioning
$\{8,\ldots,\floor{\frac n2}\}$ into $t$ modulo $n$ difference triples and a
collection of connection sets for circulant graphs such that the circulant
graphs can be decomposed into Hamilton cycles (or Hamilton cycles and a
perfect matching) using the results in Section
\ref{hamcycledecompsofcayleygraphs}. Whenever we need to consider
$q=1$, one of the given connection sets has cardinality 4 and is a
difference quadruple, which means that the corresponding circulant
graph has a $(4^n)$-decomposition. Suitable partitions are given in the
following tables, noting that for $t=1$ we have $n\not\in\{25,26\}$ and for
$t=2$ we have $n\neq 31$..

\noindent
{\bf t=1:}
$$
\begin{array}{|c|c|c|}
\hline
n&{\rm modulo\ } n&{\rm connection\ sets}\\
&{\rm difference\ triples}&\\
\hline
\hline
27&\{8,9,10\}&\{11,12\},\{13\}\\
\hline

28&\{8,9,11\}&\{10,12,13\}\\
\hline

29&\{8,9,12\}&\{10,11\},\{13,14\}\\
\hline

30&\{8,9,13\}&\{10,11\},\{12,14,15\}\\
\hline

31&\{8,9,14\}&\{10,11\},\{12,13\},\{15\}\\
\hline

32&\{8,9,15\}&\{10,11\},\{12,13,14\}\\
\hline

33&\{10,11,12\}&\{8,9,13,14\},\{15,16\}\\
\hline

34&\{8,10,16\}&\{9,11,13,15\},\{12,14,17\}\\
\hline

35&\{10,12,13\}&\{8,9,14,15\},\{11\},\{16,17\}\\
\hline

36&\{10,12,14\}&\{8,9,15,16\},\{11,13\},\{17\}\\
\hline

38&\{10,12,16\}&\{8,9,14,15\},\{11,13\},\{17,18\}\\
\hline

39&\{10,13,16\}&\{8,9,11,12\},\{14,15\},\{17,18\},\{19\}\\
\hline

40&\{10,14,16\}&\{8,9,11,12\},\{13,15\},\{17,18\},\{19\}\\
\hline
\end{array}
$$

\noindent
{\bf t=2:}
$$
\begin{array}{|c|c|c|}
\hline
n&{\rm modulo\ } n&{\rm connection\ sets}\\
&{\rm difference\ triples}&\\
\hline
\hline
32&\{8,10,14\},\{9,11,12\}&\{13,15\}\\
\hline

33&\{8,10,15\},\{9,11,13\}&\{12,14\},\{16\}\\
\hline

34&\{8,12,14\},\{10,11,13\}&\{9\},\{15,16\}\\
\hline

35&\{8,11,16\},\{10,12,13\}&\{9,14\},\{15,17\}\\
\hline

36&\{8,13,15\},\{10,12,14\}&\{9,11\},\{16,17\}\\
\hline

37&\{8,10,18\},\{11,12,14\}&\{9,13\},\{15,16\},\{17\}\\
\hline

38&\{8,10,18\},\{11,12,15\}&\{9\},\{13,14\},\{16,17\}\\
\hline

39&\{8,10,18\},\{11,12,16\}&\{9,14\},\{13,15,17,19\}\\
\hline

40&\{8,10,18\},\{11,12,17\}&\{9,19\},\{13,14,15,16\}\\
\hline

41&\{8,10,18\},\{12,14,15\}&\{9,11\},\{13,16,17,20\},\{19\}\\
\hline

42&\{8,10,18\},\{12,14,16\}&\{9,11,13,15\},\{17\},\{19,20\}\\
\hline

44&\{8,9,17\},\{10,14,20\}&\{11,12\},\{13,15\},\{16,18,19,21\}\\
\hline

45&\{8,10,18\},\{14,15,16\}&\{9,11,17,19\},\{12,13\},\{20,21\},\{22\}\\
\hline

\end{array}
$$

\vspace{0.5cm}

We now deal with $n\geq 6t+30$  and $(n,t)$ not covered earlier. This implies
$4q+h\geq 8$. Define $S_t$ by $S_t= \{8,\ldots,3t+14\}\cup\{3t+16\}$. The
following table gives a partition $\pi_t$ of $S_t$ into difference triples and
two difference quadruples $Q'_t$ and $Q'_t$ such that $Q_t$ and $Q'_t$ can
each be partitioned into two pairs of relatively prime integers.

$$
\begin{array}{|c|c|}
\hline
t&\pi_t \\
\hline  \hline
1&\{\{8,11,19\},\{9,10,12,13\},\{14,15,16,17\}\}\\
\hline
2&\{\{8,10,18\},\{9,11,20\},\{12,15,19,22\},\{13,14,16,17\}\}\\
\hline
\end{array}
$$

\vspace{0.5cm}

Thus, each of $\lan{Q_t}$ and $\lan{Q'_t}$ can be decomposed into two
connected $4$-regular Cayley graphs, which in turn can be decomposed into
Hamilton cycles using Theorem \ref{BFMTheorem}. It follows that there is a
$(3^{tn},4^{2n})$-decomposition, a $(3^{tn},4^n,n^4)$-decomposition and a
$(3^{tn},n^8)$-decomposition of $\lan{S_t}$. If $q=0$, then we use the
$(3^{tn},n^8)$-decomposition of $\lan{S_t}$, if $q=1$, then we use the
$(3^{tn},4^n,n^4)$-decomposition of $\lan{S_t}$, and if $q\geq2$, then we use
the $(3^{tn},4^{2n})$-decomposition of $\lan{S_t}$.

This leaves us needing an $(n^{h-8})$-decomposition of
$K_n-\lan{\{1,2,3,4,5,6,7\}\cup S_t}$ when $q=0$, a $(n^{h-4})$-decomposition
of $K_n-\lan{\{1,2,3,4,5,6,7\}\cup S_t}$ when $q=1$, and a
$(4^{(q-2)n},n^h)$-decomposition of $K_n-\lan{\{1,2,3,4,5,6,7\}\cup S_t}$ when
$q\geq2$. Note that $K_n-\lan{\{1,2,3,4,5,6,7\}\cup S_t}$ is isomorphic to
$\lan{\{3t+15\} \cup \{3t+17,\ldots,\floor{\frac n2}\}}$.

When the required number of Hamilton cycles (that is, $h-8$ when $q=0$, $h-4$
when $q=1$ and $h$ when $q\geq2$) is at least $2$, the required decomposition
exists by Lemma \ref{hamsand4swithgap}. So we need to consider only the cases
where
\begin{itemize}
\item $q=0$ and $h\in\{8,9\}$; and
\item $q=1$ and $h\in\{4,5\}$.
\end{itemize}

Since $3t+4q+h=\floor{\frac{n-1}{2}}-7$, and since we have already dealt with the cases where
$(n,t)\in\{(38,1),(39,1),(40,1),(44,2),(45,2)\}$, this leaves only
$$(n,t,h) \in \{(37,1,8),(43,2,8),(46,2,9)\}$$
when $q=0$, and
$$(n,t,h) \in \{(37,1,4),(43,2,4),(46,2,5)$$
when $q=1$.

In both cases ($q=0$ and $q=1$ respectively), when the required number of Hamilton cycles
($h-8$ and $h-4$ respectively) is $0$ we have $n$ odd, and when the
required number of Hamilton cycles is $1$ we have $n\equiv 10\md{12}$.
Thus, the decomposition exists by Lemma \ref{hamsand4swithgap}.

\vspace{0.5cm}

\noindent{\bf Case D2:}\quad Suppose that $t\in\{3,4,5,6\}$. Redefine $S_t$ by
$S_t=\{8,\ldots,3t+13\}$.

The case $6t+19 \leq n \leq 6t+26$ is dealt with first. Since $h\geq 2$, it follows from
$3t+4q+h=\floor{\frac{n-1}2}-7$ that in each of these cases we have
$q=0$.
The required decompositions are obtained by partitioning
$\{8,\ldots,\floor{\frac n2}\}$ into $t$ modulo $n$ difference triples and a
collection of connection sets for circulant graphs such that the circulant
graphs can be decomposed into Hamilton cycles (or Hamilton cycles and a
perfect matching) using the results in Section
\ref{hamcycledecompsofcayleygraphs}.
Suitable partitions are given in the
following tables.

\noindent
{\bf t=3:}
$$
\begin{array}{|c|c|c|}
\hline
n&{\rm modulo\ } n&{\rm connection\ sets}\\
&{\rm difference\ triples}&\\
\hline
\hline
37&\{8,10,18\},\{9,13,15\},\{11,12,14\}&\{16,17\}\\
\hline

38&\{8,10,18\},\{9,13,16\},\{11,12,15\}&\{14,17\}\\
\hline

39&\{8,10,18\},\{9,13,17\},\{11,12,16\}&\{14,15\},\{19\}\\
\hline

40&\{8,10,18\},\{9,12,19\},\{11,13,16\}&\{14,15\},\{17\}\\
\hline

41&\{8,10,18\},\{9,11,20\},\{12,13,16\}&\{14,15\},\{17,19\}\\
\hline

42&\{8,10,18\},\{9,11,20\},\{12,13,17\}&\{14,15\},\{16,19\}\\
\hline

43&\{8,10,18\},\{9,11,20\},\{12,14,17\}&\{13,15\},\{16,19\},\{21\}\\
\hline

44&\{8,10,18\},\{9,11,20\},\{13,14,17\}&\{12,15\},\{16,19\},\{21\}\\
\hline

\end{array}
$$

\noindent
{\bf t=4:}
$$
\begin{array}{|c|c|c|}
\hline
n&{\rm modulo\ } n&{\rm connection\ sets}\\
&{\rm difference\ triples}&\\
\hline
\hline
43&\{8,10,18\},\{9,12,21\},\{11,15,17\},\{13,14,16\}&\{19,20\}\\
\hline

44&\{8,10,18\},\{9,12,21\},\{11,14,19\},\{13,15,16\}&\{17,20\}\\
\hline

45&\{8,10,18\},\{9,12,21\},\{11,14,20\},\{13,15,17\}&\{16,19\},\{22\}\\
\hline

46&\{8,10,18\},\{9,12,21\},\{11,15,20\},\{13,14,19\}&\{16,17\},\{22,23\}\\
\hline

47&\{8,10,18\},\{9,12,21\},\{11,14,22\},\{13,15,19\}&\{16,17\},\{20,23\}\\
\hline

48&\{8,10,18\},\{9,12,21\},\{11,14,23\},\{13,15,20\}&\{16,17\},\{19,22\}\\
\hline

49&\{8,10,18\},\{9,12,21\},\{11,15,23\},\{13,14,22\}&\{16,17\},\{19,20\},\{24\}\\
\hline

50&\{8,10,18\},\{9,12,21\},\{11,17,22\},\{13,14,23\}&\{15,16\},\{19,20\},\{24,25\}\\
\hline

\end{array}
$$

\noindent
{\bf t=5:}
$$
\begin{array}{|c|c|c|}
\hline
n&{\rm modulo\ } n&{\rm connection\ sets}\\
&{\rm difference\ triples}&\\
\hline
\hline
49&\{8,12,20\},\{9,14,23\},\{10,11,21\}, &\{22,24\}\\
&\{13,17,19\},\{15,16,18\}&\\
\hline

50&\{8,16,24\},\{9,11,20\},\{10,12,22\}, &\{17,23\}\\
&\{13,18,19\},\{14,15,21\}&\\
\hline

51&\{8,16,24\},\{9,11,20\},\{10,12,22\}, &\{17,21\},\{25\}\\
&\{13,15,23\},\{14,18,19\}&\\
\hline

52&\{8,16,24\},\{9,11,20\},\{10,12,22\}, &\{17,21\},\{23\}\\
&\{13,14,25\},\{15,18,19\}&\\
\hline

53&\{8,16,24\},\{9,11,20\},\{10,12,22\}, &\{17,19\},\{23,26\}\\
&\{13,15,25\},\{14,18,21\}&\\
\hline

54&\{8,16,24\},\{9,11,20\},\{10,12,22\}, &\{17,19\},\{21,26\}\\
&\{13,18,23\},\{14,15,25\}&\\
\hline

55&\{8,16,24\},\{9,11,20\},\{10,12,22\}, &\{15,19\},\{21,26\},\{27\}\\
&\{13,17,25\},\{14,18,23\}&\\
\hline

56&\{8,16,24\},\{9,11,20\},\{10,12,22\}, &\{15,18\},\{21,25\},\{27\}\\
&\{13,17,26\},\{14,19,23\}&\\
\hline

\end{array}
$$

\noindent
{\bf t=6:}
$$
\begin{array}{|c|c|c|}
\hline
n&{\rm modulo\ } n&{\rm connection\ sets}\\
&{\rm difference\ triples}&\\
\hline
\hline
55&\{8,21,26\},\{10,20,25\},\{12,19,24\}, &\{9,11\}\\
&\{13,15,27\},\{14,18,23\},\{16,17,22\}&\\
\hline

56&\{8,9,17\},\{10,19,27\},\{11,22,23\}, &\{15,21\}\\
&\{12,20,24\},\{13,18,25\},\{14,16,26\}&\\
\hline

57&\{8,15,23\},\{9,16,25\},\{10,17,27\}, &\{21,22\},\{28\}\\
&\{11,13,24\},\{12,14,26\},\{18,19,20\}&\\
\hline

58&\{8,15,23\},\{9,16,25\},\{10,18,28\}, &\{20,21\},\{27\}\\
&\{11,13,24\},\{12,14,26\},\{17,19,22\}&\\
\hline

59&\{8,15,23\},\{9,16,25\},\{10,18,28\}, &\{19,21\},\{27,29\}\\
&\{11,13,24\},\{12,14,26\},\{17,20,22\}&\\
\hline

60&\{8,15,23\},\{9,16,25\},\{10,18,28\}, &\{19,20\},\{27,29\}\\
&\{11,13,24\},\{12,14,26\},\{17,21,22\}&\\
\hline

61&\{8,15,23\},\{9,16,25\},\{10,18,28\}, &\{17,21\},\{27,29\},\{30\}\\
&\{11,13,24\},\{12,14,26\},\{19,20,22\}&\\
\hline

62&\{8,15,23\},\{9,16,25\},\{10,18,28\}, &\{17,20\},\{27\},\{29,30\}\\
&\{11,13,24\},\{12,14,26\},\{19,21,22\}&\\
\hline

\end{array}
$$

We now deal with $n\geq 6t+27$. This implies
$4q+h\geq 6$. Define $S_t$ by $S_t= \{8,\ldots,3t+13\}$. The following table
gives a partition $\pi_t$ of $S_t$ into difference triples, a set $R_t$ of
two relatively prime integers, and a difference
quadruple $Q_t$ such that $Q_t$ can be partitioned into two
pairs of relatively prime integers.

$$
\begin{array}{|c|c|}
\hline
t&\pi_t \\
\hline \hline
3&\{\{8,12,20\},\{9,13,22\},\{10,11,21\},\{14,15,16,17\},\{18,19\}\}\\
\hline
4&\{\{8,12,20\},\{9,14,23\},\{10,15,25\},\{11,13,24\},\{16,17,18,19\},\{21,22\}\}\\
\hline
5&\{\{8,14,22\},\{9,15,24\},\{10,16,26\},\{11,17,28\},\{12,13,25\},\\
&\{18,19,20,21\},\{23,27\}\}\\
\hline
6&\{\{8,16,24\},\{9,17,26\},\{10,18,28\},\{11,19,30\},\{12,13,25\},\{14,15,29\}\\
&\{20,21,22,23\},\{27,31\}\}\\
\hline
\end{array}
$$

\vspace{0.5cm}

Thus, $\lan{R_t}$ is a connected 4-regular Cayley graph, and so can be
decomposed into 2 Hamilton cycles using Theorem \ref{BFMTheorem}. Additionally
$\lan{Q_t}$ can be decomposed into two connected $4$-regular Cayley graphs,
which in turn can be decomposed into Hamilton cycles using Theorem
\ref{BFMTheorem}. It follows that there is a $(3^{tn},4^n,n^2)$-decomposition,
and a
$(3^{tn},n^6)$-decomposition of $\lan{S_t}$. If $q=0$, then we use the
$(3^{tn},n^6)$-decomposition of $\lan{S_t}$, and if $q\geq1$, then we use
the $(3^{tn},4^n,n^2)$-decomposition of $\lan{S_t}$.

This leaves us needing an
$(n^{h-6})$-decomposition of $K_n-\lan{\{1,2,3,4,5,6,7\}\cup S_t}$ when $q=0$,
and a
$(4^{(q-1)n},n^{h-2})$-decomposition of $K_n-\lan{\{1,2,3,4,5,6,7\}\cup S_t}$ when
$q\geq2$. Note that $K_n-\lan{\{1,2,3,4,5,6,7\}\cup S_t}$ is isomorphic to
$\lan{\{3t+14,\ldots,\floor{\frac n2}\}}$, so the decomposition exists by
Lemma \ref{hamsand4snogap}.

\vspace{0.5cm}

\noindent{\bf Case D3:}\quad Suppose that $t\in\{7,8,9,10\}$.
Redefine $S_t$ by
$S_t=\{8,\ldots,3t+11\}$ when $t\in\{7,8\}$ and $S_t=\{8,\ldots,3t+10\} \cup
\{3t+12\}$ when $t\in\{9,10\}$.

The case $6t+19 \leq n \leq 6t+22$ is dealt with first, along with the cases $(n,t)\in\{(78,9),(80,9),$ $(84,10)\}$. Since $h\geq 2$, it follows from
$3t+4q+h=\floor{\frac{n-1}2}-7$ that in each of these cases we have
$q=0$.
The required decompositions are obtained by partitioning
$\{8,\ldots,\floor{\frac n2}\}$ into $t$ modulo $n$ difference triples and a
collection of connection sets for circulant graphs such that the circulant
graphs can be decomposed into Hamilton cycles (or Hamilton cycles and a
perfect matching) using the results in Section
\ref{hamcycledecompsofcayleygraphs}.
Suitable partitions are given in the
following tables.

\noindent
{\bf t=7:}
$$
\begin{array}{|c|c|c|}
\hline
n&{\rm modulo\ } n&{\rm connection\ sets}\\
&{\rm difference\ triples}&\\
\hline
\hline
61&\{8,18,26\},\{9,19,28\},\{10,20,30\},\{11,16,27\}, &\{22,24\}\\
&\{12,13,25\},\{14,15,29\},\{17,21,23\}&\\
\hline

62&\{8,18,26\},\{9,19,28\},\{10,20,30\},\{11,16,27\}, &\{21,24\}\\
&\{12,13,25\},\{14,15,29\},\{17,22,23\}&\\
\hline

63&\{8,18,26\},\{9,19,28\},\{10,20,30\},\{11,16,27\}, &\{21,23\},\{31\}\\
&\{12,13,25\},\{14,15,29\},\{17,22,24\}&\\
\hline

64&\{8,18,26\},\{9,19,28\},\{10,20,30\},\{11,16,27\}, &\{21,22\},\{31\}\\
&\{12,13,25\},\{14,15,29\},\{17,23,24\}&\\
\hline

\end{array}
$$

\noindent
{\bf t=8:}
$$
\begin{array}{|c|c|c|}
\hline
n&{\rm modulo\ } n&{\rm connection\ sets}\\
&{\rm difference\ triples}&\\
\hline
\hline
67&\{8,18,26\},\{9,20,29\},\{10,21,31\},\{11,22,33\}, &\{24,30\}\\
&\{12,16,28\},\{13,14,27\},\{15,17,32\},\{19,23,25\}&\\
\hline

68&\{8,18,26\},\{9,20,29\},\{10,21,31\},\{11,22,33\}, &\{23,30\}\\
&\{12,16,28\},\{13,14,27\},\{15,17,32\},\{19,24,25\}&\\
\hline

69&\{8,18,26\},\{9,20,29\},\{10,21,31\},\{11,25,33\}, &\{19,30\},\{34\}\\
&\{12,16,28\},\{13,14,27\},\{15,17,32\},\{22,23,24\}&\\
\hline

70&\{8,18,26\},\{9,19,28\},\{10,20,30\},\{11,21,32\}, &\{24,34,35\},\{33\}\\
&\{12,15,27\},\{13,16,29\},\{14,17,31\},\{22,23,25\}&\\
\hline

\end{array}
$$

\noindent
{\bf t=9:}
$$
\begin{array}{|c|c|c|}
\hline
n&{\rm modulo\ } n&{\rm connection\ sets}\\
&{\rm difference\ triples}&\\
\hline
\hline
73&\{8,20,28\},\{9,21,30\},\{10,23,33\},\{11,24,35\}, &\{29,36\}\\
&\{12,15,27\},\{13,19,32\},\{14,17,31\},\{16,18,34\},&\\
&\{22,25,26\}&\\
\hline

74&\{8,24,32\},\{9,22,31\},\{10,19,29\},\{11,16,27\}, &\{30,36,37\}\\
&\{12,21,33\},\{13,15,28\},\{14,20,34\},\{17,18,35\},&\\
&\{23,25,26\}&\\
\hline

75&\{8,21,29\},\{9,22,31\},\{10,23,33\},\{11,24,35\}, &\{26,27\},\{28\}\\
&\{12,25,37\},\{13,17,30\},\{14,18,32\},\{15,19,34\},&\\
&\{16,20,36\}&\\
\hline

76&\{8,20,28\},\{9,21,30\},\{10,23,33\},\{11,24,35\}, &\{26,36,37\}\\
&\{12,15,27\},\{13,19,32\},\{14,17,31\},\{16,18,34\},&\\
&\{22,25,29\}&\\
\hline

78&\{8,21,29\},\{9,22,31\},\{10,23,33\},\{11,24,35\}, &\{26,27\},\{28,38,39\}\\
&\{12,25,37\},\{13,17,30\},\{14,18,32\},\{15,19,34\},&\\
&\{16,20,36\}&\\
\hline

80&\{8,21,29\},\{9,22,31\},\{10,23,33\},\{11,24,35\}, &\{26,27\},\{28,38,39\}\\
&\{12,25,37\},\{13,17,30\},\{14,18,32\},\{15,19,34\},&\\
&\{16,20,36\}&\\
\hline

\end{array}
$$

\noindent
{\bf t=10:}
$$
\begin{array}{|c|c|c|}
\hline
n&{\rm modulo\ } n&{\rm connection\ sets}\\
&{\rm difference\ triples}&\\
\hline
\hline
79&\{8,21,29\},\{9,22,31\},\{10,23,33\},\{11,26,37\}, &\{36,39\}\\
&\{12,20,32\},\{13,25,38\},\{14,16,30\},\{15,19,34\},&\\
&\{17,18,35\},\{24,27,28\}&\\
\hline

80&\{8,21,29\},\{9,22,31\},\{10,23,33\},\{11,24,35\}, &\{37,39\}\\
&\{12,26,38\},\{13,17,30\},\{14,18,32\},\{15,19,34\},&\\
&\{16,20,36\},\{25,27,28\}&\\
\hline

81&\{8,20,28\},\{9,22,31\},\{10,23,33\},\{11,24,35\}, &\{36,39,40\}\\
&\{12,26,38\},\{13,17,30\},\{14,18,32\},\{15,19,34\},&\\
&\{16,21,37\},\{25,27,29\}&\\
\hline

82&\{8,22,30\},\{9,23,32\},\{10,24,34\},\{11,25,36\}, &\{28,29\},\{35\}\\
&\{12,26,38\},\{13,27,40\},\{14,19,33\},\{15,16,31\},&\\
&\{17,20,37\},\{18,21,39\}&\\
\hline

84&\{8,22,30\},\{9,23,32\},\{10,24,34\},\{11,25,36\}, &\{28,29\},\{35,41\}\\
&\{12,26,38\},\{13,27,40\},\{14,19,33\},\{15,16,31\},&\\
&\{17,20,37\},\{18,21,39\}&\\
\hline

\end{array}
$$

We now deal with $n\geq 6t+23$ and $(n,t)$ not covered earlier. This implies
$4q+h\geq 4$.
Define $S_t$ by
$S_t=\{8,\ldots,3t+11\}$ when $t\in\{7,8\}$ and $S_t=\{8,\ldots,3t+10\} \cup
\{3t+12\}$ when $t\in\{9,10\}$.
The following table
gives a partition $\pi_t$ of $S_t$ into difference triples and a difference
quadruple $Q_t$ such that $Q_t$ can be partitioned into two
pairs of relatively prime integers.

$$
\begin{array}{|c|c|}
\hline
t&\pi_t \\
\hline \hline
7&\{\{8,18,26\},\{9,19,28\},\{10,20,30\},\{11,21,32\},\{12,15,27\},\{13,16,29\},\\
&\{14,17,31\},\{22,23,24,25\}\}\\
\hline
8&\{\{8,16,24\},\{9,20,29\},\{10,21,31\},\{11,22,33\},\{12,23,35\},\{13,17,30\},\\
&\{14,18,32\},\{15,19,34\},\{25,26,27,28\}\}\\
\hline
9&\{\{8,20,28\},\{9,21,30\},\{10,25,35\},\{11,18,29\},\{12,19,31\},\{13,24,37\},\\
&\{14,22,36\},\{15,17,32\},\{16,23,39\},\{26,27,33,34\}\}\\
\hline
10&\{\{8,22,30\},\{9,23,32\},\{10,24,34\},\{11,25,36\},\{12,21,33\},\{13,18,31\},\\
&\{14,26,40\},\{15,27,42\},\{16,19,35\},\{17,20,37\},\{28,29,38,39\}\}\\
\hline
\end{array}
$$

\vspace{0.5cm}

Thus,
$\lan{Q_t}$ can be decomposed into two connected $4$-regular Cayley graphs,
which in turn can be decomposed into Hamilton cycles using Theorem
\ref{BFMTheorem}. It follows that there is a $(3^{tn},4^n)$-decomposition,
and a
$(3^{tn},n^4)$-decomposition of $\lan{S_t}$. If $q=0$, then we use the
$(3^{tn},n^4)$-decomposition of $\lan{S_t}$, and if $q\geq1$, then we use
the $(3^{tn},4^n)$-decomposition of $\lan{S_t}$.

This leaves us needing an
$(n^{h-4})$-decomposition of $K_n-\lan{\{1,2,3,4,5,6,7\}\cup S_t}$ when $q=0$,
and a
$(4^{(q-1)n},n^h)$-decomposition of $K_n-\lan{\{1,2,3,4,5,6,7\}\cup S_t}$ when
$q\geq1$. Note that $K_n-\lan{\{1,2,3,4,6,7\}\cup S_t}$ is
isomorphic to
\begin{itemize}
    \item
$\lan{\{3t+12,\ldots,\floor{\frac n2}\}}$ when $t \in \{7,8\}$; and
    \item
$\lan{\{3t+11\} \cup \{3t+13,\ldots,\floor{\frac n2}\}}$ when $t \in \{9,10\}$.
\end{itemize}

When $t \in \{7,8\}$ the required decomposition exists by Lemma
\ref{hamsand4snogap}. When $t \in \{9,10\}$ and the required number of
Hamilton cycles (that is, $h-4$ when $q=0$ and $h$ when $q\geq1$) is at least
$2$, the required decomposition exists by Lemma \ref{hamsand4swithgap}. So we
need to consider only the cases where $q=0$, $h\in\{4,5\}$ and $t\in\{9,10\}$.

Since $3t+4q+h=\floor{\frac{n-1}{2}}-7$, and since we have already dealt with
the cases where $(n,t)\in\{(78,t),(80,9),(84,10)\}$, this leaves us with the
cases where
$$(n,t,h) \in \{(77,9,4),(79,9,5),(83,10,4),(85,10,5),(86,10,5)\}.$$

In the cases $(n,t,h) \in \{(77,9,4),(83,10,4)\}$ we have that $h-4$ (the
required number of Hamilton cycles) is $0$ and $n$ is odd, and in the case
$(n,t,h) \in \{(79,9,5),(85,10,5),(86,10,5)\}$ we have that $h-4$ (the
required number of Hamilton cycles) is $1$ and $n \equiv 1,2,7 \md{12}$. So in
all these cases the required decompositions exist by Lemma
\ref{hamsand4swithgap}.

\vspace{0.5cm}

\noindent{\bf Case D4:}\quad Suppose that $t\in\{11,12,13,14\}$.
Define $S_t$ by
$S_t=\{8,\ldots,3t+9\}$

The following table
gives a partition $\pi_t$ of $S_t$ into difference triples and a set
$Q_t$ such that $Q_t$ can be partitioned into two
pairs of relatively prime integers.

$$
\begin{array}{|c|c|}
\hline
t&\pi_t \\
\hline \hline
11&\{\{8,24,32\},\{9,25,34\},\{10,26,36\},\{11,27,38\},\{12,28,40\},\{13,29,42\},\\
&\{14,19,33\},\{15,20,35\},\{16,21,37\},\{17,22,39\},\{18,23,41\},\{30,31\}\}\\
\hline
12&\{\{8,26,34\},\{9,27,36\},\{10,28,38\},\{11,24,35\},\{12,29,41\},\{13,30,43\},\\
&\{14,31,45\},\{15,18,33\},\{16,21,37\},\{17,23,40\},\{19,25,44\},\{20,22,42\},\\
&\{32,39\}\}\\
\hline
13&\{\{8,26,34\},\{9,27,36\},\{10,28,38\},\{11,29,40\},\{12,30,42\},\{13,31,44\},\\
&\{14,32,46\},\{15,33,48\},\{16,21,37\},\{17,22,39\},\{18,23,41\},\{19,24,43\},\\
&\{20,25,45\},\{35,47\}\}\\
\hline
14&\{\{8,28,36\},\{9,29,38\},\{10,30,40\},\{11,31,42\},\{12,32,44\},\{13,33,46\},\\
&\{14,34,48\},\{15,35,50\},\{16,23,39\},\{17,24,41\},\{18,19,37\},\{20,25,40\},\\
&\{21,26,47\},\{22,27,49\},\{43,51\}\}\\
\hline
\end{array}
$$

\vspace{0.5cm}

Thus,
$\lan{Q_t}$ is a connected $4$-regular Cayley graph,
which in turn can be decomposed into two Hamilton cycles using Theorem
\ref{BFMTheorem}. It follows that there is a $(3^{tn},n^2)$-decomposition
of $\lan{S_t}$, which leaves us needing an
$(n^{h-2})$-decomposition of $K_n-\lan{\{1,2,3,4,5,6,7\}\cup S_t}$. Note that $K_n-\lan{\{1,2,3,4,6,7\}\cup S_t}$ is
isomorphic to $\lan{\{3t+10,\ldots,\floor{\frac n2}\}}$ so the required decomposition exists by Lemma
\ref{hamsand4snogap}.

\noindent{\bf Case D5:}\quad Suppose that $t\geq 15$. Redefine $S_t$ by
$S_t=\{8,\ldots,3t+7\}$ when $t\equiv 0,3\md 4$, and $S_t=\{8,\ldots,3t+6\}
\cup \{3t+8\}$ when $t\equiv 1,2\md 4$. We now obtain a
$(3^{tn})$-decomposition of $\lan{S_t}$.

For $t\equiv 0,3\md 4$ (respectively $t\equiv
1,2\md 4$), we can obtain a $(3^{t)n})$-decomposition of $\lan{S_t}$ by using a
Langford sequence (respectively hooked Langford sequence) of order $t$ and
defect $8$, which exists since $t \geq 15$, to partition $S_t$ into difference
triples (see \cite{Sha,Sim}). So we have a $(3^{tn})$-decomposition of
$\lan{S_t}$, and require a $(4^{qn},n^h)$-decomposition of
$K_n-\lan{\{1,2,3,4,5,6,7\}\cup S_t}$. Since $K_n-\lan{\{1,2,3,4,5,6,7\}\cup
S_t}$ is isomorphic to
\begin{itemize}
 \item $\lan{\{3t+8,\ldots,\floor{\frac n2}\}}$ when $t\equiv 0,3\md 4$; and
 \item $\lan{\{3t+7\}\cup\{3t+9,\ldots,\floor{\frac n2}\}}$ when $t\equiv 1,2\md 4$;
\end{itemize}
this decomposition exists by Lemma \ref{hamsand4snogap} or
\ref{hamsand4swithgap}. \qed

\subsection{Proof of Lemma \ref{5sandHamsKn-S}}\label{5sHamsProof}

We now prove Lemma \ref{5sandHamsKn-S}, which we restate here for convenience.
\vspace{0.5cm}

\noindent{\bf Lemma \ref{5sandHamsKn-S}}
{\em
If $S\in\{\{1,2,3,4\},\{1,2,3,4,6\},\{1,2,3,4,5,7\},\{1,2,3,4,5,6,7\},\{1,2,3,4,5,6,7,8\}\}$
and $n\geq 2\max(S)+1$, $r\geq 0$ and $h\geq 2$ are integers satisfying
$5r+h=\lfloor\frac{n-1}2\rfloor-|S|$, then there is a $(5^{rn},n^h)$-decomposition of
$K_n-\lan{S}$.

}
\vspace{0.5cm}
\proof We give the proof for each

$$S\in\{\{1,2,3,4\},\{1,2,3,4,5,7\},\{1,2,3,4,5,6,7\},\{1,2,3,4,5,6,7,8\}\}$$
separately.

\noindent{\bf Case A: ($S = \{1,2,3,4\}$)} \quad

The conditions $h\geq 2$ and $5r+h=\lfloor\frac{n-1}2\rfloor-4$ imply $n\geq 10r+13$. If $r=0$,
then the result follows immediately by Lemma \ref{hamdecompxtonover2withgapnodd} ($n$ odd) or Lemma
\ref{hamdecompxtonover2withgapneven} ($n$ even). Thus, we assume $r\geq 1$.

Define $S_r$ by $S_r=\{5,\ldots,5r+4\}$ when $r\equiv 0,3\md 4$, and
$S_r=\{5,\ldots,5r+3\} \cup \{5r+5\}$ when $r\equiv 1,2\md 4$. We now obtain a
$(5^{rn})$-decomposition $\lan{S_r}$ by partitioning $S_r$ into difference
quintuples. We have constructed such a partition $\pi_r$ of $S$ for $1\leq
r\leq 6$, shown in Table \ref{tab:quints_1234} in the appendix, and thus assume
$r \geq 7$.

For $r\equiv 0,3\md 4$, we take a Langford sequence of order $r$ and defect
$4$, which exists since $r\geq 7$ (see \cite{Sha,Sim}), and use it to partition
$\{4,\ldots,3r+3\}$ into $r$ difference triples. We then add $1$ to each
element of each of these triples to obtain a partition of $\{5,\ldots,3r+4\}$
into $r$ triples of the form $\{a,b,c\}$ where $a+b=c+1$. It is easy to
construct the required partition of $\{5,\ldots,5r+4\}$ into difference
quintuples from this by partitioning $\{3r+5,\ldots,5r+4\}$ into $r$ pairs of
consecutive integers.

For $r\equiv 1,2\md 4$, we take a hooked Langford sequence of order $r$ and
defect $4$, which exists since $r\geq 7$ (see \cite{Sha,Sim}), and use it to
partition $\{4,\ldots,3r+2\}\cup\{3r+4\}$ into $r$ difference triples. We then
add $1$ to each element of each of these triples, except the element $3r+4$, to
obtain a partition of $\{5,\ldots,3r+4\}$ into $r-1$ triples of the form
$\{a,b,c\}$ where $a+b=c+1$, and one triple of the form $\{a,b,c\}$ where
$a+b=c+2$. It is easy to construct the required partition of
$\{5,\ldots,5r+3\}\cup\{5r+5\}$ into difference quintuples from this by
partitioning $\{3r+5,\ldots,5r+3\}\cup\{5r+5\}$ into $r-1$ pairs of consecutive
integers, and the pair $\{5r+3,5r+5\}$. The pair $\{5r+3,5r+5\}$ combines with
the triple of the form $\{a,b,c\}$ where $a+b=c+2$ to form a difference
quintuple.

So we have a $(5^{rn})$-decomposition of $\lan{S_r}$, and require an
$(n^h)$-decomposition of $K_n-\lan{S_r \cup \{1,2,3,4\}}$. Note that
$K_n-\lan{S_r \cup \{1,2,3,4\}}$ is isomorphic to
\begin{itemize}
\item $\lan{\{5r+5,\ldots,\floor{\frac n2}\}}$ when $r\equiv 0,3\md 4$; and
\item $\lan{\{5r+4\}\cup\{5r+6,\ldots,\floor{\frac n2}\}}$ when $r\equiv 1,2\md 4$.
\end{itemize}
If $r\equiv 0,3\md 4$, then the decomposition exists by Lemma \ref{hamdecompxtonover2}. If $r\equiv
1,2\md 4$, then the decomposition exists by Lemma \ref{hamdecompxtonover2withgapnodd} ($n$ odd) or
Lemma \ref{hamdecompxtonover2withgapneven} ($n$ even).

\noindent{\bf Case B: ($S = \{1,2,3,4,6\}$)} \quad

The conditions $h\geq 2$ and $5r+h=\lfloor\frac{n-1}2\rfloor-5$ imply $n\geq
10r+15$. If $r=0$, then the result follows immediately by Lemma
\ref{hamdecompxtonover2withgapnodd} ($n$ odd) or Lemma
\ref{hamdecompxtonover2withgapneven} ($n$ even). Thus, we assume $r\geq 1$.

Define $S_r$ by $S_r=\{5\} \cup \{7,\ldots,5r+5\}$ when $r\equiv 2,3\md 4$, and
$S_r=\{5\} \cup \{7,\ldots,5r+4\} \cup \{5r+6\}$ when $r\equiv 0,1\md 4$. We
now obtain a $(5^{rn})$-decomposition $\lan{S_r}$ by partitioning $S_r$ into
difference quintuples. We have constructed such a partition of $S$ for $1\leq
r\leq 30$ with the aid of a computer, shown in Table \ref{tab:quints_12346} in
the appendix, and thus assume $r\geq 31$.

For $r\geq 31$ we first partition $S_r$ into $\{5\} \cup \{7,\ldots,15\}=S_2$
and $S_r\setminus S_2$. We have already noted that $S_2$ can be partitioned
into difference quintuples, so we only need to partition $S_r\setminus S_2$
into difference quintuples. Note that $S_r\setminus S_2=\{16,\ldots,5r+5\}$
when $r\equiv 2,3\md 4$, and $S_r\setminus S_2=\{16,\ldots,5r+4\}\cup\{5r+6\}$
when $r\equiv 0,1\md 4$.

For $r\equiv 2,3\md 4$, we take a Langford sequence of order $r-2$ and defect
$15$, which exists since $r\geq 31$ (see \cite{Sha,Sim}), and use it to
partition $\{15,\ldots,3r+8\}$ into $r-2$ difference triples. We then add $1$
to each element of each of these triples to obtain a partition of
$\{16,\ldots,3r+9\}$ into $r-2$ triples of the form $\{a,b,c\}$ where
$a+b=c+1$. It is easy to construct the required partition of
$\{16,\ldots,5r+5\}$ into difference quintuples from this by partitioning
$\{3r+10,\ldots,5r+5\}$ into $r-2$ pairs of consecutive integers.

For $r\equiv 0,1\md 4$, we take a hooked Langford sequence of order $r-2$ and
defect $15$, which exists since $r\geq 31$ (see \cite{Sha,Sim}), and use it to
partition $\{15,\ldots,3r+7\}\cup\{3r+9\}$ into $r-2$ difference triples. We
then add $1$ to each element of each of these triples, except the element
$3r+9$, to obtain a partition of $\{16,\ldots,3r+9\}$ into $r-3$ triples of the
form $\{a,b,c\}$ where $a+b=c+1$, and one triple of the form $\{a,b,c\}$ where
$a+b=c+2$. It is easy to construct the required partition of
$\{16,\ldots,5r+4\}\cup\{5r+6\}$ into difference quintuples from this by
partitioning $\{3r+10,\ldots,5r+4\}\cup\{5r+6\}$ into $r-3$ pairs of
consecutive integers, and the pair $\{5r+4,5r+6\}$. The pair $\{5r+4,5r+6\}$
combines with the triple of the form $\{a,b,c\}$ where $a+b=c+2$ to form a
difference quintuple.

So we have a $(5^{rn})$-decomposition of $\lan{S_r}$, and require an $(n^h)$-decomposition of
$K_n-\lan{S_r \cup \{1,2,3,4,6\}}$. Note that $K_n-\lan{S_r \cup \{1,2,3,4,6\}}$ is isomorphic to
\begin{itemize}
\item $\lan{\{5r+6,\ldots,\floor{\frac n2}\}}$ when $r\equiv 2,3\md 4$; and
\item $\lan{\{5r+5\}\cup\{5r+7,\ldots,\floor{\frac n2}\}}$ when $r\equiv 0,1\md 4$.
\end{itemize}
If $r\equiv 2,3\md 4$, then the decomposition exists by Lemma \ref{hamdecompxtonover2}. If $r\equiv
0,1\md 4$, then the decomposition exists by Lemma \ref{hamdecompxtonover2withgapnodd} ($n$ odd) or
Lemma \ref{hamdecompxtonover2withgapneven} ($n$ even).

\noindent{\bf Case C: ($S = \{1,2,3,4,5,7\}$)} \quad

The conditions $h\geq 2$ and $5r+h=\lfloor\frac{n-1}2\rfloor-6$ imply $n\geq 10r+17$. If $r=0$,
then the result follows immediately by Lemma \ref{hamdecompxtonover2withgapnodd} ($n$ odd) or Lemma
\ref{hamdecompxtonover2withgapneven} ($n$ even). Thus, we assume $r\geq 1$.

If $r=1$ we define $S_1 = \{6,8,9,10,11,12,13\}$. Since $\{6,8,9,10,13\}$ is a
difference quintuple and $\lan{\{11,12\}}$ has a decomposition into Hamilton
cycles (by Lemma \ref{BFMTheorem}), there is a $(5^n,n^2)$-decomposition of $\lan{S_1}$. If
$h=2$, then we are finished. Otherwise, we use Lemma \ref{BFMTheorem} to obtain an
$(n^{h-2})$-decomposition of $\lan{14,\ldots,\floor{\frac{n}2}}$ and we are
finished.

If $r \geq 2$, define $S_r$ by $S_r=\{6\} \cup \{8,\ldots,5r+6\}$ when $r\equiv
1,2\md 4$, and $S_r=\{6\} \cup \{8,\ldots,5r+5\} \cup \{5r+7\}$ when $r\equiv
0,3\md 4$. We now obtain a $(5^{rn})$-decomposition $\lan{S_r}$ by partitioning
$S_r$ into difference quintuples. We have constructed such a partition of $S$
for $1\leq r\leq 32$ with the aid of a computer,
shown in Table \ref{tab:quints_123457} in the appendix, and thus assume $r \geq 33$.

For $r\geq 33$ we first partition $S_r$ into $\{6\} \cup \{8,\ldots,16\}=S_2$
and $S_r\setminus S_2$. We have already noted that $S_2$ can be partitioned
into difference quintuples, so we only need to partition $S_r\setminus S_2$
into difference quintuples. Note that $S_r\setminus S_2=\{17,\ldots,5r+6\}$
when $r\equiv 1,2\md 4$, and $S_r\setminus S_2=\{17,\ldots,5r+5\}\cup\{5r+7\}$
when $r\equiv 0,3\md 4$.

For $r\equiv 1,2\md 4$, we take a Langford sequence of order $r-2$ and defect $16$, which exists
since $r\geq 33$ (see \cite{Sha,Sim}), and use it to partition $\{16,\ldots,3r+9\}$ into $r-2$
difference triples. We then add $1$ to each element of each of these triples to obtain a partition
of $\{17,\ldots,3r+10\}$ into $r-2$ triples of the form $\{a,b,c\}$ where $a+b=c+1$. It is easy to
construct the required partition of $\{17,\ldots,5r+6\}$ into difference quintuples from this by
partitioning $\{3r+11,\ldots,5r+6\}$ into $r-2$ pairs of consecutive integers.

For $r\equiv 0,3\md 4$, we take a hooked Langford sequence of order $r-2$ and defect $16$, which
exists since $r\geq 33$ (see \cite{Sha,Sim}), and use it to partition $\{16,\ldots,3r+8\}\cup\{3r+10\}$
into $r-2$ difference triples. We then add $1$ to each element of each of these triples, except the
element $3r+10$, to obtain a partition of $\{17,\ldots,3r+10\}$ into $r-3$ triples of the form
$\{a,b,c\}$ where $a+b=c+1$, and one triple of the form $\{a,b,c\}$ where $a+b=c+2$. It is easy to
construct the required partition of $\{17,\ldots,5r+5\}\cup\{5r+7\}$ into difference quintuples
from this by partitioning $\{3r+11,\ldots,5r+5\}\cup\{5r+7\}$ into $r-3$ pairs of consecutive
integers, and the pair $\{5r+5,5r+7\}$. The pair $\{5r+5,5r+7\}$ combines with the triple of the
form $\{a,b,c\}$ where $a+b=c+2$ to form a difference quintuple.

So we have a $(5^{rn})$-decomposition of $\lan{S_r}$, and require an $(n^h)$-decomposition of
$K_n-\lan{S_r \cup \{1,2,3,4,5,7\}}$. Note that $K_n-\lan{S_r \cup \{1,2,3,4,5,7\}}$ is isomorphic to
\begin{itemize}
\item $\lan{\{5r+7,\ldots,\floor{\frac n2}\}}$ when $r\equiv 1,2\md 4$; and
\item $\lan{\{5r+6\}\cup\{5r+8,\ldots,\floor{\frac n2}\}}$ when $r\equiv 0,3\md 4$.
\end{itemize}
If $r\equiv 1,2\md 4$, then the decomposition exists by Lemma \ref{hamdecompxtonover2}. If $r\equiv
0,3\md 4$, then the decomposition exists by Lemma \ref{hamdecompxtonover2withgapnodd} ($n$ odd) or
Lemma \ref{hamdecompxtonover2withgapneven} ($n$ even).

\noindent{\bf Case D: ($S = \{1,2,3,4,5,6,7\}$)} \quad

The conditions $h\geq 2$ and $5r+h=\lfloor\frac{n-1}2\rfloor-7$ imply $n\geq 10r+19$. If $r=0$,
then the result follows immediately by Lemma \ref{hamdecompxtonover2withgapnodd} ($n$ odd) or Lemma
\ref{hamdecompxtonover2withgapneven} ($n$ even). Thus, we assume $r\geq 1$.

If $r=1$ we define $S_1 = \{8,9,10,11,12,13,14\}$. Since $\{8,9,10,13,14\}$ is a
difference quintuple and $\lan{\{11,12\}}$ has a decomposition into Hamilton
cycles (by Lemma \ref{BFMTheorem}), there is a $(5^n,n^2)$-decomposition of $\lan{S_1}$. If
$h=2$, then we are finished. Otherwise, we use Lemma \ref{BFMTheorem} to obtain an
$(n^{h-2})$-decomposition of $\lan{14,\ldots,\floor{\frac{n}2}}$ and we are
finished.

If $r \geq 2$, define $S_r$ by $S_r=\{8,\ldots,5r+7\}$ when $r\equiv
0,1\md 4$, and $S_r=\{8,\ldots,5r+6\} \cup \{5r+8\}$ when $r\equiv
2,3\md 4$. We now obtain a $(5^{rn})$-decomposition $\lan{S_r}$ by partitioning
$S_r$ into difference quintuples. We have constructed such a partition of $S$
for $2\leq r\leq 12$ with the aid of a computer, shown in Table
\ref{tab:quints_1234567} in the appendix, and thus we assume $r \geq 13$.

For $r\equiv 0,1\md 4$, we take a Langford sequence of order $r$ and defect
$7$, which exists since $r\geq 13$ (see \cite{Sha,Sim}), and use it to
partition $\{7,\ldots,3r+6\}$ into $r$ difference triples. We then add $1$ to
each element of each of these triples to obtain a partition of
$\{8,\ldots,3r+7\}$ into $r$ triples of the form $\{a,b,c\}$ where $a+b=c+1$.
It is easy to construct the required partition of $\{8,\ldots,5r+7\}$ into
difference quintuples from this by partitioning $\{3r+8,\ldots,5r+7\}$ into $r$
pairs of consecutive integers.

For $r\equiv 2,3\md 4$, we take a hooked Langford sequence of order $r$ and
defect $7$, which exists since $r\geq 13$ (see \cite{Sha,Sim}), and use it to
partition $\{7,\ldots,3r+5\}\cup\{3r+7\}$ into $r$ difference triples. We then
add $1$ to each element of each of these triples, except the element $3r+7$, to
obtain a partition of $\{8,\ldots,3r+7\}$ into $r-1$ triples of the form
$\{a,b,c\}$ where $a+b=c+1$, and one triple of the form $\{a,b,c\}$ where
$a+b=c+2$. It is easy to construct the required partition of
$\{8,\ldots,5r+6\}\cup\{5r+8\}$ into difference quintuples from this by
partitioning $\{3r+8,\ldots,5r+6\}\cup\{5r+8\}$ into $r-1$ pairs of consecutive
integers, and the pair $\{5r+6,5r+8\}$. The pair $\{5r+6,5r+8\}$ combines with
the triple of the form $\{a,b,c\}$ where $a+b=c+2$ to form a difference
quintuple.

So we have a $(5^{rn})$-decomposition of $\lan{S_r}$, and require an $(n^h)$-decomposition of
$K_n-\lan{S_r \cup \{1,2,3,4,5,6,7\}}$. Note that $K_n-\lan{S_r \cup \{1,2,3,4,5,6,7\}}$ is isomorphic to
\begin{itemize}
\item $\lan{\{5r+8,\ldots,\floor{\frac n2}\}}$ when $r\equiv 0,1\md 4$; and
\item $\lan{\{5r+7\}\cup\{5r+9,\ldots,\floor{\frac n2}\}}$ when $r\equiv 2,3\md 4$.
\end{itemize}
If $r\equiv 0,1\md 4$, then the decomposition exists by Lemma \ref{hamdecompxtonover2}. If $r\equiv
2,3\md 4$, then the decomposition exists by Lemma \ref{hamdecompxtonover2withgapnodd} ($n$ odd) or
Lemma \ref{hamdecompxtonover2withgapneven} ($n$ even).

\noindent{\bf Case E: ($S = \{1,2,3,4,5,6,7,8\}$)} \quad

The conditions $h\geq 2$ and $5r+h=\lfloor\frac{n-1}2\rfloor-8$ imply $n\geq 10r+21$. If $r=0$,
then the result follows immediately by Lemma \ref{hamdecompxtonover2withgapnodd} ($n$ odd) or Lemma
\ref{hamdecompxtonover2withgapneven} ($n$ even). Thus, we assume $r\geq 1$.

For $r=1$ and $n\in\{31,32,33,34\}$, we firstly have the following results.
For $n=31$, we can obtain the required decomposition by
noting that $\{9,11,13,14,15\}$ is a modulo 31 difference quintuple and that
$\lan{\{10,12\}}$ has a Hamilton cycle decomposition by Lemma \ref{BFMTheorem}.

For $n=32$, we can obtain the required decomposition by
noting that $\{10,12,13,14,15\}$ is a modulo 32 difference quintuple and that
$\lan{\{9,11\}}$ has a Hamilton cycle decomposition by Lemma \ref{BFMTheorem}.

For $n\in\{33,34\}$, we can obtain the required decomposition by
noting that $\{9,10,12,15,16\}$ is a difference quintuple, that
$\lan{\{11,14\}}$ has a Hamilton cycle decomposition by Lemma \ref{BFMTheorem},
and $\lan{\{13\}}$ is a Hamilton cycle.

For $r=1$ and $n\geq 35$, we define $S_1 = \{9,\ldots,17\}$. Since $\{9,10,11,14,16\}$ is a
difference quintuple and each of $\lan{\{12,13\}}$ and $\lan{\{15,17\}}$ has a decomposition into Hamilton
cycles (by Lemma \ref{BFMTheorem}), there is a $(5^n,n^4)$-decomposition of $\lan{S_1}$. If
$h=4$, then we are finished. Otherwise, we use Lemma \ref{BFMTheorem} to obtain an
$(n^{h-4})$-decomposition of $\lan{18,\ldots,\floor{\frac{n}2}}$ and we are
finished.

For $r=2$, we define $S_2 = \{9,\ldots,20\}$. Since $\{9,11,14,16,18\}$ and
$\{10,12,13,15,20\}$ are both
difference quintuples and $\lan{\{17,19\}}$ has a decomposition into Hamilton
cycles (by Lemma \ref{BFMTheorem}), there is a $(5^{2n},n^2)$-decomposition of $\lan{S_2}$. If
$h=2$, then we are finished. Otherwise, we use Lemma \ref{BFMTheorem} to obtain an
$(n^{h-2})$-decomposition of $\lan{18,\ldots,\floor{\frac{n}2}}$ and we are
finished.

For $r \geq 3$, define $S_r$ by $S_r=\{9,\ldots,5r+8\}$ when $r\equiv
0,3\md 4$, and $S_r=\{9,\ldots,5r+7\} \cup \{5r+9\}$ when $r\equiv
1,2\md 4$. We now obtain a $(5^{rn})$-decomposition $\lan{S_r}$ by partitioning
$S_r$ into difference quintuples. We have constructed such a partition of $S$
for $3\leq r\leq 14$ with the aid of a computer, shown in Table
\ref{tab:quints_12345678} in the appendix, and thus assume $r \geq 15$.

For $r\equiv 0,3\md 4$, we take a Langford sequence of order $r$ and defect
$8$, which exists since $r\geq 15$ (see \cite{Sha,Sim}), and use it to
partition $\{8,\ldots,3r+7\}$ into $r$ difference triples. We then add $1$ to
each element of each of these triples to obtain a partition of
$\{9,\ldots,3r+8\}$ into $r$ triples of the form $\{a,b,c\}$ where $a+b=c+1$.
It is easy to construct the required partition of $\{9,\ldots,5r+8\}$ into
difference quintuples from this by partitioning $\{3r+9,\ldots,5r+8\}$ into $r$
pairs of consecutive integers.

For $r\equiv 1,2\md 4$, we take a hooked Langford sequence of order $r$ and
defect $8$, which exists since $r\geq 15$ (see \cite{Sha,Sim}), and use it to
partition $\{8,\ldots,3r+6\}\cup\{3r+8\}$ into $r$ difference triples. We then
add $1$ to each element of each of these triples, except the element $3r+8$, to
obtain a partition of $\{9,\ldots,3r+8\}$ into $r-1$ triples of the form
$\{a,b,c\}$ where $a+b=c+1$, and one triple of the form $\{a,b,c\}$ where
$a+b=c+2$. It is easy to construct the required partition of
$\{9,\ldots,5r+7\}\cup\{5r+9\}$ into difference quintuples from this by
partitioning $\{3r+9,\ldots,5r+7\}\cup\{5r+8\}$ into $r-1$ pairs of consecutive
integers, and the pair $\{5r+7,5r+9\}$. The pair $\{5r+7,5r+9\}$ combines with
the triple of the form $\{a,b,c\}$ where $a+b=c+2$ to form a difference
quintuple.

So we have a $(5^{rn})$-decomposition of $\lan{S_r}$, and require an $(n^h)$-decomposition of
$K_n-\lan{S_r \cup \{1,2,3,4,5,6,7,8\}}$. Note that $K_n-\lan{S_r \cup \{1,2,3,4,5,6,7,8\}}$ is isomorphic to
\begin{itemize}
\item $\lan{\{5r+9,\ldots,\floor{\frac n2}\}}$ when $r\equiv 0,3\md 4$; and
\item $\lan{\{5r+8\}\cup\{5r+10,\ldots,\floor{\frac n2}\}}$ when $r\equiv 1,2\md 4$.
\end{itemize}
If $r\equiv 0,3\md 4$, then the decomposition exists by Lemma \ref{hamdecompxtonover2}. If $r\equiv
1,2\md 4$, then the decomposition exists by Lemma \ref{hamdecompxtonover2withgapnodd} ($n$ odd) or
Lemma \ref{hamdecompxtonover2withgapneven} ($n$ even). \qed

\newpage
\section{Data for Section \ref{sec:Sofsize3-8}}
\label{appendix}

\LTcapwidth=\textwidth
\begin{longtable}{|llll|}
    \hline
    $J_{3} \rightarrow  3,6$ & $J_{5} \rightarrow  4,5,6$ &
    $J_{6} \rightarrow  4^{3},6$ & $J_{7} \rightarrow  5^{3},6$ \\
    $J_{5} \rightarrow  4^{2},7$ & $J_{5} \rightarrow  3,5,7$ &
    $J_{7} \rightarrow  4,5^{2},7$ & $J_{9} \rightarrow  5^{4},7$ \\
    $J_{7} \rightarrow  4^{2},5,8$ & $J_{7} \rightarrow  3,5^{2},8$ &
    $J_{7} \rightarrow  3^{3},4,8$ & $J_{8} \rightarrow  4^{4},8$ \\
    $J_{9} \rightarrow  4,5^{3},8$ & $J_{11} \rightarrow  5^{5},8$ &
    & \\
    $J_{8} \rightarrow  5^{3},9$ & $J_{8} \rightarrow  3,4^{3},9$ &
    $J_{8} \rightarrow  3^{2},4,5,9$ & $J_{8} \rightarrow  3^{5},9$ \\
    $J_{9} \rightarrow  4^{2},5^{2},9$ & $J_{10} \rightarrow  4^{4},5,9$ &
    $J_{11} \rightarrow  4^{6},9$ &  \\
    $J_{9} \rightarrow  4^{3},5,10$ & $J_{9} \rightarrow  3,4,5^{2},10$ &
    $J_{9} \rightarrow  3^{3},4^{2},10$ & $J_{9} \rightarrow  3^{4},5,10$ \\
    $J_{10} \rightarrow  5^{4},10$ & $J_{10} \rightarrow  4^{5},10$ &
    $J_{11} \rightarrow  4^{2},5^{3},10$ & \\
    \hline
  \caption{These decompositions are required for Lemma \ref{123:Jgraphk=6-10}. The decompositions themselves are given in Table \ref{tab:123_decomps}.}
  \label{tab:123_6-10}
\end{longtable}

\begin{longtable}{|llll|}
    \hline
    $J_{8} \rightarrow  4^{2},5,11^{*}$ & $J_{8} \rightarrow  3,5^{2},11^{*}$ &
    $J_{8} \rightarrow  3^{3},4,11^{*}$ & $J_{9} \rightarrow  4^{4},11^{*}$ \\
    $J_{10} \rightarrow  4,5^{3},11^{*}$ & $J_{12} \rightarrow  5^{5},11^{*}$ &
    & \\
    $J_{9} \rightarrow  5^{3},12^{*}$ & $J_{9} \rightarrow  3,4^{3},12^{*}$ &
    $J_{9} \rightarrow  3^{2},4,5,12^{*}$ & $J_{9} \rightarrow  3^{5},12^{*}$ \\
    $J_{10} \rightarrow  4^{2},5^{2},12^{*}$ & $J_{11} \rightarrow  4^{4},5,12^{*}$ &
    $J_{12} \rightarrow  4^{6},12^{*}$ &  \\
    $J_{10} \rightarrow  4^{3},5,13^{*}$ & $J_{10} \rightarrow  3,4,5^{2},13^{*}$ &
    $J_{10} \rightarrow  3^{3},4^{2},13^{*}$ & $J_{10} \rightarrow  3^{4},5,13^{*}$ \\
    $J_{11} \rightarrow  5^{4},13^{*}$ & $J_{11} \rightarrow  4^{5},13^{*}$ &
    $J_{12} \rightarrow  4^{2},5^{3},13^{*}$ &  \\
    $J_{12} \rightarrow  4^{3},5^{2},14^{*}$ & $J_{12} \rightarrow  3,4,5^{3},14^{*}$ &
    $J_{12} \rightarrow  3^{2},4^{4},14^{*}$ & $J_{12} \rightarrow  3^{3},4^{2},5,14^{*}$ \\
    $J_{12} \rightarrow  3^{4},5^{2},14^{*}$ & $J_{12} \rightarrow  3^{6},4,14^{*}$ &
    $J_{13} \rightarrow  5^{5},14^{*}$ & $J_{13} \rightarrow  4^{5},5,14^{*}$ \\
\end{longtable}

\begin{longtable}{|c|>{\raggedright\arraybackslash}p{0.75\textwidth}|}
  \caption{Table of decompositions of $J_n^{\{1,2,3\}}$}
\endfoot 
\hline
\endhead
\label{tab:123_decomps}
$J_1 \rightarrow 3$ & $(0,1,3)$
\\ \hline
$J_{3} \rightarrow  4,5$ & 
$(0,1,4,2,3)$, $(1,2,5,3)$
\\ \hline
$J_4 \rightarrow 4^3$ & $(0,1,2,3)$, $(1,3,6,4)$, $(2,4,3,5)$
\\ \hline
$J_{5} \rightarrow  5^{3}$ & 
$(0,1,4,2,3)$, $(1,2,5,4,3)$, $(3,6,4,7,5)$
\\ \hline
$J_{3} \rightarrow  3,6$ & 
$(0,1,4,2,5,3)$, $(1,2,3)$
\\ \hline
$J_{5} \rightarrow  4,5,6$ & 
$(0,1,4,5,2,3)$, $(1,2,4,3)$, $(3,6,4,7,5)$
\\ \hline
$J_{6} \rightarrow  4^{3},6$ & 
$(0,1,4,5,2,3)$, $(1,2,4,3)$, $(3,6,8,5)$, $(4,7,5,6)$
\\ \hline
$J_{7} \rightarrow  5^{3},6$ & 
$(0,1,4,5,2,3)$, $(1,2,4,6,3)$, $(3,4,7,6,5)$, $(5,8,6,9,7)$
\\
\hline
\end{longtable}

\begin{longtable}{|llll|}
    \hline
$J_{3} \rightarrow  3^{2},6$ & $J_{4} \rightarrow  5^{2},6$ &
$J_{4} \rightarrow  4,5,7$ & $J_{4} \rightarrow  3^{3},7$ \\
$J_{5} \rightarrow  4^{3},8$ & $J_{5} \rightarrow  3,4,5,8$ &
$J_{5} \rightarrow  3^{4},8$ & $J_{6} \rightarrow  3^{2},5^{2},8$ \\
$J_{7} \rightarrow  5^{4},8$ & & & \\
    \hline
  \caption{These decompositions are required for Lemma \ref{1234:Jgraphk=6-8}.
  The decompositions themselves are given in Table \ref{tab:1234_decomps}.}
  \label{tab:1234_6-8}
\end{longtable}
\begin{longtable}{|llll|}
    \hline
    $J_{6} \rightarrow  5^{3},9^{*}$ & $J_{6} \rightarrow  3,4^{3},9^{*}$ &
    $J_{6} \rightarrow  3^{2},4,5,9^{*}$ & $J_{6} \rightarrow  3^{5},9^{*}$ \\
    $J_{7} \rightarrow  3^{3},5^{2},9^{*}$ & && \\
    $J_{7} \rightarrow  4^{2},5^{2},10^{*}$ & $J_{7} \rightarrow  3,5^{3},10^{*}$ &
    $J_{7} \rightarrow  3^{2},4^{3},10^{*}$ & $J_{7} \rightarrow  3^{3},4,5,10^{*}$ \\
    $J_{7} \rightarrow  3^{6},10^{*}$ & $J_{8} \rightarrow  3^{4},5^{2},10^{*}$ &
    $J_{10} \rightarrow  5^{6},10^{*}$ &  \\
    $J_{8} \rightarrow  4^{4},5,11^{*}$ & $J_{8} \rightarrow  3,4^{2},5^{2},11^{*}$ &
    $J_{8} \rightarrow  3^{2},5^{3},11^{*}$ & $J_{8} \rightarrow  3^{3},4^{3},11^{*}$ \\
    $J_{8} \rightarrow  3^{4},4,5,11^{*}$ & $J_{8} \rightarrow  3^{7},11^{*}$ &
    $J_{9} \rightarrow  5^{5},11^{*}$ & $J_{9} \rightarrow  3^{5},5^{2},11^{*}$ \\
    $J_{9} \rightarrow  4,5^{4},12^{*}$ & $J_{9} \rightarrow  4^{6},12^{*}$ &
    $J_{9} \rightarrow  3,4^{4},5,12^{*}$ & $J_{9} \rightarrow  3^{2},4^{2},5^{2},12^{*}$ \\
    $J_{9} \rightarrow  3^{3},5^{3},12^{*}$ & $J_{9} \rightarrow  3^{4},4^{3},12^{*}$ &
    $J_{9} \rightarrow  3^{5},4,5,12^{*}$ & $J_{9} \rightarrow  3^{8},12^{*}$ \\
\end{longtable}

\begin{longtable}{|c|>{\raggedright\arraybackslash}p{0.75\textwidth}|}
\caption{Table of decompositions of $J_n^{\{1,2,3,4\}}$}
\endfoot 
\hline
\label{tab:1234_decomps}
$J_{1} \rightarrow  4$ & 
$(0,4,2,3)$
\\ \hline
$J_{2} \rightarrow  3,5$ & 
$(0,4,1,5,3)$, $(2,3,4)$
\\ \hline
$J_{3} \rightarrow  3^{4}$ & 
$(0,4,3)$, $(1,5,4)$, $(2,6,4)$, $(2,5,3)$
\\ \hline
$J_{4} \rightarrow  3^{2},5^{2}$ & 
$(0,4,2,5,3)$, $(1,5,7,3,4)$, $(2,6,3)$, $(4,5,6)$
\\ \hline
$J_{5} \rightarrow  5^{4}$ & 
$(0,4,2,5,3)$, $(1,5,6,3,4)$, $(2,6,4,7,3)$, $(4,8,6,7,5)$
\\ \hline
$J_{3} \rightarrow  3^{2},6$ & 
$(0,4,1,5,2,3)$, $(2,6,4)$, $(3,4,5)$
\\ \hline
$J_{4} \rightarrow  5^{2},6$ & 
$(0,4,1,5,2,3)$, $(2,6,5,3,4)$, $(4,5,7,3,6)$
\\ \hline
$J_{4} \rightarrow  4,5,7$ & 
$(0,4,1,5,2,6,3)$, $(2,3,5,6,4)$, $(3,7,5,4)$
\\
\hline
\end{longtable}


\begin{longtable}{|llll|}
    \hline
    $J_{2} \rightarrow  4,6$ & $J_{3} \rightarrow  3^{3},6$ &
    $J_{4} \rightarrow  3^{2},4^{2},6$ & $J_{5} \rightarrow  3^{5},4,6$ \\
    $J_{6} \rightarrow  3^{8},6$ & && \\
    $J_{2} \rightarrow  3,7$ & $J_{3} \rightarrow  4^{2},7$ &
    $J_{4} \rightarrow  3^{3},4,7$ & $J_{5} \rightarrow  3^{6},7$ \\
    $J_{7} \rightarrow  3^{8},4,7$ & && \\
    $J_{3} \rightarrow  3,4,8$ & $J_{4} \rightarrow  4^{3},8$ &
    $J_{4} \rightarrow  3^{4},8$ & $J_{5} \rightarrow  3^{3},4^{2},8$ \\
    $J_{6} \rightarrow  3^{6},4,8$ & $J_{7} \rightarrow  3^{9},8$ &
    & \\
    $J_{3} \rightarrow  3^{2},9$ & $J_{4} \rightarrow  3,4^{2},9$ &
    $J_{5} \rightarrow  4^{4},9$ & $J_{5} \rightarrow  3^{4},4,9$ \\
    $J_{6} \rightarrow  3^{7},9$ & $J_{8} \rightarrow  3^{9},4,9$ &
    & \\
    $J_{4} \rightarrow  5^{2},10$ & $J_{4} \rightarrow  3^{2},4,10$ &
    $J_{5} \rightarrow  3,4^{3},10$ & $J_{5} \rightarrow  3^{5},10$ \\
    $J_{6} \rightarrow  4^{5},10$ & $J_{6} \rightarrow  3^{4},4^{2},10$ &
    $J_{7} \rightarrow  3^{7},4,10$ & $J_{8} \rightarrow  3^{10},10$ \\
        \hline
  \caption{These decompositions are required for Lemma \ref{12346:Jgraphk=6-10}.
  The decompositions themselves are given in Table \ref{tab:12346_decomps}.}
  \label{tab:12346_6-10}
\end{longtable}

\begin{longtable}{|llll|}
    \hline
    $J_{6} \rightarrow  4,5^{3},11^{*}$ & $J_{6} \rightarrow  3,4^{4},11^{*}$ &
    $J_{6} \rightarrow  3^{2},4^{2},5,11^{*}$ & $J_{6} \rightarrow  3^{3},5^{2},11^{*}$ \\
    $J_{6} \rightarrow  3^{5},4,11^{*}$ & $J_{7} \rightarrow  4^{6},11^{*}$ &
    $J_{7} \rightarrow  3^{4},4^{3},11^{*}$ & $J_{7} \rightarrow  3^{8},11^{*}$ \\
    $J_{8} \rightarrow  3^{3},4^{5},11^{*}$ & $J_{8} \rightarrow  3^{7},4^{2},11^{*}$ &
    $J_{9} \rightarrow  3^{10},4,11^{*}$ & $J_{10} \rightarrow  3^{13},11^{*}$ \\
    $J_{6} \rightarrow  4^{2},5^{2},12^{*}$ & $J_{6} \rightarrow  3,5^{3},12^{*}$ &
    $J_{6} \rightarrow  3^{2},4^{3},12^{*}$ & $J_{6} \rightarrow  3^{3},4,5,12^{*}$ \\
    $J_{6} \rightarrow  3^{6},12^{*}$ & $J_{7} \rightarrow  3,4^{5},12^{*}$ &
    $J_{7} \rightarrow  3^{5},4^{2},12^{*}$ & $J_{8} \rightarrow  4^{7},12^{*}$ \\
    $J_{8} \rightarrow  3^{4},4^{4},12^{*}$ & $J_{8} \rightarrow  3^{8},4,12^{*}$ &
    $J_{9} \rightarrow  3^{11},12^{*}$ & $J_{11} \rightarrow  3^{13},4,12^{*}$ \\
    $J_{8} \rightarrow  4^{3},5^{3},13^{*}$ & $J_{8} \rightarrow  3,4,5^{4},13^{*}$ &
    $J_{8} \rightarrow  3,4^{6},13^{*}$ & $J_{8} \rightarrow  3^{2},4^{4},5,13^{*}$ \\
    $J_{8} \rightarrow  3^{3},4^{2},5^{2},13^{*}$ & $J_{8} \rightarrow  3^{4},5^{3},13^{*}$ &
    $J_{8} \rightarrow  3^{5},4^{3},13^{*}$ & $J_{8} \rightarrow  3^{6},4,5,13^{*}$ \\
\end{longtable}
\newpage
\begin{longtable}{|c|>{\raggedright\arraybackslash}p{0.75\textwidth}|}
  \caption{Table of decompositions of $J_n^{\{1,2,3,4,6\}}$}
\endfoot
\hline
\label{tab:12346_decomps}
$J_{1} \rightarrow  5$ &
$(0,6,3,2,4)$
\\ \hline
$J_{3} \rightarrow  3,4^{3}$ &
$(0,6,2,4)$, $(1,7,4,5)$, $(2,8,5,3)$, $(3,4,6)$
\\ \hline
$J_{3} \rightarrow  3^{2},4,5$ &
$(0,6,4)$, $(1,7,4,3,5)$, $(3,2,6)$, $(2,8,5,4)$
\\ \hline
$J_{4} \rightarrow  4^{5}$ &
$(0,6,2,4)$, $(1,7,4,5)$, $(2,8,5,3)$, $(3,7,5,6)$, $(3,9,6,4)$
\\ \hline
$J_{4} \rightarrow  3^{4},4^{2}$ &
$(0,6,2,4)$, $(1,7,5)$, $(2,8,5,3)$, $(4,3,7)$, $(3,9,6)$, $(4,6,5)$
\\ \hline
$J_{5} \rightarrow  3^{3},4^{4}$ &
$(0,6,2,4)$, $(1,7,3,5)$, $(2,8,4,3)$, $(3,9,6)$, $(4,6,5)$, $(4,10,7)$,
$(5,7,6,8)$
\\ \hline
$J_{5} \rightarrow  3^{5},5^{2}$ &
$(0,6,7,3,4)$, $(2,8,6,9,3)$, $(2,6,4)$, $(5,4,8)$, $(1,7,5)$, $(3,5,6)$,
$(4,10,7)$
\\ \hline
$J_{5} \rightarrow  3^{7},4$ &
$(0,6,4)$, $(3,9,6,5)$, $(4,10,7)$, $(1,7,5)$, $(3,7,6)$, $(2,4,3)$, $(5,4,8)$,
$(2,8,6)$
\\
\hline
\end{longtable}

\begin{longtable}{|llll|}
    \hline
$J_{4} \rightarrow  3^{6},6$ & $J_{4} \rightarrow  3,5^{3},6$ &
$J_{4} \rightarrow  3^{3},4,5,6$ & $J_{4} \rightarrow  4^{2},5^{2},6$ \\
$J_{4} \rightarrow  3^{2},4^{3},6$ & $J_{5} \rightarrow  4,5^{4},6$ &
$J_{5} \rightarrow  3^{2},4^{2},5^{2},6$ & $J_{5} \rightarrow  3,4^{4},5,6$
\\
$J_{5} \rightarrow  4^{6},6$ & $J_{6} \rightarrow  5^{6},6$ &
$J_{6} \rightarrow  4^{5},5^{2},6$ & $J_{6} \rightarrow  3^{2},4^{6},6$ \\
$J_{7} \rightarrow  4^{9},6$ & $J_{9} \rightarrow  4^{12},6$ & & \\
        \hline
  \caption{These decompositions are required for Lemma \ref{123457:Jgraphk=6}.
  The decompositions themselves are given in Table \ref{tab:123457_decomps}.}
  \label{tab:123457_6}
\end{longtable}

\begin{longtable}{|llll|}
  \hline
  \caption{These decompositions are required for Lemma \ref{123457:Jgraphk=7-12}.
  The decompositions themselves are given in Table \ref{tab:123457_decomps}.}
  \endfoot
  \hline
  \endhead
  \label{tab:123457_7-12}
    $J_{5} \rightarrow  3^{6},5,7^{*}$ & $J_{5} \rightarrow  3,5^{4},7^{*}$ &
    $J_{5} \rightarrow  3^{3},4,5^{2},7^{*}$ & $J_{5} \rightarrow  3^{5},4^{2},7^{*}$ \\
    $J_{5} \rightarrow  4^{2},5^{3},7^{*}$ & $J_{5} \rightarrow  3^{2},4^{3},5,7^{*}$ &
    $J_{5} \rightarrow  3,4^{5},7^{*}$ & $J_{6} \rightarrow  4,5^{5},7^{*}$ \\
    $J_{6} \rightarrow  3^{2},4^{2},5^{3},7^{*}$ & $J_{6} \rightarrow  3,4^{4},5^{2},7^{*}$ &
    $J_{6} \rightarrow  4^{6},5,7^{*}$ & $J_{7} \rightarrow  5^{7},7^{*}$ \\
    $J_{7} \rightarrow  4^{5},5^{3},7^{*}$ & $J_{7} \rightarrow  3,4^{8},7^{*}$ &
    $J_{8} \rightarrow  4^{9},5,7^{*}$ & $J_{9} \rightarrow  3,4^{11},7^{*}$ \\
    $J_{10} \rightarrow  4^{12},5,7^{*}$ & && \\
    $J_{5} \rightarrow  3^{4},5^{2},8^{*}$ & $J_{5} \rightarrow  3^{6},4,8^{*}$ &
    $J_{5} \rightarrow  3,4,5^{3},8^{*}$ & $J_{5} \rightarrow  3^{3},4^{2},5,8^{*}$ \\
    $J_{5} \rightarrow  4^{3},5^{2},8^{*}$ & $J_{5} \rightarrow  3^{2},4^{4},8^{*}$ &
    $J_{6} \rightarrow  3,5^{5},8^{*}$ & $J_{6} \rightarrow  4^{2},5^{4},8^{*}$ \\
    $J_{6} \rightarrow  3^{2},4^{3},5^{2},8^{*}$ & $J_{6} \rightarrow  3,4^{5},5,8^{*}$ &
    $J_{6} \rightarrow  4^{7},8^{*}$ & $J_{7} \rightarrow  4,5^{6},8^{*}$ \\
\end{longtable}

\begin{longtable}{|c|>{\raggedright\arraybackslash}p{0.75\textwidth}|}
\caption{Table of decompositions of $J^{\{1,2,3,4,5,7\}}_n$}
\endfoot \hline
\label{tab:123457_decomps}
$J_{1} \rightarrow  3^{2}$ & 
$(0,7,5)$, $(6,3,7)$
\\ \hline
$J_{2} \rightarrow  3,4,5$ & 
$(0,7,5)$, $(3,7,8,1,6)$, $(6,8,4,7)$
\\ \hline
$J_{3} \rightarrow  4^{2},5^{2}$ & 
$(0,7,6,8,5)$, $(3,7,8,1,6)$, $(5,9,2,7)$, $(4,8,9,7)$
\\ \hline
$J_{3} \rightarrow  3,5^{3}$ & 
$(0,7,6,8,5)$, $(3,7,8,1,6)$, $(4,8,9,2,7)$, $(5,9,7)$
\\ \hline
$J_{3} \rightarrow  3^{2},4^{3}$ & 
$(0,7,8,5)$, $(6,3,7)$, $(1,8,6)$, $(4,8,9,7)$, $(5,9,2,7)$
\\ \hline
$J_{4} \rightarrow  4,5^{4}$ & 
$(0,7,6,8,5)$, $(3,7,8,1,6)$, $(5,9,2,7)$, $(8,3,10,6,9)$, $(4,8,10,9,7)$
\\ \hline
$J_{4} \rightarrow  4^{6}$ & 
$(0,7,8,5)$, $(1,8,9,6)$, $(5,9,2,7)$, $(3,10,9,7)$, $(6,8,4,7)$, $(3,8,10,6)$
\\ \hline
$J_{4} \rightarrow  3,4^{4},5$ & 
$(0,7,6,8,5)$, $(1,8,9,6)$, $(5,9,2,7)$, $(3,8,4,7)$, $(3,10,6)$, $(7,9,10,8)$
\\
\hline
\end{longtable}

\begin{longtable}{|lllll|}
  \caption{These decompositions required for Lemma \ref{1234567:Jgraphk=6-7}.
  For all these decompositions the $k$-cycle is incident upon vertices
  $\{4,5,\ldots,k+3\}$.
  The decompositions themselves are given in full in table
  \ref{tab:1234567_decomps}.}
  \endfoot \hline \label{tab:1234567_6-7}
    $J_{7} \rightarrow  3^{11},5^{2},6$ & $J_{7} \rightarrow  3^{6},5^{5},6$ &
    $J_{7} \rightarrow  3,5^{8},6$ & $J_{7} \rightarrow  3^{13},4,6$ &
    $J_{7} \rightarrow  3^{8},4,5^{3},6$ \\
    $J_{7} \rightarrow  3^{3},4,5^{6},6$ &
    $J_{7} \rightarrow  3^{10},4^{2},5,6$ & $J_{7} \rightarrow  3^{5},4^{2},5^{4},6$ &
    $J_{7} \rightarrow  4^{2},5^{7},6$ & $J_{7} \rightarrow  3^{7},4^{3},5^{2},6$ \\
    $J_{7} \rightarrow  3^{2},4^{3},5^{5},6$ & $J_{7} \rightarrow  3^{9},4^{4},6$ &
    $J_{7} \rightarrow  3^{4},4^{4},5^{3},6$ & $J_{7} \rightarrow  3^{6},4^{5},5,6$ &
    $J_{7} \rightarrow  3,4^{5},5^{4},6$ \\
    $J_{7} \rightarrow  3^{3},4^{6},5^{2},6$ &
    $J_{7} \rightarrow  3^{5},4^{7},6$ & $J_{7} \rightarrow  4^{7},5^{3},6$ &
    $J_{7} \rightarrow  3^{2},4^{8},5,6$ & $J_{7} \rightarrow  3,4^{10},6$ \\
    $J_{8} \rightarrow  3^{15},5,6$ & $J_{8} \rightarrow  3^{10},5^{4},6$ &
    $J_{8} \rightarrow  3^{5},5^{7},6$ & $J_{8} \rightarrow  5^{10},6$ &
    $J_{8} \rightarrow  4^{5},5^{6},6$ \\
    $J_{8} \rightarrow  4^{10},5^{2},6$ &
    $J_{9} \rightarrow  3^{19},6$ & $J_{9} \rightarrow  4^{13},5,6$ &
    $J_{10} \rightarrow  4^{16},6$ & $J_{12} \rightarrow  3^{26},6$ \\
    $J_{8} \rightarrow  3^{13},5^{2},7$ & $J_{8} \rightarrow  3^{8},5^{5},7$ &
    $J_{8} \rightarrow  3^{3},5^{8},7$ & $J_{8} \rightarrow  3^{15},4,7$ &
    $J_{8} \rightarrow  3^{10},4,5^{3},7$ \\
    $J_{8} \rightarrow  3^{5},4,5^{6},7$ &
    $J_{8} \rightarrow  4,5^{9},7$ & $J_{8} \rightarrow  3^{12},4^{2},5,7$ &
    $J_{8} \rightarrow  3^{7},4^{2},5^{4},7$ & $J_{8} \rightarrow  3^{2},4^{2},5^{7},7$ \\
\hline
\end{longtable}

\LTcapwidth=\textwidth
\begin{longtable}{|llll|}
\hline 
\endhead
\hline 
  \caption{These decompositions are required for Lemma
    \ref{1234567:Jgraphk=8-17}. These are sorted by the value of $k$ in the
    decomposition, and are given in table \ref{tab:1234567_decomps}.}
\endfoot
  \label{tab:1234567_8-17}
    $J_{6} \rightarrow  3^{8},5^{2},8^{*}$ & $J_{6} \rightarrow  3^{3},5^{5},8^{*}$ &
    $J_{6} \rightarrow  3^{10},4,8^{*}$ & $J_{6} \rightarrow  3^{5},4,5^{3},8^{*}$ \\
    $J_{6} \rightarrow  4,5^{6},8^{*}$ & $J_{6} \rightarrow  3^{7},4^{2},5,8^{*}$ &
    $J_{6} \rightarrow  3^{2},4^{2},5^{4},8^{*}$ & $J_{6} \rightarrow  3^{4},4^{3},5^{2},8^{*}$ \\
    $J_{6} \rightarrow  3^{6},4^{4},8^{*}$ & $J_{6} \rightarrow  3,4^{4},5^{3},8^{*}$ &
    $J_{6} \rightarrow  3^{3},4^{5},5,8^{*}$ & $J_{6} \rightarrow  4^{6},5^{2},8^{*}$ \\
    $J_{6} \rightarrow  3^{2},4^{7},8^{*}$ & $J_{7} \rightarrow  3^{12},5,8^{*}$ &
    $J_{7} \rightarrow  3^{7},5^{4},8^{*}$ & $J_{7} \rightarrow  3^{2},5^{7},8^{*}$ \\
    $J_{7} \rightarrow  4^{4},5^{5},8^{*}$ & $J_{7} \rightarrow  4^{9},5,8^{*}$ &
    $J_{8} \rightarrow  3^{16},8^{*}$ & $J_{8} \rightarrow  3,5^{9},8^{*}$ \\
    $J_{8} \rightarrow  4^{12},8^{*}$ & $J_{9} \rightarrow  5^{11},8^{*}$ &
    $J_{11} \rightarrow  3^{23},8^{*}$ &  \\
    $J_{7} \rightarrow  3^{10},5^{2},9^{*}$ & $J_{7} \rightarrow  3^{5},5^{5},9^{*}$ &
    $J_{7} \rightarrow  5^{8},9^{*}$ & $J_{7} \rightarrow  3^{12},4,9^{*}$ \\
    $J_{7} \rightarrow  3^{7},4,5^{3},9^{*}$ & $J_{7} \rightarrow  3^{2},4,5^{6},9^{*}$ &
    $J_{7} \rightarrow  3^{9},4^{2},5,9^{*}$ & $J_{7} \rightarrow  3^{4},4^{2},5^{4},9^{*}$ \\
\hline
\end{longtable}

\begin{longtable}{|c|>{\raggedright\arraybackslash}p{0.74\textwidth}|}
  \caption{Table of decompositions of $J_n^{\{1,2,3,4,5,6,7\}}$}
\endfoot \hline \label{tab:1234567_decomps}
$J_{1} \rightarrow  3,4$ & 
$(0,6,3,4)$, $(0,5,7)$
\\ \hline
$J_{2} \rightarrow  4,5^{2}$ & 
$(0,5,7,4)$, $(0,6,8,1,7)$, $(3,4,5,1,6)$
\\ \hline
$J_{2} \rightarrow  3^{3},5$ & 
$(0,7,4)$, $(1,6,8)$, $(0,6,3,4,5)$, $(5,1,7)$
\\ \hline
$J_{3} \rightarrow  4^{4},5$ & 
$(0,6,3,4)$, $(1,6,8,5)$, $(1,7,9,2,8)$, $(0,5,4,7)$, $(5,6,2,7)$
\\ \hline
$J_{3} \rightarrow  3^{2},5^{3}$ & 
$(0,5,6,3,4)$, $(0,6,1,5,7)$, $(1,7,4,5,8)$, $(6,2,8)$, $(2,7,9)$
\\ \hline
$J_{4} \rightarrow  4^{7}$ & 
$(0,7,3,4)$, $(0,6,1,5)$, $(3,10,8,6)$, $(2,8,3,9)$, $(5,8,1,7)$, $(4,5,6,7)$, $(2,7,9,6)$
\\ \hline
$J_{4} \rightarrow  3,5^{5}$ & 
$(0,7,9,3,4)$, $(5,0,6,3,7)$, $(1,7,2,6,8)$, $(1,6,7,4,5)$, $(2,8,5,6,9)$, $(3,8,10)$
\\ \hline
$J_{5} \rightarrow  5^{7}$ & 
$(0,7,9,3,4)$, $(5,0,6,3,7)$, $(1,7,4,8,5)$, $(1,6,2,7,8)$, $(6,7,10,3,8)$, $(2,8,10,4,9)$, $(4,11,9,6,5)$
\\
\hline
\end{longtable}

\newpage
\begin{longtable}{|llll|}
    \hline
    $J_{3} \rightarrow  3^{6},6$ & $J_{3} \rightarrow  3,5^{3},6$ &
    $J_{3} \rightarrow  3^{3},4,5,6$ & $J_{3} \rightarrow  4^{2},5^{2},6$ \\
    $J_{3} \rightarrow  3^{2},4^{3},6$ & $J_{5} \rightarrow  4,5^{6},6$ &
    $J_{6} \rightarrow  3^{14},6$ &  \\
    $J_{3} \rightarrow  3^{4},5,7$ & $J_{3} \rightarrow  3,4,5^{2},7$ &
    $J_{3} \rightarrow  3^{3},4^{2},7$ & $J_{3} \rightarrow  4^{3},5,7$ \\
    $J_{4} \rightarrow  5^{5},7$ & $J_{4} \rightarrow  3^{7},4,7$ &
    $J_{5} \rightarrow  3^{11},7$ & $J_{3} \rightarrow  3^{6},6$ \\
    $J_{3} \rightarrow  3,5^{3},6$ & $J_{3} \rightarrow  3^{3},4,5,6$ &
    $J_{3} \rightarrow  4^{2},5^{2},6$ & $J_{3} \rightarrow  3^{2},4^{3},6$ \\
    $J_{5} \rightarrow  4,5^{6},6$ & $J_{6} \rightarrow  3^{14},6$ &
    $J_{3} \rightarrow  3^{4},5,7$ & $J_{3} \rightarrow  3,4,5^{2},7$ \\
    $J_{3} \rightarrow  3^{3},4^{2},7$ & $J_{3} \rightarrow  4^{3},5,7$ &
    $J_{4} \rightarrow  5^{5},7$ & $J_{4} \rightarrow  3^{7},4,7$ \\
    $J_{5} \rightarrow  3^{11},7$ & && \\
    $J_{3} \rightarrow  3^{2},5^{2},8$ & $J_{3} \rightarrow  3^{4},4,8$ &
    $J_{3} \rightarrow  3,4^{2},5,8$ & $J_{3} \rightarrow  4^{4},8$ \\
    $J_{4} \rightarrow  3^{8},8$ & $J_{4} \rightarrow  4,5^{4},8$ &
    $J_{6} \rightarrow  5^{8},8$ & $J_{7} \rightarrow  3^{16},8$ \\
        \hline
  \caption{These decompositions are required for Lemma
    \ref{12345678:Jgraphk=6-8}. These have been generated such that the
    $k$-cycle is incident on some subset of the vertices
    $\{n,\ldots,n+9\}$, and are given in table
    \ref{tab:12345678_decomps}.}
  \label{tab:12345678_6-8}
\end{longtable}

\begin{longtable}{|llll|}
  \caption{These decompositions are required for Lemma
    \ref{12345678:Jgraphk=9-14}. These are sorted by the value of $k$ in the
    decomposition, and are given in table \ref{tab:12345678_decomps}.}
  \endfoot
  \hline
  \endhead
  \label{tab:12345678_9-14}
    $J_{4} \rightarrow  4^{2},5^{3},9^{*}$ & $J_{4} \rightarrow  3,5^{4},9^{*}$ &
    $J_{4} \rightarrow  3,4^{5},9^{*}$ & $J_{4} \rightarrow  3^{2},4^{3},5,9^{*}$ \\
    $J_{4} \rightarrow  3^{3},4,5^{2},9^{*}$ & $J_{4} \rightarrow  3^{5},4^{2},9^{*}$ &
    $J_{4} \rightarrow  3^{6},5,9^{*}$ & $J_{5} \rightarrow  3^{9},4,9^{*}$ \\
    $J_{6} \rightarrow  4,5^{7},9^{*}$ & $J_{6} \rightarrow  3^{13},9^{*}$ &
    $J_{8} \rightarrow  5^{11},9^{*}$ & $J_{9} \rightarrow  3^{21},9^{*}$ \\
    $J_{4} \rightarrow  4^{3},5^{2},10^{*}$ & $J_{4} \rightarrow  3,4,5^{3},10^{*}$ &
    $J_{4} \rightarrow  3^{2},4^{4},10^{*}$ & $J_{4} \rightarrow  3^{3},4^{2},5,10^{*}$ \\
    $J_{4} \rightarrow  3^{4},5^{2},10^{*}$ & $J_{4} \rightarrow  3^{6},4,10^{*}$ &
    $J_{5} \rightarrow  5^{6},10^{*}$ & $J_{5} \rightarrow  3^{10},10^{*}$ \\
    $J_{8} \rightarrow  4,5^{10},10^{*}$ & $J_{8} \rightarrow  3^{18},10^{*}$ &
    & \\
    $J_{5} \rightarrow  4,5^{5},11^{*}$ & $J_{5} \rightarrow  4^{6},5,11^{*}$ &
    $J_{5} \rightarrow  3,4^{4},5^{2},11^{*}$ & $J_{5} \rightarrow  3^{2},4^{2},5^{3},11^{*}$ \\
    $J_{5} \rightarrow  3^{3},5^{4},11^{*}$ & $J_{5} \rightarrow  3^{3},4^{5},11^{*}$ &
    $J_{5} \rightarrow  3^{4},4^{3},5,11^{*}$ & $J_{5} \rightarrow  3^{5},4,5^{2},11^{*}$ \\
\hline
\end{longtable}

\begin{longtable}{|c|>{\raggedright\arraybackslash}p{0.75\textwidth}|}
  \caption{Table of decompositions of $J_n^{\{1,2,3,4,5,6,7,8\}}$}
\endfoot \hline \label{tab:12345678_decomps}
$J_{1} \rightarrow  4^{2}$ & 
$(0,5,8,7)$, $(5,9,1,7)$
\\ \hline
$J_{1} \rightarrow  3,5$ & 
$(0,5,7)$, $(7,1,9,5,8)$
\\ \hline
$J_{2} \rightarrow  3^{4},4$ & 
$(0,5,7)$, $(1,6,9)$, $(7,1,8)$, $(6,10,2,8)$, $(5,9,8)$
\\ \hline
$J_{3} \rightarrow  4,5^{4}$ & 
$(0,5,9,8,7)$, $(1,8,6,9)$, $(1,6,10,2,7)$, $(5,8,2,9,7)$, $(9,3,11,7,10)$
\\ \hline
$J_{4} \rightarrow  3^{9},5$ & 
$(0,5,7)$, $(1,6,9)$, $(6,10,4,12,8)$, $(7,11,10)$, $(7,1,8)$, $(2,8,10)$, $(2,7,9)$, $(3,8,11)$, $(9,3,10)$, $(5,9,8)$
\\ \hline
$J_{5} \rightarrow  5^{8}$ & 
$(0,5,9,8,7)$, $(1,8,10,6,9)$, $(5,8,6,1,7)$, $(2,8,3,9,10)$, $(2,7,11,4,9)$, $(7,10,3,11,9)$, $(10,4,12,8,11)$, $(11,5,13,9,12)$
\\ \hline
$J_{5} \rightarrow  3^{12},4$ & 
$(0,5,7)$, $(1,6,8)$, $(5,13,9)$, $(6,10,4,9)$, $(7,11,10)$, $(7,2,8)$, $(1,7,9)$, $(8,12,9)$, $(4,11,12)$, $(5,11,8)$, $(3,9,11)$, $(3,8,10)$, $(2,9,10)$
\\ \hline
$J_{6} \rightarrow  3^{16}$ & 
$(0,5,7)$, $(1,6,8)$, $(5,13,9)$, $(6,12,9)$, $(3,9,11)$, $(8,2,9)$, $(1,7,9)$, $(9,4,10)$, $(6,14,10)$, $(5,12,8)$, $(7,11,8)$, $(2,7,10)$, $(3,8,10)$, $(10,5,11)$, $(4,11,12)$, $(10,13,12)$
\\
\hline
\end{longtable}

\section{Data for Section \ref{sec:123n}}
\label{appendix2}

\begin{longtable}{|c|p{0.75\textwidth}|}
  \caption{Decompositions of graphs for section \ref{sec:123n}}
  \endfoot
  \hline
  \label{tab:123n_decomps}
  $L_2 \rightarrow 3,4^+,1^H$ & $[2,4,1,0,3]$, $[3,4]$, $(1,2,3)$
  \\ \hline
  $L_2 \rightarrow 4,3^+,1^H$ & $[2,4,1,3]$, $[3,4]$, $(0,1,2,3)$
  \\ \hline
  $L_2 \rightarrow 5,2^+,1^H$ & $[2,1,3]$, $[3,4]$, $(0,1,5,2,3)$
  \\ \hline
  $L_2 \rightarrow 3^2,1^+,1^H$ & $[2,3]$, $[3,4]$, $(0,1,3)$, $(1,2,5)$
  \\ \hline
  $L_3 \rightarrow 4,5^+,2^H$ & $[3,0,1,2,5,4]$, $[5,3,4]$, $(1,3,2,4)$
  \\ \hline
  $L_3 \rightarrow 3^2,3^+,2^H$ & $[3,2,5,4]$, $[5,3,4]$, $(0,1,3)$, $(1,2,4)$
  \\ \hline
  $L_3 \rightarrow 3,4,2^+,2^H$ & $[3,2,4]$, $[5,3,4]$, $(0,1,3)$, $(1,2,5,4)$
  \\ \hline
  $L_4 \rightarrow 5,6^+,3^H$ & $[4,2,1,0,3,6,5]$, $[6,4,3,5]$, $(1,3,2,5,4)$
  \\
\hline

\end{longtable}

\newpage
The next five tables
contain all decompositions required for Lemma \ref{123:JgraphsClosekSmalln}.
We make extensive use of concatenation to obtain the results presented in these
tables. We therefore use the notation
$(G \rightarrow M) \oplus (H \rightarrow M')$ for the concatenation of the
decomposition $G \rightarrow M$ with the decomposition $H \rightarrow M'$.
Various different forms of concatenation are defined and used in Section
\ref{sec:123n}. In the table, the particular form of concatenation being used
is well-defined by the decompositions it involves. We also define the notation
$a \cdot (G \rightarrow M)$ to mean $(G \rightarrow M) \oplus (G \rightarrow M)
\oplus \cdots \oplus (G \rightarrow M)$, where $a$ is the number of copies of
the decomposition $G\rightarrow M$ involved.

\begin{longtable}{|c|p{0.74\textwidth}|}
  \caption{Table of decompositions for Lemma \ref{123:JgraphsClosekSmalln} with
$j=1$}
\endfoot
\hline
\label{tab:123n_concats_1}
$J_{7} \rightarrow 4^{2},6,7$ & $(L_2 \rightarrow 4,3^+,1^H)  \oplus (Q_5 \rightarrow 4,3^+,6^H)$
\\ \hline
$J_{7} \rightarrow 3,5,6,7$ & $(L_3 \rightarrow 5,4^+,2^H)  \oplus (Q_4 \rightarrow 3,2^+,5^H)$
\\ \hline
$J_{8} \rightarrow 4,5,7,8$ & $(L_3 \rightarrow 5,4^+,2^H)  \oplus (Q_5 \rightarrow 4,3^+,6^H)$
\\ \hline
$J_{8} \rightarrow 3^{3},7,8$ & $(5,2,0,3,1,4,6,7)$, $(5,3,6)$, $(2,3,4)$, $(9,6,8,5,4,7,10)$, $(7,8,9)$
\\ \hline
$J_{9} \rightarrow 5^{2},8,9$ & $(L_3 \rightarrow 5,4^+,2^H)  \oplus (Q_6 \rightarrow 5,4^+,7^H)$
\\ \hline
$J_{9} \rightarrow 3^{2},4,8,9$ & $(L_2 \rightarrow 3,4^+,1^H)  \oplus (P_3 \rightarrow 4,2^+,3^H) \oplus (Q_4 \rightarrow 3,2^+,5^H)$
\\ \hline
$J_{10} \rightarrow 3,4^{2},9,10$ & $(L_2 \rightarrow 4,3^+,1^H)  \oplus (P_4 \rightarrow 4,4^+,4^H) \oplus (Q_4 \rightarrow 3,2^+,5^H)$
\\ \hline
$J_{10} \rightarrow 3^{2},5,9,10$ & $(L_2 \rightarrow 3,4^+,1^H)  \oplus (P_4 \rightarrow 5,3^+,4^H) \oplus (Q_4 \rightarrow 3,2^+,5^H)$
\\
\hline
\end{longtable}

\begin{longtable}{|c|p{0.73\textwidth}|}
  \caption{Table of decompositions for Lemma \ref{123:JgraphsClosekSmalln} with
$j=2$}
\endfoot
\hline
\label{tab:123n_concats_2}
$J_{8} \rightarrow 5^{2},6,8$ & $(L_2 \rightarrow 5,2^+,1^H)  \oplus (Q_6 \rightarrow 5,4^+,7^H)$
\\ \hline
$J_{8} \rightarrow 3^{2},4,6,8$ & $(L_3 \rightarrow 3^2,3^+,2^H)  \oplus (Q_5 \rightarrow 4,3^+,6^H)$
\\ \hline
$J_{9} \rightarrow 3,4^{2},7,9$ & $(L_2 \rightarrow 4,3^+,1^H)  \oplus (P_3 \rightarrow 4,2^+,3^H) \oplus (Q_4 \rightarrow 3,2^+,5^H)$
\\ \hline
$J_{9} \rightarrow 3^{2},5,7,9$ & $(L_3 \rightarrow 3^2,3^+,2^H)  \oplus (Q_6 \rightarrow 5,4^+,7^H)$
\\ \hline
$J_{10} \rightarrow 4^{3},8,10$ & $(L_2 \rightarrow 4,3^+,1^H)  \oplus (P_3 \rightarrow 4,2^+,3^H) \oplus (Q_5 \rightarrow 4,3^+,6^H)$
\\ \hline
$J_{10} \rightarrow 3,4,5,8,10$ & $(L_2 \rightarrow 5,2^+,1^H)  \oplus (P_4 \rightarrow 4,4^+,4^H) \oplus (Q_4 \rightarrow 3,2^+,5^H)$
\\ \hline
$J_{10} \rightarrow 3^{4},8,10$ & $(8,5,2,0,3,1,4,7,6,9),\allowbreak (2,3,4),\allowbreak (5,3,6),\allowbreak (10,9,7,5,4,6,8,11),\allowbreak (8,7,10),\allowbreak (11,9,12)$
\\ \hline
$J_{11} \rightarrow 4^{2},5,9,11$ & $(L_2 \rightarrow 5,2^+,1^H)  \oplus (P_4 \rightarrow 4,4^+,4^H) \oplus (Q_5 \rightarrow 4,3^+,6^H)$
\\
\hline
\end{longtable}
\begin{longtable}{|c|p{0.72\textwidth}|}
  \caption{Table of decompositions for Lemma \ref{123:JgraphsClosekSmalln} with
$j=3$}
\endfoot
\hline
\label{tab:123n_concats_3}
$J_{9} \rightarrow 4^{3},6,9$ & $(7,6,4,1,3,0,2,5,8),\allowbreak (3,2,4,5),\allowbreak (8,6,3,4,7,10),\allowbreak (7,5,6,9),\allowbreak (10,9,8,11)$
\\ \hline
$J_{9} \rightarrow 3,4,5,6,9$ & $(L_3 \rightarrow 3,4,2^+,2^H)  \oplus (Q_5 \rightarrow 4,3^+,6^H)$
\\ \hline
$J_{9} \rightarrow 3^{4},6,9$ & $(6,7,4,1,3,0,2,5,8),\allowbreak (2,3,4),\allowbreak (5,3,6),\allowbreak (8,7,5,4,6,9),\allowbreak (9,7,10),\allowbreak (10,8,11)$
\\ \hline
$J_{10} \rightarrow 4^{2},5,7,10$ & $(L_2 \rightarrow 5,2^+,1^H)  \oplus (P_3 \rightarrow 4,2^+,3^H) \oplus (Q_5 \rightarrow 4,3^+,6^H)$
\\ \hline
$J_{10} \rightarrow 3,5^{2},7,10$ & $(L_2 \rightarrow 5,2^+,1^H)  \oplus (P_4 \rightarrow 5,3^+,4^H) \oplus (Q_4 \rightarrow 3,2^+,5^H)$
\\ \hline
$J_{10} \rightarrow 3^{3},4,7,10$ & $(L_2 \rightarrow 3^2,1^+,1^H)  \oplus (P_4 \rightarrow 4,4^+,4^H) \oplus (Q_4 \rightarrow 3,2^+,5^H)$
\\ \hline
$J_{11} \rightarrow 4,5^{2},8,11$ & $(L_2 \rightarrow 5,2^+,1^H)  \oplus (P_4 \rightarrow 5,3^+,4^H) \oplus (Q_5 \rightarrow 4,3^+,6^H)$
\\ \hline
$J_{11} \rightarrow 3^{2},4^{2},8,11$ & $(L_3 \rightarrow 3,4,2^+,2^H)  \oplus (P_4 \rightarrow 4,4^+,4^H) \oplus (Q_4 \rightarrow 3,2^+,5^H)$
\\
\hline
\end{longtable}
\begin{longtable}{|c|p{0.72\textwidth}|}
  \caption{Table of decompositions for Lemma \ref{123:JgraphsClosekSmalln} with
$j=4$}
\endfoot
\hline
\label{tab:123n_concats_4}
$J_{10} \rightarrow 4,5^{2},6,10$ & $(L_4 \rightarrow 4,5,2^+,3^H)  \oplus (Q_6 \rightarrow 5,4^+,7^H)$
\\ \hline
$J_{10} \rightarrow 3^{2},4^{2},6,10$ & $(L_5 \rightarrow 3^2,4,3^+,4^H)  \oplus (Q_5 \rightarrow 4,3^+,6^H)$
\\ \hline
$J_{10} \rightarrow 3^{3},5,6,10$ & $(L_4 \rightarrow 3^3,2^+,3^H)  \oplus (Q_6 \rightarrow 5,4^+,7^H)$
\\ \hline
$J_{11} \rightarrow 5^{3},7,11$ & $(8,5,2,0,3,1,4,7,6,9,10),\allowbreak (5,4,2,3,6),\allowbreak (7,5,3,4,6,8,9),\allowbreak (10,7,8,11,12),\allowbreak (12,9,11,10,13)$
\\ \hline
$J_{11} \rightarrow 3,4^{3},7,11$ & $(L_3 \rightarrow 3,4,2^+,2^H)  \oplus (P_3 \rightarrow 4,2^+,3^H) \oplus (Q_5 \rightarrow 4,3^+,6^H)$
\\ \hline
$J_{11} \rightarrow 3^{2},4,5,7,11$ & $(L_5 \rightarrow 3^2,4,3^+,4^H)  \oplus (Q_5 \rightarrow 4,3^+,6^H)$
\\ \hline
$J_{11} \rightarrow 3^{5},7,11$ & $(L_2 \rightarrow 3^2,1^+,1^H)  \oplus (P_5 \rightarrow 3^2,4^+,5^H) \oplus (Q_4 \rightarrow 3,2^+,5^H)$
\\ \hline
$J_{12} \rightarrow 4^{4},8,12$ & $(10,9,7,6,4,1,3,0,2,5,8,11),\allowbreak (3,2,4,5),\allowbreak (11,9,6,3,4,7,10,13),\allowbreak (6,5,7,8),\allowbreak (10,8,9,12),\allowbreak (13,12,11,14)$
\\
\hline
\end{longtable}
\begin{longtable}{|c|p{0.72\textwidth}|}
  \caption{Table of decompositions for Lemma \ref{123:JgraphsClosekSmalln} with
$j=5$}
\endfoot
\hline
\label{tab:123n_concats_5}
$J_{11} \rightarrow 4^{4},6,11$ & $(9,8,6,4,1,3,0,2,5,7,10),\allowbreak (3,2,4,5),\allowbreak (6,3,4,7),\allowbreak (9,7,8,11),\allowbreak (10,8,5,6,9,12),\allowbreak (12,11,10,13)$
\\ \hline
$J_{11} \rightarrow 3,4^{2},5,6,11$ & $(L_4 \rightarrow 4,5,2^+,3^H)  \oplus (P_3 \rightarrow 4,2^+,3^H) \oplus (Q_4 \rightarrow 3,2^+,5^H)$
\\ \hline
$J_{11} \rightarrow 3^{2},5^{2},6,11$ & $(8,5,2,0,3,1,4,7,6,9,10),\allowbreak (2,3,4),\allowbreak (5,3,6),\allowbreak (6,4,5,7,8),\allowbreak (12,11,9,7,10,13),\allowbreak (10,11,8,9,12)$
\\ \hline
$J_{11} \rightarrow 3^{4},4,6,11$ & $(L_4 \rightarrow 3^3,2^+,3^H)  \oplus (P_3 \rightarrow 4,2^+,3^H) \oplus (Q_4 \rightarrow 3,2^+,5^H)$
\\ \hline
$J_{12} \rightarrow 4^{3},5,7,12$ & $(10,9,6,7,4,1,3,0,2,5,8,11),\allowbreak (6,3,2,4,5,7,8),\allowbreak (4,3,5,6),\allowbreak (10,7,9,12),\allowbreak (11,9,8,10,13),\allowbreak (13,12,11,14)$
\\ \hline
$J_{12} \rightarrow 3,4,5^{2},7,12$ & $(L_4 \rightarrow 4,5,2^+,3^H)  \oplus (P_4 \rightarrow 5,3^+,4^H) \oplus (Q_4 \rightarrow 3,2^+,5^H)$
\\ \hline
$J_{12} \rightarrow 3^{3},4^{2},7,12$ & $(L_5 \rightarrow 4,3^2,3^+,4^H)  \oplus (P_3 \rightarrow 4,2^+,3^H) \oplus (Q_4 \rightarrow 3,2^+,5^H)$
\\ \hline
$J_{12} \rightarrow 3^{4},5,7,12$ & $(L_4 \rightarrow 3^3,2^+,3^H)  \oplus (P_4 \rightarrow 5,3^+,4^H) \oplus (Q_4 \rightarrow 3,2^+,5^H)$
\\
\hline
\end{longtable}

\newpage
\begin{longtable}{|c|>{\raggedright\arraybackslash}p{0.85\textwidth}|}
  \hline 
  $t$ & Difference triples \\
  \hline \hline
\endhead
  \caption{Difference triples for Lemma \ref{3s4sandHamsforKn-S}}
\endfoot
\label{tab:5789_triples}
$11$ &  $\{5, 24, 29\}$, $\{7, 20, 27\}$, $\{8, 23, 31\}$, $\{9, 21, 30\}$, $\{10, 26, 36\}$, $\{11, 28, 39\}$,
$\{12, 25, 37\}$, $\{13, 22, 35\}$, $\{14, 18, 32\}$, $\{15, 19, 34\}$, $\{16, 17, 33\}$ \\
\hline
$12$ &  $\{5, 21, 26\}$, $\{7, 27, 34\}$, $\{8, 23, 31\}$, $\{9, 24, 33\}$, $\{10, 28, 38\}$, $\{11, 29, 40\}$,
$\{12, 30, 42\}$, $\{13, 19, 32\}$, $\{14, 25, 39\}$, $\{15, 22, 37\}$, $\{16, 20, 36\}$, $\{17, 18, 35\}$
\\
\hline
$13$ &  $\{5, 20, 25\}$, $\{7, 28, 35\}$, $\{8, 26, 34\}$, $\{9, 29, 38\}$, $\{10, 23, 33\}$, $\{11, 31, 42\}$,
$\{12, 32, 44\}$, $\{13, 30, 43\}$, $\{14, 27, 41\}$, $\{15, 21, 36\}$, $\{16, 24, 40\}$, $\{17, 22, 39\}$,
$\{18, 19, 37\}$ \\
\hline
$14$ &  $\{5, 18, 23\}$, $\{7, 30, 37\}$, $\{8, 31, 39\}$, $\{9, 26, 35\}$, $\{10, 28, 38\}$, $\{11, 25, 36\}$,
$\{12, 33, 45\}$, $\{13, 34, 47\}$, $\{14, 32, 46\}$, $\{15, 29, 44\}$, $\{16, 27, 43\}$, $\{17, 24, 41\}$,
$\{19, 21, 40\}$, $\{20, 22, 42\}$ \\
\hline
$15$ &  $\{5, 33, 38\}$, $\{7, 13, 20\}$, $\{8, 34, 42\}$, $\{9, 28, 37\}$, $\{10, 31, 41\}$, $\{11, 29, 40\}$,
$\{12, 27, 39\}$, $\{14, 35, 49\}$, $\{15, 36, 51\}$, $\{16, 32, 48\}$, $\{17, 30, 47\}$, $\{18, 26, 44\}$,
$\{19, 24, 43\}$, $\{21, 25, 46\}$, $\{22, 23, 45\}$ \\
\hline
$16$ &  $\{5, 12, 17\}$, $\{7, 34, 41\}$, $\{8, 31, 39\}$, $\{9, 35, 44\}$, $\{10, 30, 40\}$, $\{11, 32, 43\}$,
$\{13, 29, 42\}$, $\{14, 36, 50\}$, $\{15, 37, 52\}$, $\{16, 38, 54\}$, $\{18, 33, 51\}$, $\{19, 27, 46\}$,
$\{20, 25, 45\}$, $\{21, 28, 49\}$, $\{22, 26, 48\}$, $\{23, 24, 47\}$ \\
\hline
$17$ &  $\{5, 9, 14\}$, $\{7, 34, 41\}$, $\{8, 37, 45\}$, $\{10, 36, 46\}$, $\{11, 31, 42\}$, $\{12, 32, 44\}$,
$\{13, 30, 43\}$, $\{15, 39, 54\}$, $\{16, 40, 56\}$, $\{17, 38, 55\}$, $\{18, 35, 53\}$, $\{19, 33, 52\}$,
$\{20, 28, 48\}$, $\{21, 26, 47\}$, $\{22, 29, 51\}$, $\{23, 27, 50\}$, $\{24, 25, 49\}$ \\
\hline
\end{longtable}

\begin{longtable}{|c|>{\raggedright\arraybackslash}p{0.85\textwidth}|}
  \hline 
  $t$ & Difference triples \\
  \hline \hline
\endhead
  \caption{Difference triples for Lemma \ref{3s4sandHamsforKn-S}}
\endfoot
\label{tab:689_triples}
$13$ &  $\{6, 26, 32\}$, $\{8, 27, 35\}$, $\{9, 25, 34\}$, $\{10, 31, 41\}$, $\{11, 33, 44\}$, $\{12, 24, 36\}$,
$\{13, 29, 42\}$, $\{14, 23, 37\}$, $\{15, 28, 43\}$, $\{16, 30, 46\}$, $\{17, 21, 38\}$, $\{18, 22, 40\}$,
$\{19, 20, 39\}$ \\
\hline
$14$ &  $\{6, 25, 31\}$, $\{8, 28, 36\}$, $\{9, 29, 38\}$, $\{10, 27, 37\}$, $\{11, 33, 44\}$, $\{12, 35, 47\}$,
$\{13, 32, 45\}$, $\{14, 34, 48\}$, $\{15, 24, 39\}$, $\{16, 30, 46\}$, $\{17, 26, 43\}$, $\{18, 22, 40\}$,
$\{19, 23, 42\}$, $\{20, 21, 41\}$ \\
\hline
$15$ &  $\{6, 23, 29\}$, $\{8, 30, 38\}$, $\{9, 31, 40\}$, $\{10, 36, 46\}$, $\{11, 28, 39\}$, $\{12, 35, 47\}$,
$\{13, 37, 50\}$, $\{14, 27, 41\}$, $\{15, 34, 49\}$, $\{16, 32, 48\}$, $\{17, 25, 42\}$, $\{18, 33, 51\}$,
$\{19, 26, 45\}$, $\{20, 24, 44\}$, $\{21, 22, 43\}$ \\
\hline
$16$ &  $\{6, 20, 26\}$, $\{8, 32, 40\}$, $\{9, 33, 42\}$, $\{10, 31, 41\}$, $\{11, 37, 48\}$, $\{12, 38, 50\}$,
$\{13, 30, 43\}$, $\{14, 39, 53\}$, $\{15, 29, 44\}$, $\{16, 35, 51\}$, $\{17, 28, 45\}$, $\{18, 34, 52\}$,
$\{19, 36, 55\}$, $\{21, 25, 46\}$, $\{22, 27, 49\}$, $\{23, 24, 47\}$ \\
\hline
$17$ &  $\{6, 17, 23\}$, $\{8, 34, 42\}$, $\{9, 35, 44\}$, $\{10, 33, 43\}$, $\{11, 36, 47\}$, $\{12, 41, 53\}$,
$\{13, 32, 45\}$, $\{14, 40, 54\}$, $\{15, 31, 46\}$, $\{16, 39, 55\}$, $\{18, 30, 48\}$, $\{19, 37, 56\}$,
$\{20, 38, 58\}$, $\{21, 28, 49\}$, $\{22, 29, 51\}$, $\{24, 26, 50\}$, $\{25, 27, 52\}$ \\
\hline
$18$ &  $\{6, 14, 20\}$, $\{8, 36, 44\}$, $\{9, 37, 46\}$, $\{10, 35, 45\}$, $\{11, 38, 49\}$, $\{12, 43, 55\}$,
$\{13, 34, 47\}$, $\{15, 33, 48\}$, $\{16, 41, 57\}$, $\{17, 42, 59\}$, $\{18, 40, 58\}$, $\{19, 31, 50\}$,
$\{21, 39, 60\}$, $\{22, 29, 51\}$, $\{23, 30, 53\}$, $\{24, 32, 56\}$, $\{25, 27, 52\}$, $\{26, 28, 54\}$
\\
\hline
\end{longtable}

\newpage
\begin{longtable}{|c|p{0.85\textwidth}|}
\hline
$r$ &$\pi_r$ \\
\hline \hline
\endhead
\caption{Table of partitions for Lemma 12 when $S=\{1,2,3,4\}$}
\endfoot
\label{tab:quints_1234}
1 & $\{\{5, 6, 7, 8, 10\}\}$\\
\hline
2 & $\{\{5, 7, 8, 9, 11\}$, $\{6, 10, 12, 13, 15\}\}$\\
\hline
3 & $\{\{5, 8, 11, 16, 18\}$, $\{6, 9, 12, 13, 14\}$, $\{7, 10, 15, 17, 19\}\}$\\
\hline
4 & $\{\{5, 11, 13, 21, 24\}$, $\{6, 8, 17, 19, 22\}$, $\{7, 14, 18, 20, 23\}$,  
  $\{9, 10, 12, 15, 16\}\}$\\
\hline
5 & $\{\{5, 13, 22, 23, 27\}$, $\{6, 14, 15, 16, 21\}$, $\{7, 8, 19, 20, 24\}$,  
  $\{9, 11, 17, 25, 28\}$, $\{10, 12, 18, 26, 30\}\}$\\
\hline
6 & $\{\{5, 6, 21, 22, 32\}$, $\{7, 8, 20, 23, 28\}$, $\{9, 12, 19, 31, 33\}$, 
  $\{10, 14, 26, 27, 29\}$, $\{11, 15, 16, 18, 24\}$, $\{13, 17, 25, 30, 35\}\}$\\
\hline
\end{longtable}

\begin{longtable}{|c|>{\raggedright\arraybackslash}p{0.85\textwidth}|}
\hline
$r$ &$\pi_r$ \\
\hline \hline
\endhead
\caption{Table of partitions for Lemma 12 when $S=\{1,2,3,4,6\}$}
\endfoot
\label{tab:quints_12346}
1 & $\{\{5, 7, 8, 9, 11\}\}$\\
\hline
2 & $\{\{5, 8, 11, 13, 15\}$, $\{7, 9, 10, 12, 14\}\}$\\
\hline
3 & $\{\{5, 7, 14, 18, 20\}$, $\{8, 10, 15, 16, 17\}$, $\{9, 11, 12, 13, 19\}\}$\\
\hline
4 & $\{\{5, 7, 15, 23, 26\}$, $\{8, 9, 20, 21, 24\}$, $\{10, 13, 16, 17, 22\}$, $\{11, 12, 14, 18, 19\}\}$\\
\hline
5 & $\{\{5, 11, 19, 26, 29\}$, $\{7, 9, 23, 24, 31\}$, $\{8, 16, 25, 27, 28\}$, $\{10, 13, 15, 18, 20\}$, $\{12, 14, 17, 21, 22\}\}$\\
\hline
6 & $\{\{5, 13, 21, 30, 33\}$, $\{7, 12, 25, 29, 35\}$, $\{8, 9, 22, 26, 31\}$, $\{10, 11, 27, 28, 34\}$, $\{14, 15, 18, 23, 24\}$, $\{16, 17, 19, 20, 32\}\}$\\
\hline
7 & $\{\{5, 10, 25, 30, 40\}$, $\{7, 12, 23, 34, 38\}$, $\{8, 14, 27, 32, 37\}$, $\{9, 19, 31, 36, 39\}$, $\{11, 20, 29, 33, 35\}$, $\{13, 15, 18, 22, 24\}$, $\{16, 17, 21, 26, 28\}\}$\\
\hline
8 & $\{\{5, 17, 31, 35, 44\}$, $\{7, 12, 29, 30, 40\}$, $\{8, 24, 36, 37, 41\}$, $\{9, 11, 25, 38, 43\}$, $\{10, 15, 33, 34, 42\}$, $\{13, 19, 22, 26, 28\}$, $\{14, 16, 20, 23, 27\}$, $\{18, 21, 32, 39, 46\}\}$\\
\hline
\end{longtable}

\begin{longtable}{|c|>{\raggedright\arraybackslash}p{0.85\textwidth}|}
\hline 
$r$ &$\pi_r$ \\
\hline \hline
\endhead
\caption{Table of partitions for Lemma 12 when $S=\{1,2,3,4,5,7\}$}
\endfoot
\label{tab:quints_123457}
2 & $\{\{6, 9, 12, 13, 16\}$, $\{8, 10, 11, 14, 15\}\}$\\
\hline
3 & $\{\{6, 8, 13, 18, 19\}$, $\{9, 12, 16, 17, 22\}$, $\{10, 11, 14, 15, 20\}\}$\\
\hline
4 & $\{\{6, 12, 15, 16, 19\}$, $\{8, 14, 17, 22, 27\}$, $\{9, 11, 18, 23, 25\}$, $\{10, 13, 20, 21, 24\}\}$\\
\hline
5 & $\{\{6, 8, 19, 22, 27\}$, $\{9, 17, 21, 25, 30\}$, $\{10, 11, 18, 26, 29\}$, $\{12, 14, 23, 28, 31\}$, $\{13, 15, 16, 20, 24\}\}$\\
\hline
6 & $\{\{6, 9, 21, 24, 30\}$, $\{8, 19, 25, 33, 35\}$, $\{10, 11, 22, 27, 28\}$, $\{12, 17, 26, 31, 34\}$, $\{13, 14, 16, 20, 23\}$, $\{15, 18, 29, 32, 36\}\}$\\
\hline
7 & $\{\{6, 13, 26, 27, 34\}$, $\{8, 19, 23, 36, 40\}$, $\{9, 18, 32, 33, 38\}$, $\{10, 12, 24, 29, 31\}$, $\{11, 16, 20, 22, 25\}$, $\{14, 21, 30, 37, 42\}$, $\{15, 17, 28, 35, 39\}\}$\\
\hline
8 & $\{\{6, 11, 24, 31, 38\}$, $\{8, 9, 27, 30, 40\}$, $\{10, 14, 28, 29, 33\}$, $\{12, 23, 32, 42, 45\}$, $\{13, 17, 36, 37, 43\}$, $\{15, 19, 26, 39, 47\}$, $\{16, 21, 34, 41, 44\}$, $\{18, 20, 22, 25, 35\}\}$\\
\hline
9 & $\{\{6, 11, 26, 36, 45\}$, $\{8, 19, 33, 37, 43\}$, $\{9, 12, 27, 40, 46\}$, $\{10, 13, 29, 38, 44\}$, $\{14, 18, 28, 47, 51\}$, $\{15, 21, 25, 30, 31\}$, $\{16, 22, 39, 41, 42\}$, $\{17, 34, 48, 49, 50\}$, $\{20, 23, 24, 32, 35\}\}$\\
\hline
\end{longtable}

\begin{longtable}{|c|>{\raggedright\arraybackslash}p{0.85\textwidth}|}
\hline
$r$ &$\pi_r$ \\
\hline \hline
\endhead
\caption{Table of partitions for Lemma 12 when $S=\{1,2,3,4,5,6,7\}$}
\endfoot
\label{tab:quints_1234567}
2 & $\{\{8, 11, 13, 14, 18\}$, $\{9, 10, 12, 15, 16\}\}$\\
\hline
3 & $\{\{8, 11, 14, 18, 23\}$, $\{9, 10, 16, 17, 20\}$, $\{12, 13, 15, 19, 21\}\}$\\
\hline
4 & $\{\{8, 13, 20, 23, 24\}$, $\{9, 11, 16, 21, 25\}$, $\{10, 12, 15, 18, 19\}$, $\{14, 17, 22, 26, 27\}\}$\\
\hline
5 & $\{\{8, 16, 21, 28, 31\}$, $\{9, 11, 23, 24, 27\}$, $\{10, 14, 20, 26, 30\}$, $\{12, 13, 15, 18, 22\}$, $\{17, 19, 25, 29, 32\}\}$\\
\hline
6 & $\{\{8, 11, 22, 27, 30\}$, $\{9, 12, 25, 28, 32\}$, $\{10, 18, 23, 29, 34\}$, $\{13, 15, 26, 33, 35\}$, $\{14, 17, 24, 31, 38\}$, $\{16, 19, 20, 21, 36\}\}$\\
\hline
7 & $\{\{8, 10, 27, 29, 38\}$, $\{9, 12, 26, 28, 33\}$, $\{11, 19, 24, 37, 43\}$, $\{13, 16, 23, 35, 41\}$, $\{14, 22, 32, 36, 40\}$, $\{15, 20, 21, 25, 31\}$, $\{17, 18, 30, 34, 39\}\}$\\
\hline
8 & $\{\{8, 20, 32, 34, 38\}$, $\{9, 22, 25, 41, 47\}$, $\{10, 11, 26, 30, 35\}$, $\{12, 18, 27, 43, 46\}$, $\{13, 23, 33, 42, 45\}$, $\{14, 19, 36, 37, 40\}$, $\{15, 21, 24, 29, 31\}$, $\{16, 17, 28, 39, 44\}\}$\\
\hline
9 & $\{\{8, 18, 32, 39, 45\}$, $\{9, 12, 30, 35, 44\}$, $\{10, 19, 33, 38, 42\}$, $\{11, 26, 34, 46, 49\}$, $\{13, 20, 40, 41, 48\}$, $\{14, 22, 24, 29, 31\}$, $\{15, 16, 28, 47, 50\}$, $\{17, 21, 25, 27, 36\}$, $\{23, 37, 43, 51, 52\}\}$\\
\hline
\end{longtable}

\begin{longtable}{|c|>{\raggedright\arraybackslash}p{0.85\textwidth}|}
\hline
$r$ &$\pi_r$ \\
\hline \hline
\endhead
\caption{Table of partitions for Lemma 12 when $S=\{1,2,3,4,5,6,7,8\}$}
\endfoot
\label{tab:quints_12345678}
3 & $\{\{9, 13, 15, 16, 21\}$, $\{10, 14, 17, 18, 23\}$, $\{11, 12, 19, 20, 22\}\}$\\
\hline
4 & $\{\{9, 12, 17, 18, 20\}$, $\{10, 15, 16, 19, 22\}$, $\{11, 21, 23, 27, 28\}$, $\{13, 14, 24, 25, 26\}\}$\\
\hline
5 & $\{\{9, 16, 20, 25, 30\}$, $\{10, 14, 21, 29, 32\}$, $\{11, 18, 19, 22, 26\}$, $\{12, 17, 23, 28, 34\}$, $\{13, 15, 24, 27, 31\}\}$\\
\hline
6 & $\{\{9, 16, 29, 30, 34\}$, $\{10, 11, 23, 26, 28\}$, $\{12, 17, 20, 22, 27\}$, $\{13, 19, 25, 32, 39\}$, $\{14, 15, 24, 31, 36\}$, $\{18, 21, 33, 35, 37\}\}$\\
\hline
7 & $\{\{9, 10, 24, 31, 36\}$, $\{11, 18, 23, 25, 27\}$, $\{12, 21, 29, 37, 41\}$, $\{13, 20, 28, 38, 43\}$, $\{14, 15, 26, 30, 33\}$, $\{16, 17, 32, 39, 40\}$, $\{19, 22, 34, 35, 42\}\}$\\
\hline
8 & $\{\{9, 13, 27, 38, 43\}$, $\{10, 23, 35, 37, 39\}$, $\{11, 21, 31, 32, 33\}$, $\{12, 15, 26, 44, 45\}$, $\{14, 20, 29, 41, 46\}$, $\{16, 17, 25, 28, 30\}$, $\{18, 24, 36, 42, 48\}$, $\{19, 22, 34, 40, 47\}\}$\\
\hline
9 & $\{\{9, 46, 49, 50, 54\}$, $\{10, 23, 26, 29, 30\}$, $\{11, 13, 33, 38, 47\}$, $\{12, 25, 31, 45, 51\}$, $\{14, 15, 35, 42, 48\}$, $\{16, 17, 28, 34, 39\}$, $\{18, 19, 40, 41, 44\}$, $\{20, 21, 27, 32, 36\}$, $\{22, 24, 37, 43, 52\}\}$\\
\hline
10 & $\{\{9, 27, 30, 35, 41\}$, $\{10, 25, 43, 46, 54\}$, $\{11, 13, 36, 37, 49\}$, $\{12, 14, 32, 44, 50\}$, $\{15, 20, 38, 45, 48\}$, $\{16, 28, 29, 33, 40\}$, $\{17, 19, 34, 53, 55\}$, $\{18, 21, 31, 51, 59\}$, $\{22, 26, 39, 47, 56\}$, $\{23, 24, 42, 52, 57\}\}$\\
\hline
\end{longtable}

\vspace{0.5cm} \noindent{\bf Acknowledgements:}\quad
The authors acknowledge the support of the Australian Research Council
via grants DP0770400, DE120100040, DP120100790 and DP120103067. The third author
was supported by an Australian Postgraduate Award.

\end{document}